\documentclass[11pt,twoside,spanish]{article}
 
\usepackage{srcltx,babel,float,amsfonts,amssymb,amsmath,epsfig,amscd}
 
\usepackage{graphicx,psfrag}
\usepackage[usenames]{color}

\usepackage{makeidx}

\newtheorem{definition}{Definici\'on}[subsection]

\newtheorem{theorem}[definition]{Teorema}
 \newtheorem{lemma}[definition]{Lema}
\newtheorem{corollary}[definition]{Corolario}
\newtheorem{proposition}[definition]{Proposici\'on}

\newtheorem{exercise}[definition]{Ejercicio}
\newtheorem{example}[definition]{Ejemplo}
\newtheorem{remark}[definition]{Observaci\'on}
\newtheorem{notation}[definition]{Notaci\'on}
\newtheorem{nada}[definition]{}
\newtheorem{conjecture}[definition]{Conjetura}

\newcommand{\dist}{\mbox{$\,$dist$\,$}}

 \makeindex

\begin{document}
\title{TEOR\'{I}A ERG\'{O}DICA DE LOS ATRACTORES TOPOL\'{O}GICOS Y ESTAD\'{I}STICOS}
\author{Eleonora Catsigeras
 \thanks{Instituto de Matem\'{a}tica y Estad\'{\i}stica \lq\lq Prof. Ing. Rafael Laguardia\rq\rq \  (IMERL),
 Facultad de Ingenier\'{\i}a,  Universidad de la Rep\'{u}blica,  Montevideo, Uruguay. \ \ Correo electr\'{o}nico: { \bf eleonora@fing.edu.uy } \ \ 
    La autora agradece la financiaci\'{o}n parcial de la Com. Sectorial de Investigaci\'{o}n Cient\'{\i}fica (CSIC) de la Universidad de la Rep\'{u}blica, de la Agencia Nacional de Investigaci\'{o}n e Innovaci\'{o}n (ANII) de Uruguay, de la organizaci\'{o}n de la XXVI Escuela Venezolana de Matem\'{a}tica y del Instituto Venezolano de Investigaci\'{o}n Cient\'{\i}fica (IVIC).}}
\date{8 de mayo de 2013}

\maketitle
  
  \pagestyle{myheadings}
\markboth{\textrm{Eleonora Catsigeras \ \ \ \ \ \ \ \ \ \ \ \ \ \ \ \ \ \ \ \ \ \ \ \ \ \ \ \ \ \ \ \ \ \ \ \ \ \ \ \ \ \ \ \ \ \ \ \ \ \ \ \ \   }}{ \textrm{\ \ \ \ \ \ \ \ \ \ \ \ \ \ \ \ \ \ \ \ \ \ \ Teor\'{\i}a  Erg\'{o}dica de los Atractores Topol\'{o}gicos y Estad\'{\i}sticos}}

  \abstract 
{Este texto corresponde a un curso corto introductorio   enfocado hacia las propiedades erg\'{o}dicas de los atractores topol\'{o}gicos y estad\'{\i}sticos. Est\'{a} dirigido a   estudiantes de postgrado, en particular a los participantes de la XXVI Escuela Venezolana de Matem\'{a}tica.
  Los primeros tres cap\'{\i}tulos incluyen una revisi\'{o}n  de los conceptos y propiedades b\'{a}sicas de la Teor\'{\i}a Erg\'{o}dica y de la Din\'{a}mica Topol\'{o}gica para los sistemas din\'{a}micos a tiempo discreto. No requiere del lector conocimimentos previos sobre Sistemas Din\'{a}micos.
  Los \'{u}ltimos tres cap\'{\i}tulos   desarrollan la teor\'{\i}a de atractores topol\'{o}gicos, existencia de medidas SRB, de Gibbs y f\'{\i}sicas para sistemas hiperb\'{o}licos suficientemente diferenciables, atractores erg\'{o}dicos, de Milnor y de Ilyashenko y su relaci\'{o}n con las medidas \lq\lq SRB-like\rq\rq \ o pseudo-f\'{\i}sicas, para sistemas din\'{a}micos continuos.}

 {  \tableofcontents}


\section{Din\'{a}mica medible y topol\'{o}gica}

Usaremos definiciones, notaci\'{o}n y algunos resultados b\'{a}sicos de teor\'{\i}a de la medida abs\-tracta  y de la topolog\'{\i}a en espacios m\'{e}tricos compactos, cuyos enunciados y demostraciones se pueden encontrar por ejemplo en \cite{Folland}, \cite{Rudin}  o \cite{Stein2005}.

\subsection{Din\'{a}mica de los automorfismos de medida.} \index{automorfismo! de esp. de medida}
Sea $(X, {\mathcal A})$ un espacio medible y sea $T: X \mapsto X$ una
transformaci\'{o}n medible, es decir $T^{-1}(A) \in {\mathcal A} \; \;
\forall A \in {\mathcal A}$.
\begin{definition} \em \index{sistema din\'{a}mico} {\bf Sistema din\'{a}mico discreto}

Se llama \em sistema din\'{a}mico discreto \em por iterados de $T$ hacia el futuro  a la
aplicaci\'{o}n que a cada $ n \in \mathbb{N}$ le hace corresponder la
transformaci\'{o}n $T^n: X \mapsto X$ donde $T^n := T \circ T \ldots
\circ T$ \ $  n $ veces, si $n \geq 1$ y $T^0 := id$. Es inmediato chequear que $T^n$ es medible para todo $n \in \mathbb{N}$.

Si adem\'{a}s $T $ es bi-medible (i.e. $T: X \mapsto X$ es medible,
invertible y su inversa $T^{-1}: X \mapsto X$ es medible), entonces \em el sistema din\'{a}mico discreto, \em  por iterados de $T$ hacia el futuro y el pasado,   se define como
 la aplicaci\'{o}n que a cada $n \in \mathbb{Z}$ le hace corresponder
$T^n$, donde $T^{-1}$ es la inversa de $T$ y $T^{n} := (T^{-1})^{|n|} := T^{-1}  \circ T^{-1} \circ \ldots T^{-1}$ \  $|n|$ veces si $n \leq -1$.

\end{definition}

Se observa que si  $T^{n+m} = T^n \circ T^m \ \forall \ n, m \in \mathbb{N} $, y adem\'{a}s si $T$ es bi-medible se cumple esa propiedad para todos $n, m \in \mathbb{Z}$. Esta propiedad
algebraica (cuando $T$ es bi-medible) se llama \em propiedad de grupo. \em Significa que el
sistema din\'{a}mico es una acci\'{o}n del grupo de los enteros   en el espacio $X$.

\begin{definition} \em \index{o'rbita}
Se llama \em \'{o}rbita positiva o futura $o^+(x)$ \em por el punto $x
\in X$, y $x$ se llama estado inicial de la \'{o}rbita, a la sucesi\'{o}n
$$o^+(x):= \{T^n(x)\}_{n \in N}.$$ Si $T$ es bi-medible se llama \'{o}rbita
negativa o pasada a la sucesi\'{o}n $$o^-(x):= \{T^{-n}(x)\}_{n \in N}$$ y se
llama \em \'{o}rbita $o(x)$ bilateral \em (o simplemente \'{o}rbita cuando $T$ es bi-medible) a la sucesi\'{o}n bi-infinita
$$o(x) := \{T^n(x)\}_{n \in \mathbb{Z}}.$$

Cuando $T$ es medible pero no bi-medible, llamamos simplemente \'{o}rbita a la \'{o}rbita positiva o futura.
\end{definition}

Se observa que $T^{n+m}(x) = T^n(T^m(x))$ y por lo tanto el
iterado $T^m(x)$ es el estado inicial   de la \'{o}rbita por $x$ que se
obtiene  corriendo el instante 0 al que antes era instante $m$.

\begin{definition} \em \index{punto! fijo} \index{punto! peri\'{o}dico}
 \em Punto fijo \em (o peri\'{o}dico de per\'{\i}odo 1) es $x_0$ tal que $T(x_0) = x_0$.
 \em Punto
 peri\'{o}dico \em es $x_0$ tal que existe $p \geq
 1$ tal que $T^p(x_0) = x_0$. El per\'{\i}odo es el m\'{\i}nimo $p \geq 1$
 que cumple lo anterior. Se observa, a partir de la propiedad
 de grupo, que si un punto es peri\'{o}dico
 de per\'{\i}odo $p$ entonces su \'{o}rbita est\'{a} formada por exactamente
 $p$ puntos.
\end{definition}

\begin{definition} \em \index{medida! invariante}
Sea $T$ medible en un espacio $(X, {\mathcal A})$. Una medida $\mu$ se
dice que es \em invariante por $T$, o que $T$ preserva $\mu$,   \em o se
dice tambi\'{e}n que $T$ es un \em automorfismo del espacio de medida $(X,
{\mathcal A}, \mu)$\em),    si
$$\mu (T^{-1}(A)) = \mu (A) \; \; \forall A \in {\mathcal A}.$$
\end{definition}

Se observa que pueden no existir medidas de probabilidad
invariantes para cierta $T: X \mapsto X $ transformaci\'{o}n medible
dada, como se muestra en el Ejemplo \ref{ejemplonoexistenmedidasinvariantes}. Sin embargo:

\begin{theorem} \label{medidasinvariantes} \index{medida! invariante!  teorema de existencia} \index{teorema! de existencia de! medidas invariantes} \label{teoremaExistenciaMedInvariantes} {\bf Existencia de medidas invariantes.}

 Sea $X$   un espacio m\'{e}trico compacto,   y sea $T: X \mapsto X$ continua. Entonces
existen
\em(usualmente infinitas) \em medidas de probabilidad en la sigma-\'{a}lgebra de Borel, que son
invariantes para $T$. \em

\end{theorem}

Demostraremos este teorema en la siguiente secci\'{o}n \ref{seccionpruebateoexistmedinvariantes}.

\begin{exercise}\em
Sean $(X, {\mathcal A}, \mu) $ e $  (Y, {\mathcal B}, \nu)$ u  espacios de medida y sea $T: X \mapsto Y$ medible. Se define $T^* \mu$ como la medida en $(Y, {\mathcal B})$ tal que $(T^* \mu)(B):= \mu(T^{-1}(B)) = \nu (B)$ para todo $B \in {\mathcal B}$.

(a)  Encontrar un ejemplo en que $T^*\mu = \nu$ pero   $\mu(A) \neq \nu (T(A))$ para alg\'{u}n $A \in {\mathcal A}$ tal que $T(A) \in {\mathcal B}$.

(b) Encontrar un ejemplo en que $T$ sea medible, cumpla $T^*\mu= \nu$ pero $T^{-1}$ no sea medible.

(c) Demostrar que si $T$ es medible, invertible y su inversa $T^{-1}$  es medible, entonces $T^* \mu = \nu$ si y solo si $(T^{-1})^* \nu = \mu$. Cuando se cumplen todas esas condiciones, se dice que $T$  (y por lo tanto tambi\'{e}n $T^{-1}$) es un isomorfismo de espacios de medida.
\end{exercise}

\begin{example} \em  \label{ejemplonoexistenmedidasinvariantes} \index{medida! invariante!  ejemplo de no existencia}
   Sea $T: [0,1] \mapsto [0,1]$ tal que $T(x) = x/2 $ si $x \neq
0$ y $T(0) = 1 \neq 0$. Afirmamos que:

 \em No existe
medida de probabilidad invariante por $T$.\em

\end{example}
{\em Demostraci\'{o}n: }  Si existiera, llam\'{e}mosla $\mu$. Consideremos
la partici\'{o}n de intervalo (0,1] dada por los subintervalos
disjuntos dos a dos $A_n= (1/2^{n+1}, 1/2^n]$ para $n \in \mathbb{N}$. Se
cumple $T^{-1}(A_{n}) = A_{n-1}$ para todo $n \geq 1$. Como $\mu$
es $T $ invariante se deduce $\mu (A_n) = \mu (A_0) \; \; \forall
n \geq 0$. Se tiene $\mu ((0,1] = \sum _{n \geq 0 } \mu (A_n) =
\sum _{n \geq 0 } \mu (A_0) \leq 1$. Luego $\mu (A_0) = 0$, de
donde $\mu ((0,1])= 0$ y $\mu (\{0\}) = 1$. Luego $\mu
(T^{-1}(\{0\})) = 1$ lo cual es absurdo porque $T^{-1}(\{0\}) =
\emptyset$.
\hfill $\Box$

\begin{exercise}\em \em
\em Probar que si $T:X \mapsto X$ es medible entonces $T$ preserva
la medida $\mu$ si y solo si para toda $f \in L^1 (\mu )$ se
cumple $$\int f \circ T \, d \mu = \int f \, d \mu $$
\end{exercise}
\begin{proposition}\label{proposicionSubalgebra}
Sea $T: X \mapsto X$ una transformaci\'{o}n medible en un espacio
medible $(X, {\mathcal A})$. Si $\mu$ es una medida finita o $sigma$-finita y
${\mathcal A}_0$ es un \'{a}lgebra que genera a ${\mathcal A}$ tal que $\mu
(T^{-1} (A)) = \mu (A) \; \forall A \in {\mathcal A}_0$ entonces $\mu
$ es invariante por $T$.
\end{proposition}
{\em Demostraci\'{o}n: } Sea $\nu (A) = \mu (T^{-1}(A))$ definida para
todo $A \in {\mathcal A }$. En la sub\'{a}lgebra ${\mathcal A}_0$ la premedida
$\nu $ coincide con la premedida $\mu$ (ambas restringidas a la
sub\'{a}lgebra son premedidas). Como existe una \'{u}nica extensi\'{o}n de una
premedida dada en un \'{a}lgebra a la sigma \'{a}lgebra generada, entonces
$\mu = \nu$ en $\mathcal A$. \hfill $\Box$

\begin{corollary} \label{corolarioMedidaInvarianteEnAlgebraGeneradora}
Si $T: \mathbb{R}^k \mapsto \mathbb{R}^k$ es una transformaci\'{o}n Borel me\-di\-ble y
$\mu$ es una medida $\sigma$-finita tal para todo conjunto $A$ que
sea uni\'{o}n finita  de rect\'{a}ngulos  de $R^k$ se cumple $\mu
(T^{-1}(A)) = \mu ( A )$ entonces $\mu$ es invariante por $T$ en
toda la sigma-\'{a}lgebra de Borel.
\end{corollary}

\subsection{Prueba de   existencia de medidas inva\-riantes} \label{seccionpruebateoexistmedinvariantes}

En esta secci\'{o}n $X$ denota un espacio m\'{e}trico compacto y $T: X \mapsto X$ una transformaci\'{o}n continua, a menos que se indique lo contrario.

Introduciomos algunas definiciones y resultados del An\'{a}lisis Funcional:

{\bf Notaci\'{o}n: El espacio de las funciones continuas y su dual.} \index{$C^0(X, \mathbb{R})$ espacio de! funciones reales continuas}

 Denotamos $C^0(X, \mathbb{R})$ el espacio de las funciones continuas $\psi: X \mapsto \mathbb{R}$ con la topolog\'{\i}a de la convergencia uniforme en $X$ (inducida por la norma del supremo). Es decir $$\mbox{dist} {(\psi_1, \psi_2)} :=   \|\psi_1 - \psi_2\|_0,$$
donde se denota
$$\|\psi\|_0 : = \max_{x \in X} | \psi (x)| \ \ \forall \ \psi \in C^0(X, \mathbb{R}).$$

Denotamos $C^0(X, [0,1])$ al subespacio de funciones continuas $\psi: X \mapsto [0,1]$ que toman valores no negativos ni mayores que 1, con la norma del supremo definida antes (o del m\'{a}ximo, en nuestro caso, pues $X$  es un espacio m\'{e}trico compacto).

El espacio $C^0(X, [0,1])$ es m\'{e}trico   acotado y cerrado. En efecto, el l\'{\i}mite de una sucesi\'{o}n uniformemente convergente de funciones continuas en $C^0(X, [0,1])$ es continua, y pertenece a $C^0(X, [0,1])$). Adem\'{a}s $C^0(X, [0,1])$ tiene una base numerable de abiertos, y por lo tanto existe un subconjunto numerable $\{\psi_i\}_{i \in \mathbb{N}}$ denso en $C^0(X, [0,1])$  (ver por ejemplo \cite[Proposition 4.40]{Folland})

Denotamos como $C^0(X, \mathbb{R})^*$ al dual de $C^0(X, \mathbb{R})$; es decir al conjunto (al que luego dotaremos de una topolog\'{\i}a adecuada) de todos los operadores lineales $$\Lambda: C^0(X, \mathbb{R})  \mapsto \mathbb{R}.$$

\begin{definition} {\bf El espacio ${\mathcal M}$ de las probabilidades y la to\-po\-lo\-g\'{\i}a d\'{e}bil$^*$} \em   \index{${\mathcal M}$ (espacio de medidas de pro\-ba\-bi\-lidad)} \index{espacio! de medidas de pro\-ba\-bi\-lidad} \index{topolog\'{\i}a d\'{e}bil estrella} \label{definicionTopologiaDebil*}
Sea $X$ un espacio m\'{e}trico compacto. Sea ${\mathcal M}$ el \em espacio de todas las medidas de probabilidad. \em

En ${\mathcal M}$ introducimos la siguiente topolog\'{\i}a, llamada \em topolog\'{\i}a d\'{e}bil$^*$: \em

Si $\mu_n, \mu \in {\mathcal M}$ decimos que la sucesi\'{o}n $\{\mu_n\}$ es convergente a $\mu$ en la topolog\'{\i}a d\'{e}bil$^*$,  es decir
$$\lim_{n \rightarrow + \infty} \mu_n = \mu \ \ \mbox{ en } {\mathcal M},$$
cuando:
$$\lim_{n \rightarrow + \infty} \int \psi \, d \mu_n = \int \psi \, d \mu \ \ \mbox{ en } {\mathbb{R}} \ \  \forall \ \psi \in C^0(X, \mathbb{R}).$$
\end{definition}

\begin{remark} \em    Por definici\'{o}n de convergencia  en la topolog\'{\i}a d\'{e}bil$^*$ de una sucesi\'{o}n de medidas de probabilidad, la sucesi\'{o}n converge a la medida
    $\mu$ si y solo si
  para cada funci\'{o}n continua, la sucesi\'{o}n de integrales converge
a la integral respecto a $\mu $.   Es falso que para todo $A$ boreliano
la sucesi\'{o}n de medidas de $A$ converja a $\mu (A)$.  En efecto, v\'{e}ase el ejercicio siguiente:

\end{remark}

\begin{exercise}\em
Sea $X = [0,1] $. Para cada $n \geq 0 $ sea $\mu _n $ la medida
delta de Dirac concentrada en el punto $ 1/2^n $.

a)
Probar que existe $\mu = \mbox{l\'{\i}m} _{n \rightarrow +\infty} \mu _n $ en la topolog\'{\i}a d\'{e}bil$^*$ y encontrar la medida l\'{\i}mite $\mu$.

b)
Encontrar $A \subset [0,1] $ boreliano tal que no existe
$\mbox{l\'{\i}m} _{n \rightarrow +\infty} \mu _n (A) $.

Sugerencia: $ A = \{ 1/2^{2j}: j \geq 0 \} $.

 c)
Encontrar $B \subset [0,1] $ boreliano tal que  existe
$\mbox{l\'{\i}m} _{n \rightarrow +\infty} (\mu _n (B)) \neq \mu (B)$.

\end{exercise}

\begin{remark} \em  Observemos que, para cada $\mu \in {\mathcal M}$, el operador $\Lambda_{\mu}: C^0(X, \mathbb{R}) \mapsto \mathbb{R}$ definido por:
$$\Lambda_{\mu} (\psi) := \int \psi \, d \mu $$
es lineal, positivo (es decir $\Lambda_{\mu}(\psi) \geq 0 $ si $\psi \geq 0$) y acotado   (es decir existe $k >0$ tal que $|\Lambda_{\mu} (\psi)| \leq k $ para toda $\psi\in C^0(X, \mathbb{R})$ tal que $\|\psi\|_0\leq 1$). Adem\'{a}s $\Lambda_{\mu}(\psi)=1$ si $\psi: X \mapsto \mathbb{R}$ es la funci\'{o}n constante igual a 1.

El siguiente teorema establece el rec\'{\i}proco de la propiedad observada arriba, y es un resultado cl\'{a}sico de la Teor\'{\i}a Abstracta de la Medida y del An\'{a}lisis Funcional:
\end{remark}

{\bf Teorema de Representaci\'{o}n de Riesz} \index{teorema! Representaci\'{o}n de Riesz}

Sea $X$ un espacio m\'{e}trico compacto.

\em Para todo operador lineal $\Lambda: C^0(X, \mathbb{R}) \mapsto \mathbb{R}$ que sea positivo y acotado, existe y es \'{u}nica una medida finita $\mu $ \em (de Borel y positiva) \em tal que \em
$$\Lambda (\psi) = \int \psi \, d \mu \ \ \forall \ \psi \in C^0(X, \mathbb{R})$$

\em Adem\'{a}s si $\Lambda (1) = 1$ entonces  $\mu $ es una probabilidad. Es decir, $\mu(X)= 1$, \'{o}, usando nuestra notaci\'{o}n, $\mu \in {\mathcal M}$. \em

\vspace{.3cm}

Una demostraci\'{o}n del Teorema de Representaci\'{o}n de Riesz puede encontrarse, por ejemplo, en \cite[Teorema 2.3.1]{Rudin}

Debido al Teorema de Representaci\'{o}n de Riesz, el espacio ${\mathcal M}$ se puede identificar con el espacio de los operadores  lineales de $C^0(X, \mathbb{R})$ que son positivos, acotados y que valen 1 para la funci\'{o}n constante igual a 1. (Recordemos que el espacio de los operadores lineales de $C^0(X, \mathbb{R})$ se llama dual de $C^0(X, \mathbb{R})$ y se denota como $C^0(X, \mathbb{R})^*$). En el An\'{a}lisis Funcional se definen diversas topolog\'{\i}as en el dual de un espacio de funciones. Una de ellas es la llamada \em topolog\'{\i}a d\'{e}bil$^*$, \em que es la topolog\'{\i}a de la convergencia \em  punto a punto, \em definida como sigue:
$$\lim_{n \rightarrow + \infty} \Lambda_n = \Lambda \mbox{ en } C^0(X, \mathbb{R})^*$$
si y solo si
$$ \lim_{n \rightarrow + \infty} \Lambda_n \psi = \Lambda \psi \mbox{ en } \mathbb{R} \ \ \forall \ \psi \in C^0(X, \mathbb{R}). $$

De las definiciones anteriores, deducimos que la topolog\'{\i}a d\'{e}bil estrella en ${\mathcal M}$ es la topolog\'{\i}a inducida por la topolog\'{\i}a d\'{e}bil estrella (o de la convergencia punto a punto) en el dual $C^0(X, \mathbb{R})^*$ del espacio funcional de las funciones continuas. La topolog\'{\i}a d\'{e}bil$^*$ puede definirse como la topolog\'{\i}a producto de las definidas por la convergencia de    los valores num\'{e}ricos $\Lambda_n(\psi)$ que toman los operadores $\Lambda_n$ para cada $\psi \in C^0(X, \mathbb{R})$ fija.

El siguiente teorema es cl\'{a}sico del An\'{a}lisis Funcional, y es una consecuencia del teorema de Tichonov (ver por ejemplo \cite[Teorema 4.43]{Folland})que establece, bajo ciertas hip\'{o}tesis, la compacidad de la topolog\'{\i}a producto:

\begin{theorem}\label{theoremTichonov}   \index{teorema!Tichonov}
{\bf (Corolario del Teorema de Tichonov)}   \label{teoremaCompacidadEspacioProbabilidades0}
Si $X$ es un espacio m\'{e}trico compacto, entonces para toda constante $k >0$ el subconjunto de    los operadores lineales acotados por $k$ es compacto en el espacio dual   $C^0(X, \mathbb{R})^*$ con la topolog\'{\i}a d\'{e}bil estrella.

\end{theorem}

 Como caso particular, observemos que el espacio ${\mathcal M}$ de las probabilidades de Borel en ${ X}$ (via el teorema de Representaci\'{o}n de Riesz) es un subconjunto \em cerrado \em  del espacio de los operadores lineales acotados   con la topolog\'{\i}a d\'{e}bil estrella (es decir, si $\Lambda_n (1)= 1$ y $\Lambda_n \rightarrow \Lambda$, entonces $\Lambda (1) = 1$).

 M\'{a}s detalladamente:

\begin{theorem} {\bf Compacidad y metrizabilidad del espacio de probabilidades}
\label{teoremaCompacidadEspacioProbabilidades} \index{espacio! de medidas de pro\-ba\-bi\-lidad} \index{teorema! de compacidad de ${\mathcal M}$} \index{teorema! de metrizabilidad de ${\mathcal M}$}
\index{m\'{e}trica en ${\mathcal M}$}

Sea $X$ un espacio m\'{e}trico compacto. Sea ${\mathcal M}^{1} $ el espacio de todas las medidas $\mu$ \em (de Borel, positivas y finitas)\em  tales que $\mu(X) \leq 1$. Sea en ${\mathcal M}^{1} $ la topolog\'{\i}a d\'{e}bil$^*$, definida   en \em \ref{definicionTopologiaDebil*}. \em
Entonces:

 {\bf (a)} ${\mathcal M}^{1} $ es compacto.

  {\bf (b)} ${\mathcal M}^{1}$ es metrizable. \em (Es decir, existe una m\'{e}trica $\mbox{dist}$ que induce la topolog\'{\i}a d\'{e}bil$^*$; esto es una distancia en ${\mathcal M}^{1}\times {\mathcal M}^{1}$ tal que $$\lim_n \mbox{dist}(\mu_n, \mu) = 0 \mbox{ si y solo si } \lim_n\mu_n = \mu$$ con la topolog\'{\i}a d\'{e}bil$^*$). \em

  {\bf (c)} Si $\{\psi_i\}_{i \in \mathbb{N}} \subset C^0(X, [0,1])$ es un conjunto numerable denso de funciones, entonces la siguiente m\'{e}trica induce la topolog\'{\i}a d\'{e}bil$^*$ en ${\mathcal M}^1$:

  \begin{equation}\label{equationDistanciaMedidas}  \mbox{dist}(\mu, \nu) := \sum_{i=1}^{+ \infty} \frac{1}{2^i} \ \left| \int \psi_i \, d \mu \ - \ \int \psi_i \, d \nu \right| \ \ \forall \ \mu, \nu \in {\mathcal M}^1. \end{equation}

{\bf (d)}   ${\mathcal M}^{1}$ es secuencialmente compacto. \em (Es decir, toda sucesi\'{o}n de medidas en ${\mathcal M}^{1}$ tiene alguna subsucesi\'{o}n convergente).

\vspace{.3cm}

Consecuencia: Siendo ${\mathcal M} = \{ \mu \in {\mathcal M}^1: \mu(X)= 1\}$ cerrado en ${\mathcal M}^1$, se deduce que

{\bf (e)} \em El espacio de probabilidades ${\mathcal M}$ con la topolog\'{\i}a d\'{e}bil$^*$ es compacto, metrizable y secuencialmente compacto. \em
\end{theorem}

 \begin{exercise}\em
 Demostrar el Teorema \ref{teoremaCompacidadEspacioProbabilidades} como consecuencia del Corolario \ref{teoremaCompacidadEspacioProbabilidades0}, identificando el espacio ${\mathcal M}^1$ con el dual de $C^0(X, \mathbb{R})$ via el Teorema de Representaci\'{o}n de Riesz.
 \end{exercise}

\begin{proposition}\label{coincidencia}
Sea $\{ \psi_i: i \geq 1 \}$ un conjunto numerable denso en
   $C^0(X, [0,1])$.

 Dos medidas $\mu _1 $ y $\mu _2 $ en ${\mathcal M}^1(X)$ coinciden
si para todo $i \geq 1 $ se cumple
$$\int \psi_i \, d \mu _1 = \int \psi_i \, d \mu _2 $$
\end{proposition}

{\em Demostraci\'{o}n: }
Por la unicidad de la medida del teorema de Riesz alcanza probar que
para toda $\psi \in C^0(X, \mathbb{R})$ se cumple:
\begin{eqnarray} \label{denso}
\int \psi \, d \mu _1 = \int \psi \, d \mu _2
\end{eqnarray}
La igualdad (\ref{denso}) vale obviamente para $\psi$ id\'{e}nticamente nula.
Si $\psi$ no es id\'{e}nticamente nula, basta demostrar (\ref{denso}) para
$\psi / \|\psi\|_0 $, donde $\|\psi\|_0 = \max_{x \in X} |\psi(x)|$. Entonces supongamos que $\|\psi \|_0 = 1$. Cualquier funci\'{o}n real $\psi$ puede escribirse como $\psi = \psi^+ - \psi^-$, donde $\psi^+ = \max\{\psi, 0\}, \ \psi^- = - \min\{\psi, 0\}$. Observemos que $\psi^+, \psi^- \in C^0([0,1])$. Si demostramos la igualdad (\ref{denso}) para  $\psi^+ $ y $\psi^-$, entonces vale tambi\'{e}n para $\psi$. Basta entonces probar la igualdad (\ref{denso}) para toda $\psi \in C^0(X, [0,1])$. Por la densidad
de las funciones $\{\psi_i\}$ en $C^0(X, [0,1]) $, existe una sucesi\'{o}n $\psi_{i_n}$
convergente uniformemente (es decir con la norma del supremo en $C^0(X, \mathbb{R})$),
a la funci\'{o}n $\psi$. Por lo tanto, converge tambi\'{e}n puntualmente
y est\'{a} uniformemente acotada por 1. Cada $\psi_{i_n}$ verifica la
igualdad (\ref{denso}). Luego, por el teorema de convergencia dominada, $\psi$ tambi\'{e}n
cumple (\ref{denso}).
\hfill $\Box$

\begin{remark} \em

Sea $X$ un espacio m\'{e}trico compacto no vac\'{\i}o, sea ${\mathcal B}$ la sigma-\'{a}lgebra de Borel, y sea ${\mathcal M}$ el espacio de todas las medidas de probabilidad de Borel en $X$ con la topolog\'{\i}a d\'{e}bil$^*$.

El espacio ${\mathcal M}$ es no vac\'{\i}o.  Por ejemplo, si elegimos un punto $x \in X$ entonces $  \delta_x \in {\mathcal M}$, donde $\delta_x$ es la probabilidad Delta de Dirac soportada en el punto $x \in X$. Esto es $\delta_x$  es la probabilidad que a cada boreliano $B \subset X$ le asigna $\delta_x(B) = 1$ si $x \in B$, y $\delta_x (B)= 0$ si $ x \not \in B$.

\end{remark}

Ahora agreguemos una din\'{a}mica  en $X$:

\begin{definition} {\bf El pull back $T^*: {\mathcal M} \mapsto {\mathcal M}$} \index{$T^*$ pull back} \index{pull back}

Sea $(X, {\mathcal A})$  un espacio medible,   sea $T: X \mapsto X$ medible y sea ${\mathcal M}$ el conjunto de las medidas de probabilidad en $(X, \mathbb{A})$. Definimos el siguiente operador en ${\mathcal M}$, llamado pull back del mapa $T$:

$$T^*:{\mathcal M} \mapsto {\mathcal M} \ \ \ \ \ \ \ \ (T^* \mu)(B) = \mu (T^{-1}(B)) \ \ \forall \ \mu \in {\mathcal M} \ \ \forall\  B \in {\mathcal B}.$$

\end{definition}

Es inmediato verificar que $\mu$ es invariante por $T$  si y solo si $T^* \mu = \mu$, es decir, las medidas invariantes por $T$ son los puntos fijos por $T^*$ en el espacio ${\mathcal M}$.

\begin{exercise}\em \label{ejercicioInvarianciaDeMedidasIntegral}

Probar que para toda $\psi \in L^1(\mu)$   se cumple:
$$\int \psi \, d T^* \mu = \int \psi \circ T \, d \mu.$$

Sugerencia: Chequear primero para las funciones caracter\'{\i}sticas $\chi_B$ de los borelianos, luego para las combinaciones lineales de las funciones caracter\'{\i}sticas (funciones simples), y luego para las funciones medibles no negativas, usando el Teorema de convergencia mon\'{o}tona. Finalmente probar la igualdad para toda $\psi \in L^1(\mu)$ separando $\psi$ en parte real e imaginaria, y cada $\psi$ real en su parte positiva y negativa.
\end{exercise}

\begin{proposition} \label{proposicionContinuidadT*} Sea $X$ espacio m\'{e}trico compacto, y ${\mathcal M}$ el espacio de las proba\-bi\-li\-dades de Borel en $X$ con la topolog\'{\i}a d\'{e}bil$^*$. \index{continuidad! del operador pull back}
Si $T: X \mapsto X$ es continuo, entonces $T^*: {\mathcal M} \mapsto {\mathcal M}$ es continuo en ${\mathcal M}$.
\end{proposition}
{\em Demostraci\'{o}n: }
Basta chequear que si $\lim_n \mu_n = \mu $ en ${\mathcal M}$ entonces $\lim_n T^* \mu_n = T^* \mu$. En efecto, como $\lim_n \mu_n = \mu$, entonces, para toda $\psi \in C^0(X, \mathbb{R}) $:
$$\lim_{n \rightarrow + \infty}  \int \psi \, d \mu_n = \int \psi \, d \mu$$
En particular, siendo $T: X \mapsto X$ continuo, la igualdad anterior se cumple para $\psi \circ T$. Luego deducimos que
$$\lim_{n \rightarrow + \infty}  \int \psi \circ T \, d \mu_n = \int \psi \circ T \, d \mu.$$
Usando el resultado del Ejercicio \ref{ejercicioInvarianciaDeMedidasIntegral}, deducimos que
$$\lim_{n \rightarrow + \infty} \int \psi \, d T^* \mu_n  = \int \psi \, d T^* \mu \ \ \forall \ \psi \in C^0(X, \mathbb{R}).$$
Por lo tanto la sucesi\'{o}n de medidas $T^* \mu_n$ converge en la topolog\'{\i}a d\'{e}bil$^*$ de ${\mathcal M}$, a la medida $T^* \mu$ como quer\'{\i}amos demostrar.
\hfill $\Box$

Ahora probaremos el Teorema \ref{teoremaExistenciaMedInvariantes}, utilizando el llamado  \em
procedimiento de Bogliubov-Krylov \em \cite{Bogliubov-Krylov}. \index{Bogliubov-Krylov! pro\-cedi\-miento de} Este procedimiento parte de cualquier medida de probabilidad en el espacio $X$, toma  promedios aritm\'{e}ticos de los iterados del operador pull back $T^*$ de esta medida hasta tiempo $n$, y finalmente una subsucesi\'{o}n convergente en la topolog\'{\i}a d\'{e}bil estrella de esos promedios. Se obtienen   medidas invariantes por $T$ (bajo la hip\'{o}tesis de que $T: X \mapsto X$ es continuo). El procedimiento de Bogliubov-Krylov permite \lq\lq fabricar\rq\rq \ medidas de probabilidad invariantes, usando como \lq\lq semilla\rq\rq cualquier medida de probabilidad.

\vspace{.5cm}

\newpage
{\bf Demostraci\'{o}n del Teorema \ref{teoremaExistenciaMedInvariantes} de existencia de medidas invariantes:}
{\em Demostraci\'{o}n: }

Elijamos una medida de probabilidad de Borel cualquiera   $\rho \in {\mathcal M}$. Construyamos para cada $1 \leq n \in \mathbb{N}$, la siguiente probabilidad:
$$\mu_n = \frac{1}{n} \sum_{j= 0}^{n-1} (T^*)^j \rho$$

Es inmediato probar, a partir de la definici\'{o}n del operador $T^*$ (que cumple $T^* \mu (B) = \mu(T^{-1}(B))$ para todo boreliano $B$), que $$T^* \mu_n = \frac{1}{n} \sum_{j= 0}^{n-1} (T^*)^{j+1} \rho$$

Como el espacio ${\mathcal M}$ de las probabilidades de Borel es secuencialmente compacto con la topolog\'{\i}a d\'{e}bil$^*$, existe una subsucesi\'{o}n $\{\mu_{n_i}\}_{i \in \mathbb{N}}$ ( con $ \lim_i n_i = + \infty$), que es convergente en ${\mathcal M}$. Llamemos $\mu \in {\mathcal M}$ a su l\'{\i}mite; es decir:
$$\mu = \lim_{i \rightarrow + \infty} \frac{1}{n_i} \sum_{j= 0}^{n_i-1} (T^*)^j \rho$$
Basta demostrar ahora que $T^* \mu = \mu$, es decir $\mu$ es una probabilidad $T$-invariante.

Usando la continuidad del operador $T^*$ (Proposici\'{o}n \ref{proposicionContinuidadT*}), deducimos que
$$T^* \mu = \lim_{i \rightarrow + \infty} T^* \left (\frac{1}{n_i} \sum_{j= 0}^{n_i-1} (T^*)^j \rho \right ) = \lim_{i \rightarrow + \infty} \frac{1}{n_i} \sum_{j= 0}^{n_i-1} (T^*)^{j+1} \rho$$

Integrando cada $\psi \in C^0(X, \mathbb{R})$ res\-pec\-to de la medida $T^* \mu$, y luego res\-pec\-to de la medida $\mu$,   aplicando la definici\'{o}n de l\'{\i}mite de la sucesi\'{o}n de medidas $\mu_{n_i}$ con la topolog\'{\i}a d\'{e}bil$^*$, y la continuidad del operador $T^*$ en ${\cal M}$, obtenemos:
$$\int \psi \, dT^* \mu   - \int \psi \, d\mu   = \lim_{i \rightarrow + \infty}  \left ( \int \psi \, d \, T^* \mu_{ n_i} - \int \psi \, d \mu_{ n_i}\right )= $$

$$ \lim_{i \rightarrow + \infty} \frac{1}{n_i} \left ( \int \psi \, d \, (T^*)^{n_i}\rho  - \int \psi \, d \rho    \right )  $$
 Entonces $$\left |\int \psi \, d\, T^* \mu   - \int \psi \, d \mu \right | \leq \lim _{i \rightarrow + \infty} \frac{1}{n_i} \|\psi\|_0 \left (\rho(X)  + \rho(T^{-n_i}(X))  \right ) =$$ $$ \lim_{n \rightarrow + \infty} \frac{2}{n} \|\psi\|_0 = 0.$$
En la \'{u}ltima igualdad tuvimos  en cuenta que $\rho$ es una probabilidad. Obtuvimos que
$$\left |\int \psi \, d \, T^* \mu   - \int \psi \, d \mu \right | = 0 \ \ \forall \ \psi \in C^0(X, \mathbb{R}).$$ Por lo tanto los operadores lineales $\psi \mapsto \int \psi \, d T^* \mu$ y $\psi \mapsto \int \psi \, d \mu$ son el mismo. Por la unicidad de la medida en el teorema de Riesz deducimos que $T^* \mu = \mu$, como quer\'{\i}amos demostrar.
\hfill $\Box$

\begin{exercise}\em
Sea $(X, {\mathcal A})$ un espacio medible y sea ${\mathcal M}$ el conjunto de todas las medidas de probabilidad en $(X, {\mathcal A})$.
Suponga que existe $\mu \in {\mathcal M}(X)$ tal que $T^*\mu (A)
\leq 2 \mu (A) $ para todo boreliano $A\subset X $.

a)
Probar que $ 2 \mu - T^* \mu \in {\mathcal M}(X) $.

b)
Si $X$ es un espacio m\'{e}trico compacto, ${\mathcal A}$ es la sigma-\'{a}lgebra de Borel y si $T$ es continua, probar que dado $\mu _0 \in {\mathcal M}(X)$ existe
$\mu \in {\mathcal M}(X)$ tal que $2 \mu  - T^* \mu = \mu _0 $. Sugerencia:
Para todo $\mu \in {\mathcal M}(X) $ definir $G(\mu ) = 1/2 \cdot (T^* \mu +
\mu _0 ), \; \; \; \mu _n = 1/n \cdot \sum _{j=0} ^{n-1} G^j(\mu _0 ) $
y tomar una subsucesi\'{o}n convergente en ${\mathcal M}(X) $. Probar que $G$ es continuo en ${\mathcal M}(X)$. Observar que $G$ conmuta con las combinaciones lineales finitas \em convexas \em de probabilidades. Es decir si   $\mu = \sum_{i= 1}^k \lambda_i \nu_i$, donde $\nu_i \in {\mathcal M}(X)$, $ \ \  0 \leq \lambda_i \leq 1$ y $\sum_{i= 1}^k \lambda_i = 1$, entonces $G(\nu) = \sum_{i= 1}^k \lambda_ i G(\nu_i) $.

\end{exercise}

\subsection{Ejemplos y puntos peri\'{o}dicos hiperb\'{o}licos}
\index{punto! peri\'{o}dico! hiperb\'{o}lico}
\begin{example} \em  \index{rotaci\'{o}n! irracional} \index{rotaci\'{o}n! racional}
Sea en $S^1$ (el c\'{\i}rculo) la rotaci\'{o}n $T(e^{2 \pi i x}) = e^{(2 \pi i) (x + a)}
$, donde $a$ es una constante real. Si $a$ es racional,  $T$ se
llama \em rotaci\'{o}n racional \em  en el c\'{\i}rculo, y si $a$ es
irracional,  $T$ se llama \em rotaci\'{o}n irracional. \em   A trav\'{e}s de la
proyecci\'{o}n $\Pi: \mathbb{R} \mapsto S^1$ dada por $\Pi (x) = e ^{2 \pi i
x}$, la medida de Lebesgue $m$ en $\mathbb{R}$ induce una medida
$m_{\sim}$ en $S^1$ dada por $m_{\sim }(A) = m (\Pi ^{-1} (A) \cap
[0,1])$. Esta medida $m_{\sim }$ se llama medida de Lebesgue en el
c\'{\i}rculo. Como $m$ en $\mathbb{R}$ es invariante por traslaciones, es f\'{a}cil
probar que $m _{\sim} $ en el c\'{\i}rculo $S^1$ es invariante por las rotaciones.
\end{example}

\begin{exercise}\em \em
\em Probar que para la rotaci\'{o}n racional en el c\'{\i}rculo todos los
puntos son peri\'{o}dicos con el mismo per\'{\i}odo. Probar que la rotaci\'{o}n
irracional en el c\'{\i}rculo no tiene puntos peri\'{o}dicos. Probar que la
medida de Lebesgue en el c\'{\i}rculo es invariante por las rotaciones.
\end{exercise}

\begin{remark} \em
Aunque no es inmediato, se puede probar que  \em
todas las \'{o}rbitas
de la rotaci\'{o}n irracional del c\'{\i}rculo son densas. \em
(Lo probaremos en  \S\ref{remarkRotacionIrracionalOrbitasDensas}).
\end{remark}

\begin{exercise}\em \em {\bf Tent map} \index{tent map}
\em Sea el intervalo [0,1] dotado de la sigma \'{a}lgebra de Borel.
Sea $ T:[0,1] \mapsto [0,1]$ dada por $T(x) = 2 x$ si $x \in
[0,1/2]$ y $T(x) = 2-2x$ si $x \in [1/2,1]$. Probar que $T$ preserva
la medida de Lebesgue en el intervalo. (Sugerencia: graficar $T$ y
probar que la preimagen de un intervalo $I$ tiene la misma medida
que $I$. Usar el corolario \ref{corolarioMedidaInvarianteEnAlgebraGeneradora}.
\end{exercise}

\begin{definition} \em \label{definicionatractoryrepulsor} \index{atractor! peri\'{o}dico} \index{repulsor peri\'{o}dico}
Sea $T:X \mapsto X$ Borel medible en un espacio m\'{e}trico $X$.
Decimos que un punto peri\'{o}dico $x_0$ de per\'{\i}odo $p$ es \em un
atractor \em si existe un entorno $V$ de $x_0$ invariante hacia
delante por $T^p $ (es decir $T^p(V) \subset V$) y tal que
$$\dist (T^n(x_0), T^n(y))_{n \rightarrow + \infty} \rightarrow 0
\; \; \forall y \in V $$ Cuando $T$ es invertible con inversa
medible decimos que un punto peri\'{o}dico $x_0$ de per\'{\i}odo $p$ es \em
un repulsor \em si existe un entorno $V$ de $x_0$ invariante hacia
atr\'{a}s por $T^p $ (es decir $T^{-p}(V) \subset V$) y tal que
$$\dist (T^n(x_0), T^n(y))_{n \rightarrow - \infty} \rightarrow 0
\; \; \forall y \in V $$

\end{definition}

\begin{proposition} \label{claimatractoreshiperbolicos} \index{punto! peri\'{o}dico! atractor}
 Sea $f: S^1 \mapsto S^1$ un difeomorfismo; \em es decir $f$ es de clase $C^1$ (i.e. derivable con derivada continua), invertible (i.e. biyectiva; existe la transformaci\'{o}n inversa $f^{-1}: S^1 \mapsto S^1$), y su inversa $f^{-1}$ es tambi\'{e}n de clase $C^1$) \em

 Supongamos que el difeomorfismo $f: S^1 \mapsto S^1$ preserva   la orientaci\'{o}n (i.e. $f'>0$). Sea
$x_0$   un punto peri\'{o}dico de per\'{\i}odo $p$ tal que la derivada $
(f^p) '(x_0)$ es menor que 1. Entonces $x_0$ es un atractor. An\'{a}logamente, si
$(f^p) ' (x_0)$ es mayor que 1 entonces $x_0$ es un repulsor. \em
\end{proposition}
{\em Demostraci\'{o}n: }
La segunda afirmaci\'{o}n se obtiene de la primera usando $f^{-p}$ en
lugar de $f^p$. Demostremos la primera afirmaci\'{o}n renombrando
como $f$ a la transformaci\'{o}n $f^p$. Entonces $x_0$ es fijo. Graf\'{\i}quese $f(x)$ para $x \in S^1 \approx [0,1]$
del intervalo $[0,1]$ en s\'{\i} mismo, en el que se ha identificado el
0 con el 1 en el punto $x_0$. La gr\'{a}fica de $f $ corta
a la diagonal por lo menos en el punto $0 \sim 1= x_0$
Gr\'{a}ficamente, los iterados futuros de  $y$ en un entorno de $x_0 $
suficientemente peque\~{n}o, se obtienen trazando la vertical de
abscisa $ y$, cort\'{a}ndola con la gr\'{a}fica de $f $, trazando luego la
horizontal por ese punto, cort\'{a}ndola con la diagonal, trazando la
vertical por ese punto, cort\'{a}ndola con la gr\'{a}fica de $f$, y as\'{\i}
sucesivamente (ver por ejemplo \cite[Figure on page 19]{Jost}). Si la funci\'{o}n es continua con
derivada continua, y
la derivada $f'(x_0) = a >0 $ en el punto fijo $x_0$ es menor que uno, entonces los
sucesivos puntos en la gr\'{a}fica de $f$ obtenidos por la
construcci\'{o}n anterior, tienden mon\'{o}tonamente al punto fijo $x_0$. En efecto, por la definici\'{o}n de diferenciabilidad y de derivada: $\|f(y) - f(x_0)\| \leq (a + (1-a)/2) \,  \|y - x_0\|$ para todo $y$ suficientemente cercano a $x_0$, digamos $\|y - x_0\| < \delta$. Es decir $\|f(y) - x_0\| \leq b \|y - x_0\|$ donde $0 < b = a + (1 - a)/2 = (1 + a)/2 < 1$. Luego,
$f(y)$ tambi\'{e}n cumple $\|f(y) -x_0\| < \delta$. Se puede aplicar, por inducci\'{o}n, la desigualdad anterior a todos los iterados futuros $f^n(y)$ (es decir para todo $n \in \mathbb{N}$). Obtenemos $\|f^n(y) - x_0\| \leq b^n \|y - x_0\|$. Siendo $0 < b < 1$, deducimos que $f^n(y)$ converge mon\'{o}tonamente a $x_0$.
\hfill $\Box$

\begin{definition} \label{definicionhiperbolico} \em \index{punto! peri\'{o}dico! hiperb\'{o}lico}
 Un punto peri\'{o}dico  $x_0$ con per\'{\i}odo $p$ de un difeormofismo $f: X \mapsto X$ en una variedad diferenciable $X$ se dice
 \em hiperb\'{o}lico \em  si los valores propios (complejos) de la derivada $df^p_{x_0}$ de $f^p$
 en $x_0$, tienen todos m\'{o}dulo diferente de 1. Se recuerda que la derivada $df^p_{x_0}$ es  una transformaci\'{o}n lineal de $\mathbb{R}^m$ en $ \mathbb{R}^m$, donde $m $ es la dimensi\'{o}n de la variedad $X$.
\end{definition}

{\bf Consecuencia: } Si $x_0$ es un punto peri\'{o}dico hiperb\'{o}lico de
un difeomorfismo $f$ de clase $C^1$ del c\'{\i}rculo $S^1$,  entonces es un atractor si  $|(f^p)'(x_0)| < 1$,  y  es
un repulsor si $|(f^p)'(x_0) |
> 1$.  (Siendo $x_0$ hiperb\'{o}lico, sabemos que $|f^p)'(x_0)| \neq 1$ por definici\'{o}n, as\'{\i} que los dos casos anteriores son los \'{u}nicos posibles). En efecto, si $f$ preserva al orientaci\'{o}n del c\'{\i}rculo, aplicamos la Proposici\'{o}n \ref{claimatractoreshiperbolicos}, y si $f$ invierte la orientaci\'{o}n, aplicamos la misma Proposici\'{o}n a $f^2$ para deducir que $x_0$ es un punto fijo atractor (resp. repulsor) de $f^2$. Es f\'{a}cil probar que si $x_0$ es un punto fijo de $f$ que es atractor (resp. repulsor) para $f^2$,  entonces tambi\'{e}n   es   atractor (resp. repulsor) para $f$.

Para un difeomorfismo $f$ en el c\'{\i}rculo $S^1$, y m\'{a}s en general para un mapa de clase $C^1$ en una variedad de dimensi\'{o}n 1, un punto peri\'{o}dico hiperb\'{o}lico $x_0$ se llama \em pozo \em si $|(f^p)'(x_0)| < 1$, y se llama \em fuente \em si $|(f^p)'(x_0)| >1$. Generalizando este resultado cuando la variedad tiene dimensi\'{o}n mayor que uno, adoptamos la siguiente definici\'{o}n:

\begin{definition} \em {\bf Pozos, fuentes y sillas} \index{pozo} \index{fuente} \index{silla}

Sea $f: X  \mapsto X$ un difeomorfismo en una variedad diferenciable $X$. Sea $x_0$ un punto peri\'{o}dico hiperb\'{o}lico para $f$ de per\'{\i}odo $p$ (i.e. los valores propios de $df^p_{x_0}$ tienen todos m\'{o}dulo diferente de 1). El punto $x_0$, y tambi\'{e}n la \'{o}rbita (finita) de $x_0$, se llama \em pozo \em si los valores propios de $df^p_{x_0}$ tienen todos m\'{o}dulo menor que 1. Se llama \em fuente \em si todos tienen m\'{o}dulo mayor que 1. Y se llama \em silla \em si alguno tiene m\'{o}dulo mayor que 1 y alg\'{u}n otro m\'{o}dulo menor que 1. (Obs\'{e}rvese que las sillas solo pueden existir si la dimensi\'{o}n de la variedad es 2 o mayor).
\end{definition}

\begin{exercise}\em
(a) Encontrar   ejemplo  de un difeomorfismo $f: \mathbb{R}^2 \mapsto \mathbb{R}^2$ con un punto fijo hiperb\'{o}lico tipo silla, otro ejemplo con un pozo y otro con una fuente. (b) Encontrar un ejemplo de $f: \mathbb{R}^2 \mapsto {\mathbb{R}^2}$ que tenga exactamente tres puntos fijos, sean los tres hiperb\'{o}licos, uno tipo fuente, otro pozo y otro silla. (c) Demostrar que para cualquier difeomorfismo $f: M \mapsto M$, los pozos son atractores, las fuentes son repulsores, y las sillas no son atractores ni repulsores.
\end{exercise}

\begin{exercise}\em
Encontrar un ejemplo  de difeomorfismo  $f: \mathbb{R}^2 \mapsto \mathbb{R}^2 $ que tenga  un punto fijo atractor que no sea pozo (i.e. que no sea hiperb\'{o}lico), otro que tenga un punto fijo repulsor que no sea fuente (i.e. que no sea hiperb\'{o}lico) y otro que tenga un punto fijo no hiperb\'{o}lico que no sea ni fuente ni pozo pero que todas las \'{o}rbitas futuras en un entorno cualquiera de $x_0$ suficientemente peque\~{n}o, o bien tiendan a $x_0$ o bien se salgan del entorno.
\end{exercise}

\vspace{.2cm}

\index{exponentes de Lyapunov! de punto peri\'{o}dico}
 {\bf Exponente de Lyapounov negativo significa contracci\'{o}n exponencial:} \em Si $f: S^1 \mapsto S^1$
es un difeomorfismo y $x_0$ es un punto fijo atractor hiperb\'{o}lico \em (i.e.  $|f'(x_0)|<1$, es decir $x_0$ es un pozo), \em entonces
$$\lim _{n \rightarrow + \infty} \frac{\log \mbox{\em dist} (f^n(x), x_0)}{n} = -\lambda <0
\; \; \forall x \mbox{ en alg\'{u}n entorno de } x_0$$ donde $-\lambda
= \log |f'(x_0)| < 0$ se llama exponente de Lyapounov en $x_0$. \em

{\begin{exercise}\em \em  \em
 Demostrar la afirmaci\'{o}n anterior. Sugerencia:
Probar que
$$\frac{\dist(f^{n+1}(x), x_0)}{\dist(f^{n}(x), x_0)}\rightarrow
e^{-\lambda}< 1$$
\end{exercise}
 {\bf Interpretaci\'{o}n:}   La distancia de
$f^n(x)$ al atractor hiperb\'{o}lico se contrae exponencialmente con
coeficiente asint\'{o}ticamente igual a $e$ elevado al exponente de Lyapounov $-\lambda < 0$.

\vspace{.2cm} {\bf Exponente de Lyapounov positivo significa
dilataci\'{o}n exponencial:} \em  Si $f: S^1 \mapsto S^1$ es un difeomorfismo  y $x_0$ es un punto fijo repulsor hiperb\'{o}lico \em (i.e.
$|f'(x_0)|>1$, es decir $x_0$ es una fuente) \em  entonces
$$\lim _{n \rightarrow - \infty} \frac{\log \mbox{\em dist} (f^n(x), x_0)}{n} = \sigma > 0
\; \; \forall x \mbox{ en alg\'{u}n entorno de } x_0$$ donde $\sigma =
\log |f'(x_0)| > 0$ se llama exponente de Lyapounov en $x_0$. \em

\begin{exercise}\em \em
\em
 Demostrar la afirmaci\'{o}n anterior. Sugerencia: Aplicar lo ya probado a $f^{-1}$ y la f\'{o}rmula de derivada de la funci\'{o}n inversa para deducir que
$$\lim_{n \rightarrow - \infty} \frac{\dist(f^{n+1}(x), x_0)}{\dist(f^{n}(x), x_0)}\rightarrow
e^{\sigma}>1$$
\end{exercise}
 {\bf Interpretaci\'{o}n:}   La distancia de
$f^n(x)$ al repulsor hiperb\'{o}lico (al crecer $n$ y mientras $f^n(x)$ est\'{e} en un entorno peque\~{n}o del repulsor) se dilata exponencialmente con
coeficiente asint\'{o}ticamente  igual a $e$ elevado al exponente de Lyapounov $\sigma > 0$.  

\vspace{.2cm}

\begin{exercise}\em {\bf Flujo polo norte-polo sur} \em  \label{ejerciciopolonortepolosur} \index{flujo polo norte-polo sur}
\em
 Llamaremos \em secci\'{o}n de  flujo polo norte-polo
 sur \em en el c\'{\i}rculo $S^1$, a un
difeomorfismo $f: S^1 \mapsto S^1$ de clase $C^1$, que preserva la
orientaci\'{o}n, y tal que existen solo 2 puntos fijos
$N$ y $S$,
son ambos hiperb\'{o}licos, $N$ repulsor y  $S$ atractor.
Graficar $f$ en
$[0,1]_{mod 1}$ tomando $0 \sim 1 = N$. Demostrar que todas las
\'{o}rbitas excepto $N$ y $S$ son mon\'{o}tonas y convergen a $S$.
(Sugerencia: ver prueba de la afirmaci\'{o}n
\ref{claimatractoreshiperbolicos}.) Demostrar que las \'{u}nicas
medidas de probabilidad invariantes son las combinaciones lineales
convexas de $\delta _N$ y $\delta _S$. Sugerencia: considerar una
partici\'{o}n numerable del intervalo $(0, S)$ formado por $A_n =
[x_{n}, x_{n+1})$ para $n \in \mathbb{Z}$ donde $x_0$ se elige cualquiera
en el intervalo abierto $(0,S)$ y $x_n := f^n(x_0)$ para todo $n \in \mathbb{Z}$. Usando argumento similar a la prueba del ejemplo
\ref{ejemplonoexistenmedidasinvariantes}, probar
que $\mu ((0,S)) =
0$. An\'{a}logamente probar que $\mu ((S, 1)) = 0$, de donde $\mu
(\{N,S\}) = 1$.

\end{exercise}

\begin{definition} \em \index{difeomorfismos! Morse-Smale} {\bf Difeomorfismo Morse-Smale en $S^1$.}
Un difeomorfismo $f: S^1 \mapsto S^1$ se dice   \em
Morse-Smale \em si preserva la orientaci\'{o}n y existen exactamente una cantidad finita de
puntos peri\'{o}dicos (todos del mismo per\'{\i}odo) y son todos ellos
hiperb\'{o}licos.
\end{definition}

\begin{exercise}\em \em
\em Probar que en un difeomorfismo $f$ Morse-Smale  en el c\'{\i}rculo las \'{u}nicas medidas invariantes son las
combinaciones lineales convexas de las medidas $$\frac{\delta
_{x_0} + \delta_{f(x_0)} + \ldots + \delta _{f^{p-1}(x_0)}}{p}$$
donde $x_0$ es un punto peri\'{o}dico de per\'{\i}odo $p$. Sugerencia:
Graficar $f^p$ en $S^1= [0,1]/\sim$ donde $0 \sim 1$ es un punto
peri\'{o}dico de per\'{\i}odo $p$. Probar que los atractores y los
repulsores se alternan. Probar que para toda medida invariante el
arco entre un repulsor y un atractor consecutivos tiene medida
cero,  usando el procedimiento
 del ejercicio \ref{ejerciciopolonortepolosur}.
\end{exercise}


\subsection{Din\'{a}mica topol\'{o}gica} \index{din\'{a}mica topol\'{o}gica}

\begin{definition} \label{definicionrecurrencia} \index{recurrencia! topol\'{o}gica} \em {\bf
Recurrencia topol\'{o}gica}. Sea $T: X  \mapsto X$ una transformaci\'{o}n
Borel medible en un espacio topol\'{o}gico $X$. Sea $x \in X$. Se
llama \em omega-l\'{\i}mite de $x$ \em al conjunto $$\omega (x)= \{ y
\in X: \exists n_j \rightarrow + \infty \mbox{ tal que }
T^{n_j}(x) \mapsto y    \}$$ \index{omega-l\'{\i}mite} \index{$\omega(x)$ omega-l\'{\i}mite} \index{alfa-l\'{\i}mite}  Cuando $T$ es bi-medible (i.e. $T$ es medible, invertible y con inversa medible) se llama \em
alfa-l\'{\i}mite de $x$ \em al conjunto
$$\alpha (x)= \{ y \in X: \exists n_j \rightarrow - \infty \mbox{
tal que } T^{n_j}(x) \mapsto y \}$$ Un punto $x$ se dice  \em
recurrente \em si
$$x \in \omega (x)$$  \index{punto! recurrente}
Dicho de otra forma $x$ es recurrente si existe una subsucesi\'{o}n
$n_j \rightarrow + \infty$ tal que $T^{n_j}(x) \rightarrow x$.
Luego para todo entorno $V$ de $x$ existe una subsucesi\'{o}n $n_j \in
\mathbb{N}$ tal que $T^{n_j}(V) \cap V \neq \emptyset$.
\end{definition}

\begin{exercise}\em Sea $X$ un espacio m\'{e}trico compacto y sea $T: X \mapsto $ continua. Probar que:

(a) $\omega(x)$ es compacto no vac\'{\i}o para todo $x \in X$.

(b) $\omega(T^n(x)) = \omega(x)$ para todo $n \in \mathbb{N}$, es decir el conjunto $\omega(x)$ depende de la \'{o}rbita por $x$ y no de qu\'{e} punto en la \'{o}rbita de $x$ se elija.

(c) $T (\omega(x)) = \omega(x) \subset T^{-1}(\omega(x))$ para todo $x \in \mathbb{N}$. Es decir $\omega(x)$ es un conjunto invariante por $T$ hacia el futuro.

(d) Si adem\'{a}s $T$ es un homeormorfismo (i.e. $T$ es continua, invertible y con inversa $T^{-1}$ continua) probar que:

 (i) $\alpha(x)$ es compacto no vac\'{\i}o para todo $x \in X$, $\alpha(T^n(x))= \alpha (x)$ para todo $n \in \mathbb{Z}$ y $T(\alpha(x)) = \alpha(x) = T^{-1}(\alpha(x))$.

(ii) Si $x \in \omega(x)$ entonces $\alpha(x) \subset \omega(x)$. Sugerencia: Si $x \in \omega(x)$ entonces $T^{-n}(x) \in \omega(x)$ para todo $n \in \mathbb{Z}$, luego $\alpha(x) \subset \overline{\omega(x)} = \omega(x)$.

 (iii) $x \in \omega(x)$ si y solo si $\omega(x) = \overline{o(x)}$, donde $o(x) = \{T^n(x): \ n \in \mathbb{Z}\}$.

\end{exercise}

\begin{definition} \em {\bf Conjunto no errante.} \index{conjunto! no errante}
Sea $T: X  \mapsto X$ una transformaci\'{o}n Borel medible en un
espacio topol\'{o}gico $X$. Un punto $x \in X$ es \em no errante \em
si para todo entorno $V$ de $x$ existe una sucesi\'{o}n $n_j
\rightarrow + \infty$ tal que $T^{n_j}(V) \cap V \neq \emptyset$.
El conjunto de los puntos no errantes de $T$ (que puede ser vac\'{\i}o)
se denota como $\Omega (T)$ y se llama \em conjunto no errante.
\em
\end{definition}

\begin{exercise}\em \em
\em Sea $X$ un espacio topol\'{o}gico de Hausdorff y sea $T: X \mapsto X$ Borel me\-dible.

(a) Probar que el conjunto de los puntos recurrentes est\'{a} contenido en
el conjunto no errante $\Omega (T)$  (la inclusi\'{o}n opuesta no es
necesariamente cierta, como se ver\'{a} m\'{a}s adelante).

(b) Sea $\mu $ una
medida de probabilidad invariante por $T$. Si $X$ tiene base
numerable de abiertos probar (sin  usar\- el enunciado del teorema de
recurrencia de Poincar\'{e} que viene m\'{a}s adelante) que $\mu (\Omega
(T))= 1$, es decir: casi todo punto es no errante para cualquier
medida de probabilidad invariante por $T$.

Sugerencia: Probar  que
para todo $V$ de la base de abiertos que tenga medida positiva la
sucesi\'{o}n de conjuntos medibles $T^{-n}(V), \; n \in \mathbb{N}$ no puede
ser de conjuntos disjuntos dos a dos a partir de un cierto $n_0$
en adelante. Deducir que existe una subsucesi\'{o}n $n_j \rightarrow +
\infty$ tal que $T^{-n_j}(V) \cap V \neq \emptyset$ y esto implica
$V \cap T^{n_j}(V) \neq \emptyset$. Un punto es errante (no es no
errante) si est\'{a} contenido en alg\'{u}n $V$ abierto tal que no cumple
lo anterior. Deducir que los puntos errantes forman un conjunto de
medida nula.
\end{exercise}

\begin{definition} \label{definitionTransitividad} \index{transitividad}
\em {\bf Transitividad topol\'{o}gica} Sea $X$ un espacio topol\'{o}gico y
$T: X \mapsto X$ Borel medible. Se dice que $T$ es \em transitiva
\em si dados dos abiertos $U$ y $V$ no vac\'{\i}os,  existe $n \geq 1$
tal que $T^n(U) \cap V \neq \emptyset$.

Sup\'{o}ngase que  $X$ es de Hausdorff sin puntos aislados.  Es f\'{a}cil
probar que \em si existe una \'{o}rbita positiva
 densa entonces $T$ es transitiva. \em Y si adem\'{a}s, $T$ es continua y el espacio topol\'{o}gico
  tiene base numerable de abiertos y es
  de Baire (esto es: toda intersecci\'{o}n numerable de abiertos densos es densa)
 entonces $T$ \em es transitiva si y solo si existe una \'{o}rbita
 positiva
 densa. \em

 La transitividad significa que para cualquier abierto $U$, por
 peque\~{n}o que sea, \em los iterados positivos de $U$ transitan por todo el espacio
  desde el punto de vista topol\'{o}gico \em (es decir, por todos los abiertos del espacio).

\end{definition}

\begin{remark} \em  \label{remark00}
  Se observa que para dos conjuntos cualesquiera $U$ y $V$:
$$T^n(U) \cap V \neq \emptyset \mbox{ si y solo si
 } T^{-n}(V) \cap U \neq \emptyset$$
 Es f\'{a}cil ver que si $T$ es continua y transitiva, entonces dados
 dos abiertos $U$ y $V$ no vac\'{\i}os, existe $n_j \rightarrow +
 \infty$ tal que $T^{-n_j}(V) \cap U \neq \emptyset$. En
 particular, tomando $U=V \ni x$, se deduce que:

 \em Si $T: X \mapsto X$ es continua y transitiva entonces $\Omega (T) =
 X$. \em

\end{remark}

\begin{exercise}\em
Probar las afirmaciones de la Observaci\'{o}n \ref{remark00} y las que est\'{a}n inmediatamente despu\'{e}s de la definici\'{o}n de transitividad en \ref{definitionTransitividad}.
\end{exercise}

\subsection{Recurrencia y Lema de Poincar\'{e}} \index{teorema! lema de Poincar\'{e}} \index{recurrencia}

     Sea $T: X \mapsto X $ Borel medible en un
espacio m\'{e}trico compacto $X$. Recordando  la Definici\'{o}n \ref{definicionrecurrencia} y teniendo en cuenta la compacidad secuencial de $X$, un punto $x$ es recurrente si y solo si \em para todo
entorno $V$ de $x$ existen infinitos iterados  hacia el futuro de
$x$ en $V$. \em Es decir, existe $n_j \rightarrow + \infty$ tal que
$T^{n_j}(x)\in V$.

\begin{definition}
\em Sea $T: X \mapsto X$ medible en un espacio medible $(X, {\mathcal A})$.
Sea $E \in {\mathcal A}$. Un punto $x \in E$ \em vuelve
infinitas veces a $E$ \em si existen infinitos iterados hacia el
futuro de $x$ en $E$. Mejor dicho: existe $n_j \rightarrow +
\infty$ tal que $T^{n_j}(x)\in E$ para todo $j \in \mathbb{N}$.
\end{definition}

Los siguientes dos teoremas, llamados Lemas de
Recurrencia de Poincar\'{e}, se encuentran  por ejemplo  en \cite[pag. 32-35]{Mane} (ver tambi\'{e}n \cite{ManeIngles}). En
\cite[\S 1.4]{Walters} se encuentra tambi\'{e}n la versi\'{o}n medible siguiente:

\begin{theorem} \label{teoPoincare} \index{teorema! lema de Poincar\'{e}}
 {\bf Lema de Recurrencia de Poincar\'{e}. Versi\'{o}n me\-dible}

Sea $T:X \mapsto X$ medible que preserva una medida de
probabilidad $\mu$. Sea $E$ un conjunto medible tal que $\mu
(E)>0$. Entonces $\mu$-c.t.p de $E$ vuelve infinitas veces a $E$.
\end{theorem}

{\em Demostraci\'{o}n: }
Sea $F_N := \bigcap_{n \geq N} T^{-n}(X \setminus E)$ y sea $F:= \bigcup_{N \geq 0} F_N$. Por construcci\'{o}n $x \in F$ si y solo si  $T^n(x)   \not \in E$ para todo $n $ suficientemente grande. Esto ocurre si y solo si la \'{o}rbita futura $\{T^n(x)\}_{n \geq 0}$ de $x$ no pasa infinitas veces por $E$. Basta probar entonces que $\mu(F \cap E)= 0$.

Por construcci\'{o}n $T^{-1}(F_N) = F_{N+1}$.   Como $\mu$ es $T$-invariante, tenemos $\mu (F_{N+1}) = \mu (F_N)$ para todo $N \geq 0$. Luego $\mu (F_N) = \mu (F_0) $ para todo   $N \geq 0$. Siendo $F_ {N+1} \supset F_N$ para todo $N \geq 0 $, entonces $\mu (F) =\lim _{N \rightarrow + \infty} \mu(F_N)  = \mu(F_0) $ y $F \supset F_0 $. Luego   $\mu (F \setminus F_0) = 0$. Como $E \cap F_0 = \emptyset$, deducimos que $E \cap F \subset F \setminus F_0$, de donde $\mu (E \cap F) = 0$, como quer\'{\i}amos demostrar.
 \hfill $\Box$

\begin{theorem}
  {\bf Lema de Recurrencia de Poincar\'{e}. Versi\'{o}n to\-po\-l\'{o}\-gica} \index{teorema! lema de Poincar\'{e}}

Sea $T: X \mapsto X$ Borel-medible en un espacio topol\'{o}gico $X$
con base numerable. Si $T$ preserva una medida de probabilidad
$\mu$, entonces $\mu$-c.t.p. es recurrente (es decir $x \in \omega
(x) \; \; \mu$-c.t.p. $x \in X$).

\end{theorem}

{\em Demostraci\'{o}n: } Sea $\{V_j\}_{j \in \mathbb{N}}$ una base de
abiertos. Por \ref{definicionrecurrencia}: $x \not \in \omega(x)$
si y solo si  $x \in \bigcup _{j \in \mathbb{N}} A_j$ para alg\'{u}n $j \in \mathbb{N}$, donde $A_j = \{x \in
V_j: x \mbox{ no vuelve infinitas veces a } V_j \}$. Por
\ref{teoPoincare} $\mu (A_j) = 0$ Luego, la uni\'{o}n numerable de los
conjuntos $A_j$, que coincide con el conjunto de los puntos no
recurrentes, tiene $\mu$ medida nula. $\Box$

\begin{exercise}\em \em  \em
Sea $(X, {\mathcal A})$ un espacio medible. Sea $T: X \mapsto X$
me\-di\-ble que preserva una probabilidad $\mu$. Sea $E \in {\mathcal A}$
tal que $\mu (E) >0$. Probar que
$$\sum_{n \in \mathbb{N}} \chi _E (T^n(x))$$ diverge $\mu$- c.t.p $ x \in
E$
\end{exercise}

\begin{exercise}\em \em  \em
Sea $T: X \mapsto X$ Borel medible en el espacio topol\'{o}gico $X$
compacto , y preservando una medida de probabilidad $\mu$. Sea
$\mbox{ supp }(\mu )$ el soporte compacto de $\mu$ (i.e. el
m\'{\i}nimo compacto  con medida $\mu$ igual a 1). Probar que
$\emptyset \neq \mbox{supp}(\mu )\subset \overline {\mbox
{Rec}(T)}$ siendo $\mbox{Rec}(T)$ el conjunto de los puntos
recurrentes de $T$.
\end{exercise}

\begin{theorem} \index{teorema! Hopf}
\em {\bf Teorema de Hopf.} \em Sea $T: \mathbb{R}^k \mapsto \mathbb{R}^k$ Borel
bimedible (medible, invertible con inversa medible) que preserva
la medida de Lebesgue $m$. Entonces casi todo punto de $\mathbb{R}^n$ o
bien es recurrente o bien tiene omega l\'{\i}mite vac\'{\i}o.
\end{theorem}

\begin{exercise}\em \em  \em
Demostrar el teorema de Hopf enunciado antes. Su\-ge\-rencia: $\mathbb{R}^k =
\bigcup_{i \in \mathbb{N}} X_i$ donde $X_i$ es una bola abierta de radio $r_i
\rightarrow + \infty$ creciente con $i$. Sea $$\widetilde {X}_i = \{x
\in X_i: T^j(x) \in X_i \mbox{ para infinitos valores positivos de
} j \} .$$ Sea $\widetilde  T _i : \widetilde  X_ i \mapsto \widetilde  X_i$ la
transformaci\'{o}n que a cada $x \in \widetilde  X_i$ hace corresponder el
primer retorno a $ X_i$: $\widetilde  T(x) = T^j (x) \in \widetilde
X_i$ para el m\'{\i}nimo $j = j(x)$ natural positivo tal que $T^j(x) \in X_i$. Probar que $m(\widetilde  T_i \widetilde  X_i) = m (\widetilde  X_i)$;
luego c.t.p. de $\widetilde  X_i$ est\'{a} en la imagen de $\widetilde  T_i$.
Probar que $\widetilde  T _i $ preserva $m$. Aplicar el teorema de
recurrencia de Poincar\'{e} para deducir que $m$-c.t.p de $\widetilde  X_i$
es recurrente. Probar que c.t.p. de $X_i$ o bien es recurrente o
bien su omega l\'{\i}mite no intersecta a $X_i$.
\end{exercise}

\begin{remark} \em
 En \cite{FrantzEnciclopedia}, y en la bibliograf\'{\i}a all\'{\i} citada, se rese\~{n}an varios otros resultados sobre recurrencia, adem\'{a}s de los lemas b\'{a}sicos de recurrencia de Poincar\'{e}. Algunos de estos resultados miden, en relaci\'{o}n a potencias de $n$, la frecuencia con la que \'{o}rbita futura del punto recurrente   $x$ se acerca a $x$, la vinculan con las medidas erg\'{o}dicas y con la entrop\'{\i}a m\'{e}trica del sistema (la entrop\'{\i}a m\'{e}trica es, gruesamente hablando,   una medici\'{o}n, ponderada seg\'{u}n una probabilidad invariante $\mu$, del desorden espacial que produce $f$ al ser iterada).
\end{remark}

\subsection{Ergodicidad} \index{ergodicidad} \index{medida! erg\'{o}dica}
\index{transformaci\'{o}n! erg\'{o}dica}
\begin{definition} \label{ergodicidad} \label{definitionErgodicidadI}
\em {\bf Ergodicidad I}

Sea $(X, {\mathcal A}, \mu ) $ un espacio de medida de probabilidad y
$T: X \mapsto X$
 medible que preserva $\mu $. Se dice que $T$ es \em erg\'{o}dica \em
 respecto a la medida $\mu$,
   o que $\mu$ es una \em medida erg\'{o}dica
\em para $T$,
  si dados dos
conjuntos medibles con medida positiva $U$ y $V$ existe $n \geq 1
$ tal que $T^n(U) \cap V \neq \emptyset$.
\end{definition}

{\bf Nota: } Observar que la ergodicidad es la versi\'{o}n en el contexto medible  de la transitividad topol\'{o}gica.  Se resalta  que por definici\'{o}n, si una medida es erg\'{o}dica para $T$, entonces es  $T$-invariante. No se define ergodicidad de medidas no invariantes.

\begin{definition} \label{definitionErgodicidadII} \index{ergodicidad} \index{medida! erg\'{o}dica} \index{transformaci\'{o}n! erg\'{o}dica}
\label{ergodicidad2} \em {\bf Ergodicidad II }

Sea $(X, {\mathcal A}, \mu )$ un espacio de medida de probabilidad y
$T: X \mapsto X$
 medible que preserva $\mu $. Se dice que \em $T$ es
 erg\'{o}dica  \em respecto a la medida $\mu$,
   o que $\mu$ es una \em medida erg\'{o}dica
\em para $T$,  si todo conjunto medible $A$ que sea invariante por
$T$ (es decir $T^{-1}(A) = A$) tiene o bien medida nula o bien
medida 1.
\end{definition}

\begin{theorem} \label{teoremaDefinicionesErgodicidad} \index{equivalencia de definiciones! de ergodicidad}
$T$ es erg\'{o}dica seg\'{u}n la definici\'{o}n I para la medida de
probabilidad $\mu$ si y solo si es erg\'{o}dica seg\'{u}n la definici\'{o}n
II.
\end{theorem}

{\em Demostraci\'{o}n: } Supongamos que no se cumple la definici\'{o}n II. Entonces, existe un conjunto  $A$   invariante por $T$ que tiene medida
positiva distinta de 1. Luego el complemento $A^c$ de $A$ es
tambi\'{e}n invariante por $T$ y tiene medida positiva. Como $A = T^{-1}(A)$, toda
\'{o}rbita positiva con estado inicial en $A$ est\'{a}
contenida en $A$.
Deducimos que no existe $n \in \mathbb{N}$ tal que $T^n(A)$ intersecta a $A^c$. Concluimos que
$T $ no  es erg\'{o}dica seg\'{u}n la definici\'{o}n I.

Rec\'{\i}procamente, suponemos por hip\'{o}tesis que se cumple la definici\'{o}n II. Todo conjunto $A$ invariante por $T$ tiene
medida o bien nula o bien 1. Supongamos por absurdo
que no se cumple
la definici\'{o}n I. Entonces existen  conjuntos medibles $U$ y $V$
con medida positiva tales que $\bigcup _{n \in \mathbb{N}} T^{-n}(V) \cap U =
\emptyset$. Entonces el conjunto $$A= \bigcap _{N \in \mathbb{N}} \big (\bigcup
_{n \geq N} T^{-n}(V)\big)$$ es medible, invariante por $T$ (verificar que $T^{-1}(A) = A$) y tiene
medida positiva (verificar que $\mu(A) \geq \mu (V) $ usando que $\mu$ es medida de probabilidad $T$-invariante),  pero $A$ no intersecta a $U$ que tambi\'{e}n tiene
medida positiva. De lo anterior se deduce que $A$ no puede tener
medida 1, con lo que encontramos un conjunto  invariante
que tiene medida positiva menor que 1,  contradiciendo la
hip\'{o}tesis. \hfill $\Box$

\begin{exercise}\em Sea $T: X \mapsto X$ que preserva una medida de probabilidad $\mu$.
(a) Probar que   $\mu$ es erg\'{o}dica para $T$ si y solo si todo conjunto $A$ medible invariante para el futuro (i.e. $A \subset T^{-1}(A) $) cumple $\mu(A)= 0$ \'{o} $\mu(A)= 1$). (b) Probar que $\mu$ es erg\'{o}dica para $T$ si y solo si todo conjunto $A$ medible invariante para el pasado (i.e. $T^{-1}(A) \subset A$) cumple $\mu(A)= 0$ \'{o} $\mu(A)= 1$. Sugerencias: (a) Considerar $B= \bigcup_{n \geq 0} T^{-n}(A)  $,  probar que  $\mu(B) = \mu(A)$ y que $B$ es $T$-invariante. (b) Considerar el complemento de $A$.
\end{exercise}

\vspace{.3cm}

\index{medida! positiva sobre abiertos}

 {\bf Definici\'{o}n: } Una medida de Borel $\mu$ en un espacio
 topol\'{o}gico se dice \em positiva sobre abiertos \em si $\mu (V) >0$
  para todo abierto $V$ no vac\'{\i}o.
\begin{exercise}\em \em  \em
Sea $X$ un espacio topol\'{o}gico y $T: X \mapsto X$ Borel medible que
preserva una medida de probabilidad $\mu$ que es positiva sobre
abiertos. Probar que si $T$ es erg\'{o}dica respecto de $\mu$ entonces
$T$ es transitiva (topol\'{o}gicamente).

\end{exercise}
\begin{exercise}\em \em \index{difeomorfismos! Morse-Smale} \index{transitividad}
 \em Probar que los difeomorfismos Morse-Smale en el c\'{\i}rculo tienen
 medidas erg\'{o}dicas pero no son transitivos topol\'{o}gicamente.
\end{exercise}

\begin{remark} \em  M\'{a}s adelante demostraremos que  \em la medida de Lebesgue es
 erg\'{o}dica para la rotaci\'{o}n irracional del c\'{\i}rculo \em (Teoremas \ref{teoremaErgodicidadRotacionIrracional} y \ref{teorema1}).
 Luego, como la medida de Lebesgue es positiva sobre abiertos, la rotaci\'{o}n irracional  es transitiva topol\'{o}gicamente. Entonces existe alguna \'{o}rbita
 densa. Es f\'{a}cil ver, usando que la rotaci\'{o}n en el c\'{\i}rculo conserva
 las distancias,  que al existir una \'{o}rbita densa, \em todas las
 \'{o}rbitas son densas. \em

 Tambi\'{e}n probaremos que el \em tent map $T$ en el intervalo es erg\'{o}dico respecto a la
 medida de Lebesgue. \em Luego $T$ es topol\'{o}gicamente transitivo y existe \'{o}rbita densa. Sin embargo no todas las \'{o}rbitas en el futuro por $T$ son densas: en efecto, existen \'{o}rbitas peri\'{o}dicas (que, como \'{o}rbitas en el futuro, son conjuntos finitos, y por lo tanto no son densos).
 \end{remark}

\begin{nada} \em  \index{promedio! temporal}
{\bf Promedios temporales asint\'{o}ticos de Birkhoff}

\end{nada} El Teorema Erg\'{o}dico de Birkhoff-Khinchin, que enunciaremos completamente m\'{a}s adelante ({Teorema \ref{theoremBirkhoff}}), establece que si $T: X \mapsto X$ es una transformaci\'{o}n medible que tiene medidas de probabilidad $T$-invariantes, entonces  para toda probabilidad   $\mu $ que sea $T$-invariante  y para toda funci\'{o}n $\psi \in L^1(\mu) $, existe $\mu$-c.t.p. $x \in X$ el siguiente l\'{\i}mite:
 $$ \widetilde \psi (x):= \lim_{n \rightarrow + \infty} \frac{1}{n} \sum_{j= 0}^{n-1} \psi \circ T^j(x).$$

\begin{theorem} {\bf Ergodicidad III.} \label{TeoremaErgodicidadIII}  \index{ergodicidad} \index{medida! erg\'{o}dica} \index{equivalencia de definiciones! de ergodicidad} \index{transformaci\'{o}n! erg\'{o}dica} \label{proposition_mu_ergodicapromediobirkhoff1}
Sea $T: X \mapsto X$ medible en un espacio medible $(X, {\mathcal A})$ y sea $\mu$ una probabilidad invariante por $T$. Las siguientes afirmaciones son equivalentes:

{\bf (i)} $\mu$ es  erg\'{o}dica para $T$

{\bf (ii)} Toda funci\'{o}n  $\psi : X \mapsto \mathbb{R}$    que sea medible e invariante por $T$  \em (es decir $\psi(x) = \psi (T(x))  \ \forall \, x \in X$), \em  es constante $\mu$-c.t.p.

{\bf (iii)} Para toda funci\'{o}n medible y acotada $\psi: X \mapsto \mathbb{R}$ existe el siguiente l\'{\i}mite \em $\mu$-c.t.p \em y es igual a $\int \psi \, d \mu$: \em
\begin{equation} \label{equationBirkhoffErgodica} \lim_{n \rightarrow + \infty} \frac{1}{n} \sum_{j= 0}^{n-1} \psi \circ T^j (x) = \int \psi \, d \mu \ \ \ \mu-\mbox{c.t.p.} \ x \in X \end{equation}

\end{theorem}

Probaremos el Teorema \ref{proposition_mu_ergodicapromediobirkhoff1} en el par\'{a}grafo \ref{paragraphPrueba00}.


\begin{remark} \em
{\bf Hip\'{o}tesis de Bolzmann de la Mec\'{a}nica Estad\'{\i}stica:} \index{hip\'{o}tesis de Bolzmann} \index{ergodicidad} \index{medida! erg\'{o}dica} \index{transformaci\'{o}n! erg\'{o}dica}
Antes de demostrar el Teorema \ref{proposition_mu_ergodicapromediobirkhoff1}, interpretaremos el significado de la afirmaci\'{o}n en la parte (iii) para fundamentar su relevancia.
Ella es un caso particular del Teorema Erg\'{o}dico de Birkhoff, que veremos m\'{a}s adelante, que  establece que para toda medida $\mu$ erg\'{o}dica se cumple la igualdad (\ref{equationBirkhoffErgodica}), no solo cuando $\psi$ es medible y acotada, sino tambi\'{e}n para toda $\psi \in L^1(\mu)$.

La igualdad (\ref{equationBirkhoffErgodica}) posee un significado relevante en la teor\'{\i}a erg\'{o}dica, pues afirma que      \em el promedio  espacial \em de cada  funci\'{o}n $\psi$ con respecto a la probabilidad $\mu$ en el espacio $X$  (i.e. el valor esperado $\int \psi \, d \mu$ de cada variable aleatoria $\psi$) coincide con \em el promedio temporal asint\'{o}tico \em de los valores observados de $\psi$ a lo largo de   $\mu$-casi toda \'{o}rbita. Este promedio temporal asint\'{o}tico es el l\'{\i}mite cuando $n \rightarrow + \infty$ de los promedios temporales $(1/n) \sum _{j= 0}^{n-1} \psi (T^j(x))$  de los valores de $\psi$,    observados a lo largo del pedazo finito de \'{o}rbita desde el instante 0 hasta el instante $n-1$. Salvo casos excepcionales, es falso que el l\'{\i}mite de los promedios temporales exista para todos los puntos $x \in M$. Adem\'{a}s, aunque para toda medida invariante $\mu$ ese l\'{\i}mite existe $\mu$-c.t.p. (Teorema erg\'{o}dico de Birkhoff),      es falso en general (salvo cuando $\mu$ es erg\'{o}dica) que coincida  con el promedio espacial de $\psi$ respecto a la probabilidad $\mu$. Por lo tanto las medidas erg\'{o}dicas $\mu$ para $T$ tienen un significado estad\'{\i}stico relevante, pues permite estimar el promedio temporal a largo plazo (esto es el promedio estad\'{\i}stico de las   series de observaciones $\psi(T^j(x))$  a largo plazo, llam\'{e}smole por ejemplo el \lq\lq clima\rq \rq)  en los sistemas determin\'{\i}sticos, para $\mu$-casi todo estado inicial $x \in X$, calculando el valor esperado de $\psi$ respecto a la probabilidad $\mu$. Sin embargo, aplicar la igualdad (\ref{equationBirkhoffErgodica}) para hacer esa estimaci\'{o}n, puede ser  muy err\'{o}neo cuando $\mu$ no es erg\'{o}dica,  o cuando   $\mu$ no es invariante por $T$, o cuando el estado inicial $x$ no pertenece al conjunto de $\mu$-probabilidad igual a 1.

La igualdad (\ref{equationBirkhoffErgodica}) de los promedios temporales asint\'{o}ticos con el promedio espacial (o sea el valor esperado) es lo que en la Mec\'{a}nica Estad\'{\i}stica se llama Hip\'{o}tesis de Boltzmann. Es una hip\'{o}tesis importante para demostrar propiedades de la din\'{a}mica de sistemas formado por una cantidad finita pero muy grande de part\'{\i}culas que evolucionan determin\'{\i}sticamente en el tiempo (por iteraciones de un mapa $T$ medible) preservando una medida de probabilidad dada  $\mu$   en el espacio $X$ de todos los estados posibles (llamado \lq\lq espacio de fases\rq\rq). Para los sistemas llamados conservativos esta medida de probabilidad $\mu$ es la medida de Lebesgue, normalizada para que $\mu(X)= 1$. Aplicando el Teorema \ref{proposition_mu_ergodicapromediobirkhoff1}, la hip\'{o}tesis de Boltzmann se traduce  en la hip\'{o}tesis de ergodicidad de esa pro\-ba\-bi\-lidad  $\mu$.

\end{remark}

\begin{nada} \em
\label{paragraphPrueba00}
{\bf Prueba del Teorema \ref{proposition_mu_ergodicapromediobirkhoff1}:}
\end{nada}

{\em Demostraci\'{o}n: }

{\bf (i) $\Rightarrow$ (ii): }

Sea $a \in \mathbb{R}$ y sea $A_a = \{x \in X: \psi(x) \leq a\} \subset X$. Como $\psi$ es medible, $A_a$ es medible . Siendo $\psi (T(x))= \psi (x)$ para todo $x \in X$, tenemos $T^{-1}(A_a) = A_a$. Como $\mu$ es erg\'{o}dica, entonces $\mu(A_a)$ vale 0 o 1. Considere la funci\'{o}n $g: \mathbb{R} \mapsto \{0,1\}$ definida como $g(a):= \mu (A_a)$. Por construcci\'{o}n $A_a \subset A_b$ si $a < b$. Luego $g(a) \leq g(b)$ si $a < b$. Entonces $g$ es creciente  y   solo puede tomar valores 0 \'{o} 1. Sea  \begin{equation}
\label{eqn11} k= \inf \{a \in \mathbb{R}: g(a)= 1\} = \sup \{ a \in \mathbb{R}: g(a)= 0\} \in [-\infty, + \infty].\end{equation} Probemos que $k \in \mathbb{R}$. En efecto $\psi(x) \in \mathbb{R}$ para todo $x \in M$; entonces $\emptyset=  \bigcap_{n \in \mathbb{N}} A_{-n}, \ \ X = \bigcup_{n \in \mathbb{N} }A_n$. Como $A_n \subset A_{n+1}$ para todo $n \in \mathbb{N}$, entonces $\lim_{n \rightarrow + \infty} \mu (A_n) = 1$ y $\lim _{n \rightarrow - \infty} \mu(A_{n}) = 0$. Esto implica que existe $n_0 \in \mathbb{N}$ tal que   $g(n_0) = 1$  y $g(-n_0) = 0$. Por lo tanto $-n_0 \leq k \leq n_0 $. Ahora probemos que $\mu(A_k)= 1$. De (\ref{eqn11}) deducimos que
 $0 =\mu(A_{k-(1/n)}) \leq \mu (A_k) \leq \mu (A_{k+(1/n)}) \ \forall \ 1 \leq n \in \mathbb{N}$
 Adem\'{a}s $A_{k} = \bigcap_{n= 1}^{\infty} A_{k + (1/n)} = X $, y $A_{k + (1/{n+1})} \subset A_{k + (1/n)}$. Luego $\mu (A_k)= \lim_{n \rightarrow + \infty}\mu(A_{k + (1/n)}) = 1$.  Hemos probado que $\mu(A_k) = 1$.
 Por otra parte $ B_k:= \{x \in X: \psi(x) \geq k\}  = \bigcap _{n= 1}^{\infty} (X \setminus A_{k-(1/n)})$, donde $\mu (X \setminus A_{k-(1/n})) = 1 $. Como $A_{k-(1/{n+1})} \supset A_{k -(1/n)}$, obtenemos $\mu(B_k) = 1$. Por lo tanto $\mu (A_k \cap B_k) = 1$, es decir para $\mu-$c.t.p. $x \in M$ se cumple $\psi(x) = k$.

 {\bf ii) $\Rightarrow$ iii):}

 Como $\mu$ es invariante por $T$ tenemos
 $\int \psi \circ T^j \, d \mu = \int \psi \, d \mu \ \ \forall \ j \in \mathbb{N}.$
 Luego $ \int \frac{1}{n} \sum_{j= 0}^{n-1} \psi \circ T^j \, d \mu = \int \psi \, d \mu \ \ \forall \ \ n \geq 1.$
 Tomando l\'{\i}mite cuando $n \rightarrow + \infty$ y aplicando el Teorema de Convergencia Dominada, deducimos que
 $$\int \lim _{n \rightarrow  + \infty} \frac{1}{n} \sum_{j= 0}^{n-1} \psi \circ T^j \, d\mu = \lim_{n \rightarrow + \infty} \int \frac{1}{n} \sum_{j= 0}^{n-1} \psi \circ T^j \, d \mu = \int \psi \, d \mu  $$

   Consideremos la siguiente funci\'{o}n  definida $\mu$-c.t.p: $$\Phi(x)  :=   \lim_{n \rightarrow + \infty} \frac{1}{n} \sum_{j = 0}^{n-1} \psi \circ T^j(x)    \mu-\mbox{c.t.p. } x \in X. $$
   Aqu\'{\i}, en la hip\'{o}tesis de existencia $\mu$-c.t.p. de esta funci\'{o}n $\Phi$, es decir en la hip\'{o}tesis de existencia $\mu$-c.t.p. del l\'{\i}mite de los promedios temporales, estamos aplicando el Teorema Erg\'{o}dico de Birkhoff (como si ya estuviera demostrado). Entonces tenemos
   $$\int \Phi \, d \mu  = \int \psi \, d \mu.$$

   Para tener definida la funci\'{o}n $\Phi$ en todo punto $x \in X$, tomemos $$\Phi (x):= 0 \ \forall \ x \in X \mbox { tal que } \not \exists \lim _{n \rightarrow + \infty} \frac{1}{n} \sum_{j = 0}^{n-1} \psi \circ T^j(x). $$

   Afirmamos que $\Phi$ es invariantes por $T$. Basta chequear que dado $\epsilon >0 $ existe $N$ tal que para todo $n \geq N$ se cumple   $$\Big | \frac{1}{n} \sum _{j= 0}^{n-1} \psi \circ T^j(T(x)) -    \frac{1}{n} \sum _{j= 0}^{n-1} \psi \circ T^j( x)\Big| =$$ $$ =\Big |\frac{1}{n} \psi(T^{n+1}(x)) -  \psi(x) \Big| < \epsilon \ \ \ \forall \ x \in X.$$
   En efecto la desigualdad anterior se verifica para todo $n$ suficientemente grande porque, por hip\'{o}tesis, $0 \leq |\psi(x)| \leq k$ para todo $x \in X$, para cierta constante $k >0$.
   Entonces, para toda sucesi\'{o}n $n_j \rightarrow + \infty$ tal que exista el l\'{\i}mite del promedio temporal hasta $n_j$ de $\psi$, con estado inicial $x$, tambi\'{e}n existe ese l\'{\i}mite con estado inicial $T(x)$, y ambos  l\'{\i}mites coinciden. Deducimos que $$\Phi(x) = \Phi(T(x)) \ \ \forall \ x \in X.$$ Es decir, $\Phi$ es funci\'{o}n real  invariante  por $T$.  Aplicando la hip\'{o}tesis (ii) tenemos que $\Phi$ es $\mu$-c.t.p. constante, igual  a cierto n\'{u}mero real $K$. Entonces $$\int \psi \, d \mu = \int \Phi \, d \mu = \int K \, d \mu  = K.$$
   Deducimos que $K=   \int \psi \, d \mu$, es decir $\Phi(x) = \int \psi \, d \mu$ para $\mu$-c.t.p. $x \in X$. Esto implica la  igualdad (\ref{equationBirkhoffErgodica}), probando   (iii).

   {\bf (iii) $\Rightarrow $ (i)}

    Sea $A \subset X$ medible e invariante por $T$, y consideremos la funci\'{o}n caracter\'{\i}stica $\chi_A$. Es una funci\'{o}n medible y acotada, y como $A$ es invariante por $T$, para todo $x \in X$ tenemos $\chi_A \circ T^j(x) = \chi_{T^{-j}(A)}(x) = \chi_A(x) \in \{0,1\}$. Entonces $$\widetilde  \chi_A   = \lim_{n \rightarrow +\infty} \sum_{j= 0}^{n-1} \chi_A (T^{j}(x)) = \chi _{A} (x) \in \{0,1\} \ \forall \ x \in X.$$ En particular vale la igualdad anterior para $\mu$-c.t.p. $x \in X$. Como  $$\widetilde (\chi_A)(x) = \int \chi_A \, d \mu = \mu (A) \ \ \mu-\mbox{c.t.p. } x \in X   $$
    deducimos que $\mu(A) \in \{0,1\}$. Por lo tanto $\mu$ es erg\'{o}dica, terminando de probar (i).
\hfill $\Box$

\subsection{Existencia de medidas erg\'{o}dicas} \label{sectionExistenciaMedidasErgodicas}

\begin{theorem} \index{teorema! de existencia de! medidas erg\'{o}dicas}
{\bf Existencia de medidas erg\'{o}dicas}    \label{teoremaExistenciaMedErgodicas} \index{teorema! de medidas erg\'{o}dicas! existencia}
Sea $X$ un espacio m\'{e}trico compacto y sea $T: X\mapsto X$ continua. Entonces existen medidas de probabilidad   erg\'{o}dicas para $T$. Adem\'{a}s, toda medida $T$-invariante  es el l\'{\i}mite en la topolog\'{\i}a d\'{e}bil$^*$ de una sucesi\'{o}n de medidas que son combinaciones lineales finitas de medidas erg\'{o}dicas.
\end{theorem}

Demostraremos el teorema \ref{teoremaExistenciaMedErgodicas} de existencia de medidas erg\'{o}dicas al final de esta secci\'{o}n, en el par\'{a}grafo   \ref{paragraphPruebaTeoremaExistMedidasErgodicas}.

\newpage

\begin{nada} \em
{\bf Singularidad mutua  y continuidad absoluta} \index{medida! absolutamente continua} \index{medida! mutuamente singulares}
\end{nada}

Recordamos que dos medidas de probabilidad $\mu$ y $\nu$ se dicen mutuamente singulares $\mu \perp \nu$ cuando existe alg\'{u}n conjunto medible $A \subset X$ tal que $\mu(A)= 1$ y $\nu(A)= 0$. Entonces $\mu (X \setminus A)= 0$ y $\nu(X \setminus A)= 1$ y la relaci\'{o}n mutuamente singular es sim\'{e}trica.

Dadas dos medidas de probabilidad $\mu$ y $\nu$ se dice que $\mu$ es absolutamente continua respecto de $\nu$, y se denota $\mu \ll \nu$, cuando para todo conjunto medible $A$  tal que $\nu(A)= 0$ se cumple $\mu (A) = 0$.

Se dice que dos medidas $\mu$ y $\nu$ son equivalentes, y se denota $\mu \sim \nu $, cuando $\mu \ll \nu$ y $\nu \ll \mu$.

Se observa que si $\mu \ll \nu$ entonces $\mu \not \perp \nu$ (el rec\'{\i}proco es falso).

El siguiente teorema es cl\'{a}sico en la Teor\'{\i}a abstracta de la Medida (en particular en la  Teor\'{\i}a de Probabilidades):

{\bf Teorema de Descomposici\'{o}n de Lebesgue-Radon-Nikodym. } \index{teorema! Lebesgue-Radon-Nikodym}

\em Dadas dos medidas de probabilidad $\mu$ y $\nu$ existen   dos probabilidades $\mu_1$ y $\mu_2$, y un \'{u}nico real $t \in [0,1]$, tales que
$$\mu = t \mu_1 + (1-t) \mu_2, \ \ \mu_1 \ll \nu, \ \ \mu_2 \perp \nu.$$
Si adem\'{a}s $ t \neq 0, 1$ entonces $\mu_1$ y $\mu_2$ son \'{u}nicas. \em

El enunciado cl\'{a}sico de este Teorema, establece que    \em existen   \'{u}nicas las medidas finitas $t \mu_1$ y $(1-t) \mu_2$ \em (posiblemente alguna de ellas es cero) tales que sumadas dan $\mu$, siendo $t \mu_1 \ll \nu$ y $(1-t) \mu_2 \perp \nu$.

La demostraci\'{o}n del Teorema de Descomposici\'{o}n de Lebesgue-Radon-Nikodym se encuentra por ejemplo en \cite[Theorem 3.8]{Folland} \'{o} en \cite[Teorema 6.2.3]{Rudin}.

Del teorema de Lebesgue-Radon-Nikodym se deduce que  en el caso particular  $\mu \ll \nu$, se cumple $t= 1$, $\mu_1 = \mu$, y $\mu_2$ es cualquiera. An\'{a}logamente, si $\mu \perp \nu$, entonces $t= 0$, $\mu_2= \mu$ y $\mu_1$ es cualquiera. En el caso que $\mu \not \ll \nu, \ \mu \not \perp \nu$, es \'{u}nica la pareja  $(\mu_1, \mu_2)$ de probabilidades en la descomposici\'{o}n de Radon-Nikodym.

Volvamos ahora a las propiedades de las medidas de probabilidad erg\'{o}dicas para una transformaci\'{o}n $T$:
\begin{theorem} \label{teoremaErgodicaMutuamSingular} {\bf Singularidad mutua de medidas erg\'{o}dicas} \index{medida! erg\'{o}dica} \index{teorema! de medidas erg\'{o}dicas! singularidad mutua} \index{transformaci\'{o}n! erg\'{o}dica} \index{ergodicidad}

 Sea $T: X \mapsto X$ medible en el espacio medible  $(X, {\mathcal A})$.

 {\em (a) } Si existen dos medidas de probabilidad diferentes $\mu$ y $\nu$, ambas  erg\'{o}dicas   para la transformaci\'{o}n $T$, entonces $\mu \perp \nu$.

{\em (b) } Si $\mu$ y $\nu$ son   medidas de probabilidad erg\'{o}dicas para $T$ y si $\mu \ll \nu$, entonces $\mu= \nu$.
\end{theorem}
{\em Demostraci\'{o}n: }
Afirmamos que, si  $\mu$ y $\nu$ son ambas erg\'{o}dicas para $T$, y si para todo conjunto $A$ invariante por $T$ se cumple $\mu(A) = \nu(A)$, entonces $\mu = \nu$. En efecto, sea $B$ cualquier conjunto medible, y denotemos $\chi_B$ a la funci\'{o}n caracter\'{\i}stica de $B$.   Aplicando la propiedad (iii) del Teorema \ref{proposition_mu_ergodicapromediobirkhoff1}  a la funci\'{o}n $\psi= \chi_B$, sabemos que existe un conjunto $A_1$, que es $T$-invariante, tal que $\mu(A_1) = 1$, y se cumple la igualdad (\ref{equationBirkhoffErgodica}) para la medida $\mu$ y para todo $x \in A_1$. An\'{a}logamente, existe un conjunto $A_2$, que es $T$-invariante, tal que $\nu(A_2) = 1$, y se cumple la igualdad (\ref{equationBirkhoffErgodica}) para la medida $\nu$ y para todo $x \in A_2$. Por hip\'{o}tesis $\nu(A_1) = \mu(A_1)= 1 = \nu(A_2) = \mu(A_2)$. Entonces $\nu(A_1 \cap A_2) = \mu(A_1 \cap A_2) = 1$ y para todos los puntos $x \in A_1 \cap A_2$ se cumple  la igualdad
   $$\mu(B) = \int \chi_B \, d \mu = \lim _{n \rightarrow + \infty} \frac{1}{n} \sum_{j= 0}^{n-1} \chi_B \circ T^j(x) = \int \chi_B \, d \nu = \nu (B)$$
   Deducimos entonces  que $\mu(B) = \nu(B)$ para todo conjunto medible $B$, de donde $\mu= \nu$. Hemos terminado de probar la afirmaci\'{o}n del principio.

   Demostremos ahora la parte (a) del teorema \ref{teoremaErgodicaMutuamSingular}.  Debido a la afirmaci\'{o}n reci\'{e}n demostrada, como $\nu \neq \mu$  y son ambas erg\'{o}dicas, existe un conjunto $A$ que es $T$-invariante, tal que $\mu(A)   \neq \nu (A)$. Por definici\'{o}n de ergodicidad $\mu(A), \nu (A) \in \{0,1\}$. Luego (eventualmente sustituyendo $A$ por su complemente en caso necesario), tenemos $\mu(A)= 1$ y $\nu(A)= 0$, de donde $\mu \perp \nu$, demostrando (a).

   La parte (b) es una consecuencia inmediata de (a), pues si $\mu \ll \nu$ entonces $\mu \not \perp \nu$. Por lo tanto, como $\mu$ y $\nu$ son erg\'{o}dicas y    $\mu  \not \perp  \nu$, aplicando la parte (a) deducimos que  $\mu = \nu$, como quer\'{\i}amos probar.
\hfill $\Box$

\begin{nada} \em
 \label{paragraphPruebaTeoremaExistMedidasErgodicas}   \index{teorema! de existencia de! medidas erg\'{o}dicas}
{\bf Demostraci\'{o}n del Teorema \ref{teoremaExistenciaMedErgodicas}  (existencia de medidas er\-g\'{o}\-dicas)}
\end{nada}

El teorema \ref{teoremaExistenciaMedErgodicas} es  un corolario inmediato del   teorema \ref{teoremaconvexhullmedidasergodicas} que demostraremos a continuaci\'{o}n.  Las definiciones de \em puntos extremales de un conjunto compacto y convexo \em en un espacio vectorial topol\'{o}gico, y de \em envolvente compacta convexa, \em se emplean en el enunciado siguiente, y se incluyen abajo del mismo.

\begin{theorem} \label{teorema2_5} \label{teoremaErgodicasExtremales}
\label{teoremaconvexhullmedidasergodicas} \index{teorema! de medidas erg\'{o}dicas! extremalidad}
\index{medida! extremal}
\index{medida! erg\'{o}dica} \index{equivalencia de definiciones! de ergodicidad} \index{transformaci\'{o}n! erg\'{o}dica} \index{ergodicidad}

  Sea $T: X \mapsto X$ continua en el espacio m\'{e}trico compacto $X$. Sea ${\mathcal M}$ el conjunto de medidas de probabilidad borelianas en $X$, con la topolog\'{\i}a d\'{e}bil$^*$. Sea ${\mathcal M}_ T \subset {\mathcal M}$ el conjunto de probabilidades invariantes por $T$ y sea ${\mathcal E}_T \subset {\mathcal M}_T$ el conjunto de las probabilidades (invariantes) erg\'{o}dicas para $T$.

Entonces

 {\bf (a)} ${\mathcal M}_T$ es compacto y convexo.

 {\bf (b)} ${\mathcal E}_T$ coincide con  el conjunto de puntos extremales de ${\mathcal M}_T$.

  {\bf (c)}   ${\mathcal M}_T$ coincide con la envolvente convexa compacta de ${\mathcal E}_T$.

\vspace{.5cm}

\em Demostraremos el Teorema \ref{teoremaExistenciaMedErgodicas} m\'{a}s adelante en esta secci\'{o}n. Para poder demostrarlo, necesitamos definir convexidad,  envolvente convexa, puntos extremales y ver las propiedades de estos conceptos:

\end{theorem}

{\bf Combinaciones convexas  y puntos extremales} \index{combinaci\'{o}n convexa}

Recordemos las siguientes definiciones y el teorema de Krein-Milman del An\'{a}lisis Funcional:

$\bullet$ Sea ${\mathcal A}$  s un conjunto no vac\'{\i}o de un espacio vectorial topol\'{o}gico. Se llama \em   combinaci\'{o}n convexa \em de puntos de ${\mathcal A}$ a cualquier punto del espacio que pueda escribirse como una combinaci\'{o}n lineal finita $$t_1 a_1 + t_2 a_2 + \ldots t_k a_k, $$  tal que $a_i \in {\mathcal A}$, \    $0 \leq t_i \leq 1$ para todo $1 \leq i \leq k$,   y $\sum_{i= 1}^k t_i= 1$.
Denotamos como $ec({\mathcal A})$ al conjunto de todas las combinaciones convexas de puntos de ${\mathcal A}$, llamado envolvente convexa de ${\mathcal A}$. Observar que ${\mathcal A} \subset ec({\mathcal A})$.

$\bullet$ El conjunto ${\mathcal A}$ se dice \em convexo \em si contiene a todas las combinaciones convexas de sus puntos, es decir ${\mathcal A} = ec({\mathcal A})$. Es inmediato deducir, por inducci\'{o}n en $k \geq 2$, que un conjunto ${\mathcal A}$ es convexo, si y solo si $$ta + (1-t) b \in {\mathcal A} \ \ \ \ \forall \   t \in [0,1], \ \ \forall \ a,b \in {\mathcal A}.$$

$\bullet$ Se llama \em envolvente convexa cerrada \em de ${\mathcal A}$  a $\overline {ec({\mathcal A})}$, donde $\overline{\  \ \ }$   indica la clausura (o adherencia).
\index{envolvente convexa}

$\bullet$  Si ${\mathcal A}$ es tal que $\overline {ec({\mathcal A})}$ es compacto, llamaremos a este \'{u}ltimo conjunto    \em la envolvente convexa compacta \em de ${\mathcal A}$.

$\bullet$ Si ${\mathcal K}    $ es un subconjunto   no vac\'{\i}o   de un espacio vectorial topol\'{o}gico, se llama \em punto extremal \em de ${\mathcal K}$ (cuando existe)  a un punto $a \in {\mathcal K}$ tal que las \'{u}nicas combinaciones convexas $t  b + (1-t) c = a $, con $0 \leq t \leq 1 \ \ b, \  c \in {\mathcal K}$, son aquellas para las cuales $b = a$ \'{o} $c= a$.

\vspace{.2cm}

$\bullet$ {\bf Teorema de Krein-Milman}

\em Todo compacto no vac\'{\i}o y convexo en un espacio vectorial topol\'{o}gico  contiene puntos extremales y coincide con la envolvente convexa compacta del conjunto de sus puntos extremales. \em

La demostraci\'{o}n del Teorema de Krein-Milman puede encontrarse por ejemplo en \cite[Teorema 3.21]{RudinFuncional}.

\begin{exercise}\em Sean $x_1, x_2 \in \mathbb{R}^2$. Chequear que toda combinaci\'{o}n convexa de $\{x_1, x_2\}$ est\'{a} en el segmento con extremos $x_1, x_2$ y rec\'{\i}procamente. Para ilustrar el teorema de Krein-Milman, considerar un pol\'{\i}gono regular $K$ en $\mathbb{R}^2$. ($K$ es la uni\'{o}n del interior del pol\'{\i}gono con su borde). Entonces $K$ es compacto. Chequear con argumentos geom\'{e}tricos  que: (a) el pol\'{\i}gono $K$ es   convexo; (b) los puntos extremales de $K$ son los v\'{e}rtices del pol\'{\i}gono; (c) todo punto de $K$ es una combinaci\'{o}n convexa de sus v\'{e}rtices.  En el ejemplo del pol\'{\i}gono,   la cantidad de puntos extremales es finita. Para ilustrar el teorema de Krein-Milman cuando la cantidad de puntos extremales es infinita, considerar una circunferencia $S^1 \subset \mathbb{R}^2$ y el compacto $K$ que es la uni\'{o}n de la circunferencia $S^1$ (borde de $K$) con la regi\'{o}n acotada encerrada por $S^1$ (interior de $K$). Entonces $K$ es compacto. Chequear  con argumentos geom\'{e}tricos  que: (d) $K$ es   convexo; (e) el conjunto de   puntos extremales de $K$ es la circunferencia $S^1$; (f)   todo punto de $K$ es combinaci\'{o}n convexa de dos puntos extremales. \end{exercise}

{\bf Demostraci\'{o}n del Teorema \ref{teoremaconvexhullmedidasergodicas}}

{\em Demostraci\'{o}n: }
{\bf (a)} Es inmediato chequear que ${\mathcal M}_T$ es convexo: en efecto,  si $\mu_1, \mu_2 \in {\mathcal M}_T$ entonces  para todo conjunto boreliano $B$ se cumple $$\mu_i(T^{-1}(B))= \mu_i(B) \ \ \ i= 1,2.$$ Luego, cualquiera sea $0 \leq t \leq 1$ se tiene $$[t\mu_1 + (1-t) \mu_2] (T^{-1}(B)) =[t\mu_1 + (1-t) \mu_2] ( B),$$ de donde $t \mu_1 + (1-t) \mu_2 \in {\mathcal M}_T$, probando que ${\mathcal M}_T$ es convexo.

Ahora probemos que ${\mathcal M}_T$ es compacto.
Recordemos que ${\mathcal M}$ es compacto con la topolog\'{\i}a d\'{e}bil$^*$ (esto fue demostrado en la secci\'{o}n \ref{seccionpruebateoexistmedinvariantes}). Luego, para demostrar que ${\mathcal M}_T$ es compacto, basta probar   que ${\mathcal M}_T$ es cerrado en ${\mathcal M}$. Sea $\mu_n \rightarrow \mu$ en ${\mathcal M}$ tales que $\mu_n \in {\mathcal M}_T$. Para deducir que $\mu \in{\mathcal M}_T$ basta recordar que   es continuo el operador $T^*:{\mathcal M} \mapsto {\mathcal M}$,   definido como $[T^*\mu] (B) = \mu(T^{-1}(B))$ para todo boreliano $B$. (La continuidad de $T^*$ fue demostrada en la secci\'{o}n \ref{seccionpruebateoexistmedinvariantes}). Entonces $$T^* \mu = \lim_ n T^*\mu_n = \lim \mu_n = \mu,$$ de donde $\mu \in {\mathcal M}_T$. Esto termina la prueba de que ${\mathcal M}_T$ es cerrado en el espacio m\'{e}trico compacto ${\mathcal M}$, y por lo tanto ${\mathcal M}_T$ es compacto.

{\bf (b)} Probemos que $\mu$ es erg\'{o}dica para $T$ si y solo si es punto extremal de ${\mathcal M}_T$. Primero asumamos que $\mu$ no es erg\'{o}dica. Entonces existe un boreliano $T$-invariante   $A \subset X$ tal que $0 <   \mu(A) < 1$. Sea $A^c= X \setminus A$ y  consid\'{e}rense las probabilidades definidas   para todo boreliano $B$ por: $$\mu_1 (B) := \frac{\mu(B \cap A)}{\mu(A)}, \ \ \ \ \ \mu_2(B):= \frac{\mu(B \cap A^c)}{\mu(A^c)}.$$ Por construcci\'{o}n, tenemos $$\mu = t \mu_1 + (1-t) \mu_2, \ \ \ \ \mbox{ donde }  t:= \mu(A).$$ N\'{o}tese que $\mu_1 \neq \mu_2$, porque $\mu_1(A) = 1$ y $ \mu_2(A) = 0$. Adem\'{a}s $0 <t < 1$, de donde $$\mu \neq \mu_1, \ \ \mu \neq \mu_2.$$ Por lo tanto $\mu$ no es punto extremal de ${\mathcal M}_T$. Hemos probado que si $\mu$ no es erg\'{o}dica, entonces $\mu$ no es punto extremal de ${\mathcal M}_T$.

Ahora probemos el rec\'{\i}proco. Asumamos que $\mu$ es erg\'{o}dica y probemos que $\mu$ es punto extremal de ${\mathcal M}_T$. Sea   $\mu = t \mu_1 + (1-t) \mu_2$,  donde $ 0 \leq t \leq 1$ y $\mu_1, \mu_2 \in {\mathcal M}_T$. Por un lado, si $t= 0$ \'{o} $t= 1$ entonces $\mu = \mu_1$ \'{o}   $\mu = \mu_2$. Por otro lado, si $0 < t < 1$ entonces $\mu_1 \ll \mu$ y $\mu_2 \ll \mu$ (pues si $\mu(B)= 0$ siendo la suma de dos sumandos no negativos, cada sumando debe ser cero, de d\'{o}nde $0=\mu_1(B)=\mu_2(B)$). Sea $A$ un conjunto $T$-invariante. Por hip\'{o}tesis $\mu$ es erg\'{o}dica, entonces $\mu(A)= 0$ \'{o} $\mu(A^c)= 0$. Como $\mu_i \ll \mu$  para $i= 1,2$, deducimos que   $\mu_i(A) = 0$ \'{o} $\mu_i(A)= 1$. Luego $\mu_i$ es erg\'{o}dica. Por lo demostrado en el Teorema \ref{teoremaErgodicaMutuamSingular}, $\mu_i \ll \mu$ siendo $\mu_i$ y $\mu$ erg\'{o}dicas, implica que $\mu_i = \mu$. En este caso tenemos entonces $\mu= \mu_1 = \mu_2$. Hemos probado, en todos los casos, que las \'{u}nicas combinaciones convexas de $\mu$ son aquellas para las cuales $\mu = \mu_1$ \'{o} $\mu = \mu_2$. Entonces por definici\'{o}n, $\mu $ es extremal  como quer\'{\i}amos demostrar.

{\bf (c)} Es consecuencia directa de a) y b), y del  Teorema de Krein Milman. \hfill $\Box$

\vspace{.3cm}

{\bf Fin de la prueba del Teorema \ref{teoremaExistenciaMedErgodicas}:

Existencia de medidas erg\'{o}dicas}

{\em Demostraci\'{o}n: }
 Por la parte c) del Teorema \ref{teoremaconvexhullmedidasergodicas}, $${\mathcal M}_T = \overline{e.c.{\mathcal E}_T},$$ donde ${\mathcal E}_T$ denota el conjunto de medidas erg\'{o}dicas para $T$ y ${\mathcal M}_T$ el de todas las medidas $T$-invariantes. Por el Teorema \ref{teoremaExistenciaMedInvariantes}, ${\mathcal M}_T \neq \emptyset$. Luego ${\mathcal E}_T \neq \emptyset$. Para probar la \'{u}ltima afirmaci\'{o}n del Teorema \ref{teoremaExistenciaMedInvariantes} aplicamos la Definici\'{o}n  de la clausura  $\overline{e.c.{\mathcal E}_T}$ de la envolvente convexa del conjunto ${\mathcal E}_T$ de probabilidades erg\'{o}dicas. Recordamos que la envolvente convexa $ {e.c.{\mathcal E}_T}$ es el conjunto de todas las medidas que se obtienen como combinaciones lineales convexas de medidas erg\'{o}dicas. Como toda medida invariante $\mu$ est\'{a} en $\overline{e.c.{\mathcal E}_T}$, entonces $\mu$ se puede aproximar, tanto como uno desee en la topolog\'{\i}a d\'{e}bil$^*$ del espacio de probabilidades,  por medidas que son por combinaciones lineales finitas y convexas de medidas erg\'{o}dicas.
\hfill $\Box$

\begin{exercise}\em
Sea $T: X \mapsto X$ una transformaci\'{o}n continua en un espacio m\'{e}trico compacto.

(I) Asuma que existe una medida $\mu$ invariante por $T$ y positiva sobre abiertos (no necesariamente erg\'{o}dica).

(a) Probar que para cada abierto $V$ no vac\'{\i}o, existe una medida erg\'{o}dica $\nu_V $ tal que $\nu_V(V) >0$.

Sugerencia: Por absurdo, si $\nu(V)= 0$ para toda medida erg\'{o}dica $\nu$, entonces $\rho(V)= 0$ para toda medida de probabilidad $\rho$ que sea combinaci\'{o}n convexa finita de erg\'{o}dicas. Probar que $\rho(V)= 0$ para toda medida $\rho$ en la envolvente convexa compacta de las erg\'{o}dicas (con la topolog\'{\i}a d\'{e}bil$^*$ del espacio ${\mathcal M}$ de probabilidades). Usar el Teorema \ref{teoremaconvexhullmedidasergodicas}  para concluir que $\mu(V)= 0$, contradiciendo la hip\'{o}tesis.

(b) Usando el Lema de Recurrencia de Poincar\'{e} demostrar que el conjunto de los puntos recurrentes es denso en $X$.

(c) Concluir que $\Omega(T) = X$.

 (II) Asuma que para todo abierto $V$ existe una medida de probabilidad invariante $\mu_V$ (no necesariamente erg\'{o}dica), tal que $\mu_V(V) >0$.

 (d) Demostrar que existe una medida de probabilidad invariante $\nu_V$ erg\'{o}dica tal que $\nu_V(V) >0$.

 (e) Demostrar que existe una medida de probabilidad invariante $\rho$ positiva sobre abiertos. Sugerencia: Tomar una base numerable de abiertos $\{V_i\}_{i \geq 1}$, demostrar que las medidas  $\rho_n := \sum_{i= 1}^n (1/2^i) \mu_{V_i}$   (que no son probabilidades, pero son finitas) satisfacen $0 <\rho_n (X) \leq 1$, son $T$-invariantes y existe $\rho= \lim^*_{n \rightarrow + \infty} \rho_n$ en la topolog\'{\i}a d\'{e}bil$^*$ del espacio ${\mathcal M}^1$ de medidas de probabilidad finitas uniformemente acotadas por 1.

\end{exercise}


\section{Teoremas Erg\'{o}dicos.}
\label{chapterTeoBirkhoff}

Las hip\'{o}tesis generales para este cap\'{\i}tulo  son las siguientes:

 $(X, {\mathcal A})$ es un espacio medible, $T: X \mapsto X$ es una transformaci\'{o}n
 medible que preserva una medida de probabilidad $\mu $, y  $f: X \mapsto {\mathbb{R}}$ es una
 funci\'{o}n  medible.

Consideremos el siguiente resultado: \index{medida! invariante}

\em Sea $T: X \mapsto X$ medible en el espacio m\'{e}trico compacto $X$. Sea $\mu$ una medida de probabilidad en la sigma-\'{a}lgebra de Borel. Entonces $\mu$ es $T$-invariante si y solo si para toda $f: X \mapsto \mathbb{R}$ continua se cumple
\begin{equation} \label{eqn27}\int f \,d \mu  = \int f \circ T \, d \mu\end{equation}
Adem\'{a}s,    vale la igualdad \em (\ref{eqn27}) \em para toda $f$ continua si y solo si vale   para toda $f \in L^p(\mu)$, cualquiera sea el natural $p \geq 1$. \em

\begin{exercise}\em
Probar la afirmaci\'{o}n anterior.
\end{exercise}

\subsection{Teorema Erg\'{o}dico de Birkhoff-Khinchin  - Enunciado y Corolarios}

\begin{definition} {\bf Promedios temporales }  \index{promedio! temporal} \label{definitionPromediosBirkhoff}   \index{promedio! de Birkhoff}
\em
Sea $(X, {\mathcal A})$ un espacio medible.
 Sea $T: X \mapsto X$ medible que preserva una
medida de pro\-ba\-bi\-lidad $\mu$ (no necesariamente erg\'{o}dica). Sea $f
\in L^p (\mu)$ para $1 \leq p \in \mathbb{N}$.

Se denota con $\widetilde {f}^+ $ o simplemente con $\widetilde f $ al
l\'{\i}mite de los llamados \em promedios orbitales (o promedios temporales o promedios de Birkhoff) hacia
el futuro \em de $f$, en los puntos $x \in X$ donde exista, esto
es:
$$\widetilde f^+ (x)= \lim _{n \rightarrow + \infty} \frac{1}{n} (f (x)+ f \circ T (x) +
 \ldots f \circ T^{n-1})$$

Si adem\'{a}s $T$ es invertible, se denota con $\widetilde f^-  $ al
l\'{\i}mite de los promedios orbitales (o temporales o  de Birkhoff) hacia el pasado, en
los puntos $x$ donde exista. Esto es:
$$\widetilde f^- (x) = \lim _{n \rightarrow + \infty} \frac{1}{n} (f (x)+
 f \circ T^{-1}(x) + \ldots f \circ T ^{-(n-1)}(x))$$
\end{definition}

\begin{theorem}
  {\bf Teorema erg\'{o}dico de Birkhoff-Khinchin}  \index{teorema! Birkhoff-Khinchin} \label{theoremBirkhoff} \index{teorema! erg\'{o}dico}

Sea $(X, {\mathcal A})$ un espacio medible.
Si $T: X \mapsto X$ es medible que preserva una medida de
probabilidad $\mu$ entonces:

\begin{itemize}
\item [a) ] Para toda $f \in L^1(\mu )$ existe $\widetilde f (x)= \lim _{n \rightarrow \infty}
\frac{1}{n} \sum _{j= 0}^{n-1} f \circ T ^j (x)$ $\mu$-c.t.p. en $x
\in X$.

\item[b) ] $\widetilde  f$ es $T$-invariante, es decir : $\widetilde f \circ
T = \widetilde f \; \mu$-c.t.p. M\'{a}s precisamente, para todo $x \in X$
existe $\widetilde f (x)$ si y solo si existe $\widetilde f
(T(x))$ y en ese caso $\widetilde f(x) = \widetilde f(T(x))$.

\item[c) ]Para todo natural $p \geq 1$, si $f \in L^p (\mu ) \subset L^1(\mu)$,
entonces $\widetilde f \in L^p(\mu )$ y la convergencia es tambi\'{e}n
en $L^p(\mu )$.
\item [d) ] $\int  f \, d \mu = \int \widetilde f \, d \mu$
\end{itemize}
\end{theorem}

   La demostraci\'{o}n de Birkhoff  del Teorema \ref{theoremBirkhoff} se encuentra en \cite{Birkhoff}. Otra demostraci\'{o}n, diferente de la de Birkhoff,    puede encontrarse  por ejemplo,  en
  \cite[proof of Theorem 1.14, pag. 38-39]{Walters}, en
  \cite[Theorem 2.1.5]{Keller} o en
  \cite[p\'{a}g. 114-122]{Mane} {Mane} (ver tambi\'{e}n \cite{ManeIngles}).

  El Teorema Erg\'{o}dico de Birkhoff-Khinchin es un caso particular del llamado {\bf Teorema Erg\'{o}dico Subaditivo de Kingmann} \cite{Kingmann}, que establece la convergencia de la sucesi\'{o}n $\{f_n / {n}\}_n$ donde $f_n$, en vez de ser necesariamente una suma de Birkhoff, es \em una sucesi\'{o}n subaditiva de funciones \em en $L^1(\mu)$. M\'{a}s precisamente, se asume por hip\'{o}tesis, que $$f_{n+m} \leq f_n + f_m \circ T^n \ \ \ \forall \ \ n, m \in \mathbb{N}^+,$$ donde $f_1 \in L^1(\mu)$. En particular, $f_n = \sum_{j= 0}^{n-1} f_1 \circ T^j$ es un ejemplo de sucesi\'{o}n subaditiva.

  El Teorema Erg\'{o}dico Subaditivo de Kingmann establece que, para toda medida $\mu$ que sea invariante por $T$, para toda sucesi\'{o}n $\{f_n\}_{n \geq 1}$ subaditiva de funciones reales en $L^1(\mu)$, \em la sucesi\'{o}n $\{f_n/n\}_n$ converge $\mu$-c.t.p. \em Luego, este Teorema generaliza el Teorema Erg\'{o}dico de Birkhoff-Khinchin. La demostraci\'{o}n del Teorema Subaditivo de Kingmann se encuentra en \cite{Kingmann}. Otra demostraci\'{o}n, que no usa  el Teorema Erg\'{o}dico de Birkhoff-Khinchin, y por lo tanto puede sustituir la demostraci\'{o}n de este \'{u}ltimo teorema, se encuentra en \cite{Bochi_KingmannNotes}.


\begin{corollary} \em \label{corollary1}
{\bf Igualdad de los promedios temporales  hacia el futuro y hacia el pasado. } \index{promedio! temporal}

\em
Si $T$ es medible, invertible, con inversa medible, y preserva una
medida de probabilidad $\mu$, entonces para toda $f \in L^1(\mu )$
se cumple $\widetilde f^+$ \em (promedio temporal hacia el futuro) \em y $\widetilde f^-$ \em (promedio temporal hacia el pasado) \em  existen $\mu$-c.t.p.
y son iguales $\widetilde f^+ = \widetilde f^- \;\mu$-c.t.p. \em

(Ver la definici\'{o}n de $\widetilde f^+$ y $\widetilde f^-$ en \ref{definitionPromediosBirkhoff}.)

En el siguiente Ejercicio \ref{exerciseBirkhoffPatras}, se da una gu\'{\i}a para la demostraci\'{o}n del Corolario \ref{corollary1}.
\end{corollary}

\begin{exercise}\em \label{exerciseBirkhoffPatras}

 Probar el corolario \ref{corollary1} como consecuencia del Teorema de Birkhoff-Khinchin.  Sugerencias: Para demostrar que existen $\widetilde f^+(x)$ y $\widetilde f^-(x)$ para $\mu$-c.t.p. $x \in X$, aplicar el teorema de Birkhoff a $T^{-1}$, probar que toda medida $\mu$ es invariante por $T$ si y solo s\'{\i} lo es por $T^{-1}$,  y usar que la intersecci\'{o}n de dos conjuntos con $\mu$-medida igual a $1$ tiene medida $\mu$-medida $1$. Para probar que $\widetilde f^+(x) = \widetilde f^-(x)$ para $\mu$-c.t.p. $x \in X$, para cada natural $n \geq 1$ denote $f^+_n := (1/n) \sum_{j= 0}^{n-1} f \circ T^j, \ f^-_{n} := (1/n) \sum_{j= 0}^{n-1} f \circ T^{-j}$. Por el Teorema de Birkhoff,   $f^-_n - f^+_n$ converge en $L^1(\mu)$ a $\widetilde f^- - \widetilde f^+$. Entonces para todo $\epsilon >0$ existe $n$ suficientemente grande tal que $\|\widetilde f^- - \widetilde f^+\|_{L^1} \leq \|f^-_n - f^+_n\|_{L^1} + \epsilon$. Basta probar entonces que $\int |f^-_n - f^+_n| \, d \mu$ tiende a cero cuando $n \rightarrow + \infty$. Chequear que para todo $x \in X$ se cumple la siguiente igualdad: $$f^-_{n+1}(x) - f^+_{n+1}(x) = f_{n+1}^-(x) + f^+_{n+1}(x) - 2 f^+_{n+1}(x) = $$ $$=\frac{2n+1}{n} \, f^+_{2n+1}(T^{-n}(x)) \ \ +  \frac{f(x)}{n+1} \ \  - \ \ 2 f^+_{n+1}(x)$$
 Como $f^+_{2n+1}$ converge a $\widetilde f^+$ en $L^1(\mu)$, y adem\'{a}s la medida $\mu$ y la funci\'{o}n $\widetilde f^+$ son invariantes por $T$, tenemos, para todo $n $ suficientemente grande:
  $$\int |f^+_{2n+1} - \widetilde f^+| \, d \mu < \epsilon$$
 $$\int |f^+_{2n+1} - \widetilde f^+| \, d \mu = \int |f^+_{2n+1}\circ T^{-n} - \widetilde f^+ \circ T^{-n}| \, d \mu = $$ $$=\int |f^+_{2n+1} \circ T^{-n}- \widetilde f^+| \, d \mu  < \epsilon$$
 Juntando todo lo anterior, deducir, para todo $n$ suficientemente grande, que: $$\|f^-_{n+1}  - f^+_{n+1} \|_{L^1(\mu)} \leq $$ $$\frac{2n+1}{n+1} \ \|f^+_{2n+1} \circ T^{-n} - \widetilde f^+  \|_{L^1(\mu)} \ \ + \ \ \frac{|f(x)|}{n+1} \ \ + \ \ 2 \, \| f^+_{n+1} - \widetilde f^+\|_{L^1(\mu)} $$ $$\leq 3\epsilon + \epsilon + 2\epsilon  = 6 \epsilon.$$
\end{exercise}

\begin{corollary} \label{corolarioPromediosTransitividad}
{\bf Promedios de medida de transitividad.} \index{promedio! de medida de transitividad}
Para toda medida $\mu$ invariante por $T$  y para todos los conjuntos $A$ y $B$ medibles, existe el l\'{\i}mite siguiente:
$$\tau(A, B) = \lim _{n \rightarrow + \infty} \frac{1}{n} \sum _{j=0}^{n-1}
\mu (T^{-j}(A) \cap B)$$
\end{corollary}
{\em Demostraci\'{o}n: } Denotemos $\chi_C$ a la funci\'{o}n caracter\'{\i}stica de cualquier conjunto $C$. Sea $$I_n = \frac{1}{n} \sum_{j= 0}^{n-1}
\mu(T^{-j}(A) \cap B) = $$ $$\frac{1}{n} \sum_{j= 0}^{n-1}\int
\chi_{T^{-j}(A)} \chi _B \, d \mu =  \int ( \frac{1}{n} \sum_{j=
0}^{n-1} \chi_A \circ T^j) \chi_B \, d \mu, $$ El integrando a la
derecha de la igualdad anterior est\'{a} dominado por $1 \in L^1(\mu
)$. Por el teorema de convergencia dominada $\lim I_n = \int
\widetilde \chi_A \chi _B \, d \mu \; \; \; \; \Box$

\begin{exercise}\em
Probar que el l\'{\i}mite $\tau(A, B)$ del Corolario \ref{corolarioPromediosTransitividad} verifica las siguientes desigualdades
$$\mu(A)- \sqrt{\mu(A) [1 - \mu(B)]} \leq \tau(A, B) \leq \sqrt{\mu(A) \mu(B) }.$$
Concluir que:

 $\tau(A, B) = 0$ si $\mu(A) = 0$ \'{o} $\mu(B) = 0$.

 Si $\tau(A, B) = 0$ entonces $\mu(A) + \mu(B) \leq 1$

 $\tau(A, B) = 1$ si y solo si $\mu(A) = \mu(B)= 1$.

 Sugerencia para la primera parte: En la demostraci\'{o}n del Corolario \ref{corolarioPromediosTransitividad} se prob\'{o} que $\tau(A, B) = \int \widetilde \chi_A \chi _B, d \mu$. En $L^2(\mu)$ se cumple
 $$\int f g \, d \mu \leq \|f\|_{L^2} \, \|g\|_{L^2}.$$ Para probar la desigualdad de la derecha, aplicar lo anterior a $\widetilde \chi_A$ y $\chi _B$, y recordar que $0 \leq \widetilde \chi_A \leq 1$, por lo cual $\widetilde \chi_A^2 \leq \widetilde \chi_A$. Deducir que $\|\widetilde \chi_A\|_{L^2} \leq \sqrt {\mu(A)}$ aplicando el teorema de Birkhoff. Para probar la   desigualdad de la izquierda, aplicar la desigualdad de la derecha a $A$ y $B^c$ y probar que $\tau(A, B) + \tau (A, B^c) = \mu(A)$.
\end{exercise}

\begin{corollary} \em {\bf Promedios  de sucesi\'{o}n de funciones.}
  \label{corollary3} \index{promedio! de sucesi\'{o}n de funciones}

\em Sea $T: X \mapsto X$ medible que preserva una medida de
probabilidad $\mu$. Sea $f_n \in L^1(\mu )$ una sucesi\'{o}n de
funciones, dominada por $f_0 \in L^1(\mu)$, que converge
$\mu$-c.t.p. y en $L^1(\mu )$ a $f \in L^1(\mu )$. Entonces
$$\lim _{n \rightarrow + \infty} \frac{1}{n}\sum _{j=0}^{n-1} f_{j} \circ T^j = \widetilde f \ \ \
\mu \mbox{-c.t.p. y en }L^1(\mu ). $$
\end{corollary}
 \begin{exercise}\em
   Demostrar el corolario \ref{corollary3}. Sugerencia: Basta
 probarlo para $f_n \geq 0; f_n \rightarrow 0 \; c.t.p.$ Sea $G_k (x) = \sup _{n \geq k}
 \{f_n(x)\}$. Entonces $G_k \rightarrow 0\;  c.t.p.$, y por
 convergencia dominada $\|G_k\|_{L^1} \rightarrow 0$. Sea $\widetilde
 G_k = \lim _{n \rightarrow \infty}(1/n) \sum _{j= 0}^{n-1} G_k \circ
 T^j$ que existe $\mu$-c.t.p. por el teorema de Birkhoff. La sucesi\'{o}n $\widetilde G_k$ es decreciente con $k$, por lo que tiene l\'{\i}mite, y por el
 lema de Fatou $$0 \leq \int \lim \widetilde G_k \, d \mu \leq \lim
 \int \widetilde G_k \, d \mu = \lim \int G_k \, d \mu = 0.$$ Luego
 $\widetilde G_k \rightarrow 0 \; c.t.p.$. Finalmente, usar que
  $$\lim \sup _{n \rightarrow \infty} \frac{1}{n}
 \sum _{j=0}^{n-1} f_j \circ T^j (x) = \lim \sup _{n \rightarrow \infty} \frac{1}{n}
 \sum _{j= k}^{n-1} f_j \circ T^j (x) \leq
   \widetilde G_k(x).$$
 \end{exercise}
\begin{definition} \em \index{tiempo medio de estad\'{\i}a} \index{promedio! de tiempo de estad\'{\i}a}
Dado $A$ conjunto medible, se denomina  \em tiempo medio de
estad\'{\i}a $\tau _A(x)$ de un punto $x \in X$ en $A$, \em al
siguiente l\'{\i}mite, cuando existe:
$$\tau _A(x) =
\lim _{n \rightarrow + \infty} \frac{1}{n} \# \{j \in \mathbb{N}: 0 \leq j
\leq n-1,\;  T^j(x) \in A\}$$
\end{definition}

\begin{exercise}\em
Probar el siguiente teorema, para todo conjunto me\-di\-ble $A \subset X$ y para toda medida $\mu$ que sea invariante por $T: X \mapsto X$:

\em El tiempo medio de estad\'{\i}a $\tau_A $ existe $\mu $-c.t.p. y adem\'{a}s
la convergencia es en $L^p(\mu )$ para todo $p \geq 1$ natural.
Adem\'{a}s $\int \tau _A \, d \mu = \mu (A)$. \em  Sugerencia:
Observar que el  tiempo medio de estad\'{\i}a en $A$ es la funci\'{o}n
$\tau_A = \widetilde \chi _A$, donde $\chi_A$ es la funci\'{o}n caracter\'{\i}stica de $A$. Aplicar el teorema erg\'{o}dico de
Birkhoff-Khinchin.
\end{exercise}

\begin{definition} \em {\bf Conjuntos de probabilidad total para $T$} \label{definitionProbabilidadTotal} \index{probabilidad! total} \index{conjunto! de probabilidad total}
Sea $T:X \mapsto X$ una transformaci\'{o}n medible tal que el conjunto ${\mathcal M}_T$ de las medidas de probabilidad $T$-invariantes es no vac\'{\i}o.  Un conjunto  medible $\Lambda \subset X$ se dice que tiene \em probabilidad total para $T$ \em si para toda $\mu \in {\mathcal M}_T$ se cumple $\mu(\Lambda)= 1$.

Entonces, la primera parte del ejercicio anterior se puede enunciar de la siguiente forma:

\em Para todo conjunto medible $A \subset X$, el tiempo medio de estad\'{\i}a $\tau_A(x)$ existe para un conjunto de puntos  con probabilidad total. Si adem\'{a}s $\mu(A) >0$, entonces $\tau_A(x)$ no nula $\mu$-c.t.p. \em (pues $\int \tau_A \, d \mu = \mu(A)$).
\end{definition}

En el Corolario \ref{CorolarioDescoErgodicaEspaciosMetricos} veremos un criterio para que un conjunto tenga probabilidad total, que requiere solo el conocimiento de las medidas erg\'{o}dicas.

\begin{exercise}\em \label{exercise4} \label{exerciseLimiteBirkhoffConjuntoEstable} {\bf Conjuntos estables e inestables.} \index{conjunto! estable} \index{conjunto! inestable}

  Sea $T: X \mapsto X$ Borel medible e invertible con inversa
medible, en un espacio m\'{e}trico compacto $X$, tal que el conjunto de medidas invariantes por $T$ es no vac\'{\i}o. Sea $f: X \mapsto \mathbb{C} $ una
funci\'{o}n compleja continua. Sea $\Lambda \subset X$ el conjunto de probabilidad total tal que existen, y son iguales entre s\'{\i},  los l\'{\i}mites  $\widetilde f ^+ (x)$ y $\widetilde f^-(x)$ de los promedios de Birkhoff hacia el futuro  y hacia el pasado, respectivamente. Sea $x_0 \in \Lambda$.  Se
definen los conjuntos estable e inestable respectivamente por el
punto $x_0$ (quiz\'{a}s se reducen solo a $\{x_0\}$):
$$W^s(x_0) := \{y \in X: \lim _{n \rightarrow + \infty} \dist (T^ny, T^n x_0) = 0\}$$
$$W^u(x_0) := \{y \in X: \lim _{n \rightarrow - \infty} \dist (T^ny, T^n x_0) = 0\}$$

  Probar que para todo $y \in W^s(x_0)$ existe $\widetilde f^+(y)$ y $\widetilde f^+(y) =
\widetilde f^+(x_0) $. Probar que para todo $y \in W^u(x_0)$ existe $\widetilde f^-(y)$ y $\widetilde f^-(y) = \widetilde f^-(x_0). $

\end{exercise}

\begin{exercise}\em {\bf Promedios de Birkhoff para funciones reales fuera de $L^1(\mu)$.}

  Probar la siguiente generalizaci\'{o}n del teorema de
Birkhoff-Khinchin:
Sea $(X, {\mathcal A})$ un espacio medible y sea
  Sea $T: X \mapsto X$ medible que preserva la probabilidad
$\mu$. Sea $f : X \mapsto \mathbb{R}$ medible.

\em Entonces, $\mu$-c.t.p. o bien $\widetilde{| f |} \; (x) = + \infty $ o bien $\widetilde f$ existe y es finito. \em

Sugerencia:
Basta probarlo para $f \geq 0$. Dado $c >0$
 sea $$X_c = \{ x \in X: \liminf_{n \rightarrow + \infty} \frac{1}{n} \sum _{j= 0}^{n-1} f \circ T^j (x) \leq c
 \}$$
 $X_c$ es $T$- invariante. Sea $f_m (x) = \min (f(x), m)$.  Usar el teorema de Birkhoff para probar que existe $\widetilde f_m$ $\mu$-c.t.p. $x \in X_c$. Sea $f|_{X_c} = \chi_{X_c} f$ donde $\chi_{X_c}$ denota la funci\'{o}n caracter\'{\i}stica de $X_c$. Probar que $f|_{X_c} \in L^1(\mu)$ usando el teorema de convergencia mon\'{o}tona y la igualdad del teorema de Birkhoff $\int f_m|_{X_c} \, d \mu  = \int \widetilde {f_m|_{X_c}}\, d \mu \leq c $. Deducir que existe el l\'{\i}mite $\widetilde f (x)$ para $\mu$-c.t.p. $x$ en $X_c$. Tomar la uni\'{o}n de los $X_c$  para todo $c \geq 1$ natural y concluir que en el complemento de esa uni\'{o}n se cumple $ \widetilde {f } = + \infty$.

\end{exercise}

\begin{exercise}\em Sea $T: X \mapsto X$ una transformaci\'{o}n medible en un espacio m\'{e}trico compacto $X$, tal que es no vac\'{\i}o el conjunto ${\mathcal M}_T$ de probabilidades invariantes por $T$.
Probar que el conjunto siguiente tiene probabilidad total para $T$ (i.e. tiene probabilidad 1 para toda   $\mu \in {\mathcal M}_T$):
$$ \big \{ x \in X: \mbox{ existe l\'{\i}m} _{n \rightarrow +\infty}
\frac{1}{n} \cdot \sum _{j=0} ^{n-1} f \circ T ^{p j} (x) \: \: \:
\forall \: f \in C^0(X, \mathbb{R}), \: \: \forall \: p \geq 1 \big \} $$
Sugerencia: ${\mathcal M} _{T^p} (X) \supset {\mathcal M} _T (X) $.
\end{exercise}

\subsection{Otras caracterizaciones de la ergodicidad}

 En esta secci\'{o}n, salvo indicaci\'{o}n en contrario, $(X, {\mathcal A})$ denota un espacio medible y $T: X \mapsto X$ una transformaci\'{o}n medible que preserva alguna medida de probabilidad $\mu$.

Se recuerda las Definiciones \ref{definitionErgodicidadI} y \ref{definitionErgodicidadII} de ergodicidad, y los Teoremas \ref{teoremaDefinicionesErgodicidad}, \ref{TeoremaErgodicidadIII} y \ref{teoremaErgodicasExtremales}, en los que dimos diferentes caracterizaciones de ergodicidad. Agregamos ahora las siguientes:

\begin{theorem} \label{teoremaergodicidad} {\bf Ergodicidad IV}

Sea $T: X \mapsto X$ medible que preserva una medida de
probabilidad $\mu$. Las siguientes propiedades son equivalentes:

 {\bf a) } $T$ es erg\'{o}dica respecto de $\mu$.

{\bf b) } Para toda $f \in L^1(\mu ) $ se cumple
$$\widetilde f (x) = \int f \, d \mu \; \; \; \mu -\mbox{c.t.p.}$$

{\bf c) } Para toda pareja de conjuntos medibles $A$ y $B$ se
cumple:
$$\lim _{n \rightarrow + \infty} \frac{1}{n} \sum _{j= 0}^{n-1} \mu (T^{-j}(A) \cap B)
= \mu (A) \mu (B).$$
\end{theorem}
 El enunciado y la prueba del Teorema \ref{teoremaergodicidad} fueron extra\'{\i}dos, con leves modificaciones, de \cite[p\'{a}gs. 130-131]{Mane} (ver tambi\'{e}n \cite{ManeIngles}).

\vspace{.2cm}

 {\bf Demostraci\'{o}n de que a) $\Rightarrow $ b) en el
Teorema \ref{teoremaergodicidad}:} Por el Teorema de Birkhoff-Khinchin $\int
\widetilde f \, d \mu = \int f \, d \mu$. Entonces basta demostrar
que $\widetilde f$ es constante $\mu-$c.t.p.

Por el teorema de Birkhoff-Khinchin $\widetilde f(x)$ existe
$\mu$-c.t.p., $\widetilde f(x)$ existe si y solo si existe
$\widetilde f(Tx)$ y en ese caso $\widetilde f(Tx) = \widetilde
f(x)$. Dicho de otra forma, el conjunto $A= \{x \in X: \widetilde
f (x) \mbox { existe } \}$ cumple $\mu(A)= 1$, $T^{-1}(A) = A$ y
para todo $x \in A: \; \; \widetilde f(x) = \widetilde f(Tx)$.

Sea $g: X \mapsto \mathbb{C} $  definida como $g(x) = \widetilde f(x)$ si $x
\in A$, $g(x) = 0$ si $x \not \in A$. Entonces para todo $x \in X$
se cumple $g \circ T (x)= g(x)$. Por la afirmaci\'{o}n demostrada
antes $g = cte\; \; \; \mu-$c.t.p. Pero por construcci\'{o}n $g=
\widetilde f \; \; \; \mu-$c.t.p., de donde se deduce que
$\widetilde f = cte \; \; \; \mu-$c.t.p. como quer\'{\i}amos. \hfill $\Box$

\vspace{.2cm}

{\bf Demostraci\'{o}n de que b) $\Rightarrow $ c) en el Teorema
\ref{teoremaergodicidad}:}

 Obs\'{e}rvese que $\mu (T^{-j}A \cap
B) = \int \chi_{T^{-j}(A)}\chi_B \, d \mu = \int (\chi_A \circ
T^j) \chi _ B \, d \mu$, de donde:
$$I_n= \frac{1}{n} \sum _{j=0}^{n-1} \mu (T^{-j}A \cap
B) = \int (\frac{1}{n} \sum _{j=0}^{n-1} \chi_A \circ T^j) \chi _
B \, d \mu$$ Por convergencia dominada (el integrando est\'{a} acotado
por la funci\'{o}n constante $1 \in L^1(\mu )$), resulta:
\begin{equation} \label{eqn1} \lim I_n = \int \lim _{n \rightarrow + \infty}
(\frac{1}{n}\sum _{j=0}^{n-1} \chi_A \circ T^j) \chi_ B \, d \mu =
\int (\widetilde \chi _A) \chi _B \, d \mu \end{equation} Por
  el teorema de Bikhoff-Khinchin, observando
que $\chi _A \in L^p(\mu )$ y $\widetilde \chi _A $ es invariante
con $T$, luego contante $\mu-c.t.p$, se tiene que
$$\widetilde \chi_A (x) = \int \widetilde \chi_A \, d \mu = \int \chi_A \, d \mu = \mu (A)\; \; \; \mu-c.t.p. x \in X$$
Sustituyendo en (\ref{eqn1}) resulta \
$\lim_n I_n = \int (\mu (A)) \chi _B \, d \mu = \mu (A) \mu (B). $ \hfill $\Box$

 \vspace{.2cm}

{\bf Demostraci\'{o}n de que c) $\Rightarrow $ a) en el Teorema
\ref{teoremaergodicidad}:} Basta demostrar que si $A, B
\subset X$ son medibles con $\mu$ medida posi\-tiva, entonces existe
$j \geq 1$ tal que $\mu (T^{-j}(A) \cap B) >0$ (cf. Definici\'{o}n \ref{definitionErgodicidadI}). Por hip\'{o}tesis:
$$\lim \frac{1}{n} \sum _{j=1}^{n-1} \mu (T^{-j}(A) \cap B) =
\lim \frac{1}{n} \sum _{j=0}^{n-1} \mu (T^{-j}(A) \cap B) = \mu
(A) \mu (B) >0
$$ Entonces  $\sum _{j=1}^{n-1} \mu (T^{-j}(A) \cap B) >0 \ \  \forall \ n $ suficientemente
grande, de donde $ \mu (T^{-j}(A) \cap B)>0 $ para alg\'{u}n $j
\geq 1$ como quer\'{\i}amos probar. \hfill $\Box$

\begin{exercise}\em {\bf Ergodicidad V.} \label{ejercicio0} \index{medida! erg\'{o}dica} \index{ergodicidad} \index{transformaci\'{o}n! erg\'{o}dica}

  Sea $T: X \mapsto X$ que preserva una probabilidad $\mu$.
Probar que:

{\bf (a) } \em $T$ erg\'{o}dica respecto de $\mu $ si y solo si para todo 
conjunto $A$ medible tal que $T^{-1}(A) \subset A$ \'{o} $T^{-1}(A) \supset A$, se cumple $\mu
(A)$ es o bien cero o bien uno. \em

{\bf (b) } \em $T$ es erg\'{o}dica respecto de $\mu$ si y solo si para toda $f \in
L^1(\mu)$ tal que \em $f \circ T \leq f\;   \mu$-c.t.p., \em se cumple \em $f =
cte \ \ \ \ \mu-$ c.t.p. \em

{\bf  (c) }  $T$ es erg\'{o}dica respecto de $\mu$ si y solo si para todos los conjuntos $A$ medibles que cumplen $T^{-1}(A) \subset A$ \'{o} $A \subset T^{-1}(A)$ se verifica $\mu(A) = 0$ \'{o} $\mu(A)= 1$. \em

{\bf (d) } \em $T$ es erg\'{o}dica respecto de $\mu$ si y solo si para todas las funciones $f \in
L^1(\mu)$ tales que \em $f \circ T \leq f\;   \mu$-c.t.p., \ \'{o} \  $f \circ T \geq f\;   \mu$-c.t.p.\em se cumple \em $f =
cte \ \ \ \ \mu-$ c.t.p.

Sugerencia para (a): Basta probarlo cuando $T^{-1}(A) \supset A$, pues en caso contrario, sustituimos $A$ por su complemento. Denote $\chi_{A}$ a la funci\'{o}n caracter\'{\i}stica de $A$. Como $T^{-1}(A) \supset A$, pruebe que  $\chi_{T^{-n}(A)}(x) = 1$ para todo $x \in   A$. Luego $\widetilde \chi_{A}(x) = 1$ para todo $x \in A$, y usando que $\chi_A$ es constante $\mu$-c.t.p.  si $\mu$ es erg\'{o}dica, se deduce $\mu(A)= 1$ \'{o} $\mu(A)= 0$.
\end{exercise}

\begin{exercise}\em \label{ejercicio513}
Sea $X$ un espacio m\'{e}trico compacto y sea $T: X \mapsto X$ Borel-medible tal que preserva una probabilidad $\mu$.
Sea $\{g_i: i\geq 1 \} $ un conjunto numerable denso en
  $C^0(X, [0,1])$.
 Probar  que son equivalentes
las afirmaciones siguientes:

i)
$\mu $ es erg\'{o}dica.

ii)
$\widetilde  f (y) = \int f \mbox{d} \mu \: \: \: \mu\mbox{-c.t.p.} \; y;
\: \: \forall \: f \in C_0 (X, \mathbb{R}) $

iii)
$\widetilde  g_i (y) = \int g_i \mbox{d} \mu \: \: \: \mu \mbox{-c.t.p.} \; y;
\: \: \forall \: i \geq 1 $

Sugerencia: Recordar que $\mu $ es erg\'{o}dica si y solo si
para toda $h \in L^1(\mu ): \widetilde  h (y ) = \int h \mbox{d} \mu $ ,  $ \mu \mbox{-c.t.p.} \; y $; y que las funciones continuas son densas en $L^1 (\mu ) $.

\end{exercise} \index{medida! erg\'{o}dica} \index{ergodicidad} \index{transformaci\'{o}n! erg\'{o}dica} \index{equivalencia de definiciones! de ergodicidad}
\begin{exercise}\em \em \label{exercise1} \em
 Sea $X$ un espacio m\'{e}trico compacto y sea $T: X \mapsto X$ Borel-medible tal que preserva una probabilidad $\mu$. Probar que $T$ es erg\'{o}dica respecto de
$\mu$ si y solo si para toda funci\'{o}n compleja $f: X \mapsto \mathbb{C} $
continua,  el l\'{\i}mite $\widetilde  f$ de los promedios de Birkhoff de
$f$ es constante $\mu$-c.t.p. Sugerencia: Usar lo probado en el ejercicio \ref{ejercicio513} y chequear, usando el Teorema de Birkhoff-Khinchin, que si $\widetilde f$ es una constante $\mu$-c.t.p., entonces esta constante es $\int f \, d \mu$.

\end{exercise}

Volvamos al caso general de un espacio medible $(X, {\mathcal A})$ con una transformaci\'{o}n medible $T: X \mapsto X$ que preserva una medida de probabilidad $\mu$.
Por el Lema de recurrencia de Poincar\'{e} (Teorema \ref{teoPoincare}) si un conjunto medible $A$ cumple $\mu(A) >0$, entonces la \'{o}rbita futura de $\mu$-c.t.p. $x \in A$ vuelve infinitas veces a $A$. Sin embargo, ese Lema no dice nada sobre la frecuencia de visita y la duraci\'{o}n de las estad\'{\i}as de la \'{o}rbita de $x$ en el conjunto $A$. Las medidas erg\'{o}dicas dan exactamente el valor de la frecuencia asint\'{o}tica en que la \'{o}rbita de $x$ pasa dentro de $A$.

\begin{definition} \em {\bf Tiempo medio de estad\'{\i}a.} \index{tiempo medio de estad\'{\i}a} \index{promedio! de tiempo de estad\'{\i}a}

 Llamamos \em tiempo medio de estad\'{\i}a \em $\tau _A (x)$ de la
\'{o}rbita por $x \in X$ en un conjunto medible $A$ a:

$$\tau _A(x) = \lim _{n \rightarrow + \infty} \frac{\#\{ j \in \mathbb{N}: 0 \leq j \leq n-1, T^j(x)
 \in A
 \}}{n} = $$ $$ \lim _{n \rightarrow + \infty} \frac{1}{n} \sum _{j= 0}^{n-1} \chi _A \circ T^j (x)=
 \widetilde \chi _A (x)$$
\end{definition}

\vspace{.3cm}

\begin{theorem} {\bf Ergodicidad VI.} \index{medida! erg\'{o}dica} \index{ergodicidad}   \index{equivalencia de definiciones! de ergodicidad} \index{transformaci\'{o}n! erg\'{o}dica}

Sea $T: X \mapsto X$ medible que preserva una medida de probabilidad $\mu$. Entonces  $\mu$ es erg\'{o}dica
para $T$ si y solo si para todo conjunto medible $A$   el
tiempo medio de estad\'{\i}a $\tau_A (x)$ es constante $\mu$-c.t.p. Adem\'{a}s, en ese caso $\tau _A(x) = \mu (A)\; \;
 \; \mu-c.t.p.$
\end{theorem}
 El enunciado y la prueba de este teorema, con leves modificaciones, fue extra\'{\i}do de \cite[p\'{a}g. 133]{Mane} (ver tambi\'{e}n \cite{ManeIngles}).

{\em Demostraci\'{o}n:} Por el teorema de Birkhoff-Khinchin
$\widetilde \chi _A \in L^1(\mu)$ y cumple $\widetilde \chi_A =
\widetilde \chi_A \circ T \; \; \mu-c.t.p.$
Si $\widetilde \chi _A = cte \; \; \mu-c.t.p.$ entonces, aplicando
nuevamente el teorema de Bikhoff-Khinchin:
\begin{equation} \label{eqn31}\widetilde \chi_A(x) =
\int \widetilde \chi _A \, d \mu = \int \chi _A \, d \mu = \mu
(A)\; \; \; \; \mu-c.t.p. \end{equation}
Si $\tau_A = \widetilde \chi _A = cte \; \; \mu-c.t.p.$ para todo
$A $ medible, tomemos en particular $A$ tal que $T^{-1}(A) = A$.
Para demostrar la ergodicidad de $\mu$ hay que probar que $\mu (A)
$ es cero o uno.
$$\widetilde \chi_A = \lim \frac{1}{n} \sum _{j= 0}^{n-1} \chi_A \circ T^j$$
Pero $\chi _A \circ T^j = \chi _{T^{-j}(A)} = \chi _A$. Luego
resulta:
\begin{equation} \label{eqn32}\widetilde \chi _A(x) = \chi _A(x) \in \{0,1\} \; \; \; \mu-\mbox{c.t.p. }x \in X \end{equation}

Por (\ref{eqn31}) y (\ref{eqn32}) se tiene $\mu (A) \in \{0,1\}$ y $\mu$ es erg\'{o}dica
como quer\'{\i}amos probar.

Rec\'{\i}procamente, si $\mu$ es erg\'{o}dica, entonces por la parte b) del
teorema \ref{teoremaergodicidad}. $\tau _A = \widetilde \chi _A = cte \; \;
\mu-c.t.p.$ \hfill $ \Box$

\begin{exercise}\em \label{exercise2}  \index{medida! erg\'{o}dica} \index{transformaci\'{o}n! erg\'{o}dica} \index{equivalencia de definiciones! de ergodicidad} \index{ergodicidad} \label{exerciseTildefLocalmenteConstante} \em
\em Sea $T:X \mapsto X$ es Borel medible en un espacio
topol\'{o}gico   $X$ conexo,   que preserva una
medida de probabilidad $\mu$ positiva sobre abiertos. Probar que:

(a) Una funci\'{o}n $g \in L^1 (\mu) $ invariante con $T$ (i.e. $g
\circ T = g \; \mu $-c.t.p.) es constante $\mu$-c.t.p. si y solo
si es localmente  constante c.t.p. (es decir: existe un
cubrimiento de $X$ por abiertos, tales que en cada abierto $V$ del
cubrimiento se cumple $g|_V = K_V $ constante $\mu$-c.t.p. de
$V$.)

(b) $T$ es erg\'{o}dica si y solo si toda funci\'{o}n $g \in L^1(\mu )$
que sea invariante con $T$  es localmente constante.
\end{exercise}

\subsection{Ergodicidad \'{U}nica} \index{medida! erg\'{o}dica} \index{transformaci\'{o}n! erg\'{o}dica} \index{ergodicidad \'{u}nica}
\index{transformaci\'{o}n! \'{u}nicamente erg\'{o}dica} \index{ergodicidad}

 \begin{definition}
\em
 La transformaci\'{o}n $T$ es \em \'{u}nicamente erg\'{o}dica  \em  si existe
una \'{u}nica medida de probabilidad que es $T$ invariante.

Por lo visto en la secci\'{o}n \ref{sectionExistenciaMedidasErgodicas} esta \'{u}nica medida
es extremal en el conjunto de las probabilidades invariantes; luego es erg\'{o}dica.

\em
\end{definition}

\begin{theorem} \label{ergouni} \label{teoremaErgodicidadUnica}
Sea $T$ continua en un espacio m\'{e}trico compacto.
Las siguientes afirmaciones son equivalentes:
\begin{description}
\item[i)]
$T$ es \'{u}nicamente erg\'{o}dica
\item[ii)]
Para toda $f \in C^0(X, \mathbb{R})$, para todo $x \in X $ existe el siguiente l\'{\i}mite y es un n\'{u}mero independiente de $x$: $$\widetilde  f := \lim_{n \rightarrow + \infty} \frac{1}{n} \sum_{j= 0}^{n-1} f \circ T^j(x).$$
\item[iii)]
Para toda $f \in C^0(X, \mathbb{R})$, la sucesi\'{o}n de funciones continuas $\{ f_n \} $
definidas por
$$ f_n  = \frac{1}{n} \sum _{j=0}^{n-1} f \circ T^j $$
converge uniformemente a una constante cuando $n \rightarrow \infty $.
\end{description}
\end{theorem}

{\em Demostraci\'{o}n:}
i) implica iii):

Sea $\mu $ la \'{u}nica medida de ${\mathcal M}_T(X)$. Probaremos
$f_n $ converge uniformemente a $\int f \, d \mu $ en $X$.
Por absurdo, supongamos que existe $\epsilon > 0 $, y una sucesi\'{o}n
$n_j \rightarrow \infty $ tal que $\sup _{x \in X} |f_{n_j} (x) - \int
f \, d \mu | \geq \epsilon $ para todo $j \geq 1$.
Como el supremo es un m\'{a}ximo, existe $x_j \in X$, donde se alcanza.

Por el teorema de Riesz, existe una medida $\mu _j$ (no necesariamente $T$- invariante) tal que:
$$ \forall \: g \in C^0(X, \mathbb{R}): \:\:\: \int g \, d \mu _j =
\frac{1}{n_j} \sum _ {i =0 }^{n_j -1} g \circ T^i (x_j) $$

Por la compacidad del espacio ${\mathcal M}(X)$, existe una subsucesi\'{o}n
convergente de esta medidas $\mu _j$. Por simplicidad seguiremos usando
la misma notaci\'{o}n para la subsucesi\'{o}n que para la sucesi\'{o}n original.

$$\mu _j \rightarrow  \nu $$ Afirmamos que $\nu $ es $T$- invariante.
En efecto:
$$ \int g \, d T^* \nu = \int g \circ T \, d \nu = \lim \int g\circ T \, d \mu _j = \lim \frac{1}{n_j} \sum _{ i=0}^{n_j -1} g \circ T^{i+1} (x_j) $$
$$ = \lim \frac{1}{n_j} \sum _{ i=0}^{n_j -1} g \circ T^{i} (x_j)
= \lim \int g \, d \mu _{n_j} = \int g \, d \nu $$

Como por hip\'{o}tesis $T$ es \'{u}nicamente erg\'{o}dica, se tiene $\mu = \nu $.
Entonces
$$ \int f \, d \mu = \int f \, d \nu = \lim f _{n_j}(x_j)$$
Pero por construcci\'{o}n
$$ |f_{n_j} (x_j) - \int f \, d \mu | \geq \epsilon \:\:\: \mbox{ para todo }
j \geq 1 $$ contradiciendo
la igualdad anterior.

 iii) implica ii) porque la convergencia uniforme  de funciones continuas en $X$ implica la convergencia
en todo punto de $X$.

 ii) implica i):
Sea $\Lambda $ el funcional lineal positivo definido por
$$\Lambda (f) = \widetilde  f $$ para toda $f \in C^0(X, \mathbb{R})$.

Para toda medida $\mu $ que sea $T$ invariante, por el teorema de
Birkhoff
$$ \int f \, d \mu = \int \widetilde  f \, d \mu = \widetilde  f = \Lambda (f) $$
Por el teorema de Riesz, existe una \'{u}nica $\mu $ que cumple
$$\int f \, d \mu = \Lambda (f) $$
Luego existe una \'{u}nica $\mu $ que es $T$- invariante.
\hfill $\Box$

 Como ejemplo, veremos en la pr\'{o}xima secci\'{o}n que la rotaci\'{o}n
irracional en el c\'{\i}rculo es \'{u}nicamente erg\'{o}dica.

\begin{definition} \em {\bf Conjuntos minimales} \index{conjunto! minimal} \label{definicionMinimalTopologico}

Sea $T$ una transformaci\'{o}n continua en un espacio m\'{e}trico
compacto $X$.

Un subconjunto $\Lambda \subset X$ es \em minimal \em (desde el punto de vista topol\'{o}gico) si es
compacto, no vac\'{\i}o, invariante hacia el futuro (es decir: $ \Lambda \subset T^{-1}(\Lambda)$), y no existe ning\'{u}n   subconjunto propio de
$\Lambda $ que sea compacto, no vac\'{\i}o e invariante hacia el futuro.

\vspace{.2cm}

Tenemos la siguiente caracterizaci\'{o}n (ver parta (a) del Ejercicio \ref{ejercicioMinimalTopologico}):

$\Lambda$ \em es minimal si y solo si \em
  $\Lambda $ es compacto, no vac\'{\i}o y $T$-invariante (es decir: $T^{-1}(\Lambda)= \Lambda$) y no contiene  subconjuntos propios que sean compactos, no vac\'{\i}os e invariantes por $T$ hacia el futuro.

  \vspace{.2cm}

  Si adem\'{a}s el mapa continuo $T$ es   invertible con inversa continua, entonces $\Lambda$ \em es minimal si y solo si \em
  $\Lambda $ es compacto, no vac\'{\i}o y $T$-invariante, y no contiene  subconjuntos propios que sean tambi\'{e}n compactos, no vac\'{\i}os y $T$-invariantes  (ver parte (b) del Ejercicio \ref{ejercicioMinimalTopologico}).

\vspace{.2cm}

Finalmente, es f\'{a}cil ver que cualquiera sea $T$ continua (no necesariamente invertible), un conjunto $\Lambda $
compacto, no vac\'{\i}o e invariante, \em es mi\-ni\-mal si y solo si  \em  todas
sus \'{o}rbitas hacia el futuro son densas en $\Lambda $. Esto es porque la clausura
de cada una de ellas es un compacto no vac\'{\i}o, invariante hacia adelante,
y contenido en $\Lambda $ (ver tambi\'{e}n parte (a) del Ejercicio \ref{ejercicioMinimalTopologico}).

\end{definition}

\begin{exercise}\em \label{ejercicioMinimalTopologico}
(a) Probar que las siguientes afirmaciones son equivalentes (y por lo tanto cualquiera de ellas puede utilizarse como definici\'{o}n de $\Lambda $ minimal):

(i)  $\Lambda \subset X$ es   compacto no vac\'{\i}o e invariante por $T$ hacia el futuro, y no contiene subconjuntos propios compactos no vac\'{\i}os e invariantes hacia el futuro.

(ii) $\Lambda \subset X$ es   compacto no vac\'{\i}o y $T$-invariante, y no contiene subconjuntos propios compactos no vac\'{\i}os e invariantes por $T$ hacia el futuro. (Sugerencia    para demostrar (i) $\Rightarrow$ (ii): probar que si $\Lambda $ cumple (i), entonces $T^{-1}(\Lambda) \setminus \Lambda$ tambi\'{e}n.)

(iii) $\Lambda \subset X$ es compacto no vac\'{\i}o, y para todo $x \in \Lambda$  la clausura  de $\{f^j(x)\}_{j \geq 0}  $ es igual a $\Lambda$.

(b) Asumir ahora que $T$ es un homeormorfismo (es decir, $T$  es continua, invertible y su $T^{-1}$ es continua). Probar que las siguientes afirmaciones son equivalentes (y por lo tanto cualquiera de ellas puede utilizarse como definici\'{o}n de $\Lambda$ minimal):

(i)     $\Lambda \subset X$ es   compacto no vac\'{\i}o e invariante por $T$ hacia el futuro, y no contiene subconjuntos propios compactos no vac\'{\i}os e invariantes hacia el futuro.

(ii) $\Lambda \subset X$ es   compacto no vac\'{\i}o y $T$-invariante, y no contiene subconjuntos propios compactos no vac\'{\i}os y $T$-invariantes.

(iii) $\Lambda \subset X$ es compacto no vac\'{\i}o, y para todo $x \in \Lambda$  la clausura  de $\{f^j(x)\}_{j \geq 0}  $ es igual a $\Lambda$.

(iv) $\Lambda \subset X$ es compacto no vac\'{\i}o   e invariante por $T$ hacia el pasado (es decir $T^{-1}(\Lambda) \subset \Lambda$), y no contiene subconjuntos propios compactos no vac\'{\i}os e invariantes hacia el pasado.

(v)   $\Lambda \subset X$ es compacto no vac\'{\i}o, y para todo $x \in \Lambda$  la clausura  de $\{f^{-j}(x)\}_{j \geq 0}  $ es igual a $\Lambda$.
\end{exercise}

Veamos como se vincula la ergodicidad \'{u}nica con los conjuntos minimales:

\begin{theorem} \index{ergodicidad \'{u}nica} \index{transformaci\'{o}n! \'{u}nicamente erg\'{o}dica} \label{theoremErgodicidadUnica->Minimal}
Si $T$ es continua en un espacio m\'{e}trico compacto $X$, y \'{u}nicamente
erg\'{o}dica, entonces existe un \'{u}nico minimal $\Lambda $, y adem\'{a}s
$\Lambda $  es el soporte de
la medida invariante por $T$.

\end{theorem}

{\bf Ejemplo de Furstenberg: } El rec\'{\i}proco del Teorema \ref{theoremErgodicidadUnica->Minimal}  es falso:
Furnstenberg en \cite{Furstenberg} (ver tambi\'{e}n
\cite[Cap\'{\i}tulo 2, \S 7, p\'{a}g. 172]{Mane} o  \cite{ManeIngles}), dio un ejemplo de transformaci\'{o}n continua en el toro que preserva la medida de Lebesgue,
para la cual todo el toro es el \'{u}nico minimal pero la medida
de Lebesgue   no es erg\'{o}dica. Luego, en el Ejemplo de Furstenberg, la transformaci\'{o}n no
es \'{u}nicamente erg\'{o}dica, pero existe un \'{u}nico minimal (y adem\'{a}s este minimal   es
el soporte de una medida $T$-invariante).

Otros ejemplos en los que se prueba   la existencia de m\'{a}s de una medida erg\'{o}dica para mapas con un \'{u}nico conjunto minimal, se encuentran en \cite{BachurinStatAttr}.

\vspace{.3cm}

{\em Demostraci\'{o}n del Teorema} \ref{theoremErgodicidadUnica->Minimal}:
Sea $\mu $ la \'{u}nica medida invariante de $T$. Sea $\Lambda $ el soporte compacto
de $ \mu $, definido como $$ \Lambda := \{ x \in X: \forall \: V \mbox{ entorno de } x \: \:
\mu (V) > 0 \}. $$  $\Lambda $  es cerrado en $X$ compacto, luego es compacto.

Veamos que $\Lambda \subset T^{-1}(\Lambda)$. Sea $x \in \Lambda$ y sea $y= T(x)$. Hay que probar que $y \in \Lambda$.  Para todo entorno
$U$ de $y$, $T^{-1}(U) $ es entorno de $x$, luego $0 < \mu (T^{-1} (U))
= \mu (U) $. Esto prueba que $ y \in \Lambda$; luego
  $\Lambda$ es invariante hacia adelante.

Sea $\Lambda _0  $ compacto, no vac\'{\i}o, invariante
hacia adelante. Sea $\widehat  T = T | _{\Lambda _{0}} : \Lambda _0 \mapsto
\Lambda _0 $. Por el teorema de existencia de medidas invariantes,
existe $ \widehat  \nu $ probabilidad que es $\widehat  T$ invariante.

Sea $\nu $ probabilidad en $X$, definida as\'{\i}:
$$ \nu (A) = \widehat  \nu (A \cap \Lambda _0 ) $$ El soporte de
$\nu $ est\'{a} contenido en $\Lambda _0 $. En efecto, si $x \not \in \Lambda_0$, entonces existe $V$, entorno de $x$ disjunto con $\Lambda_0$. Luego $\nu(V) = 0$ y  $x$ no pertenece al soporte de $\nu$.

Probemos que $\nu $ es $T$-invariante:
$$\nu (T^{-1}(A) = \widehat  \nu (T^{-1} (A) \cap \Lambda _0 ) =
\widehat  \nu \{ x \in \Lambda _0: T(x) \in A \} =$$ $$ \widehat  \nu (\widehat  T ^{-1}
(A \cap \Lambda _0)) = \widehat  \nu (A \cap \Lambda _0 ) = \nu (A). $$

Como $T$ es \'{u}nicamente erg\'{o}dica, $\nu = \mu $. Luego
$\Lambda = \mbox{sop } \nu \subset \Lambda _0 $.
Hemos probado que todo $\Lambda _0 $ compacto, no vac\'{\i}o e invariante hacia adelante
contiene a $\Lambda $.
De ello se deducen dos resultados:
Primero, si $\Lambda _0 $ adem\'{a}s est\'{a} contenido en
  $\Lambda $, entonces coincide con $\Lambda $.
Luego $\Lambda $ es minimal.
Segundo: si $\Lambda _0 $ es minimal, entonces coincide con $\Lambda $.
Luego $\Lambda $ es el \'{u}nico  minimal.
\hfill $\Box$

\begin{exercise}\em
Sea $X$ un espacio m\'{e}trico compacto y sea $T:X \mapsto X$ medible, tal que existe alguna medida invariante $\mu$.
$$ P =
\big \{ x \in X: \mbox{ existe } \lim_{n \rightarrow \infty}
\frac{1}{n} \cdot \sum _{j=0} ^{n-1} h \circ T ^{ j} (x) \
\forall \  h : X \mapsto \mathbb{R} \mbox{ medible acotada}    \big \} $$

a)
Probar que $P$ es el conjunto de los puntos peri\'{o}dicos de $T$.
Sugerencia: 1) Inventar una sucesi\'{o}n $\{a _i \}_{i \geq 0 } $ de
ceros y unos tal que no exista el l\'{\i}mite cuando $n \rightarrow \infty $
de $ 1/n \cdot \sum _{j=0}^{n-1} a _j $. (Por ejemplo: 1 uno, 1 cero,
10 unos, 10 ceros, 100 unos, 100 ceros, 1000 unos, 1000 ceros, etc).
2) Si $x $ no es peri\'{o}dico mostrar que $\widetilde  \chi _A (x) $ no existe
para $A = \{ T^i(x): a_i = 1 \}$.

b)
Sea $A \subset X$ cualquier boreliano dado. Probar que tiene
probabilidad total el conjunto
$ \{x\in X: \mbox { existe } \lim_{n\rightarrow
\infty} (1/n) \, \sum _{j=0} ^ {n-1} \chi _A (T^j (x)) \} $.

c)
Probar que si $T$ es continua, \'{u}nicamente erg\'{o}dica y existe
un conjunto minimal con infinitos puntos, entonces $P = \emptyset $.
(Un ejemplo de tal $T$ es la rotaci\'{o}n irracional del c\'{\i}rculo, como veremos a continuaci\'{o}n).

\end{exercise}

\subsection{Ergodicidad de la rotaci\'{o}n irracional} \index{rotaci\'{o}n! irracional} \index{ergodicidad! de la rotaci\'{o}n irracional} \label{sectionRotacionIrracionalErgodica}
\begin{theorem} \label{teoremaErgodicidadRotacionIrracional}
\label {irrota}
La rotaci\'{o}n irracional en el c\'{\i}rculo es \'{u}nicamente erg\'{o}dica.
Luego, es erg\'{o}dica respecto a la medida de Lebesgue.
\end{theorem}

Una prueba muy breve y cl\'{a}sica del Teorema \ref{teoremaErgodicidadRotacionIrracional} requiere la aplicaci\'{o}n de la Teor\'{\i}a Espectral para el estudio de las propiedades de ergodicidad y de mixing de los sistemas din\'{a}micos. No cubriremos ese tema en este texto, por lo que daremos otra prueba,   pedestre y m\'{a}s larga. La prueba breve que aplica la Teor\'{\i}a Espectral puede encontrarse por ejemplo en \cite[Theorem 1.8]{Walters}.

{\em Demostraci\'{o}n: }
Sea $T(x) = x+ x_0 $ (m\'{o}d. 1), para todo $x \in S^1 = \mathbb{R}/\sim (\mbox{m\'{o}d.} 1)$,  donde $x_0 \in (0,1) $ es   irracional dado.
 Aplicando el Teorema \ref{teoremaErgodicidadUnica} parte ii),  para probar que $T$ es \'{u}nicamente erg\'{o}dica, probaremos que para cada $f \in C^0(S^1, \mathbb{R})$ la sucesi\'{o}n $$f_n = \frac{1}{n} \sum_{j= 0}^{n-1} f \circ T^j$$ converge  en todo punto a una constante cuando $n \rightarrow + \infty$.

 Por el Teorema Erg\'{o}dico de Birkhoff-Khinchin para Lebesgue-casi todo punto $x_1 \in S^1$ existe $\lim_{n \rightarrow  + \infty} f_n(x_1) = \widetilde  f(x_1)$. Dado $\epsilon >0$, como $f$ es uniformemente continua en $S^1$ (porque es continua en el compacto $S^1$), existe $\delta >0$ tal que $$|x-y|< \delta \ \Rightarrow \ |f(x)- f(y)| < \epsilon.$$  Por la definici\'{o}n de la rotaci\'{o}n $T$ tenemos la siguiente igualdad m\'{o}d. $1$, para todo punto $y \in S^1$ y para todo $j \geq 0$:
  $$T^j(y) = y + j x_0  = (y-x_1) + x_1 + j x_0 = (y-x_1) + T^j(x_1). $$
  Luego,
  $|y-x_1| < \delta \ \Rightarrow \ |f(T^j(y)) - f(T^j(x_1))| < \epsilon \ \ \forall \ j \in \mathbb{N} \ \Rightarrow $
   $$|f_n(y) - f_n(x_1)| < \epsilon \ \ \forall \ n \in \mathbb{N} .$$

Sea $N \geq 1$ tal que $|f_n(x_1) - \widetilde  f (x_1)| < \epsilon$ para todo $n \geq N$. Obtenemos:
  \begin{equation} \label{equation20}|f_n(y) - \widetilde  f(x_1) | < 2 \epsilon \ \ \forall \ n \in \mathbb{N} \mbox{ si } |y-x_1| < \delta\end{equation}
   Hemos probado que  para todo $\epsilon >0$ existe $\delta >0$ (m\'{o}dulo de continuidad uniforme de $f$, y por lo tanto independiente de $x_1$) tal que se cumple la desigualdad (\ref{equation20}) para todo punto $y$ a distancia menor que $\delta$ de $x_1$. Luego $\limsup f_n (y) - \liminf f_n(y) < 2\epsilon$ si $|y- x_1| < \delta$. Como esta desigualdad vale para Lebesgue-c.t.p. $x_1 \in S^1$, dado $\epsilon >0$ y dado $y \in S^1$, podemos siempre elegir alg\'{u}n $x_1 \in S^1$  tal que $|y- x_1| < \delta$ y tal que existe $\widetilde  f (x_1)$. Entonces podemos aplicar la desigualdad (\ref{equation20}) para todo $y \in S^1$ eligiendo, para cada $y$, alg\'{u}n $x_1$ adecuado. Deducimos que   para todo $y \in S^1$ se cumple $\limsup f_n(y) - \liminf f_n(y) < 2 \epsilon$. Como $\epsilon >0$ es arbitrario,  existe $$\widetilde  f(y) = \lim f_n(y) \ \ \forall \ y \in S_1.$$ Aplicando nuevamente la desigualdad (\ref{equation20}), deducimos que  \em la funci\'{o}n $ \widetilde  f $     es con\-ti\-nua. \em
Lo anterior vale para cualquier rotaci\'{o}n del c\'{\i}rculo. Todav\'{\i}a no usamos que la rotaci\'{o}n es irracional, que aplicaremos ahora  para probar que $\widetilde  f$ es constante.
La funci\'{o}n $\widetilde  f$ es invariante por $T$, es decir $\widetilde  f = \widetilde  f \circ T$, pues el l\'{\i}mite de los promedios temporales de Birkhoff. Luego, $\widetilde  f$ toma un valor constante para cada \'{o}rbita. Entonces, para probar que $\widetilde  f$ es constante (sabiendo ya que es continua) alcanza con probar que existe una \'{o}rbita densa.
  Para probar que alguna \'{o}rbita es densa  alcanza con probar que existe $x_1 \in S^1$ tal que para todo $\delta >0$ la \'{o}rbita  $o^+(x_1) = \{ T^n (x_1): \ n \geq 0\}$ es $\delta$-densa (i.e. todo intervalo de longitud  $\delta$ contiene alg\'{u}n punto de $o^+(x_1)$). La transformaci\'{o}n $T$ preserva la medida de Lebesgue. Aplicando el Lema de Recurrencia de Poincar\'{e}, Lebesgue-c.t.p. es recurrente. Elijamos un punto recurrente $x_1$. Entonces existe $n_j \rightarrow + \infty$ tal que $$\lim_j |T^{n_j}(x_1) - x_1| = 0.$$ Fijemos $m_1  \geq 1$ tal que $|T^{m_1}(x_1)- x_1| < \delta$. Observemos que $T^{m_1}(x_1) - x_1 = m_1 x_0$ (todas las igualdades son m\'{o}dulo 1). Como $x_0 \in (0,1) \setminus \mathbb{Q}$, tenemos que $m_1 x_0 \neq 0$ (m\'{o}d. $1$).  Adem\'{a}s $T^{2m_1}(x_1)  =  T^{m_1}(x_1) + m_1x_0$, de donde
  $$|T^{2m_1}(x_1) - T^{m_1}(x_1)| =  |T^{m_1} (x_1) - x_1|=  m_1 x_0 = a \in (0, \delta) \neq 0$$
  Por inducci\'{o}n en $k \geq 0$ obtenemos
  \begin{equation}
  \label{equation21}
  |T^{(k+1)m_1}(x_1) - T^{k m_1}(x_1)| =  a\neq 0  \ \ \forall \ k \in \mathbb{N}.   \end{equation}
  Afirmamos que el conjunto $A:= \{T^{k m_1}(x_1): \    k  \in \mathbb{N}\} \subset o^+(x_1)$ es $\delta$-denso en $S^1$. Por un lado, para valores diferentes de $k \in \mathbb{N}$, los puntos  respectivos $T^{k n_1}(x_1) \in A$ son diferentes. En efecto, si $T^{k m_1}(x_1) = T^{h m_1}(x_1)$, entonces $km_1 x_0 = h m_1 x_0$, de donde $|k-h| m_1 x_0 = 0$ (m\'{o}d. $1$), donde $k,h, m_1 \in \mathbb{N}$. Como $x_0 \not \in \mathbb{Q}$ deducimos que $h= k$.   Por otro lado, dos puntos consecutivos de $A$ cumplen la igualdad (\ref{equation21}).   Entonces   $\bigcup_{k \in \mathbb{N}} (T^{km_1}(x_1) - a, T^{km_1}(x_1)] = S^1$, lo que demuestra que $A$ es $\delta$-denso en $S^1$ (pues $\delta > a$). Siendo $A \subset o^+(x_1)$, la \'{o}rbita por $x_1$ es $\delta$-densa, para todo $\delta >0$, luego es densa, terminando la demostraci\'{o}n del Teorema \ref{teoremaErgodicidadRotacionIrracional}.  \hfill $\Box$

\begin{remark} \em  \label{remarkRotacionIrracionalOrbitasDensas}
En la demostraci\'{o}n del Teorema \ref{teoremaErgodicidadRotacionIrracional} probamos que existe una \'{o}rbita $\{f^n(x_0)\}_{n \geq 0}$ densa en la rotaci\'{o}n irracional del c\'{\i}rculo $S^1= [0,1]\ |(\mbox{m\'{o}d. } 1)$ dada por $f(x) = x+ a \ (\mbox{m\'{o}d. } 1)$ donde $a $ es irracional. Esto implica que
\begin{center}

\em Toda \'{o}rbita por la rotaci\'{o}n irracional del c\'{\i}rculo es densa. \em
\end{center}

{\em Demostraci\'{o}n: }

Tenemos: $f^n(y_0) = y_0 + na = x_0 + na + (y_0 - x_0) = f^n(x_0) + (y_0 - x_0) \  (\mbox{m\'{o}d. } 1)$. Entonces, la \'{o}rbita $\{f^n(y_0)\}_{n \geq 0}$ se obtiene de la \'{o}rbita $\{f^n(x_0)\}_{n \geq 0}$ rot\'{a}ndola $y_0 - x_0 \ (\mbox{m\'{o}d. } 1)$. Como cualquier rotaci\'{o}n en el c\'{\i}rculo es un homeomorfismo, lleva un conjunto denso a un conjunto tambi\'{e}n denso. Por lo tanto, existe una \'{o}rbita densa si y solo si todas las \'{o}rbitas son densas. Ya probamos (al final de la demostraci\'{o}n del Teorema \ref{teoremaErgodicidadRotacionIrracional}), que existe una \'{o}rbita densa cuando $a$ es irracional. Concluimos que todas las \'{o}rbitas son densas. \hfill $\Box$
\end{remark}

De la prueba del Teorema \ref{teoremaErgodicidadRotacionIrracional} deducimos que las rotaciones en el c\'{\i}rculo son erg\'{o}dicas (respecto de la medida de Lebesgue) si y solo s\'{\i} son \'{u}nicamente erg\'{o}dicas, y esto ocurre si y solo si la rotaci\'{o}n es topologicamente transitiva. M\'{a}s en general, la transitividad
topol\'{o}gica es equivalente a la ergodicidad
\'{u}nica para la rotaci\'{o}n r\'{\i}gida en cualquier grupo topol\'{o}gico compacto
abeliano (ver por ejemplo \cite[page 266]{NicolPetersenEnciclopedia}).

 \begin{exercise}\em

Consid\'{e}rese el toro $k$-dimensional $$\mathbb{T}^k = (S^1)^k =[0,1]^k/(\mbox{m\'{o}d} \{1\}^k),$$ y la operaci\'{o}n de grupo $+$ $(\mbox{m\'{o}d} \{1\}^k)$. Sea $x_0 \in \mathbb{T}^k$. Sea $T$ la traslaci\'{o}n $T(x) = x + x_0  \ (\mbox{m\'{o}d} \{1\}^k) $ para todo $x \in \mathbb{T}^k$.Si $\widetilde  x_0 \in {\mathbb{R}^n}$ es un representante de $x_0 \in T^k$ y si $<,>$ denota el producto interno usual en $\mathbb{R}^n$, asuma que
  $< x_0 , m > \not \in  Z  \; \; \; \forall m \in Z^k - \{0\} $.

a)
Demostrar que $T$ es \'{u}nicamente erg\'{o}dica.

b)
Demostrar que $\forall x \: \in X $, la medida de Lebesgue es el
l\'{\i}mite cuando   $ n \rightarrow \infty $ de las medidas
$ 1/n \cdot \sum _{j=0}^{n-1} \delta _{T^j (x)} $.

c)   Probar que la medida de Lebesgue es la \'{u}nica medida
de proba\-bili\-dad en el toro invariante por todas las traslaciones, pero que no   toda traslaci\'{o}n es \'{u}nicamente erg\'{o}dica.
{ Sugerencia: las traslaciones que tienen puntos peri\'{o}dicos no  son \'{u}nicamente erg\'{o}dicas.}

d) Deducir que para toda traslaci\'{o}n erg\'{o}dica del toro, \'{e}ste es minimal y todas las \'{o}rbitas son densas.

\end{exercise}

\subsection{Transformaciones Mixing.}

\index{transformaci\'{o}n! mixing}
\index{medida! mixing}

\begin{definition}
\em \label{definicionmixing} Sea $T: X \mapsto X$ medible que
preserva una medida de probabilidad $\mu$. Se dice que $T$ \em es
mixing respecto de $\mu$,  o que $\mu$ es mixing respecto de $T$,
\em si para toda pareja $A,B$ de conjuntos medibles se cumple:
\begin{equation} \label{eqn23a}\lim _{n \rightarrow + \infty} \mu (T^{-n} A \cap B ) = \mu (A) \mu (B).\end{equation}
\end{definition}

\begin{theorem} \label{teoremaMixingImplicaErgodica}
Toda transformaci\'{o}n mixing es erg\'{o}dica.

\em (El rec\'{\i}proco es falso como veremos en el Ejemplo \ref{rotirracnoesmixing}.)
\end{theorem}
{\em Demostraci\'{o}n: } En el teorema \ref{teoremaergodicidad} se
prueba que $T$ es erg\'{o}dica si y solo si
\begin{equation} \label{eqn23b}\lim _{n \rightarrow + \infty} \frac{1}{n} \sum _{j=0}^{n-1}
\mu (T^{-j} A \cap B ) = \mu (A) \mu (B)\end{equation} Si se cumple la
igualdad (\ref{eqn23a}) de la definici\'{o}n de transformaci\'{o}n mixing, entonces se
cumple la igualdad (\ref{eqn23b}) de ergodicidad. En efecto si una sucesi\'{o}n $\{a_n\}_{n \in \mathbb{N}}$ de n\'{u}meros reales es convergente a $a \in \mathbb{R}$, entonces (aplicando la definici\'{o}n de l\'{\i}mite se puede chequear f\'{a}cilmente):   $\lim_{n \rightarrow + \infty} (1/n) \, \sum_{j=0}^{n -1} a_j = a$. \hfill $\Box $

Para probar que el
rec\'{\i}proco del Teorema \ref{teoremaMixingImplicaErgodica} es falso, basta ver  que la medida de Lebesgue para la rotaci\'{o}n irracional en el c\'{\i}rculo (que ya probamos que es
erg\'{o}dica en al Teorema \ref{teoremaErgodicidadRotacionIrracional})   no es mixing. Probaremos que no es mixing en el Ejemplo
\ref{rotirracnoesmixing} de esta secci\'{o}n.

\begin{definition}
\em Sea $T: X \mapsto X$ Borel-medible en un espacio topol\'{o}gico
$X$. Se dice que $T$ \em es topol\'{o}gicamente mixing  \em si y solo
si para toda pareja de abiertos $U$ y $V$ no vac\'{\i}os existe $n_0
\in \mathbb{N}$ tal que
$$T^{n} (U) \cap V \neq \emptyset \; \; \; \forall n \geq n_0 $$
\end{definition}

\begin{exercise}\em
Probar que \em si $T$ es Borel medible en un espacio topol\'{o}gico
$X$ y es mixing respecto a una medida de probabilidad $\mu$
positiva sobre abiertos, entonces es topol\'{o}gicamente mixing.
\end{exercise}

 No toda medida erg\'{o}dica es mixing. En efecto:

\begin{example} \em
  \label{rotirracnoesmixing}

  \em La medida de Lebesgue en el c\'{\i}rculo no es mixing para la
rotaci\'{o}n irracional. \em

 {\em Demostraci\'{o}n: } Sea la rotaci\'{o}n irracional $T:S^1
\mapsto S^1$ en el c\'{\i}rculo $S^1$, definida por $T(z) = z + \alpha \ (\mbox{m\'{o}d.} 1)$, donde
 $\alpha$ es un n\'{u}mero irracional (que puede tomarse en $(0,1)$). Para demostrar que su \'{u}nica medida de probabilidad invariante (la medida $m$ de Lebesgue) no es mixing, basta demostrar que $T$ no es
 topol\'{o}gicamente mixing. Sea
  $U \subset S^1$ un intervalo abierto de longitud
 $\epsilon:
 0< \epsilon < (1/4)   \min \{\alpha, 1 - \alpha \}$.  Probaremos que, si
  para cierto $n_0 \in \mathbb{N}$ se cumple $T^{n_0} (U) \cap U \neq
  \emptyset$,
 entonces $T^{n_0 +1}(U) \cap U = \emptyset$.  De lo contrario
 la longitud del intervalo uni\'{o}n $T^{n_0 +1}(U) \cup U \cup
 T^{n_0}(U)$ ser\'{\i}a  menor que $3 \epsilon <   \min \{\alpha, 1 - \alpha
 \}$, pero contendr\'{\i}a dos puntos $x_0$ y $T(x_0)= x_0 + \alpha$
 que distan $ \min \{\alpha, 1 - \alpha \}$. \hfill $\Box$

\end{example}

\begin{exercise}\em
Sea $T$ medible, invertible con inversa medible, que preserva una
medida de probabilidad $\mu$.  Probar que $T$ \em es mixing si y
solo si $T^{-1}$ lo es. \em (Sugerencia: para dos conjuntos $A$ y
$B $ cualesquiera $\mu (T^{-n}(A) \cap B ) = \mu (T^n(B) \cap A
).) $
\end{exercise}

\begin{example} \em
\label{exampleTentMapMixing}
Sea en $S^1 = [0,1]/\sim \mbox{m\'{o}d. } 1$ el tent map $T$ definido por $T(x)= 2 x$ si $0 \leq x \leq 1/2$, y $T(x) = 2- 2 x$ si $1/2 \leq x \leq 1$.  Probaremos el siguiente resultado:

\em La medida de Lebesgue $m$ en $S^1$ es mixing para el tent map $T$; luego es erg\'{o}dica. \em

\end{example}

{\em Demostraci\'{o}n: }
Por definici\'{o}n de mixing, debemos probar que \begin{equation}\label{eqn24}\lim_{n \rightarrow + \infty} m (T^{-n}A \cap B) B= m (A) m (B)\end{equation} para toda pareja de conjuntos medibles $A$ y $B$. Si $A = \emptyset$ \'{o} $B = \emptyset$, entonces la igualdad anterior es trivialmente igual a cero. Consideremos entonces el caso en que $A$ y $B$ son no vac\'{\i}os. Fijemos el boreliano $A$ no vac\'{\i}o. Primero   probaremos (\ref{eqn24}) cuando $B$ es un intervalo abierto; luego cuando $B$ es un abierto, despu\'{e}s para $B$ compacto, y finalmente para $B$ boreliano cualquiera.

Sea $B$ es un intervalo  abierto, con longitud $a$. Para cada $n \geq 1$, la gr\'{a}fica de $T^n$ est\'{a} compuesta por $2^n$ segmentos, con pendiente $2^n$ (en valor absoluto) cada uno. La imagen por $T^n$ de cualquier intervalo de longitud $2^{n-1}$ cubre todo el intervalo $[0,1]$ (croquizar la gr\'{a}fica). La preimagen por $T^n$ de cualquier boreliano $A$ no vac\'{\i}o est\'{a} formada por $2^n$ copias de $A$, todas homot\'{e}ticas a $A$ con raz\'{o}n $1/2^n$ y equidistantes en el segmento $[0,1]$. Luego, la intersecci\'{o}n $T^{-n}(A) \cap B$, cuando $B$ es un intervalo de longitud $a$, contiene $k = [\mbox{parte entera}(a \cdot 2^n)]-1$ de esas copias homot\'{e}ticas a $A$, y menos de $k+2$ de esas copias. Luego:
$$ ({ \mbox{parte entera}(a \cdot 2^n)}-1) \,  \frac{m(A)} {2^n} \leq m(T^{-n}(A) \cap B) $$ $$\leq {(\mbox{parte entera}(a \cdot 2^n) + 1)} \, \frac{m(A)}{2^n}.$$
De las desigualdades anteriores se deduce que $\lim _n m(T^{-n}(A) \cap B) = a \cdot m (A) = m (B) m (A)$. Hemos probado (\ref{eqn24}) cuando $B$ es un intervalo.

Ahora consideremos el caso en que $B$ es un abierto no vac\'{\i}o en el c\'{\i}rculo.  $B= \bigcup I_j$ donde $\{I_j\}$ es una colecci\'{o}n finita o infinita numerable de intervalos abiertos disjuntos dos a dos. Dividimos en dos subcasos: Si $B$ es uni\'{o}n finita de intervalos disjuntos dos a dos, o si es uni\'{o}n infinita numerable. Si $B = \bigcup_{j= 1}^N I_j$, tenemos
$$\lim_{n \rightarrow + \infty} m (T^{-n} (A) \cap B) = \lim_{n \rightarrow + \infty} \sum_{j= 1}^N m (T^{-n}(A) \cap I_j) = $$
$$= \sum_{j= 1}^N \lim_{n \rightarrow + \infty} m (T^{-n}(A) \cap I_j) = \sum_{j= 1}^N m (A) m (I_j) = m (A) \cdot m(B).$$
Si $B = \bigcup_{j= 1}^{+ \infty} I_j$, dado $\epsilon >0$ existe $N$ tal que \begin{equation} \label{eqn25a} 0 \leq m(B) - m(B_N) = m(B \setminus B_N) < \epsilon,\end{equation} donde $B_N = \bigcup_{j= 1}^N I_j$ Luego:
$$0 \leq m(T^{-n}(A) \cap B) - m(T^{-n}(A) \cap B_N) =$$ $$ =m (T^{-n}(A) \cap (B \setminus B_N))  \leq m (B \setminus B_N) < \epsilon \ \ \forall \ n \geq 0. $$ Es decir: \begin{equation} \label{eqn25b} 0 \leq m(T^{-n}(A) \cap B) - m(T^{-n}(A) \cap B_N) < \epsilon \ \ \forall \ n \geq 0\end{equation}
Por lo probado antes $\lim_{n \rightarrow + \infty} m(T^{-n}(A) \cap B_N) = m (A) m (B_N)$. Entonces, para todo $n$ suficientemente grande \begin{equation} \label{eqn25c}|m (T^{-n}(A) \cap B_N) - m(A) m(B_N)| < \epsilon.\end{equation} Reuniendo las desigualdades (\ref{eqn25a}), (\ref{eqn25b}), (\ref{eqn25c}), obtenemos, para todo $n$ suficientemente grande, la siguiente cadena de desigualdades:
$$|m (T^{-n}(A) \cap B) - m(A) m(B)| \leq |m (T^{-n}(A) \cap B) - m (T^{-n}(A) \cap B_N)| + $$ $$+ |m (T^{-n}(A) \cap B_N) - m(A) m(B_N)| + m(A) | m(B_N) - m(B)| < 3 \epsilon.$$
Lo anterior demuestra que $\lim_{n \rightarrow + \infty} m (T^{-n}(A) \cap B) = m(A) m(B)$. Hemos probado la  igualdad (\ref{eqn24}) cuando $B$ es abierto. Ahora demostremos que si la igualdad (\ref{eqn24})  se cumple para un conjunto $B$, entonces tambi\'{e}n se cumple para su complemento $B^c$. En efecto: $$m (T^{-n}(A) \cap B^c) = m (T^{-n} (A)) - m (T^{-n}(A) \cap B) =$$ $$= m (A) - m (T^{-n}(A) \cap B) \rightarrow_n m (A) - m(A) m (B) = m(A) m (B^c).$$
Entonces, como (\ref{eqn24}) vale para todos los abiertos $B$, y   es una propiedad cerrada en complementos, se cumple tambi\'{e}n para todos los compactos $B$. Ahora prob\'{e}mosla para cualquier boreliano $B$. Dado $\epsilon$, existe un compacto $K$ y un abierto $V$ tales que
$K \subset B \subset V$ y $m(V \setminus K) < \epsilon$. Luego $$m(V) = m(B) + m (V \setminus B) \leq m(B) + m (V \setminus K) < m(B) + \epsilon.$$ An\'{a}logamente
$$m(K) = m(B) - m (B \setminus K) \geq m(B) - m(V \setminus K) > m(B) - \epsilon.$$
Entonces:
$$m(T^{-n}(A) \cap B) \leq m (T^{-n}(A) \cap V) \rightarrow _n m(A) m (V) \leq m(A) m (B) + \epsilon, $$
$$m(T^{-n}(A) \cap B) \geq m (T^{-n}(A) \cap K) \rightarrow _n m(A) m (K) \geq m(A) m (B) - \epsilon. $$
Lo anterior prueba que para todo $n$ suficientemente grande
$$|m(T^{-n}(A) \cap B) - m(A) \, m (B)| < 2\epsilon, $$
de donde se deduce que la igualdad (\ref{eqn24}), como quer\'{\i}amos demostrar.
\hfill $\Box$


\subsection{Descomposici\'{o}n Erg\'{o}dica}

 El prop\'{o}sito de esta secci\'{o}n es enunciar el Teorema \ref{theoremDescoErgodicaEspaciosMetricos}, de descomposici\'{o}n o desintegraci\'{o}n erg\'{o}dica. El mismo   extiende el resultado de existencia de medidas erg\'{o}dicas para transformaciones continuas en espacios m\'{e}tricos compactos, probando c\'{o}mo se puede descomponer o desintegrar una medida invariante $\mu$ en funci\'{o}n de las que son erg\'{o}dicas.
%
%
 \begin{definition} {\bf Descomposici\'{o}n o Desintegraci\'{o}n Erg\'{o}dica } \label{definitionDescoErgodica}

 \em Sea $(X, \mathcal A)$ un espacio medible y $T: X \mapsto X$ una transformaci\'{o}n medible tal que existe alguna medida  de probabilidad $\mu$ (definida  en ${\mathcal A} $) invariante  por $T$.

 {\bf (i)} Sea $A \in {\mathcal A}$. Decimos que $\mu$ tiene \em descomposici\'{o}n o desintegraci\'{o}n erg\'{o}dica \em para el conjunto $A$ si para $\mu$-c.t.p. $x \in X$ existe una medida erg\'{o}dica $\mu_x$ tal que:

 \ \ \ \ \ {\bf (a)} La funci\'{o}n real definida $\mu$-c.t.p por $x   \mapsto \mu_x(A)$ es medible para $\mu$-c.t.p. $x \in X$. (Entonces est\'{a} en $L^1(\mu)$ pues est\'{a}   acotada por 1)

 \ \ \ \ \ {\bf (b)} $$\mu(A) = \int_{x \in X} \big(\mu_x(A)\big) \, d \mu.$$

 {\bf (ii)} Sea $h \in L^1(\mu)$. Decimos que $\mu$ tiene \em descomposici\'{o}n o desintegraci\'{o}n erg\'{o}dica \em para la funci\'{o}n $h$ si para $\mu$-c.t.p $x \in X$ existe una medida erg\'{o}dica $\mu_x$  tal que $h \in L^1(\mu_x)$ y tal que:

 \ \ \ \ \ {\bf (c)} La funci\'{o}n real definida   $\mu$-c.t.p  por $x     \mapsto \int h \, d\mu_x $ es medible  para $\mu$-c.t.p. $x \in X$ y est\'{a} en $L^1(\mu)$.

 \ \ \ \ \ {\bf (d)}
 $$\int h \, d \mu = \int_{x \in X}     \Big( \int h \, d \mu_x  \Big) \, d \mu.$$
 \end{definition}

 \begin{theorem}.
 \label{theoremDescoErgodicaEspaciosMetricos}

 {\bf Descomposici\'{o}n Erg\'{o}dica en espacios m\'{e}tricos compactos}

 Sea $X$ un espacio m\'{e}trico compacto  provisto de la sigma-\'{a}lgebra de Borel. Sea $T: X \mapsto X$ una transformaci\'{o}n medible que preserva una medida de probabilidad $\mu$.

Entonces, para todo $A \in {\mathcal A}$ y para toda $h \in L^1(\mu)$ existe descomposici\'{o}n erg\'{o}dica de     $\mu$.

M\'{a}s a\'{u}n, para $\mu$-\mbox{c.t.p. } $x \in X$ existe y es \'{u}nica una medida erg\'{o}dica $\mu_x$ \em (llamada {\bf componente erg\'{o}dica } de $\mu$ a la que pertenece $x$) \em tal que
$$\lim_{n \rightarrow + \infty} \frac{1}{n} \sum_{j= 0}^{n-1} \delta_{T^j(x)} = \mu_x,$$
donde el l\'{\i}mite es en la topolog\'{\i}a d\'{e}bil estrella del espacio ${\mathcal M}$ de probabilidades. \index{componentes erg\'{o}dicas}

 \end{theorem}

La demostraci\'{o}n del Teorema \ref{theoremDescoErgodicaEspaciosMetricos} para transformaciones continuas en espacios m\'{e}tricos se encuentra por ejemplo  en \cite[Cap.II \S]{Mane} (ver tambi\'{e}n \cite{ManeIngles}). Una generalizaci\'{o}n para transformaciones medibles que preservan una medida de probabilidad, se encuentra en en \cite[Theorem 2.3.3]{Keller}.
Ahora veamos alguna de   sus consecuencias:

Recordamos la Definici\'{o}n \ref{definitionProbabilidadTotal}: Un conjunto medible $A \subset X$  se dice que tiene \em probabilidad total \em si $\mu(A)= 1$ para toda medida de probabilidad $\mu$ en $X$ que sea invariante por $T$ (bajo la hip\'{o}tesis que existen medidas de probabilidad invariantes por $T$).

 \begin{corollary} \label{CorolarioDescoErgodicaEspaciosMetricos}
  Si $X$ es un espacio  m\'{e}trico compacto, y si $T: X \mapsto X$ es continua, entonces existen medidas erg\'{o}dicas para $T$. Adem\'{a}s un conjunto medible $A \subset X$ tiene probabilidad total  si y solo si $\nu(A)= 1$ para toda medida de probabilidad $\nu$ erg\'{o}dica   para $T$.
 \end{corollary}
{\em Demostraci\'{o}n: }
Usando el Teorema \ref{theoremDescoErgodicaEspaciosMetricos},
y la definici\'{o}n \ref{definitionDescoErgodica}, existen medidas erg\'{o}dicas para $T$ . Adem\'{a}s, para cualquier conjunto medible $A$, cualquiera sea la medida invariante $\mu$ tenemos, para el complemento $A^c$ de $A$, la siguiente igualdad:
$$  \mu(A^c) = \int \Big (\int \chi_{A^c} \, d \mu_x \Big)\, d \mu,$$
donde $\mu_x$ es una medida erg\'{o}dica, que depende del punto $x$ y est\'{a} definida para $\mu$-c.t.p. $x \in X$. Como $\chi_{A^c} \geq 0$
entonces la funci\'{o}n $ x \mapsto \int \chi_{A^c} \, d \mu_x = \mu_x(A^c)$ es no negativa. Las medidas $\mu_x$ son erg\'{o}dicas seg\'{u}n enuncia el Teorema \ref{theoremDescoErgodicaEspaciosMetricos} y la Definici\'{o}n \ref{definitionDescoErgodica}. Concluimos que  $\mu(A^c) = 0$ para toda medida invariante $\mu$, si  y solo si $\nu(A^c) = 0$ para toda medida erg\'{o}dica $\nu$.
\hfill $\Box$


\section{Din\'{a}mica diferenciable: \\ Hiperbolicidad uniforme y no uniforme} \label{chapterTeoriaPesinElementos}

El siguiente ejemplo es conocido como "Arnold's cat map" (el mapa del gato de Arnold). Esto es porque su din\'{a}mica fue representada por el matem\'{a}tico ruso Vladimir Arnold en \cite{arnold}, con un dibujo similar al de la Figura \ref{figuraManzana} de este cap\'{\i}tulo. En su dibujo (ver por ejemplo uno parecido al original  en \cite[Figure 2]{Pikovsky_Arnoldcat}), Arnold utiliza    el contorno de la figura estilizada de un gato, en lugar de una manzana como hacemos nosotros en la Figura \ref{figuraManzana} de este cap\'{\i}tulo. En su dibujo, Arnold   \lq\lq muestra\rq\rq \ el efecto de  la propiedad de mixing del mapa sobre los trazos del contorno del gato, el cual, en pocos iterados, se vuelve irreconocible.

\subsection{Ejemplo de automorfismo lineal hiperb\'{o}lico en el toro.} \index{automorfismo! lineal del toro} \label{section2111} \index{hiperbolicidad! uniforme} \index{transformaci\'{o}n! hiperb\'{o}lica! uniforme} \index{difeomorfismos! de Anosov}

Sea el toro ${\mathbb{T}}^2 = \mathbb{R}^2 /\sim$ donde la relaci\'{o}n de equivalencia
$\sim$ est\'{a} dada por:$(a,b) \sim (c,d) $ en ${\mathbb{T}}^2$ si $c-a$ y $d-b$
son enteros. (Otra notaci\'{o}n que se usa para $\mathbb{R}^2 /\sim $ es $\mathbb{R}^2
|_{mod \mathbb{Z}^2} = \mathbb{R}^2/\mathbb{Z}^2$.)

Sea $\Pi: \mathbb{R}^2 \mapsto \mathbb{Z}^2$ la proyecci\'{o}n del espacio de
recubrimiento $\mathbb{R}^2$ del toro definida por $\Pi (a,b) = (a,b)_{mod \
\mathbb{Z}^2}$ donde esto \'{u}ltimo indica la clase de equivalencia de $(a,b)
\in \mathbb{R}^2$.

Sea $f: {\mathbb{T}}^2 \mapsto {\mathbb{T}}^2$ dada por $$f (x) = \Pi ( \left (%
\begin{array}{cc}
  2 & 1 \\
  1 & 1 \\
\end{array}%
\right) (\Pi ^{-1}(x))  )  \ \ \ \forall \ x \in \mathbb{T}^2.$$ Llamaremos a $f$ automorfismo lineal hiperb\'{o}lico de matriz $\left(%
\begin{array}{cc}
  2 & 1 \\
  1 & 1 \\
\end{array}%
\right)$ en el toro, o simplemente \lq \lq $\left (
\begin{array}{cc}
  2 & 1 \\
  1 & 1 \\
\end{array} \right )$ en el toro".
Abusando de la notaci\'{o}n, a un punto $x \in {\mathbb{T}}^2$ lo denotaremos con
cualquier representante suyo $(a,b) \in \mathbb{R}^2$: $(a,b) \in
\Pi^{-1}(\{x\})$.

\begin{exercise}\em
Probar que $(0,0)$ es el \'{u}nico punto fijo por la transformaci\'{o}n $f = \left(%
\begin{array}{cc}
  2  & 1     \\
  1  & 1 \\
\end{array}%
\right)$ en el toro. Probar que
$$\{(1/5, 2/5), (4/5, 3/5)\} \ \ \ \mbox{ y } \ \ \
  \{(3/5, 1/5), (2/5, 4/5)\}$$  son dos \'{o}rbitas peri\'{o}dicas por $f$ y
que son las \'{u}nicas de per\'{\i}odo $2$. Probar que
$$\{(1/2, 1/2), (1/2,
1), (1, 1/2)\} $$  es una \'{o}rbita peri\'{o}dica de per\'{\i}odo $3$.

\end{exercise} La topolog\'{\i}a en el toro es el cociente de la topolog\'{\i}a usual en
$\mathbb{R}^2$. Es metrizable y la m\'{e}trica est\'{a} dada por $$\dist (x,y) = $$ $$
\min \{ \sqrt{(c-a)^2 + (d-b)^2}: {(a,b) \in \Pi^{-1}(x), (c,d)
\in \Pi ^{-1}(y)} \}.$$

La medida de Lebesgue en el toro es la medida de Borel $\widetilde  m$
definida por $\widetilde  m (B) = m (\Pi ^{-1} B \cap [0,1]^2)$ donde
$m$ es la medida de Lebesgue en $\mathbb{R}^2$. Se observa que la medida
de Lebesgue en el toro es una medida de probabilidad. Donde no d\'{e}
lugar a confusi\'{o}n renombraremos como $m$ a la medida de Lebesgue
en el toro.

\begin{proposition} \index{medida! de Lebesgue}
La medida de Lebesgue $m$ en el toro es invariante por la
transformaci\'{o}n $f = \left(%
\begin{array}{cc}
  2  & 1     \\
  1  & 1 \\
\end{array}%
\right)$.
\end{proposition}

{\em Demostraci\'{o}n: }

Denotamos $A =\left(%
\begin{array}{cc}
  2  & 1     \\
  1  & 1 \\
\end{array}%
\right)$ a la matriz de coeficientes enteros que define la transformaci\'{o}n $f$ en el toro. Observamos que $\mbox{det}(A)= 1$.

Sea $B$ un boreliano en el toro. $$m (f^{-1}
(B)) = \int \chi _{f^{-1}(B)} (x) \; dm(x) = \int \chi _B \circ
f(x)\, dm(x).$$ Haciendo el cambio de variables lineal e invertible
$z= f(x)$ en la integral anterior resulta $$m (f^{-1} (B)) =  \int
\chi _B (z) J(z) \, dm(z)$$  donde $J(z) = |\mbox{det} df^{-1}(z)| = |\mbox{det}
df(x)|^{-1}$ es el jacobiano del cambio de variables $z =f (x)$. \index{jacobiano}
En nuestro caso una parametrizaci\'{o}n local de la superficie del toro $\mathbb{T}^2$ est\'{a} dada por $\Pi|_{B_{\delta}(\Pi^{-1}(x))}$, donde $ { B_{\delta}(\Pi^{-1}(x))}$ es la bola abierta en $\mathbb{R}^2$ de radio $\delta >0$ (suficientemente peque\~{n}o) y centro en un punto denotado como $\Pi^{-1}(x) \in \mathbb{R}^2$ que se proyecta por $\Pi: \mathbb{R}^2 \mapsto \mathbb{T}^2$ en el punto $x \in \mathbb{T}^2$. Calculando la derivada $df(x)$ y el Jacobiano con las coordenadas en esa carta local,  resulta $T_x\mathbb{T}^2 \sim \mathbb{R}^2$, y  $$J(z) = (\mbox{det} A)^{-1} = 1 \; \; \forall z \in
{\mathbb{T}}^2.$$ Entonces
 $m (f^{-1} (B)) = \int \chi _B (z)\, dm(z) = m(B) $. \hfill $ \Box$ 

\vspace{.2cm}

{\bf Din\'{a}mica del ejemplo de automorfismo lineal hiperb\'{o}lico en el toro.}

Estudiemos la din\'{a}mica por iterados de la transformaci\'{o}n lineal
$F: \mathbb{R}^2 \mapsto \mathbb{R}^2$ que tiene como matriz
$A = \left(%
\begin{array}{cc}
  2  & 1     \\
  1  & 1 \\
\end{array}%
\right)$.

$F$ se llama levantado de $f$ a $\mathbb{R}^2$ y $f$ se llama proyecci\'{o}n
de $F$ en el toro. La din\'{a}mica de $F$ est\'{a} relacionada fuertemente
con la din\'{a}mica de su proyecci\'{o}n $f$ en el toro.

En efecto la proyecci\'{o}n $\Pi: \mathbb{R}^2 \mapsto {\mathbb{T}}^2$ \em cuando
restringida a un entorno suficientemente peque\~{n}o del origen
$(0,0)$ \em  es un homeomorfismo sobre su imagen.

$\Pi$ transforma la din\'{a}mica de $F$ en la de $f$. M\'{a}s precisamente
$$\Pi \circ F = f \circ \Pi$$

Luego, aplicando $\Pi: \mathbb{R}^2 \mapsto {\mathbb{T}}^2$ a una \'{o}rbita de $F$ se
obtiene una \'{o}rbita de $f$. Y la preimagen por $\Pi$ de una \'{o}rbita
por $f$ es una infinidad numerable de \'{o}rbitas por $F$ tal que una
se obtiene de otra traslad\'{a}ndola en $\mathbb{R}^2$ seg\'{u}n un vector de
coordenadas enteras.

 El comportamiento topol\'{o}gico local de las
\'{o}rbitas de $F$ en un entorno suficientemente peque\~{n}o $V$ del
origen es el mismo (a menos del homeomorfismo $\Pi|_V$) que el de
las \'{o}rbitas de $f$ en el abierto $\Pi (V)$.

\vspace{.2cm}

{\bf Variedades invariantes.} \index{variedad invariante! estable}  \index{variedad invariante! inestable}

 Los valores propios de la matriz $A$
son $$ \sigma := (3+\sqrt5)/2
>1, \ \ \ 0 < \lambda:= (3-\sqrt5)/2 <1.$$ Las direcciones propias respectivas
tienen pendientes irracionales, la primera positiva y la segunda
negativa. Se deduce que las dos rectas que pasan por el origen y
tienen direcciones seg\'{u}n los vectores propios de la matriz $A$ son
invariantes por $F$ en el plano. Entonces sus proyecciones en el
toro son curvas invariantes por $f$ y se cortan transversalmente
en el origen (y tambi\'{e}n se cortan transversalmente en todos sus otros puntos de intersecci\'{o}n en el
toro).

$F$ restringida a la recta $r_1$ que tiene direcci\'{o}n propia de
valor propio mayor que 1, expande las distancias exponencialmente
con tasa $\log (3+\sqrt5) /2 >0$. Es decir, contrae distancias hacia el pasado como $\sigma ^{-n} = e^{-n \log \sigma}$: $$\frac{\mbox{dist} (F^{-n}(a,b), (0,0))}{\mbox{dist}((a,b), (0,0))} =e^{-n \log \sigma}   \ \forall \ n \in \mathbb{N}, \ \ \forall \ (a,b) \in r_1 \setminus (0,0).$$ Proyectando $r_1$ en el toro se
obtiene una curva $W^u(0,0) = \Pi (r_1)$ inmersa en el toro, que se llama \em
variedad inestable por $(0,0)$. \em

$W^u(0,0)$  pasa por el origen, \em  es densa en el toro  \em
(esto se puede demostrar usando que la pendiente de $r_1$ en el
plano es irracional y usando que la rotaci\'{o}n irracional en el
c\'{\i}rculo es densa en el c\'{\i}rculo), \em  es invariante por $f$ y
cumple: \em
 $$W^u(0,0) = \{y \in {\mathbb{T}}^2: \lim_{n \rightarrow - \infty} f^n(y) =   (0,0)\}$$
Esto \'{u}ltimo se debe a que $r_1 = \{(a,b) \in \mathbb{R}^2: \lim_{n
\rightarrow - \infty} F^n(a,b) = (0,0) \}$. Observamos que la subvariedad $W^u(0,0) = \Pi (r_1)$ es inmersa y densa en  $\mathbb{T}^2$, pero \em no es subvariedad encajada en ${\mathbb{T}^2}$. \em Es decir, la topolog\'{\i}a que se define a lo largo de   la subvariedad $W^u(0,0)$ no es la inducida por su inclusi\'{o}n en $\mathbb{R}^2$. Los abiertos en $W^u(0,0)$ est\'{a}n generados por  los arcos abiertos (homeomorfos a intervalos abiertos en la recta real). Estos no se obtienen como intersecci\'{o}n de un abierto en $\mathbb{T}^2$ con $W^u(0,0)$ pues  cualquier abierto en $\mathbb{T}^2$ cortado con $W^u(0,0)$ contiene una infinidad de arcos conexos, debido a la densidad de $W^u(0,0)$.

{\bf Variedad inestable local:} \index{variedad invariante! local} En este ejemplo, existe $\epsilon >0$ suficientemente peque\~{n}o, tal que, denotando $W^u_{loc}(0,0)$ (variedad inestable local) a la componente conexa de $W^u(0,0)$ intersecada con la bola de centro $(0,0)$ y radio $\epsilon$ en el toro $\mathbb{T}^2$, se tiene:
$$\frac{\mbox{dist} (f^{-n}(y), (0,0))}{\mbox{dist}(y, (0,0))} =e^{-n \log \sigma}   \ \forall \ n \in \mathbb{N}, \ \ \forall \ y  \in W^u_{loc}(0,0). $$
Esta igualdad se obtiene porque la bola de centro $(0,0)$ y radio $\epsilon >0$ en el toro $\mathbb{T}^2$ es difeomorfa por un preimagen de $\Pi$, con la bola de centro $(0,0)$ y radio $\epsilon >0$ en $\mathbb{R}^2$, si $\epsilon >0$ es suficientemente peque\~{n}o.

An\'{a}logamente $F$ restringida a la recta $r_2$ que tiene direcci\'{o}n
propia de valor propio menor que 1, contrae las distancias
exponencialmente con tasa $\log (3-\sqrt5) /2 <0$. Es decir, contrae distancias hacia el futuro como $\lambda ^{n} = e^{n \log \lambda}$: $$\frac{\mbox{dist} (F^{n}(a,b), (0,0))}{\mbox{dist}((a,b), (0,0))} =e^{n \log \lambda}   \ \forall \ n \in \mathbb{N}, \ \ \forall \ (a,b) \in r_2 \setminus (0,0).$$ Proyectando
$r_2$ en el toro se obtiene una curva $W^s(0,0) = \Pi (r_1)$ inmersa en el
toro, que  se llama \em variedad estable por $(0,0)$.  \em

$W^s(0,0)$ pasa por el origen, \em es densa en el toro \em (porque
la pendiente de $r_2$ en el plano es irracional), \em es
invariante por $f$, no encajada en $\mathbb{R}^2$  \em y cumple:
 $$W^s(0,0) = \{y \in {\mathbb{T}}^2: \lim_{n \rightarrow + \infty}f^n(y) = (0,0)\}$$
Esto \'{u}ltimo se debe a que $r_2 = \{(a,b) \in \mathbb{R}^2:  \lim_{n
\rightarrow + \infty} F^n(a,b)= (0,0) \}$.

Denotando $W^s_{loc}(0,0)$ (variedad estable local) a la componente conexa de $W^s(0,0)$ intersecada con la bola de centro $(0,0)$ y radio $\epsilon$ en el toro $\mathbb{T}^2$, se tiene:
$$\frac{\mbox{dist} (f^{n}(y), (0,0))}{\mbox{dist}(y, (0,0))} =e^{n \log \lambda}   \ \forall \ n \in \mathbb{N}, \ \ \forall \ y  \in W^s_{loc}(0,0). $$

{\bf Variedades invariantes por cualquier punto:}
El argumento anterior puede aplicarse a cualquier punto peri\'{o}dico $x$. Deducimos que todos los puntos peri\'{o}dicos son hiperb\'{o}licos tipo silla (tienen un valor propio mayor que uno y otro positivo menor que uno).

En general, para cualquier punto $x \in \mathbb{T}^2$, aunque no sea peri\'{o}dico, la variedad estable $W^s(x)$ y la variedad inestable $W^u(x)$ se definen como la proyecciones sobre el toro de las rectas en $\mathbb{R}^2$ que pasan por $\Pi^{-1}(x)$, seg\'{u}n las direcciones de los vectores propios de la matriz $A$ (que son las mismas que las direcciones de las rectas $r_1$ y $r_2$ en $\mathbb{R}^2$ que pasan por el origen). Argumentando como m\'{a}s arriba se tiene $$\lim_{n \rightarrow + \infty} \mbox{dist}(f^n(y), f^n(x))= 0 \ \ \forall \ y \in W^s(x)$$
$$\lim_{n \rightarrow + \infty} \mbox{dist}(f^{-n}(y), f^{-n}(x))= 0 \ \ \forall \ y \in W^u(x),$$
y el acercamiento a cero de esas distancias se realiza como $\lambda^n$ o $\sigma^{-n}$ respectivamente.  Las variedad inestable por un punto  $x$ no es invariante si el punto $x$ no es   fijo por $f$, pero  su imagen por $f$ es la variedad inestable por el punto $f(x)$ (esto se chequea inmediatamente de la construcci\'{o}n de $W^u(x) $ y $W^u(f(x))$ como las im\'{a}genes por $\Pi$ de las rectas paralelas a $r_1$ por $\Pi^{-1}(x)$ y $\Pi^{-1}(f(x)) = F(\Pi^{-1}(x))$ respectivamente).

 {\bf Foliaciones invariantes:} \index{foliaci\'{o}n! invariante! estable} \index{foliaci\'{o}n! invariante! inestable}

 La  familia  de todas las variedades inestables, forman en el toro $\mathbb{T}^2$ lo que se llama una \em foliaci\'{o}n \em invariante: pues cada subvariedad de la foliaci\'{o}n (llamada hoja), al aplicarle $f$ se transforma en otra hoja de la foliaci\'{o}n. An\'{a}logamente, la foliaci\'{o}n formada por las variedades estables, es invariante.

 {\begin{figure} [h]

 \vspace{-.3cm}

\begin{center}\includegraphics[scale=.6]{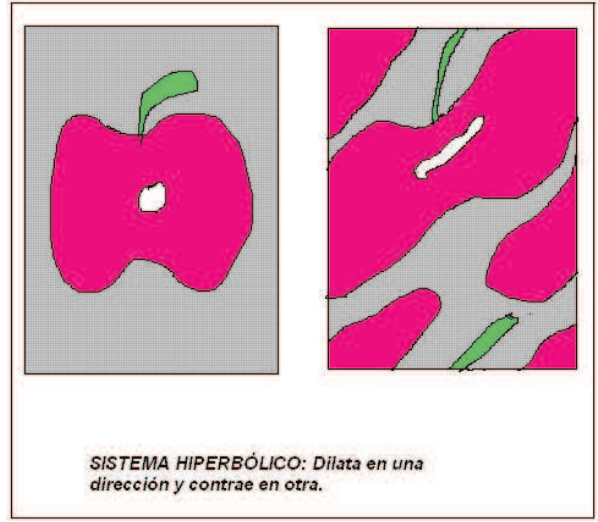}

\vspace{0cm}

\caption{\label{figuraManzana} Deformaci\'{o}n producida por el automorfismo lineal hiperb\'{o}lico en el toro.}
\end{center}
\vspace{0cm}
\end{figure}}

 {\bf Interpretaci\'{o}n gr\'{a}fica del automorfismo lineal hiperb\'{o}lico:} La deformaci\'{o}n que produce $f$ en este ejemplo de automorfismo lineal hiperb\'{o}lico en el toro $\mathbb{T}^2$, est\'{a} representado en la figura \ref{figuraManzana}. Esa figura es una modificaci\'{o}n de la conocida  llamada   \lq\lq gato de Arnold\rq\rq (que en vez de una manzana,   deforma la imagen de un gato, figura creada por Arnold, en la d\'{e}cada de 1960 para ilustrar la deformaci\'{o}n hiperb\'{o}lica de un automorfismo lineal hiperb\'{o}lico en el toro $\mathbb{T}^2$). Al iterar sucesivas veces $f$, la figura representada, se estira a lo largo de la foliaci\'{o}n inestable y se contrae a lo largo de la foliaci\'{o}n estable. Como $f$ es invertible, los pedazos que se obtienen de identificar $0$ con $1$ en vertical y horizontal, no se intersecan (provienen de subconjuntos disjuntos antes de aplicar $f$).

\begin{definition}
{\bf Exponentes de Lyapunov en puntos fijos hi\-per\-b\'{o}\-licos.} \index{exponentes de Lyapunov}
\em
Los logaritmos de los m\'{o}dulos de los valores propios de $df^p (x_0)$ dividido $p$, cuando $x_0$  es un punto peri\'{o}dico de per\'{\i}odo $p$, se llaman \em
exponentes de Lyapunov \em en $x_0$. En el cap\'{\i}tulo \ref{chapterTeoriaPesinElementos}  definiremos con   generalidad los exponentes de Lyapunov de   cualquier sistema  din\'{a}mico  diferenciable, para casi todas sus \'{o}rbitas (no necesariamente puntos fijos ni \'{o}rbitas peri\'{o}dicas).

En el ejemplo del $\left(%
\begin{array}{cc}
  2 & 1 \\
  1 & 1 \\
\end{array}%
\right)$   $f$ es lineal (mejor dicho $F$, dada por la matriz $A$, es lineal), por lo tanto la matriz $A$ es, en una base adecuada de $T_{(0,0)} {\mathbb{T}^2}$), la derivada $df (0,0): T_{(0,0)} \mathbb{T}^2 \sim \mathbb{R}^{2}  \mapsto T_{(0,0)} \mathbb{T}^2 \sim \mathbb{R}^2$. Los logaritmos de los valores propios de $A = df(0,0)$
son dos, uno positivo $\chi^+= \log \sigma$ y el otro negativo $\chi^- = \log \lambda$. (El origen es un punto
fijo hiperb\'{o}lico tipo silla, pues $0 <\lambda < 1 < \sigma$).

El exponente de Lyapunov $\chi^+ := \log \mu =  \log (3+\sqrt5) /2 >0$ es \em la tasa  exponencial \em
\em   de dilataci\'{o}n \em al aplicar $df_{(0,0)}$, de las normas de los vectores en el subespacio propio $U$, que corresponde al valor propio $\mu$  de $df_{(0,0)}$. Precisamente:
$$\lim _{n \rightarrow + \infty} \frac { \log (\|df^n_{(0,0)} u\|/ \|u\|) }{n} = \log \mu = \chi^+ >0 \; \
\forall \  u \in U \setminus \{{\bf 0} \}\subset T_xM.$$
Observar que para $n < 0$ tambi\'{e}n se cumple la misma igualdad anterior; es decir:
$$\lim _{n \rightarrow - \infty} \frac { \log (\|df^n_{(0,0)} u\|/ \|u\|) }{n} = \log \mu = \chi^+ \; \
\forall \  u \in U \setminus \{{\bf 0} \}\subset T_xM.$$

Adem\'{a}s, el exponente de Lyapunov positivo $\chi^+  $  es la tasa  exponencial   de dilataci\'{o}n de las distancias a lo largo de la variedad inestable local por el punto fijo. M\'{a}s precisamente
$$\lim _{n \rightarrow - \infty} \frac { \log \dist (f^n(y), (0,0))}{n} = \chi^+ = \log \mu >0\; \;
\forall y \in W^u(0,0)$$  \index{exponentes de Lyapunov! no nulos}
\index{exponentes de Lyapunov! positivos}

An\'{a}logamente, el exponente de Lyapunov $\chi^-:= \log \lambda  =  \log (3-\sqrt5) /2 <0$
negativo, se interpreta como \em la tasa \em  \em exponencial de contracci\'{o}n por $f$ \em en un entorno
suficientemente peque\~{n}o del punto fijo a lo largo de la variedad
estable por ese punto. M\'{a}s precisamente
$$\lim _{n \rightarrow + \infty} \frac { \log \dist (f^n(y), (0,0))}{n} = \chi^- = \log \lambda < 0 \; \;
\forall y \in W^s(0,0).$$

\end{definition}

\begin{exercise}\em \label{exercise2111transitivo}
Sea $f$ la transformaci\'{o}n $\left(%
\begin{array}{cc}
  2 & 1 \\
  1 & 1 \\
\end{array}%
\right)$ en el toro. Probar que $m$- c.t.p es recurrente pero que
hay puntos no recurrentes (Sugerencia: $y \in W^s(0,0) \setminus
\{(0,0)\}$ no es recurrente). Probar que dados dos abiertos $U$ y
$V$ cualesquiera en el toro,  existe una subsucesi\'{o}n $n_j$ de
naturales tales que $f^{n_j}U \cap V \neq \emptyset$. (Sugerencia:
Usar que $W^u((0,0))$ es invariante por $f$ y pru\'{e}bese que para
cualquier arco compacto $K \subset W^u(0,0)$ la uni\'{o}n de sus
iterados hacia el futuro $\bigcup _{n \in \mathbb{N}} f^n(K) $ es densa en el
toro.)

Deducir que $f$ es transitiva. \index{automorfismo!
lineal del toro} \index{transformaci\'{o}n! hiperb\'{o}lica! uniforme}
\index{difeomorfismos! de Anosov}

Probar que todo punto es no errante. Concluir que si bien todo
punto recurrente es no errante, no necesariamente todo punto no
errante es recurrente.
\end{exercise}

{\bf Observaci\'{o}n: } \index{automorfismo!
lineal del toro} \index{transformaci\'{o}n! erg\'{o}dica} \index{ergodicidad} \index{medida! erg\'{o}dica} \index{medida! de Lebesgue} En el cap\'{\i}tulo \ref{chapterAtractores}, Corolario \ref{corolarioAnosovTransitivomedidaLebesgue}     probaremos que la transformaci\'{o}n $$f =
\left(%
\begin{array}{cc}
  2 & 1 \\
  1 & 1 \\
\end{array}%
\right)$$ en el toro es \em erg\'{o}dica respecto a la medida de
Lebesgue.\em

\begin{exercise}\em
Sea el automorfismo lineal hiperb\'{o}lico $f= \left(
                                        \begin{array}{cc}
                                          2 & 1 \\
                                          1 & 1 \\
                                        \end{array}
                                      \right)
$ del toro $\mathbb{T}^2$. Sea $P$ el conjunto de puntos peri\'{o}dicos por $f$. Probar que $P \neq \emptyset $, pero
$\mu (P) = 0 $ para toda $\mu$ invariante y erg\'{o}dica para $f$ que sea positiva sobre abiertos.
\end{exercise}

\subsection{Difeomorfismos de Anosov e hiperbolicidad uniforme}

En esta secci\'{o}n asumiremos que $M$ es una variedad diferenciable, compacta y conexa, provista de una estructura riemanniana (i.e. existe  un producto interno $< \cdot, \cdot >:TM \times TM \mapsto \mathbb{R}$ diferenciable, definido en el fibrado tangente $TM$. Por lo tanto, para todo $x \in M$ y para todo $v \in T_xM$, est\'{a} definida la norma $\| v \| := \sqrt {< v, v>_x}$, determinada por la m\'{e}trica riemanniana en $M$).

Sea $f: M \mapsto M$ un difeomorfismo, lo que se denota $f \in \mbox{Diff }^1(M)$ y significa que   $f$ es de clase $C^1$, invertible, y su inversa $f^{-1}: M \mapsto M$ tambi\'{e}n es de clase $C^1$. Si adem\'{a}s $f$ y $f^{-1}$ fueran de clase $C^r$ para alg\'{u}n natural $r > 1$, se denota $f \in \mbox{Diff }^r(M)$ (para lo cual se requiere que $M$ sea una variedad de clase $C^r$ por lo menos).
En   esta secci\'{o}n, asumiremos que     $f \in \mbox{Diff }^1(M)$  e indicaremos expresamente cuando adem\'{a}s $f \in \mbox{Diff }^r(M)$ para alg\'{u}n $r > 1$.

\begin{definition} \em \label{definicionAnosov} {\bf Difeomorfismos de Anosov \cite{Anosov}}

\index{hiperbolicidad! uniforme}
\index{transformaci\'{o}n! hiperb\'{o}lica! uniforme}
\index{difeomorfismos! de Anosov}
\index{splitting! hiperb\'{o}lico}

  $f \in \mbox{Diff }^1(M)$ se llama  \em difeomorfismo de Anosov \em si existe una descomposici\'{o}n   (llamada \em splitting\em) $S \oplus U = TM$ del fibrado tangente $TM$, no trivial (i.e. $0 <\mbox{dim}(S) < \mbox{dim}(M)$), que es invariante por $df$ (i.e. $df(x) S_x = S_{f(x)}, \ \ df(x) U_x = U_{f(x)}$), y existen constantes $C >0$ y $0 < \lambda < 1 < \sigma$, tales que, para todo $x \in M$:
\begin{equation}
 \label{equationAnosovStable}
 \|df^n(x) s\| \leq C \lambda^n \|s\| \ \ \forall \ n \geq 0, \ \ \forall \  s \in S_x,\end{equation}
\begin{equation} \index{subespacio! inestable} \index{fibrado! inestable}
 \label{equationAnosovUnstable} \|df^{n}(x) u\| \geq C^{-1} \sigma^n \|u\| \ \ \forall \ n \geq 0, \ \ \forall \  u \in U_x.\end{equation}
Para cada punto $x \in M$, los subespacios $S_x$ y $U_x$ se llaman \em estable e inestable \em respectivamente, en $x$. Los subfibrados $S$ y $U$ (formados por los subespacios $S_x$ y $U_x$   al variar $x \in M$), se llaman fibrados estable e inestable  respectivamente. La constante $ 0 <\lambda < 1$ se llama  tasa o coeficiente de contracci\'{o}n (uniforme) en el futuro a lo largo del fibrado estable, y la constante $\sigma > 1$, tasa o coeficiente de dilataci\'{o}n (uniforme) en el futuro a lo largo del fibrado inestable.

{\bf Nota: } En la definici\'{o}n de difeomorfismo de Anosov, la constante ${C} $  no depende de $x$, as\'{\i} como tampoco dependen de $x$ los coeficientes $\lambda$ y $\sigma$ de contracci\'{o}n y dilataci\'{o}n respectivamente.  Siendo  la variedad $M$ compacta,  si se cambia la m\'{e}trica Riemanniana, la norma de los vectores en cada subespacio tangente se cambia por una equivalente. Luego, se modifica la constante $C$ por otra constante $C'$ (que tampoco depende de $x$) manteni\'{e}ndose las mismas tasas $\lambda$ y $\sigma$ de contracci\'{o}n y dilataci\'{o}n respectivamente. Esto permite que la definici\'{o}n de difeomorfismo de Anosov sea intr\'{\i}nseca al difeomorfismo, y no dependa de la m\'{e}trica riemanniana elegida. Se puede probar que para todo $f$ de Anosov existe una m\'{e}trica riemanniana (llamada \em m\'{e}trica adaptada a    $f$\em) para la cual la constante $C$ puede tomarse igual a 1).
\end{definition}

\begin{exercise}\em \index{automorfismo!
lineal del toro}
Probar que el autormorfismo lineal $f =  \left(%
\begin{array}{cc}
  2 & 1 \\
  1 & 1 \\
\end{array}%
\right)$     en el toro $\mathbb{T}^2$ es un difeomorfismo de Anosov. Sugerencia: Tomar $S_x$ y $U_x$ las direcciones propias de la matriz $A$ que define a $f$,  $\lambda$ y $\mu$ los valores propios respectivos y $C= 1$.
\end{exercise}

\begin{exercise}\em
Sea $f$ un difeomorfismo de Anosov. Probar que la inversa $0 <\sigma^{-1} < 1$ de la tasa $\sigma$ de dilataci\'{o}n hacia el futuro a lo largo del fibrado inestable, es una tasa de contracci\'{o}n hacia el pasado  a lo largo del mismo fibrado. Precisamente:
\begin{equation}
 \label{equationAnosovUnstableBB} \|df^{-n}(x) u\| \leq C  \sigma^{-n} \|u\| \ \ \forall \ n \geq 0, \ \ \forall \  u \in U_x, \ \ \ \forall \ x \in M.\end{equation}

 An\'{a}logamente, probar que la inversa   $\lambda^{-1} >1 $ de la tasa contracci\'{o}n $\lambda$ hacia el futuro  a lo largo del fibrado estable, es una tasa de dilataci\'{o}n hacia el futuro a lo largo del mismo fibrado estable.

Deducir que   $f: M \mapsto M$ es un difeomorfismo de Anosov, si y solo si $f^{-1}$ tambi\'{e}n lo es.

\end{exercise}

\begin{exercise}\em \label{ejercicioUnicidadSplittingHiperbolico}
\index{splitting! hiperb\'{o}lico}
Mostrar que el splitting $S_x \oplus U_x$ es \'{u}nico.
Sugerencia: Probar, para todas las direcciones $[s] \subset S_x$ y $[v]  \subset T_xM \setminus S_x$, las siguientes desigualdades:

 $\limsup_{n \rightarrow + \infty} \frac{1}{n}{\log \|df^n s\|}  < 0 < \limsup_{n \rightarrow + \infty} \frac{1}{n}{\log \|df^n v\|}. $
\end{exercise}
\begin{exercise}\em \label{ejercicioContinuidadSplittingHiperbolico}
Probar que las aplicaciones $x \mapsto S_x$ y $x \mapsto U_x$ son continuas y deducir que si $M$ es conexa, entonces $\mbox{dim}S_x$ y $\mbox{dim}U_x$ son constantes.

Sugerencia: Considerar una  sucesi\'{o}n $\{x_n\}$ convergente a $x$ y tomar una subsucesi\'{o}n de ella tales que $S_{x_n}$ y $U_{x_n}$ tengan dimensiones constantes con $n$ y sean convergentes. Mostrar que  los subespacios $\lim_n S_{x_n} $ y $\lim_n U_{x_n}$  satisfacen las desigualdades de hiperbolicidad uniforme, y   forman un splitting de $T_xM$. Finalmente, usar la unicidad del splitting hiperb\'{o}lico en el punto $x \in M$.
\end{exercise}

\begin{exercise}\em
Probar que en la Definici\'{o}n \ref{definicionAnosov}, la condici\'{o}n $0 <\mbox{dim}(S_x) < \mbox{dim}(M)$ es redundante.

Sugerencia: Por absurdo, si $\mbox{dim}(S_x)= \mbox{dim}(M)$, elegir $N \geq 1$ tal que $C \lambda^N < 1/4$. Usando la definici\'{o}n de diferenciabilidad de $f$ en el punto $x$, probar que   existe $\delta_x >0$   tal  que: $$\mbox{dist}(x,y) < \delta_x \ \Rightarrow \ \mbox{dist}(f(x), f(y)) \leq \frac{ \mbox{dist}(x,y)}{2}.$$
Cubrir $M$ con una cantidad $k$ finita de bolas abiertas $\{B_{\delta_i}(x_i)\}_{1 \leq i \leq k} $ de centros $x_1, \ldots, x_i, \ldots x_k$ y radios $\delta_i := \delta_{x_i}$.
Probar que para todo  $y \in M$, y para todo $n \in \mathbb{N}$ existe $x_i$ (que depende de $y$ y de $n$) tal que $\mbox{dist}(f^{nN}(x_i), f^{nN}(y)) < \epsilon.$ Deducir que $f^{n N}(M)$  se puede cubrir con $k$ bolas de radio $\epsilon$. Siendo $f^{n N}$ un difeomorfismo, toda la variedad $M$ se puede cubrir con $k$ bolas de radio $\epsilon >0$, siendo $k$ constante y $\epsilon >0$ arbitrario. Deducir que $\mbox{diam}(M) < k\epsilon$ para todo $\epsilon >0$, lo cual es una contradicci\'{o}n pues el di\'{a}metro de $M$ es positivo.
\end{exercise}

\begin{exercise}\em \label{ejercicioAnosovExponentesLyapunov} \index{exponentes de Lyapunov}
Sea $f: M \mapsto M$ un difeomorfismo de Anosov con tasa $ \lambda  < 1$ de contracci\'{o}n a lo largo del fibrado estable $S$, y tasa $\sigma > 1$ de dilataci\'{o}n a lo largo del fibrado inestable $U$. Probar que para todo $x \in M$:
$$\limsup_{n \rightarrow + \infty} \frac{\log \|df^n(x) s\|}{n} \leq \log \lambda < 0 \ \forall  s \in S_x  $$
$$\limsup_{n \rightarrow - \infty} \frac{\log \|df^n(x) s\|}{-n} \leq \log \lambda < 0 \ \forall  s \in S_x  $$
$$\liminf_{n \rightarrow + \infty} \frac{\log \|df^n(x) u\|}{n} \geq \log \sigma > 0 \ \forall  u \in U_x  $$
$$\liminf_{n \rightarrow - \infty} \frac{\log \|df^n(x) u\|}{-n} \geq \log \sigma > 0 \ \forall  u \in U_x  $$
\end{exercise}

{\bf Exponentes de Lyapunov para difeomorfismos de Anosov: } \index{exponentes de Lyapunov}
\index{teorema! Oseledets}

En las pr\'{o}ximas secciones     enunciaremos el teorema de Oseledets que establece que los l\'{\i}mites superior e inferior de las desigualdades de arriba, son l\'{\i}mites que existen para $\mu$-casi todo punto $x \in M$ (donde $\mu$ es cualquier medida de probabilidad invariante por $f$). Esos l\'{\i}mites se llaman \em exponentes de Lyapunov \em  (en el futuro) para la \'{o}rbita de $x \in M$.

 Por lo tanto, admitiendo el teorema de Oseledets, de las desigualdades de arriba se deducen los  siguientes enunciados, que satisface todo difeomorfismo de Anosov $f$   para toda medida de probabilidad $\mu$ que sea $f$-invariante:

{\bf (a)} \em  Los exponentes de Lyapunov   de la \'{o}rbita por $\mu$- casi todo punto $x \in M$, correspondientes a las direcciones $u   $ del subespacio inestable, son positivos y mayores o iguales que el logaritmo del coeficiente de dilataci\'{o}n $\sigma > 1$ del difeomorfismo de Anosov. \em

\index{exponentes de Lyapunov! no nulos} \index{exponentes de Lyapunov! positivos}

{\bf (b)} \em  Los exponentes de Lyapunov   de la \'{o}rbita por $\mu$- casi todo punto $x \in M$, correspondientes a las direcciones $s   $ del subespacio  estable, son negativos y menores iguales que el logaritmo del coeficiente de contracci\'{o}n $\lambda < 1$ del difeomorfismo de Anosov. \em

De los enunciados (a) y (b), teniendo en cuenta que por definici\'{o}n de difeomorfismo de Anosov, el espacio tangente $T_xM$ es la suma directa de los subespacios estable $S_x$ e intestable $U_x$, se deduce el siguiente enunciado  (para una demostraci\'{o}n del mismo usar el teorema de Oseledets, que probaremos en el pr\'{o}ximo cap\'{\i}tulo, y ver Ejercicio  \ref{ejercicioExponentesLyapunov}):

{\bf (c)} \em Los exponentes de Lyapunov para un difeomorfismo de Anosov no son cero, y est\'{a}n \lq\lq bounded away from zero\rq\rq \ \em (i.e. est\'{a}n fuera de un entorno de cero).

\begin{remark} \em .
\index{transitividad}

  {\bf Sobre la transitividad de los difeomorfismos de Anosov.}
  Una conjetura cuya demostraci\'{o}n es a\'{u}n un problema abierto es la siguiente:

{\bf Conjetura: } Los difeomorfismos de Anosov en variedades compactas y conexas son transitivos.

Se conocen pruebas parciales: si la variedad $M$ donde act\'{u}a el difeomorfismo de Anosov es un toro $\mathbb{T}^n$, entonces $f$ es transitivo. Otro caso conocido: si la dimensi\'{o}n del fibrado inestable o la del estable es uno (de un difeomorfismo $f$ de Anosov), entonces $f$ es transitivo.  En gene\-ral, no se conocen ejemplos de difeomorfismos $f$ de Anosov que no sean transitivos.
\end{remark}

\begin{remark} \em .

  {\bf Sobre la ergodicidad de los difeomorfismos de Anosov.} \index{ergodicidad} \index{difeomorfismos! de Anosov}
  \index{transformaci\'{o}n! erg\'{o}dica}

Sea $M$ una variedad compacta y conexa de dimensi\'{o}n $n \geq 2$. Si $f \in \mbox{Diff }^2(M)$ es un difeomorfismo de Anosov  (de clase $C^2$), si $f$ es transitivo y si $f$ preserva la medida de Lebesgue $m$, entonces $m$ es erg\'{o}dica.    Demostraremos este resultado en el Corolario \ref{corolarioAnosovTransitivomedidaLebesgue}, al definir y estudiar las medidas invariantes de probabilidad llamadas SRB (Sinai-Ruelle-Bowen)  para los difeomorfismos de Anosov.

\end{remark}

\subsection{Conjuntos uniformemente hiperb\'{o}licos} \index{hiperbolicidad! uniforme}
\index{transformaci\'{o}n! hiperb\'{o}lica! uniforme}
\index{conjunto! hiperb\'{o}lico! uniforme}

\begin{definition}
 \em  \label{definicionHiperbolicidadUniforme}
  Sea $f: M \mapsto M$  y sea $\Lambda \subset M$ un subconjunto compacto e invariante (i.e. $f^{-1} (\Lambda) = \Lambda$).
El conjunto $\Lambda $ se llama \em uniformemente hiperb\'{o}lico \em  (en breve, hiperb\'{o}lico) si para todo $x \in \Lambda$ existe un splitting $S_x \oplus U_x = T_xM$ del espacio tangente $T_xM$ a $M$ en $x$, que dependen continuamente de $x \in \Lambda$, y constantes $C >0$ y  $ 0 < \lambda < 1 < \sigma$, que verifican las desigualdades (\ref{equationAnosovStable}) y (\ref{equationAnosovUnstable}) para todo $x \in \Lambda$. \index{splitting! hiperb\'{o}lico}

\index{subespacio! inestable}
\index{fibrado! inestable}

{\bf Nota: } La condici\'{o}n de dependencia continua de $S_x$ y $U_x$ al variar $x \in \Lambda$ es redundante. Se puede demostrar esta continuidad a partir de las invariancia de esos subespacios, de las desigualdades de hiperbolicidad uniforme y de la compacidad del conjunto $\Lambda$ (ver Ejercicios \ref{ejercicioUnicidadSplittingHiperbolico} y \ref{ejercicioContinuidadSplittingHiperbolico} y generalizarlos sustituyendo $M$ por $\Lambda$).

{\bf Sobre las dimensiones de   $S_x$ y $U_x$:} Son complementarias  (su suma es $\mbox{dim} M$), pero no necesariamente son ambas mayores que cero. Es decir, alguno de los dos subespacios puede tener dimensi\'{o}n cero, y el otro coincidir con $T_xM$. Las dimensiones pueden depender de $x \in \Lambda$. La dependencia continua de los subespacios invariantes implica que las dimensiones de $S_x$ y $U_x$ sean localmente constantes. Como $\Lambda$ es compacto, si existen  varios subconjuntos de $\Lambda$ para los cuales $S_x$ y $U_x$ tienen dimensiones diferentes, entonces son dos a dos aislados.

\end{definition}

De las definiciones anteriores, es inmediato deducir que $f: M \mapsto M$ es un difeomorfismo de Anosov, si y solo si   la variedad $M$ es un conjunto uniformemente hiperb\'{o}lico (con dimensiones estable e inestable no nulas).

\index{difeomorfismos! de Anosov}

\begin{exercise}\em \label{exercisepozosillafuente} \index{pozo} \index{fuente} \index{silla} \index{punto! peri\'{o}dico! hiperb\'{o}lico}
Sea $f: M \mapsto M$. Probar que   un punto peri\'{o}dico $x$ de per\'{\i}odo $p$ (i.e. $f^p(x)= x, \ f^j(x) \neq x \ \forall \ 1 \leq j < p $ es hiperb\'{o}lico (i.e. los valores propios de $df^p$ tienen m\'{o}dulo diferente de 1), si y solo si su \'{o}rbita (finita) $\{x, f(x), \ldots, f^{p-1}(x)\}$ es un conjunto uniformemente hiperb\'{o}lico. Probar que si $x$ es punto silla (es decir, $df^p(x)$ tiene valores propios con m\'{o}dulo menor que uno y tambi\'{e}n con m\'{o}dulo mayor que uno), entonces los subespacios propios de $df^p(x)$ son los subespacios estable $S_x$ e inestable $U_x$, respectivamente. Probar que si $x$ es un pozo  hiperb\'{o}lico (i.e. $df^p(x)$ tiene todos sus valores propios con m\'{o}dulo menor que uno), entonces $S_x = T_xM$ y $U_x = {\bf 0}$. Enunciar y probar resultado dual si $x$ es una fuente hiperb\'{o}lica (i.e. $df^p(x)$ tiene todos sus valores propios con m\'{o}dulo mayor que uno).
\end{exercise}

 \subsection{Ejemplo: Herradura de Smale lineal}  \index{herradura de Smale} \label{sectionHerraduraSmale}
Un ejemplo paradigm\'{a}tico de conjunto hiperb\'{o}lico (transitivo) que no es toda la variedad, es el siguiente, debido a Smale \cite{Smale} (ver tambi\'{e}n, por ejemplo, \cite[pages 97-98]{Jost}):

\begin{definition} \label{definicionHerradura}
\em {\bf Herradura de Smale lineal.}

Se llama \em herradura de Smale lineal   \em  (en
dimensi\'{o}n 2 y con 2 patas)  a un difeomorfismo  $T: Q= [0,1]^2 \subset \mathbb{R}^2 \mapsto T(Q) \subset \mathbb{R}^2$     que satisface las siguientes condiciones (ver Figura \ref{figuraHerraduraSmale}):
\begin{itemize}
\item [a) ] $T(Q)\cap Q = Q_0 \cup Q_1$ donde $Q= [0,1]^2, $

$ Q_0=[1/5,2/5]\times [0,1], \; Q_1 =
[3/5,4/5]\times [0,1] $.

 \item[b) ] $T^{-1}(Q_0) = [0,1] \times
[1/5,2/5], \; T^{-1}(Q_1) = [0,1] \times [3/5,4/5] $.

\item [c) ] Para $j=0,1$, la restricci\'{o}n $T|_{T^{-1}(Q_j)} (x,y) $
 es una transformaci\'{o}n af\'{\i}n en $(x,y)$ con
direcciones propias (1,0) y (0,1) y valores propios $\lambda$ y
$\mu$ reales tales que $|\lambda | = 1/5$ y $|\mu | = 5$
respectivamente. Por ejemplo:
$$T|_{T^{-1}(Q_0)} (x,y) = ((1/5)(x+1), 5 y -1 ),$$ $$
T|_{T^{-1}(Q_1)} (x,y) = ((-1/5)(x-4), -5 y +4 )$$
\end{itemize}

Para comprender c\'{o}mo es la transformaci\'{o}n $T$, ver la figura \ref{figuraHerraduraSmale} e imaginar $T$ como la composici\'{o}n de dos transformaciones: Primero considerar una transformaci\'{o}n af\'{\i}n que lleva el cuadrado $Q= [0,1]^2$ a un rect\'{a}ngulo 5 veces m\'{a}s alto y 5 veces menos ancho que el cuadrado $Q $ (contrae en horizontal y dilata en vertical). Despu\'{e}s aplicar una transformaci\'{o}n que \rq\rq dobla\rq\rq en forma biyectiva al rect\'{a}ngulo alto y flaco, d\'{a}ndole forma de herradura, superponiendo esta sobre el cuadrado $Q$ en $Q_0$ y $Q_1$ (sin deformar $Q_0$ ni $Q_1$).

 Se observa que hemos restringido la definici\'{o}n eligiendo valores num\'{e}ricos fijos para $|\lambda|$ y $|\mu|$, iguales a $1/5$ y $5$ respectivamente. Sin embargo, si se toman otros valores num\'{e}ricos, $0 < |\lambda| < 1/2$ y $2 < |\sigma|$, y se definen (coherentemente con esos nuevos valores num\'{e}ricos) los rect\'{a}ngulos compactos $Q_0 $ y $Q_1$ disjuntos, como en la Figura \ref{figuraHerraduraSmale}, y tales $T(Q) \cap Q = Q_0 \cup Q_1$; y si las  transformaciones $T|_{T^{-1}(Q_0)}$ y $T|_{T^{-1}(Q_1)}$ son afines con valores propios $ \pm\lambda $ y $\pm \sigma$, entonces $T$ tambi\'{e}n se llama \em herradura de Smale lineal.\em

\end{definition}


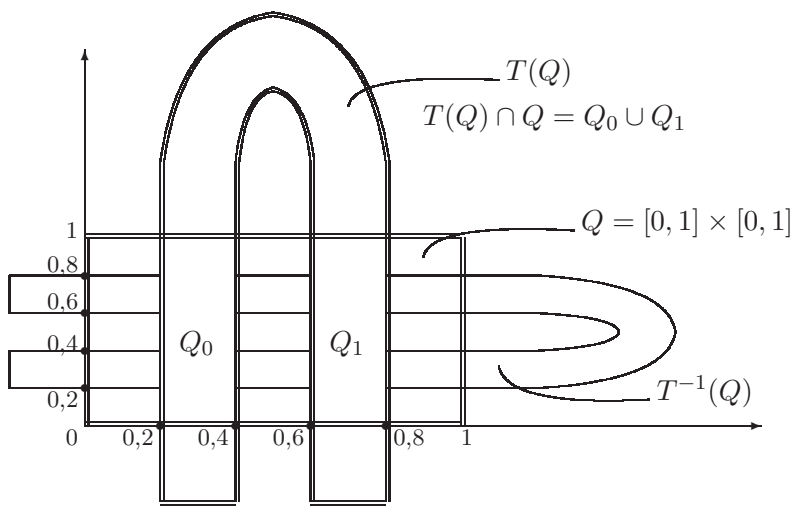
\begin{figure}[h]
\setlength{\unitlength}{0.5truecm}

\begin{picture}(18,12)(0,0)

\put(2.5,0.5){\footnotesize $0$}
\put(3,1){\vector(3,0){18}}

\put(3,1){\vector(0,3){10}}

\put(3,6){\line(3,0){10}}
\put(13,1){\line(0,3){5}}

\put(3,6.1){\line(3,0){10}}
\put(13.1,1){\line(0,3){5}}
 \put(3,1.1){\line(3,0){10}}
\put(3.1,1){\line(0,3){5}}

\put(5,1){\circle*{.2}}
\put(7,1){\circle*{.2}}
\put(9,1){\circle*{.2}}
\put(11,1){\circle*{.2}}

\put(4.0,0.5){\footnotesize $0.2$}
\put(6,0.5){\footnotesize $0.4$}
\put(8,0.5){\footnotesize $0.6$}
\put(11.2,0.5){\footnotesize $0.8$}
\put(13,0.5){\footnotesize $1$}

\put (5.5,3) {$Q_0$}
\put (9.5,3) {$Q_1$}

\put(5,-1){\line(0,1){9}}
\put(7,-1){\line(0,1){9}}
\put(9,-1){\line(0,1){9}}
\put(11,-1){\line(0,1){9}}

\put (5,-1){\line(1,0){2}}
\put (9,-1){\line(1,0){2}}

\put(5.1,-1){\line(0,1){9}}
\put(7.1,-1){\line(0,1){9}}
\put(9.1,-1){\line(0,1){9}}
\put(11.1,-1){\line(0,1){9}}

\put (5,-1.1){\line(1,0){2}}
\put (9,-1.1){\line(1,0){2}}

{\qbezier(5,8)(5.5,11.6)(8,12)}
{\qbezier(11.1,8)(10.6,11.6)(8,12)}
{\qbezier(5.1,8)(5.6,11.6)(8,11.9)}
{\qbezier(11,8)(10.5,11.6)(8,11.9)}

{\qbezier(7,8)(7.2,9.8)(8,10)}
{\qbezier(9.1,8)(8.9,9.8)(8,10)}
{\qbezier(7.1,8)(7.3,9.8)(8,9.9)}
{\qbezier(9,8)(8.8,9.8)(8,9.9)}

\put(3,2){\circle*{.2}}
\put(3,3){\circle*{.2}}
\put(3,4){\circle*{.2}}
\put(3,5){\circle*{.2}}

\put(2,1.5){\footnotesize $0.2$}
\put(2,3){\footnotesize $0.4$}
\put(2,4.1){\footnotesize $0.6$}
\put(2,5.1){\footnotesize $0.8$}
\put(2.5,6){\footnotesize $1$}

\put(1,2){\line(1,0){4}}
\put(1,3){\line(1,0){4}}
\put(1,4){\line(1,0){4}}
\put(1,5){\line(1,0){4}}

\put (1,2){\line(0,1){1}}
\put (1,4){\line(0,1){1}}

\put(7,2){\line(1,0){2}}
\put(7,3){\line(1,0){2}}
\put(7,4){\line(1,0){2}}
\put(7,5){\line(1,0){2}}

\put(11,2){\line(1,0){4}}
\put(11,3){\line(1,0){4}}
\put(11,4){\line(1,0){4}}
\put(11,5){\line(1,0){4}}

{\qbezier(15,2)(18.5,2.2)(18.7,3.5)}
{\qbezier(15,5)(18.5,4.8)(18.7,3.5)}

{\qbezier(15,3)(17.1,3.1)(17.2,3.5)}
{\qbezier(15,4)(17.1,3.9)(17.2,3.5)}

{\qbezier(12,5.5)(12.2,6.4)(16,6.2)}
 \put(16.2, 6.2){$Q = [0,1] \times [0,1]$}

 {\qbezier(10,9.5)(10.2,10.4)(14,10.2)}
 \put(14.2, 10.2){$T(Q)$}

 \put (12, 9) {$  T(Q) \cap Q = Q_0 \cup Q_1$}

{\qbezier(14,2.6)(14.2,1.5)(18,1.7)}
 \put(18.2, 1.7){$T^{-1}(Q)$}

\end{picture}

\vspace{.5cm}

\caption{\small Herradura de Smale}
\index{herradura de Smale}
\label{figuraHerraduraSmale}

\vspace{-.4cm}
\end{figure}


\begin{exercise}\em
  Sea $T: \mathbb{R}^2 \mapsto \mathbb{R}^2$ una herradura de Smale lineal. Sea $Q = [0,1]^2$.

(a) Dibujar esquem\'{a}ticamente $T(T(Q) \cap Q) \subset T(Q)$, $T^2(Q)$ y $\cap_{n= 0}^N T^n(Q)$ para $N= 2, 3$.

(b)   Dibujar  esquem\'{a}ticamente)
$T^{-1}(Q) \cap Q, $  y $\cap_{n=0}^{N} T^{-n} Q $ para $N= 2, 3$.

(c) Dibujar  esquem\'{a}ticamente)
 el conjunto \lq\lq estable\rq\rq \ $  W^s \cap Q $ de todos los puntos de $Q$ cuyas
\'{o}rbitas futuras permanecen en $Q$ para todos los iterados $n
\geq 0$.   (Sugerencia: Ver que  $W^s \cap Q =   \cap_{n= 0}^{+ \infty}T^{-n}(Q)$.)

(d) Calcular el {\bf exponente de Lyapunov  negativo (tasa exponencial de contracci\'{o}n hacia el futuro)}
de la herradura de Smale: $\lim _{n \rightarrow + \infty} (1/n)
\log \|DT^n _{x} (1,0) \|$ para todo punto $x \in \bigcap_{n \in
\mathbb{N}}T^{-n}(Q)$.

(d)   Definir el conjunto \lq\lq inestable\rq\rq \ $W^u \cap Q$ y el exponente de
Lyapunov  positivo (tasa exponencial de dilataci\'{o}n hacia el futuro.) Calcularlo.
\end{exercise}

\begin{definition} \em \label{definicionMaximalHerraduraSmale} Sea $T: Q =[0,1]^2 \subset \mathbb{R}^2 \mapsto \mathbb{R}^2$
es la herradura de Smale lineal.  Se llama \em conjunto invariante maximal de $T$ en $Q$ \em a
$$\Lambda : = \bigcap _{n \in \mathbb{Z}} T^{-n}Q.$$ Obs\'{e}rvese que por la propiedad de intersecciones finitas no vac\'{\i}as
de compactos, el conjunto $\Lambda$ es compacto no vac\'{\i}o. \index{conjunto! maximal invariante} \index{maximal invariante}

\begin{exercise}\em \label{ejercicioHerraduraSmaleHiperbolica}
Probar que el maximal invariante $\Lambda$ de una herradura de Smale lineal, es un conjunto hiperb\'{o}lico. Probar que $\Lambda$ es el producto cartesiano de dos conjuntos de Cantor en el intervalo $[0,1]$.
\end{exercise}

 \end{definition}

\subsection{Variedades invariantes de conjuntos hi\-per\-b\'{o}licos} \index{variedad invariante! estable}
\index{variedad invariante! inestable}
\index{hiperbolicidad! uniforme}
\index{transformaci\'{o}n! hiperb\'{o}lica! uniforme}
\index{conjunto! hiperb\'{o}lico! uniforme}

A continuaci\'{o}n enunciamos  algunos teoremas sobre din\'{a}mica diferenciable, que generalizan propiedades de   los difeomorfismos de Anosov  en una variedad compacta $M$, y en particular algunos de los resultados vistos en el ejemplo del automorfismo lineal hiperb\'{o}lico $f =  \left(%
\begin{array}{cc}
  2 & 1 \\
  1 & 1 \\
\end{array}%
\right)$  en el toro $\mathbb{T}^2$.

\begin{theorem} {\bf Variedades invariantes para $\Lambda$ unif. hiperb\'{o}lico} \index{teorema! de existencia de! variedades invariantes}
\label{teoremavariedadesinvariantesAnosov}
Sea    $f: M \mapsto M$ un difeomorfismo y $\Lambda \subset M$ un conjunto  invariante, compacto y   uniformemente hiperb\'{o}lico. Para   $x \in M$ denotemos $U_x, S_x$ los subespacios inestable y estable respectivamente.   Entonces, para  cada $x \in \Lambda$ existen y son \'{u}nicas:

 {\bf A)} Una subvariedad  conexa $W^s(x)$, $C^1$-inmersa en $M$ \em (pero no  necesariamente \lq\lq embedded\rq\rq, i.e.   no  encajada     en $M$, ni necesariamente compacta), \em llamada variedad estable por $x$, tal que: \em
\begin{equation} \label{equationvarestable1} x \in W^s(x), \ \ \ f (W^s(x)) = W^s(f (x)),
 \ \ \ T_x (W^s(x)) = S_x. \end{equation} \em

{\bf B)} Una subvariedad $C^1$ conexa $W^u(x)$, $C^1$-inmersa en $M$ \em (pero  no necesariamente encajada en $M$ ni compacta), \em llamada variedad inestable por $x$, tal que: \em
\begin{equation} \label{equationvarinestable1}x \in W^u(x), \ \ \ f (W^u(x)) = W^u(f (x)), \ \ \ T_x (W^u(x)) = U_x.\end{equation}  \em

tales que: \em \begin{equation}
 \label{equationvarestable2} \lim _{n \rightarrow + \infty} \mbox{dist}(f^n(y), f^n(x)) = 0 \ \ \Leftrightarrow \ \  \ y \in W^s(x), \end{equation}
 \begin{equation} \label{equationvarinestable2}\lim _{n \rightarrow - \infty} \mbox{dist}(f^n(y), f^n(x)) = 0 \ \ \Leftrightarrow \ \  \ y \in W^u(x).\end{equation}    \em

 \em Si adem\'{a}s $f$ es un difeomorfismo $C^r$ para alg\'{u}n    $r > 1$ entonces $W^s_x$ y $W^u_x$ son subvariedades de clase $C^r$.
\end{theorem}

Una prueba del teorema de existencia de variedades estable e inestable puede encontrarse en \cite{HirschPughShub}.

\begin{remark} \em

\index{subespacio! inestable}
\index{fibrado! inestable}
\index{variedad invariante! estable}
\index{variedad invariante! inestable}

Sea   $f: M \mapsto M$ un difeomorfismo de Anosov con coeficiente $ \lambda  < 1$ de contracci\'{o}n en el fibrado estable $S$, y coeficiente $\sigma > 1$ de dilataci\'{o}n en el fibrado inestable $U$.
Admitiendo la existencia de $C^1$-subvariedades    $W^s(x)$ y $W^u(x)$ conexas, que cumplen las condiciones (\ref{equationvarestable1}) y (\ref{equationvarinestable1}) respectivamente, y sabiendo que los subespacios $S_x $ y $ U_x$ dependen continuamente de $x$, se puede demostrar  que
\begin{equation} \label{equationdistanciaestable}   \limsup _{n \rightarrow + \infty} \frac{\log \mbox{dist}(f^n(y), f^n(x))}{n} \leq  \log \lambda < 0 \ \ \forall \ y \in W^s(x); \end{equation}
\begin{equation} \label{equationdistanciainestable}   \liminf _{n \rightarrow - \infty} \frac{\log \mbox{dist}(f^n(y), f^n(x))}{n} \geq  \log \sigma >0 \ \ \forall \ y \in W^u(x).\end{equation}

\end{remark}

\begin{remark} \em  \index{difeomorfismos! de Anosov}
\index{teorema! Franks}

Volvamos a los difeomorfismos de Anosov, como caso particular de sistema uniformemente hiperb\'{o}lico:

El Teorema o Lema de Franks \cite{Franks} (ver tambi\'{e}n \cite{Manning} para difeomorfismos de Anosov en el toro n-dimensional, para todo $n \geq 2$)    establece que:

{\bf Teorema o Lema de Franks} \em Las \'{u}nicas superficies  compactas, conexas, sin borde y orientables
que soportan un difeomorfismo de Anosov $f$ son homeomorfas al
toro ${\mathbb{T}}^2$, y  $f$ es  conjugado a un automorfismo lineal hiperb\'{o}lico. \em

\index{homeomorfismos conjugados} \index{conjugaci\'{o}n} \index{difeomorfismos! conjugados}

La demostraci\'{o}n del Lema de Franks se encuentra en \cite{Franks}.   En \cite{Manning}, Manning generaliza la \'{u}ltima parte del Lema de Franks, probando que  los difeomorfismos de Anosov en el toro de dimensi\'{o}n $n$ son conjugados a automorfismos lineales. Recientemente, en \cite{Hammerlindl-LemadeFranks} se obtiene el mismo resultado que Manning, pero para difeomorfismos de Anosov en nilmanifols de dimensi\'{o}n 3, y no solo para el toro.

Nota: Dos homeomorfismos $f: X \mapsto X$ y $g: Y \mapsto Y$, en
respectivos espacios topol\'{o}gicos $X$ e $Y$,  se dice que son \em
conjugados o topol\'{o}gicamente equivalentes \em si existe un
homeormorfismo $h: X \mapsto Y$ llamado \em conjugaci\'{o}n \em tal
que $g \circ h= h \circ f$. La conjugaci\'{o}n implica que cada una de las
propiedades de la din\'{a}mica topol\'{o}gica de $f$ (transitividad,
conjunto no errante, omega y alfa l\'{\i}mites, recurrencia) se
satisface  para $f$ si y solo si se satisface  para $g$.

Debido al teorema de Franks los automorfismos lineales en el toro
${\mathbb{T}}^2$ son el paradigma de los difeomorfismos de Anosov en
superficies, y entre ellos en particular el $\left(%
\begin{array}{cc}
  2 & 1 \\
  1 & 1 \\
\end{array}%
\right)$ en el toro ${\mathbb{T}}^2$.

En cambio,   conjuntos hiperb\'{o}licos como la herradura de Smale, o topol\'{o}gicamente conjugadas a ella, se pueden construir en un abierto homeormorfo a un cuadrado de $\mathbb{R}^2$, en cualquier superficie
(variedad de dimensi\'{o}n dos).

Se puede generalizar la herradura de Smale a dimensiones mayores que dos, tomando un cubo  $n$-dimensional en el rol de cuadrado $Q$, y eligiendo dimensiones complementarias de contracci\'{o}n uniforme y de dilataci\'{o}n uniforme.
\end{remark}

\subsection{Expansividad o caos topol\'{o}gico.} \index{expansividad} \index{caos}  \index{transformaci\'{o}n! ca\'{o}tica} \index{transformaci\'{o}n! expansiva}

Se han dado varias definiciones precisas de \em caos \em en la literatura matem\'{a}tica. Pero estas definiciones no son equivalentes entre s\'{\i} (ver por ejemplo \cite{buzziEnciclopedia}). Entonces cuando uno se refiere a \em sistema ca\'{o}tico \em deber\'{\i}a siempre incluir la definici\'{o}n que est\'{a} utilizando, con precisi\'{o}n, y observar que los resultados que se obtengan con esa definici\'{o}n, no son necesariamente ciertos si se hubiese adoptado otra. Seg\'{u}n sea el objetivo de investigaci\'{o}n de quien estudia el sistema din\'{a}mico (por ejemplo, su objetivo puede ser estudiar la din\'{a}mica topol\'{o}gica de un conjunto abierto y denso de \'{o}rbitas, o solo estudiar la de $\mu$-casi toda \'{o}rbita cuando $\mu$ es invariante, o la de Lebesgue-casi toda \'{o}rbita, cuando la medida de Lebegue no es invariante, etc), adopta una u otra definici\'{o}n.

 El estudio de las propiedades topol\'{o}gicas de los difeomorfismos
linea\-les en el toro ${\mathbb{T}}^2$   es paradigm\'{a}tico para estudiar m\'{a}s en general,
los sistemas ca\'{o}ticos desde el punto de vista topol\'{o}gico. En efecto:

\begin{definition} \em
Sea $T: X \mapsto X$ un homeomorfismo en un espacio m\'{e}trico
compacto $X$. Se dice que $T$ es \em ca\'{o}tico \em (topol\'{o}gicamente)
o \em expansivo\em, si existe una constante $\alpha >0$ (llamada
\em constante de expansividad) \em  tal que
$$\sup_{n \in \mathbb{Z}} \dist (T^n(x), T^n(z)) > \alpha \; \; \forall x \neq y \in X$$
\end{definition}

{\bf Interpretaci\'{o}n: } El caos topol\'{o}gico o expansividad es una versi\'{o}n de \lq\lq hiperbolicidad topol\'{o}gica.\rq\rq \ Significa que la \em
din\'{a}mica es $\alpha$-sensible a las condiciones iniciales: \em dos
\'{o}rbitas con estados iniciales $x \neq y $ diferentes se separan
m\'{a}s que la distancia $\alpha $, hacia el futuro \'{o} hacia el pasado. Podr\'{a} haber parejas de \'{o}rbitas, con estados iniciales pr\'{o}ximos $x \neq y$, que se separan solo para el futuro pero no para el pasado, o al rev\'{e}s. Generalmente (por ejemplo en un difeomorfismo de Anosov), para la mayor parte de las parejas de estados iniciales $x    \neq y$, sus \'{o}rbitas se separan en el futuro y en el pasado.

La separaci\'{o}n hacia el futuro o hacia el pasado, m\'{a}s que una constante uniforme $\alpha$, de dos \'{o}rbitas con estados iniciales $x \neq y$ (por m\'{a}s que $x $ e $y $ est\'{e}n tan cerca entre s\'{\i} como se desee)  implica  lo siguiente:

Si uno aproxima con error
$\epsilon
>0$, aunque sea arbitrariamente peque\~{n}o,   el estado inicial $x$ sustituy\'{e}ndolo por $y$ tal que $0 <\dist (x,y)
\leq \epsilon$, entonces el estado del sistema $T^n(y)$ en alg\'{u}n
instante $n \in \mathbb{Z}$,  se modificar\'{a} m\'{a}s que $\alpha >0 $ respecto
del  estado $T^n(x)$ que tendr\'{\i}a si no se hubiese cometido el
error en el estado inicial. La sensibilidad a las condiciones
iniciales, o expansividad, o caos topol\'{o}gico,  se llama tambi\'{e}n
\lq \lq efecto mariposa": la leve modificaci\'{o}n del estado inicial
producido por el aleteo de una mariposa, por m\'{a}s leve que esta sea (es decir por m\'{a}s peque\~{n}o que sea la diferencia $\epsilon >0$ provocada en ese estado inicial) produce una modificaci\'{o}n
\lq\lq dr\'{a}stica\rq\rq \ (es decir mayor que una constante uniforme $\alpha >0$) en el estado del sistema en otro instante).
\begin{exercise}\em \index{rotaci\'{o}n! racional} \index{rotaci\'{o}n! irracional}
Demostrar que   la rotaci\'{o}n del c\'{\i}rculo   (racional o
irracional)  no es expansiva (i.e. no es topol\'{o}gicamente ca\'{o}tica). \em
\end{exercise}

\begin{exercise}\em \index{tent map}
Demostrar que el tent map $f: [0,1] \mapsto [0,1]$ es expansivo para el futuro, esto es, existe una constante $\alpha >0$ tal que   $$\mbox{Si \ dist}(f^j(x), f^j(y)) < \alpha \ \forall \ j \in \mathbb{N} \ \mbox{ entonces } \ x= y.$$
Sugerencia: Probar que la derivada $|(f^n)'|$ tiende uniformemente a $+ \infty$ con $n$ y usar $\alpha= 1/4$.
\end{exercise}
\begin{exercise}\em  \index{automorfismo!
lineal del toro} Demostrar que la transformaci\'{o}n $f = \left(%
\begin{array}{cc}
  2 & 1 \\
  1 & 1 \\
\end{array}%
\right)$ en el toro ${\mathbb{T}}^2$ es   expansiva.
\end{exercise}

\begin{remark} \em

En \cite{Lewowicz-Expansivos}, Lewowicz   demostr\'{o} los siguientes resultados:

\vspace{.2cm}   \index{expansividad} \index{transformaci\'{o}n! expansiva} \index{teorema! Lewowicz}

{\bf Teorema de Lewowicz} \em Las \'{u}nicas superficies \em (variedades de dimensi\'{o}n dos) \em compactas y conexas donde existen homeomorfismos expansivos, son homeomorfas al toro $\mathbb{T}^2$. \em

\vspace{.2cm}

\em  Todos los homeomorfismos
 expansivos en el toro ${\mathbb{T}}^2$ son
conjugados a  un Anosov \em (y por el teorema de Franks son conjugados
a un difeomorfismo de Anosov lineal).   Por eso resulta   paradigm\'{a}tico
estudiar los difeomorfismos lineales, y el $ \left(%
\begin{array}{cc}
  2 & 1 \\
  1 & 1 \\
\end{array}%
\right)$ en particular.

\end{remark}

\subsection{Foliaciones invariantes para dif. de Anosov.} \index{difeomorfismos! de Anosov} \index{foliaci\'{o}n! invariante! estable}
\index{foliaci\'{o}n! invariante! inestable}

 \label{sectionFoliacionesInvariantesDifeosAnosov}

La unicidad de las variedades invariantes para los difeomorfismos de Anosov, implica que la siguiente colecci\'{o}n no numerable de subconjuntos (subvariedades) $\{W^s_x\}_{x \in M}$, sea una partici\'{o}n de $M$. En efecto:

\begin{proposition}
 Sea $f: M \mapsto M$ un difeomorfismo de Anosov, y sean $x, y \in M$.
Entonces o bien $  W^s_x \cap W^s_y = \emptyset$, o bien $W^s_x = W^s_y$.

An\'{a}logamente para $W^u_x$ y $W^u_y$.
\end{proposition}
{\em Demostraci\'{o}n: } Debido a (\ref{equationvarestable2}) y (\ref{equationvarinestable2}),   $y \in W^s_x$ si y solo si $x \in W^s_y$.
Supongamos que existe $z \in W^s_x \cap W^s_y$.  Sea $x' \in W^s_x$. Usando nuevamente   (\ref{equationvarestable2}) y (\ref{equationvarinestable2}) y la propiedad triangular, deducimos que $$\lim_{ n \rightarrow + \infty}  \mbox{dist}(f^n(y), f^n(x')) = 0.$$  En efecto: $$\mbox{dist}(f^n(y), f^n(x') \leq $$ $$ \mbox{dist} (f^n(y), f^n(z))+\mbox{dist} (f^n(z), f^n(x))+\mbox{dist} (f^n(x), f^n(x')),$$ y estos tres sumandos tienden a cero cuando $n \rightarrow + \infty$.
Luego $x' \in W^s_y$ para todo $x' \in W^s_x$, probando que $W^s_x \subset W^s_y$. Sim\'{e}tricamente, $W^s_ y \subset W^s_x$. Hemos demostrado que si existe $z \in W^s_x \cap W^s_y$ entonces $W^s_y = W^s_x$.
\hfill $\Box$

\begin{definition} \em  {\bf Foliaciones estable e inestable}.\index{foliaci\'{o}n! trivializaci\'{o}n de}

Se llama \em foliaci\'{o}n estable \em a la partici\'{o}n $\{W^s_x\}_{x \in M}$ de $M$ en variedades estables. An\'{a}logamente, se llama \em foliaci\'{o}n inestable \em a la partici\'{o}n $\{W^u_x\}_{x \in M}$ de $M$ en variedades inestables.

{\bf Nota sobre la definici\'{o}n geom\'{e}trica de \lq\lq foliaci\'{o}n\rq\rq: } Por definici\'{o}n, una foliaci\'{o}n  tiene  una estructura geom\'{e}trica precisa que transciende a la mera partici\'{o}n del espacio $M$ en subvariedades inmersas   disjuntas  dos a dos todas de la misma dimensi\'{o}n $1 \leq k_1 < \mbox{dim} M  $. Estas variedades inmersas se llaman hojas de la foliaci\'{o}n. La estructura geom\'{e}trica de foliaci\'{o}n consiste  en la existencia de un atlas de cartas locales $\xi$ de $M$ que son \lq\lq trivializadoras\rq\rq \  de las hojas   de la foliaci\'{o}n: i.e. la imagen de cada $\xi$ es el producto cartesiano $B_{k_1} \times B_{k_2}$, donde $B_{k_i}$ es la bola unitaria de $R^{k_i}$ ($k_1 + k_2 = \mbox{dim}M$), tal que transforma la  componentes conexas locales  de cada hoja, en las secciones $B_{k_1} \times \{v_2\}$ (con $v_2 \in B_{k_2}$ fijo).

\index{foliaci\'{o}n! din\'{a}micamente definida}

En general, para las foliaciones din\'{a}micamente definidas (por ejemplo, las foliaciones estable e inestable de un $f$ de Anosov), las cartas locales trivializadoras  $\xi$ existen, pero son solo homeormofismos sobre sus im\'{a}genes  (no son siquiera $C^1$). Como   las hojas de la foliaci\'{o}n son subvariedades inmersas de clase $C^{1}$, entonces las \em restricciones \em de las   cartas trivializadoras   a las hojas locales, son $C^{1}$.
\end{definition}
\begin{definition} \em {\bf Regularidad de   foliaciones invariantes.} \index{foliaci\'{o}n! regular}
\index{foliaci\'{o}n! de clase $C^0$}

Un partici\'{o}n ${\mathcal W} := \{W _x\}_{x \in M}$ del espacio $M$ en subvariedades diferenciables inmersas en $M$ y disjuntas dos a dos, todas de la misma dimensi\'{o}n, se dice que es \em una foliaci\'{o}n invariante $C^0$ \em (desde el punto de vista din\'{a}mico) si:

 (1) $f (W_x) = W_{f(x)} \ \ \forall \ x \in M$ (invariancia),

 (2) cada hoja $W_x$ es de clase $C^1$, y

 (3)  la aplicaci\'{o}n $x \in M \mapsto T_x W_x$ que lleva cada punto $x$ en el subespacio tangente $T_xW_x $ en $x$ a la hoja $W_x$, es continua (es decir, el subespacio tangente a la hoja var\'{\i}a en forma $C^0$ con $x \in M$).

{\bf Nota: } La condici\'{o}n (2) de que cada hoja sea de clase $C^1$, implica que su subespacio tangente var\'{\i}e en forma $C^0$ con $x$ \em variando a lo largo de la hoja respectiva. \em Pero no implica que ese subespacio tangente \em var\'{\i}e adem\'{a}s en forma $C^0$ con $x$ en todas las direcciones transversales a la hoja. \em  Por lo tanto la condici\'{o}n (3) no es redundante: es m\'{a}s fuerte que la condici\'{o}n (2).

An\'{a}logamente, para todo $r \in \mathbb{N}$, la partici\'{o}n ${\mathcal W}$ es \em una foliaci\'{o}n inva\-rian\-te $C^{r}$ \em (desde el punto de vista din\'{a}mico) si cada hoja $W_x \in {\mathcal W}$   es de clase $C^{r+1}$ y la aplicaci\'{o}n $x \in M \mapsto T_xW_x$ es de clase $C^{r}$ . Para que esta definici\'{o}n sea aplicable, hay que asumir que la variedad ambiente $M$ es de clase $C^{r+1}$ por lo menos.

\index{foliaci\'{o}n! H\"{o}lder continua}
Sea $0 <\alpha < 1$. Una foliaci\'{o}n invariante ${\mathcal W}$ es $\alpha$-H\"{o}lder continua y se denota de clase $C^{\alpha}$, si  cada hoja $W_x \in {\mathcal W}$ es de clase $C^1$ y la aplicaci\'{o}n $x \in M \mapsto T_xW_x$ satisface:
\begin{equation}
\label{eqn101}
\mbox{dist}(W_x, W_y) \leq K \, [\mbox{dist}(x,y)]^{\alpha}\end{equation}
para cierta constante $K$.
Se dice que la foliaci\'{o}n es Lipschitz y se denota como $C^{Lip} $, si es $1$-H\"{o}lder continua (es decir vale la desigualdad (\ref{eqn101}) para $\alpha= 1$.
\index{foliaci\'{o}n! Lipschitz}

\end{definition}

Debido al Teorema \ref{teoremavariedadesinvariantesAnosov}, tenemos el siguiente resultado:
 \begin{corollary} \index{foliaci\'{o}n! invariante! estable} \index{foliaci\'{o}n! invariante! inestable}

 Si $f$ es difeomorfismo de Anosov de clase $C^1$, entonces las foliaciones estable e inestable son de clase $C^0$ .
  \end{corollary}
  {\em Demostraci\'{o}n: }
  $T_x W^s_x = S_x, \ T_x W^u_x = U_x$ y los subespacios $S_x$ y $U_x$ dependen continuamente de $x$. \hfill $\Box$

Uno podr\'{\i}a esperar que si $f$ es difeomorfismo de Anosov de clase $C^{r+1}$ con $r \geq 0$, entonces las foliaciones invariantes estable e inestables sean foliaciones de clase $C^{r}$. Este \'{u}ltimo resultado es FALSO. En gene\-ral, cada hoja es $C^{r+1}$. Pero  no se puede  mejorar casi nada la regularidad $C^0$ de las foliaciones invariantes  simplemente aumentando la regularidad de $f$.

\begin{theorem} {\bf H\"{o}lder continuidad de foliaciones invariantes para difeomorfismos de Anosov.} \index{teorema! de existencia de! foliaciones invariantes} \index{foliaci\'{o}n! H\"{o}lder continua} \index{foliaci\'{o}n! regular}

Si $f: M \mapsto M$ es un difeomorfismo de Anosov de clase $C^k$ con $k > 2$, entonces las variedades estables e inestables de $f$ son subvariedades $C^k$, y las foliaciones estables e inestables que forman son $\alpha-$H\"{o}lder continuas para cierto $0 < \alpha < 1$.

\end{theorem}
Una prueba de este teorema se puede encontrar en \cite[Theorem 2.3.1, pag. 48]{BarreiraPesin}.



\subsection{Exponentes de Lyapunov}

A lo largo de las pr\'{o}ximas secciones asumiremos las siguientes hip\'{o}tesis:

$\bullet$ $M$ es una variedad diferenciable, compacta, conexa y provista de una estructura riemanniana.

$\bullet$ $f: M \mapsto M$ es de clase $C^1$ (es decir $f$ es diferenciable y su derivada $df$ es continua). El mapa $f$ no es necesariamente invertible.

\begin{notation} \em 
\label{notationDiferenciable}
.

$\bullet $ $f: M \mapsto M$ es un difeomorfismo   si de clase $C^1$, invertible y adem\'{a}s su inversa tambi\'{e}n es de clase $C^1$. Denotamos $f \in \mbox{Diff }^1(M)$.

$\bullet$ Si adem\'{a}s $f$  y su inversa  son de clase $C^r$ para alg\'{u}n $r \geq 2$ (para lo cual la variedad $M$ tambi\'{e}n debe ser de clase $C^r$ por lo menos), indicaremos $f \in \mbox{Diff }^r (M)$.

$\bullet$ Denotaremos $f \in \mbox{Diff }^{ 1 + \alpha }$ cuando $f \in \mbox{Diff }^1(M)$ y adem\'{a}s $df: TM \mapsto TM$ es $\alpha$-H\"{o}lder continua para cierta constante $0 < \alpha < 1$; i.e.  existen constantes $\delta >0$ y $K >0$ tales que
$$ \| df_x - df_y\| \leq K (\mbox{dist}(x,y) )^{\alpha}\ \ \forall \ x,y \in M \mbox{ such that } \mbox{dist}(x,y) < \delta.$$
En el primer miembro de la desigualdad anterior,  $\| A \|$ denota la norma de la transformaci\'{o}n lineal $A \in L(R^k)$ donde $k = \mbox{dim}(M)$, es decir $\|A\| = \max\{\|Av\|: v \in \mathbb{R}^k, \|v\| = 1\}$.

$\bullet$ Denotaremos $f \in \mbox{Diff }^{1 + Lip}(M)$ si $f \in \mbox{Diff } ^1(M)$ y adem\'{a}s $df$ es Lipschitz, i.e.  existen constantes $\delta >0$ y $K >0$ que satisfacen la desigualdad de $\alpha$-H\"{o}lder continuidad con $\alpha = 1$.

\end{notation}

%

\begin{definition} {\bf Exponentes de Lyapunov}  \index{punto! d\'{e}bilmente regular} \index{punto! regular}  \index{exponentes de Lyapunov! de puntos regulares} \index{exponentes de Lyapunov} \index{regularidad! de puntos}
 \index{regularidad! d\'{e}bil} \label{definicionExponentesLyapunov}

\em Sea $f \in \mbox{Diff }^1(M)$.
Un punto $x \in M$ se llama \em d\'{e}bilmente regular \em si   existen   los siguiente  l\'{\i}mites  para todo $v \in T_xM \setminus \{{\bf 0}\}$:
$$\lim_{n \rightarrow + \infty} \frac{\log (\|df_x^n(v)\|)}{n} ; \ \ \ \ \  \lim_{n \rightarrow + \infty} \frac{\log (\|df_{x}^{-n}  (v)\|)}{-n}.$$
Estos dos l\'{\i}mites (n\'{u}meros reales), se llaman \em exponentes de Lyapunov en el futuro y en el pasado respectivamente, de la \'{o}rbita por $x$ correspondientes a la direcci\'{o}n $[v]$. \em

No se definen los exponentes de Lyapunov en los puntos no regulares.

Se puede definir tambi\'{e}n puntos d\'{e}bilmente regulares y exponentes de Lyapunov (en el futuro) para $f \in C^1(M)$ aunque $f$ no sea  un difeomorfismo.
\vspace{.3cm}

M\'{a}s adelante veremos el Teorema de Oseledets, que demuestra, entre otros resultados, lo siguiente:
\index{teorema! Oseledets}

\begin{center}

\em Para toda   medida de probabilidad $\mu$ que sea $f$-invariante, $\mu$-casi todo punto $x \in M$ es regular. \em

\end{center}

Dicho de otra forma, el conjunto de los puntos no regulares tiene medida nula para toda medida de probabilidad $\mu$ invariante por $f$.

{\bf Nota sobre el concepto de regularidad:} \index{punto! regular}
\index{punto! Lyapunov regular} \index{regularidad! de puntos} \index{regularidad! Lyapunov} En la literatura, suele llamarse punto regular a aquellos puntos $x$ que cumplen una condici\'{o}n m\'{a}s fuerte que la que hemos adoptado nosotros en la Definici\'{o}n \ref{definicionExponentesLyapunov}. En efecto, se definen condiciones adicionales de  existencia de subespacios invariantes tales que para todo vector en ellos el exponente de Lyapunov en el futuro coincide con el exponente de Lyapunov en el pasado. A los puntos regulares que cumplen esa condici\'{o}n m\'{a}s fuerte, los llamaremos Lyapunov-regulares (ver Definici\'{o}n \ref{definitionLyapunovRegulares}). El teorema de Oseledets establece que $\mu$-casi todo punto $x \in M$   no solo es d\'{e}bilmente regular seg\'{u}n nuestra definici\'{o}n \ref{ejercicioAnosovExponentesLyapunov}, sino tambi\'{e}n Lyapunov-regular seg\'{u}n la Definici\'{o}n \ref{definitionLyapunovRegulares}.
\end{definition}

Los exponentes de Lyapunov dependen de la \'{o}rbita  de $x$, pero no dependen de cu\'{a}l punto se tome en la misma \'{o}rbita. En efecto:

\begin{exercise}\em
Probar que si $x$ es d\'{e}bilmente regular, entonces para todo $k \in \mathbb{Z}$ fijo, el   punto $y = f^k(x)$  (es decir, cualquier punto en la \'{o}rbita   de $x$)   es tambi\'{e}n d\'{e}bilmente regular y que el conjunto de los exponentes de Lyapunov  en $y$ coincide con el de   los de $x$.

\end{exercise}

El siguiente ejercicio tiene como objetivo  adelantarse  al enunciado del teorema de Oseledets (que enunciaremos al final de esta secci\'{o}n).

\begin{exercise}\em \label{ejercicioExponentesLyapunov0}  \index{exponentes de Lyapunov! de puntos regulares} \index{exponentes de Lyapunov} \index{punto! d\'{e}bilmente regular}
\index{punto! regular}
\index{regularidad! de puntos}
 \index{regularidad! d\'{e}bil}

Sea  $x$ un punto  d\'{e}bilmente regular.

(a) Sean   ${\bf 0} \neq v \in T_xM$, $\chi$ el exponente de Lyapunov en el futuro (o en el pasado) de la \'{o}rbita de $x$ correspondiente a una direcci\'{o}n $u \neq 0$. Sea $[u] \subset T_xM$ el subespacio de dimensi\'{o}n 1 generado por $u$. Demostrar que para todo $0 \neq u' = [u]$ el exponente de Lyapunov en el futuro (o en el pasado respectivamente) correspondiente a $u$ es el mismo que el de $u'$.

(b) Sean $\chi^+_u$ y $\chi^+_v$ los exponentes de Lyapunov en el futuro de dos direcciones l.i. $0 \neq u, v \in T_xM$. Sean $\chi^-_u$ y $\chi^-_v$ los exponentes de Lyapunov en el pasado de esas dos mismas direcciones $u$ y $v$. Asuma $\chi^+_u \neq \chi^+_v, \ \chi^-_u \neq \chi^-v$. Sea $w =   u +   v$. Probar que $\chi^+_w = \max \{\chi^+_u, \chi^+_v\}$ y $\chi^-_w = \min \{\chi^-_u, \chi^-_v\}$.

Sugerencia:  Asuma $\chi^+_u >  \chi^+_v$.   Use el primer l\'{\i}mite de  \ref{definicionExponentesLyapunov} con los vectores $u$ y $v$, para probar que para todo $\epsilon >0$ existe   $N \in \mathbb{N}$  tal que:
$$\|df^n(w)\| \geq \|df^n(u) \| -\| df^n(v)\| \geq $$
$$e^{n(\chi_u ^+ - \epsilon)}( \|u\| - e^{n(\chi^+_v - \chi^+_u + 2 \epsilon)}\|v\| ) \ \ \forall \ n \geq N. $$
Fije $\epsilon < (\chi^+_u - \chi^+_v)/2$. Tome logaritmo, divida entre $n$ y luego $n \rightarrow + \infty$.

(c) Como en la parte (b) asuma ahora $\chi_u^+ = \chi_v^+$ y $\chi_u ^- = \chi_v ^-$. Pruebe que $\chi_w ^+ \leq \chi_u^+= \chi_v^+, \ \chi_w ^- \geq \chi_u^-= \chi_v^- $.

(d) Sea una base $u_1, \ldots, u_k$ de $T_xM$. Asuma que los exponentes de Lyapunov en el futuro $\chi^+_{i}$ y en el pasado $\chi^-_i$ en las direcciones $u_i$, cumplen $\chi^+_i \neq \chi^+ _j $ para todo $i \neq j$. Probar que para todo $0 \neq u = \sum_{i= 1}^k b_i u_i  \in T_xM$,  se cumple:
$$\chi^+_u = \max_{1 \leq i \leq k; \ b_i \neq 0} \chi^+_i.$$
$$\chi^-_u = \min_{1 \leq i \leq k; \ b_i \neq 0} \chi^-_i.$$
Sugerencia: usar inducci\'{o}n en la cantidad de coeficientes $b_i \neq 0$, y la parte (b).
\end{exercise}

\begin{exercise}\em \label{ejercicioExponentesLyapunov}
Sea $x$ un punto d\'{e}bilmente regular.

    Para cada $k \geq 0$ sean $E^1_{f^k(x)}$ y $E^2_{f^k(x)}$ dos subespacios L.I. (de dimensiones no nulas) de $T_{f^k(x)}M$, invariantes por $f$ (es decir $ E^i_{f^{k+1}(x)} = df_{f^k(x)} E^i_{f^k(x)} $ para todo $k \geq 0$ y para $i = 1, 2$). Asuma la siguiente hip\'{o}tesis:

 {\bf Hip\'{o}tesis I} Para cada  $i$     coinciden  los exponentes de Lyapunov  en el futuro y en el pasado entre s\'{\i} y en todas las direcciones del subespacio $E^i_x$. Ll\'{a}melo  $\chi^i$.  Adem\'{a}s para todo $0 \neq v \in E^1_x \oplus E^2_x$, el exponente de Lyapunov en el futuro de $v$ es mayor o igual que el exponente de Lyapunov en el pasado de $v$.

 (a) Probar que

 (i) Para todo $k \geq 0$  y para cada $i= 1, 2$, los exponentes de Lyapunov en el futuro y en el pasado en $f^k(x)$, correspondientes a la  direcci\'{o}n $E^i_{f^k(x)}$, coinciden entre s\'{\i} y con $\chi^i$.

(ii)   Probar que para todo $0 \neq v \in E^1_{f^k(x)} \oplus E^2_{f^k(x)} $ tal que $ v \not \in   E^1_{f^k(x)},   E^2_{f^k(x)}$, el exponente de Lyapunov $\chi^+_v(x)$ de $[v]$ en el futuro  es igual a $\max\{\chi^1, \chi^2\}$, y el exponente de Lyapunov $\chi_v ^-(x)$ de $[v]$ en el pasado, es igual a $\min\{\chi^1, \chi^2\}$.

 (iii) Deducir que si $\chi_1 = \chi_2 = \chi $, entonces el subespacio $E_{f^k(x)}= E^1_{f^k(x)} \oplus E^2_{f^k(x)}$ es invariante, de dimensi\'{o}n mayor que los subespacios $E^1_{f^k(x)} $ y $ E^2_{f^k(x)}$  que lo generaron, y los exponentes de Lyapunov en el futuro y en el pasado para todos los vectores ${\bf 0} \neq v \in E$ coinciden con $\chi$.

(b) Extender (enunciar y demostrar) los resultados de la  parte  a), asu\-mien\-do que existe un splitting $T_x M= E^1_x \oplus E^2_x \oplus \ldots \oplus E^{h}_x$ que verifica la Hip\'{o}tesis I.

(c) Deducir, asumiendo (b), que el conjunto de exponentes de Lyapunov  diferentes en cualquier punto regular que cumpla la Hip\'{o}tesis I es finito, y   menor o igual que    dim($M$).
\end{exercise}

 El ejercicio anterior motiva preguntarse cu\'{a}ndo se cumple  la  hip\'{o}tesis asumida  sobre la existencia de los subespacios $E^i_x$ invariantes para los cuales los exponentes   de Lyapunov en el futuro y en el pasado existen y   son iguales entre s\'{\i}, y   diferentes para diferentes valores de $i$.    Esta pregunta motiva la siguiente definici\'{o}n de regularidad:

 \begin{definition} \em \label{definitionLyapunovRegulares} {\bf Puntos Lyapunov-regulares} \index{punto! regular}
\index{regularidad! de puntos}
\index{regularidad! Lyapunov}
\index{punto! Lyapunov regular}

Un punto $x \in M$ se llama \em Lyapunov regular \em si
existe un splitting    del espacio tangente $$T_xM = E^1_x \oplus E^2_x \oplus E^{k(x)}_x $$  \index{splitting! de Oseledets} (que puede reducirse como caso particular a $k(x)= 1, \ \ T_x M= E^1_x$), tal  que existen y coinciden entre s\'{\i} \em{ los exponentes de Lyapunov   en el futuro y en el pasado $\chi_i(x)$} \em \index{exponentes de Lyapunov! de puntos regulares}
\index{exponentes de Lyapunov} en toda direcci\'{o}n $ [v] \subset E^i_x$  y adem\'{a}s $$\chi_1(x) < \chi_2(x) <\ldots <\chi_h(x) (x).$$

  En otras palabras:
$$\lim_{n \rightarrow + \infty} \frac{\log (\|df_x^n(v)\|)}{n} =  \lim_{n \rightarrow - \infty} \frac{\log (\|df_{x}^{n}  (v)\|)}{n} = \chi_i(x) \ \ \forall \  \{{\bf 0}\}\neq [v] \subset E^i_x.$$

Los subespacios $E^i_x$ en un punto Lyapunov regular, se llaman \em subespacios de Osele\-dets. \em \index{subespacio! de Oseledets}
\end{definition}

 {\bf Nota: } El Ejercicio \ref{ejercicioExponentesLyapunov} muestra que todo punto Lyapunov-regular es regular en el sentido de la Definici\'{o}n \ref{definicionExponentesLyapunov}. El rec\'{\i}proco es falso (ver el ejemplo del Ejercicio \ref{ejerciciopolonortepolosur}).
 \vspace{.3cm}

 En el enunciado del siguiente ejercicio  se establece que la regularidad Lyapunov, el splitting y los exponentes de Lyapunov dependen de la \'{o}rbita  y no de qu\'{e} punto $x$ en cada \'{o}rbita se elija.

 \begin{exercise}\em Sea $x$ un punto Lyapunov-regular.

 (a) Demostrar que  el splitting $T_xM = \bigoplus_{i= 1}^{h(x)} E^i(x)$ y el conjunto de exponentes de Lyapunov, son \'{u}nicos.

 (b) Demostrar que para todo $k \in \mathbb{Z}$ el punto $f^k(x)$ tambi\'{e}n es Lyapunov regular, y adem\'{a}s $E^i_{f^k(x)} = df^k E^i(x) \forall \ 1 \leq i \leq h(f^k(x)) = h(x), $
 y $\chi_i(x) = \chi_i(f^k(x))$.

 \end{exercise}

 {\bf Pregunta: } ?`Existen abundantes puntos Lyapunov regulares?

 El siguiente resultado  es un teorema fundamental en la Teor\'{\i}a Erg\'{o}dica Diferenciable y responde afirmativamente a la pregunta anterior, por lo menos desde el punto de vista medible, y con respecto a   las medidas invariantes:

\begin{theorem}
{\bf  de Oseledets } \label{theoremOseledecs} \index{teorema! Oseledets}

Sea $M$ una variedad compacta riemanniana de dimensi\'{o}n finita. Sea \em $f \in \mbox{Diff }^1(M)$.  \em  Entonces

{\bf   (a) } El conjunto $R$ de puntos Lyapunov regulares para $f$ es medible.

{\bf   (b) } Las funciones que a cada punto $x \in R$   asignan los exponentes de Lyapunov son medibles.

{\bf   (c) } El conjunto $R$ tiene probabilidad total \em (para toda medida de probabilidad $f$-invariante $\mu$, se cumple $\mu(R)= 1$).

\end{theorem}

 La demostraci\'{o}n de V.I. Oseledets del Teorema \ref{theoremOseledecs} se encuentra en \cite{Oseledecs} para difeomorfismos que preservan la medida de Lebesgue erg\'{o}dica,  y en \cite{Oseledecs2}, en general.
Otras demostraciones   del teorema de Oseledets pueden  encontrarse por ejemplo en \cite{BarreiraPesin}, en \cite[Cap. IV \S 10]{Mane} (ver tambi\'{e}n \cite{ManeIngles} y \cite{VianaOseledets}) y en \cite{Raghunathan}.

Generalizaciones del Teorema de Oseledets, llamados Teoremas Erg\'{o}dicos Multiplicativos, enuncian resultados en los cuales $df_x$ es sustituido por una aplicaci\'{o}n lineal que depende de $x$ en un espacio  vectorial finito dimensional, o incluso por un cociclo en ciertos espacios de Banach infinito dimensionales. Por ejemplo, en \cite{MargulisMultiplicativeErgodicTheorem}, se prueba un Teorema Erg\'{o}dico Multiplicativo que generaliza   al Teorema de Oseledets a ciertos cociclos en espacios de Banach uniformemente convexos.

\subsection{Hiperbolicidad no uniforme} \label{sectionHiperbNoUniforme}

{\bf Interpretaci\'{o}n de los exponentes de Lyapunov no nulos: } \index{exponentes de Lyapunov} \index{exponentes de Lyapunov! no nulos}

Veremos por qu\'{e} los exponentes de Lyapunov, cuando no son nulos, se interpretan como \em la tasa exponencial asint\'{o}tica de crecimiento (dila\-taci\'{o}n) o decrecimiento (contracci\'{o}n) \em en el futuro   de la norma de los vectores en el subespacio tangente, por iteraci\'{o}n del diferencial, es decir al aplicar  $df^{n}$
En efecto, supongamos que un exponente de Lyapunov en el futuro y en el pasado, para la direcci\'{o}n $[s] \subset T_xM$ de la \'{o}rbita por $x$, es $\chi(x) <0$. Entonces, de la definici\'{o}n del l\'{\i}mite en \ref{definicionExponentesLyapunov}, para todo $\epsilon >0$ existe $N \geq 1$:
 $$\|df^n_x(s)\| \leq e^{n(\chi  + \epsilon)} \|s\| \ \forall \ n \geq N$$
Luego,   para ciertos n\'{u}meros reales   $C(x) >0$,  $\lambda(x) = e^{\chi(x) + \epsilon} <1$ (si se fija $0 < \epsilon < -\chi(x)$), se cumple
\begin{equation}\label{equationChi+}\|df^n(s)\| \leq C(x) [\lambda(x)] ^n \|s\| \ \ \forall \ n \geq 0, \ \ \mbox{ donde }0 < \lambda(x) <1,\end{equation} 
(Demostramos con detalle la existencia de $C(x) >0$ y la desigualdad  anterior   en el Lema \ref{lemaHiperbolicidadNoUniforme}.)
La   desigualdad  anterior  significa   que, en la direcci\'{o}n $[s]$, el diferencial $n$-\'{e}simo contrae exponencialmente las  normas de los vectores hacia el futuro  con coeficiente $0<\lambda = e^{\chi + \epsilon} < 1$. Este coeficiente tiende a $e^{\chi}$ cuando $\epsilon \rightarrow 0 ^+$ (y por lo tanto cuando $n \rightarrow + \infty$). Entonces un exponente de Lyapunov negativo $\chi$ es asint\'{o}ticamente igual a $\log \lambda$.      Decimos as\'{\i} que \em un exponente de Lyapunov negativo $\chi(x)$ es la tasa exponencial asint\'{o}tica de contracci\'{o}n en el futuro (o de dilataci\'{o}n hacia el pasado) por la derivada de $f^n$ en la direcci\'{o}n $[s]$.\em

 An\'{a}logamente, si para la misma \'{o}rbita de $x$  existe alguna direcci\'{o}n $[u] \subset T_xM  $, para la cual el exponente de Lyapunov en el pasado y en el futuro es $\chi(x) > 0$, entonces existen un n\'{u}mero real $C(x) >0$  y un valor real $\sigma(x) = e^{\chi(x) - \epsilon} > 1$ (si se fija $0 < \epsilon <  \chi(x)$) tales que:
 \begin{equation}\label{equationChi-}\|df^{-n}(u)\| \leq C(x) [\sigma(x)] ^{-n} \|u\| \ \ \forall \ n \geq 0, \ \ \mbox{ donde }  \sigma(x) >1,\end{equation} 
 (Demostraremos la existencia de la constante $C(x)>0$ y la  desigualdad  anterior  en el Lema \ref{lemaHiperbolicidadNoUniforme}.) Siendo  $\sigma $ asint\'{o}ticamente igual a $ e^{\chi}$, decimos que \em un exponente de Lyapunov positivo $\chi(x)$ es la tasa exponencial asint\'{o}tica  de constracci\'{o}n hacia el pasado (o de dilataci\'{o}n hacia el futuro)  por la derivada de $f^{-n}$ en la direcci\'{o}n $[u]$.\em

 \vspace{.3cm}

Observemos las similitudes y diferencias entre las desigualdades (\ref{equationChi+}) y (\ref{equationChi-}) y las de  la  Definici\'{o}n \ref{definicionHiperbolicidadUniforme} de hiperbolicidad uniforme  (desigualdades (\ref{equationAnosovStable}) y (\ref{equationAnosovUnstable})) Las similitudes  justifican la siguiente definici\'{o}n:

\newpage

\begin{definition} \label{definitionHiperbolicidadNoUniforme} \index{conjunto! hiperb\'{o}lico! no uniforme} \index{transformaci\'{o}n! hiperb\'{o}lica! no uniforme} \index{hiperbolicidad! no uniforme} \index{no uniformemente hiperb\'{o}lico}
{\bf Hiperbolicidad no uniforme} \em
Sea $\Lambda \subset M$ un conjunto medible $f$-invariante: $f^{-1}(\Lambda)= \Lambda$   ($\Lambda$ no es necesariamente compacto).

Decimos que \em $f$ es no uniformemente hiperb\'{o}lico en $\Lambda$, \em (o que $\Lambda$ es un conjunto no uniformemente hiperb\'{o}lico para $f$), si para todo punto $x \in \Lambda$ existe  un splitting $T_xM = S_x \oplus U_x$ que depende mediblemente de $x \in \Lambda$ y que es $df$-invariante, i.e. \index{splitting! hiperb\'{o}lico} \index{fibrado! inestable} \index{subespacio! inestable}

\noindent $\bullet$
   $dfS_x = S_{f(x)}, \ \ df U_x = U_{f(x)} \ \ \forall \ x \in \Lambda,$

   y existen n\'{u}meros reales $C(x)>0$ y $0 \leq \lambda (x) < 1 < \sigma(x)$, que dependen mediblemente de $x \in \Lambda$, tales que

\noindent $\bullet$ $\lambda(f(x)) = \lambda(x), \ \ \sigma(f(x)) = \sigma (x)$,

\noindent $\bullet$    se verifican las dos desigualdades (\ref{equationChi+}) y (\ref{equationChi-}); es decir, para todo $n \geq 0$ y para todos $u \in U_x$ y $s \in S_x$:
 $$\|df^n s \|  \leq C(x) \lambda(x)^n \|s\| , \ \ \ \  \|df^{-n} u\| \leq C(x) \sigma(x)^{-n} \|u\| $$ 


\begin{remark} \em
\label{remarkSplittingMedibleHiperbolico} {\bf Hiperbolicidad No Uniforme.}
A diferencia de la Definici\'{o}n \ref{definicionHiperbolicidadUniforme}  de hiperbolicidad uniforme en compactos,  la hiperbolicidad no uniforme no implica la continuidad del splitting $S_x \oplus U_x$ al variar $x \in \Lambda$. Pero s\'{\i} exige, por definici\'{o}n, que el splitting sea medible.

\end{remark}

{\bf Nota: }
Si $C(x)$, $\lambda(x)$ y $\sigma(x)$ son constantes independientes de $x \in \Lambda$, y adem\'{a}s $\Lambda$ es compacto, se dice que $f$ es uniformemente hiperb\'{o}lico en $\Lambda$, o que $\Lambda$ es un conjunto uniformemente hiperb\'{o}lico para $f$, de acuerdo con la definici\'{o}n \ref{definicionHiperbolicidadUniforme}.
\end{definition}

\begin{lemma}
\label{lemaHiperbolicidadNoUniforme} Sea $\Lambda$ un conjunto $f$-invariante, medible tal que  todo $x \in \Lambda$ es un punto Lyapunov regular cuyo splitting    $E^1_x \oplus E^2_x \oplus \ldots \oplus E^r_x = T_xM $   depende  mediblemente de $x $ \em ($r= r(x)$ tambi\'{e}n), \em y cuyos exponentes de Lyapunov respectivos
 $\chi_1(x) < \chi_2(x) < \ldots < \chi_i(x) < \ldots \chi_r(x)$ son {\bf \em todos diferentes de cero}  y   dependen medi\-blemente de $x$. Entonces $\Lambda$ es un conjunto no uniformemente hiperb\'{o}lico. \em

En extenso, existe   un splitting medible $T_xM= S_x \oplus U_x$ invariante por $df$ y funciones medibles $C(x) >0$ y $0 < \lambda(x) < 1 < \sigma(x)$, tales que $\lambda(f(x))= \lambda(x)$ y $\sigma(f(x)) = \sigma(x)$ y
tales que se verifican las   desigualdades    (\ref{equationChi+}) y (\ref{equationChi-})     para todo $n \geq 0$, para todo $u \in U_x$ y todo $s \in S_x$, y para todo $x \in \Lambda$.
\end{lemma}
{\em Demostraci\'{o}n: }
Sean:
\hfill $\alpha := \max\{\chi_i(x): 1 \leq i \leq r, \ \chi_i(x) <0\} < 0,$

\hfill $\beta := \min\{\chi_i(x): 1 \leq i \leq r, \ \chi_i(x) >0\} >0.$

Como $\alpha$ y $\beta$ son m\'{a}ximo y m\'{\i}nimo de funciones medibles, son medibles. Fijemos un valor constante $\epsilon >0$ suficientemente peque\~{n}o tal que
$\alpha +   \epsilon <0, $ $ \beta -   \epsilon > 0.$  
Denotamos $$0 < \lambda = \lambda(x) := e^{\displaystyle \alpha +  \epsilon} < 1, \ \ \ \sigma = \sigma(x) := e^{\displaystyle \beta -  \epsilon} > 1.$$
Como $\lambda(x)$ y $\sigma(x)$ son composici\'{o}n de funciones continuas con funciones medibles, son   medibles.
Por construcci\'{o}n, como los exponentes de Lyapunov son los mismos para $x$ que para $f(x)$, tenemos
$\lambda(x)= \lambda (f(x)), $ $ \sigma(x) = \sigma(f(x)).$
Sean
$$S_x := \oplus \{E_x^i(x): 1 \leq i \leq r, \ \chi_i(x) <0 \},$$
$$U_x := \oplus \{E_x^i(x): 1 \leq i \leq r, \ \chi_i(x) >0 \},$$
donde $\oplus_{i= 1}^r E_x^i(x) = T_xM$ es el splitting en los subespacios de Oselecs del espacio tangente en el punto $x$. Estos subespacios existen por la Definici\'{o}n \ref{definitionLyapunovRegulares} de punto Lyapunov regular; y  por hip\'{o}tesis dependen mediblemente de $x$. Las funciones $\chi_i(x) $ son medibles. Luego las preimagenes por $\chi_i: T\Lambda \mapsto \mathbb{R}$ de ${\mathbb{R}^+}$ y de $\mathbb{R}^-$ son conjuntos medibles.  Finalmente, la suma directa de un conjunto finito de subfibrados medibles, es un subfibrado medible. Concluimos que $S_x $ y $U_x$ son subfibrados medibles de $T \Lambda$.

Por construcci\'{o}n, como los subespacios de Oseledets son invariantes por $df$, tenemos
$$df S_x = S_{f(x)}, \ \ \ df U_x = U_{f(x)}.$$

Aplicando el resultado del Ejercicio \ref{ejercicioExponentesLyapunov}, tenemos:

Para todo $s \in S_x$ el exponente de Lyapunov en el futuro (y tambi\'{e}n en el pasado) es menor o igual que $\alpha < 0$. Para todo $u \in U_x$ el exponente de Lyapunov en el pasado (y tambi\'{e}n en el futuro) es mayor o igual que $\beta >0$. En extenso:
\begin{equation}\label{eqn35a}\chi_i^+(x,s) = \lim_{n \rightarrow + \infty} \frac {\log \displaystyle  { \|df^n(s)\|}}{n} \leq \alpha < \alpha + \epsilon =  \log \lambda   < 0 \ \ \forall \ 0 \leq s \in S_x,\end{equation}
\begin{equation}\label{eqn35b}\chi_i^-(x,u) = \lim_{m \rightarrow - \infty} \frac{\log \displaystyle {\|df^m(u)\|} }{m} \geq \beta > \beta -  \!\epsilon =  \log \sigma  > \!0 \; \forall \; 0 \leq u \in U_x.\end{equation}

Probemos que para cada $0 \leq s \in S_x$ y para cada $0 \leq u \in U_x$, existen los siguientes n\'{u}meros reales $H(x,s) \ , K(x,s) \ >0$:
\begin{equation}\label{eqn35ff}H(x, s) :=     \sup_{n \geq 0} \frac{\|df_x^n(s)\|}{ \ \lambda^n   \, \|s\| \ }, \ \ \ \ \ \   K(x, u) :=  \  \sup_{n \geq 0} \frac{\|df_x^{-n}(u)\|}{\ \sigma^{-n} \, \|u\| \ }.\end{equation}

En efecto, fijemos $x, u, s$.  En las igualdades  (\ref{eqn35a}) y (\ref{eqn35b}), aplicamos la definici\'{o}n de l\'{\i}mite,  multiplicamos por $n$ (con $|n|$ suficientemente grande), y aplicamos la exponencial. (Hay que cuidar que cuando $n$ es negativo, al multiplicar por $n$ se invierten el sentido de las desigualdades). Concluimos que  existe $N= N(x,u,s) \geq 0$ tal que
$${\|df^n(s)\|}\ / \ {  \lambda^{n}    \, \|s\| \ } < 1, \ {\|df^{-n}(u)\|}\ / \ {\ \sigma^{-n }\, \|u\| \ } < 1,$$
para todo $n \geq N \geq 0$. Entonces, el supremo que define $H(x,s)$, as\'{\i} como el supremo que define $K(x,s)$, existe y es un n\'{u}mero  real no negativo  (porque  para $|n| \geq N$, todos los cocientes cuyo supremo buscamos son menores que 1; y para $0 \leq |n| \leq N$ los cocientes a maximizar son positivos y una cantidad finita).
Afirmamos que existen los n\'{u}meros reales $K(x), H(x) >0$, definidos por:
\begin{equation}\label{eqn35fh} H(x) :=    \sup \big\{ H(x,s) : \ \ s \in S_x, \ \|s\| = 1 \big \},\end{equation}
\begin{equation}\label{eqn35fj}K(x) :=    \sup \big\{ K(x,u) : \ \ u \in U_x, \ \|u\| = 1 \big \} \end{equation}

Probaremos que existe  $K(x)$ real (la prueba de que existe $H(x)$ real es similar). Tomemos una base   $u_1, \ldots, u_{k}$ de $U_x$, donde $k = \mbox{dim}(U_x)$, formada por  vectores $u_i$ de norma 1, que se encuentran todos en los subespacios de Oseledets, seg\'{u}n la definici\'{o}n \ref{definitionLyapunovRegulares} de punto Lyapunov regular. Entonces, dado   $u \in U_x$ se puede escribir:
$ u  = \sum_{i= 1}^{k_2} b_i u_i.$
Si $  \|u\| = 1$, entonces existe $M(x)$ tal que $  0 \leq |b_i| \leq M(x) \ \forall \ 1 \leq i \leq k.$  En efecto, fijada la base, cada  $|b_i| $ es una funci\'{o}n real continua del vector $u \in U_x$ (pues es la norma de la proyecci\'{o}n ortogonal del vector $u$ sobre el subespacio generado por $u_i$). Por el teorema de Weierstrass, la funci\'{o}n continua $u_i$ tiene un m\'{a}ximo $M_i$ en el subconjunto compacto $\{u \in U_x: \|u\| = 1\} \subset T_xM$. Entonces basta tomar $M(x) := \max_{i= 1}^{k} M_i$.

  Tenemos:
  \begin{equation}\label{eqn35y}\frac{\|df^{-n}(u)\|}{\sigma^{-n}} \leq     \sum_{i= 1}^{k} |b_i|    \frac{\|df^n(u_i)\|}{{\sigma^{-n}}} \leq   M(x) \sum_{i= 1}^{k} K(x, u_i) =: K_1(x)   \ \ \forall \ n \geq 0. \end{equation}

 Para la primera desigualdad usamos la propiedad triangular de la norma. Para la \'{u}ltima desigualdad, usamos que $|b_i| \leq M(x)$ y la definici\'{o}n de $K(x, u_i)$.
Esto prueba que existe  $K_1(x) < + \infty$ definido  por la  igualdad   (\ref{eqn35y}). Como la   desigualdad a la izquierda en (\ref{eqn35y}) vale para todo $u \in U_x$ con $\|u\|= 1$, entonces el supremo $K(x)$   definido en  (\ref{eqn35fj}) cumple $K(x) \leq K_1(x) < + \infty$. An\'{a}logamente se prueba que existe $H(x) < + \infty$ definido por (\ref{eqn35fh}).

 De las definiciones de los n\'{u}meros $H(x,s), K(x,u), H(x), K(x)$ en las igualdades (\ref{eqn35ff}),   (\ref{eqn35fh}) y (\ref{eqn35fj}), definiendo $C(x) = \max\{K(x), H(x)\}$, deducimos:
 $$ \|df^n(s) \| \leq     C(x)   \lambda(x)^n \ \|s\|   \ \forall n \geq 0, \ \forall \ s \in S_x,$$
 $$\|df^{-n}(u) \| \leq     C(x)   \sigma(x)^{-n} \ \|u\| \  \forall n \geq 0, \ \forall \ u \in U_x.$$
 En efecto, estas desigualdades valen  cuando $\|s\| = 1$  y $\|u\|= 1$ por la definici\'{o}n de $H(x)$ y $K(x)$. Entonces valen  para todo $s \in S_x$ y para todo $u \in U_x$ por la linealidad de $df^n$.
Esto termina de probar el Lema \ref{lemaHiperbolicidadNoUniforme}.
 \hfill $\Box$

\subsection{Regi\'{o}n de Pesin y medidas hiperb\'{o}licas}
 En esta secci\'{o}n $f\in \mbox{Diff }^1(M)$ y $M$ es una variedad compacta y riemanniana.

\begin{definition}
\label{definicionRegionDePesin} \em \index{exponentes de Lyapunov}  \index{exponentes de Lyapunov! no nulos}   \index{Pesin! regi\'{o}n de} \index{regi\'{o}n de Pesin}
  La \em   regi\'{o}n de Pesin $P_f \subset M$ \em es el conjunto de los puntos $x \in M$ Lyapunov-regulares tales que los exponentes de Lyapunov $\chi^1_x < \chi^2_x \ldots \chi^{h(x)}_x$ son todos diferentes de cero.
\end{definition}

Observar que por el Teorema \ref{theoremOseledecs} de Oseledets, la regi\'{o}n de Pesin es medible.
Para algunos difeomorfismos la regi\'{o}n de Pesin   puede ser vac\'{\i}a. Por ejemplo, trivialmente, si $f$ es la identidad $P_f = \emptyset$. Otro ejemplo: las rotaciones de la esfera $S^2$ ($f$ es la rotaci\'{o}n de la esfera, de \'{a}ngulo constante alrededor de un di\'{a}metro de $S^2$, llamado eje polo norte-polo sur): $P_f = \emptyset$.

\begin{exercise}\em
(a) Construir un difeomorfismo $f: S^1 \mapsto S^1$ en el c\'{\i}rculo $S^1$ que preserve la orientaci\'{o}n, tal que   $P_f = \emptyset$ y tal que en todo abierto $V \subset S^1$ la derivada $f'$ no sea id\'{e}nticamente igual a 1.

(b)Idem en la esfera $S^2$  con la condici\'{o}n de que en todo abierto $V \subset S^2$ la derivada $df$ no es id\'{e}nticamente igual a la identidad $Id$, ni a $-Id$.

(c) Construir un ejemplo en el c\'{\i}rculo $S^1$ que cumplan las condiciones de la parte (a) y adem\'{a}s tal que    el conjunto de puntos Lyapunov regulares sea finito.

(d) ?`Existen ejemplos en el c\'{\i}rculo que cumplan las condiciones de la parte (a) y adem\'{a}s tal que el conjunto de los puntos Lyapunov-regulares sea infinito?

(e) Construir un difeomorfismo $f: S^1 \mapsto S^1$ tal que la regi\'{o}n de Pesin $P_f \neq \emptyset$ pero que no coincida con el conjunto de todos los puntos Lyapunov-regulares.

(f) Probar que para todo difeomorfismo $f$ del c\'{\i}rculo, o bien la regi\'{o}n de Pesin $P_f$ es vac\'{\i}a o bien es finita o bien es infinita numerable, y encontrar ejemplos de los tres casos (sugerencia: probar que todo $x \in P_f$ es aislado en $P_f$).

(g) Demostrar que existen difeomorfismos $f: S^2 \mapsto S^2$ en la esfera tal que la regi\'{o}n de Pesin $P_f$ es no numerable (sugerencia: la Herradura de Smale definida en \ref{definicionHerradura}).

\end{exercise}



 \begin{definition}
 \em {\bf (Medida hiperb\'{o}lica)} \label{definitionMedidaHiperbolica} \index{medida! hiperb\'{o}lica}
 \index{hiperbolicidad! en regi\'{o}n de Pesin}

  Una medida de probabilidad $f$-invariante $\mu$ se dice \em hiperb\'{o}lica \em si $\mu(P_f) = 1$, donde $P_f$ es la regi\'{o}n de Pesin. En otras palabras, $\mu$  es hiperb\'{o}lica si y solo si los exponentes de Lyapunov son diferentes de cero $\mu$-c.t.p.
 \end{definition}

Para demostrar el siguiente resultado, utilizaremos   el Teorema \ref{theoremOseledecs} de Oseledets.

 \begin{theorem}
 {\bf Medidas hiperb\'{o}licas y conjuntos no uniformemente  hi\-per\-b\'{o}\-li\-cos} \label{theoremMuHiperbolicaNoUnifHiperbolico}
\index{conjunto! hiperb\'{o}lico! no uniforme} \index{hiperbolicidad! no uniforme} \index{medida! hiperb\'{o}lica}  \index{teorema! de hiperbolicidad no uniforme}

 {\bf (a)} Si $\mu$ es medida de probabilidad hiperb\'{o}lica entonces existe un conjunto invariante  $\Lambda$ (no necesariamente compacto)    tal que $\mu(\Lambda)= 1$  y $f$  es hiperb\'{o}lica \em (unif. o no unif.)   \em  en $\Lambda$.

 {\bf (b)} Si $\Lambda$ es un conjunto invariante tal que $f$ es
 \em (unif. o no unif.) \em    hiperb\'{o}lica en $\Lambda$, y si $\mu$ es una probabilidad $f$-invariante tal que $\mu(\Lambda)= 1$, entonces $\mu$ es medida hiperb\'{o}lica.

 {\bf (c)} Si $\Lambda$ es un conjunto invariante y compacto tal que $f$ es  \em (unif. o no unif.) \em    hiperb\'{o}lica en $\Lambda$, entonces existen medidas de probabilidad   $\mu$ tales que $\mu(\Lambda)= 1$. Luego, por la parte (b), todas estas medidas son hiperb\'{o}licas. En particular existen medidas de probabilidad hiperb\'{o}licas y erg\'{o}dicas soportadas en $\Lambda$.

 \end{theorem}
 {\em Demostraci\'{o}n: }
 {\bf (a)} Sea $\Lambda$ la regi\'{o}n de Pesin. Por definici\'{o}n de puntos regulares, $\Lambda $ es $f$-invariante, y por el Teorema \ref{theoremOseledecs} de Oseledets, $\Lambda$ es medible. Por definici\'{o}n de medida hiperb\'{o}lica $\mu(\Lambda)= 1$. Por el Teorema \ref{theoremOseledecs} de Oseledets,   los exponentes de Lyapunov y el splitting en subespacios correspondientes, son funciones medibles.  Debido  al Lema \ref{lemaHiperbolicidadNoUniforme}, como los exponentes de Lyapunov de todo punto $x \in \Lambda$ son no nulos,  $\Lambda$ es (unif. o no unif.) hiperb\'{o}lico.

 {\bf (b)}   De la  desigualdad  (\ref{equationChi+}), tomando logaritmo, dividiendo entre $n$ y haciendo $n \rightarrow + \infty$, se deduce que para $\mu$-c.t.p.   $x \in \Lambda$, y para toda direcci\'{o}n $[s] \in S_x$, los exponentes de Lyapunov  \em hacia el futuro  \em (seg\'{u}n Definici\'{o}n \ref{definicionExponentesLyapunov}) son menores que $\log \lambda(x) < 0$. An\'{a}logamente, de la desigualdad (\ref{equationChi-}), tomando logaritmo, dividiendo entre $-n < 0$ (se invierte el sentido de la desigualdad), y haciendo $n \rightarrow + \infty$, deducimos que para toda direcci\'{o}n $[u] \in U_x$ los exponentes de Lyapunov  \em hacia el pasado  \em   son mayores que $\log \sigma(x) > 0$.

 Por el teorema de Oseledets, $\mu$-casi todo punto es Lyapunov regular. Entonces para $\mu$-casi todo punto existen los subespacios de Oseledets para los cuales los exponentes de Lyapunov hacia el pasado son iguales a los exponentes de Lyapunov hacia el futuro. A priori, hay tres casos:  el subespacio de Oseledets $E_x^i$ est\'{a} contenido en $S_x$, o est\'{a} contenido en $U_x$, o ninguna de las dos cosas. En el primer caso, el exponente de Lyapunov hacia el futuro y hacia el pasado en $E_x^i$ es negativo  (porque es negativo hacia el futuro por estar contenido en $S_x$ y coincide con el exponente hacia el pasado por ser un subespacio de Oseledets). En el segundo caso, el Lyapunov hacia el futuro y hacia el pasado en $E_x^i$ es positivo (porque es positivo hacia el pasado por estar contenido en $U_x$ y porque coincide con el exponente hacia el futuro por ser un subespacio de Oseledets). Probemos que el tercer caso es vac\'{\i}o; es decir todo subespacio de Oseledets est\'{a} contenido en $S_x$ o en $U_x$. Por absurdo, sea una direcci\'{o}n $[v] \in E_x^i, \  \ [v] \not \in S_x, U_x$. Como $T_xM = U_x \oplus S_x$, entonces   $v=   u +   s$ con $0 \neq u \in U_x$, \ $0 \neq s \in S_x$. Como $\chi^+_s <0$, aplicando lo probado en el Ejercicio \ref{ejercicioExponentesLyapunov}, el exponente de Lyapunov $\chi_i$ hacia el futuro en $E_x^i$ es menor o igual que $\chi^+_s <0$. Entonces es negativo. An\'{a}logamente, como $\chi^-_u >0$ (porque $u \in U_x$ y en $U_x$ los exponentes de Lyapunov hacia el pasado son positivos), aplicamos lo probado en el Ejercicio \ref{ejercicioExponentesLyapunov} y deducimos que el exponente de Lyapunov $\chi_i$ hacia el pasado en $E_x^i$ es mayor o igual que $\chi^-_u  >0$. Entonces es positivo. Pero en un subespacio de Oseledecs, el exponente de Lyapunov hacia el futuro y hacia el pasado coincide. No puede ser negativo y positivo a la vez. Entonces no hay subespacio de Oseledets que no est\'{e} incluido en $S_x$ o en $U_x$.

 Concluimos que en todos los subespacios de Oseledets, los exponentes de Lyapunov son no nulos.    Por lo probado en el ejercicio \ref{ejercicioExponentesLyapunov}, todos los exponentes de Lyapunov en cualquier direcci\'{o}n, hacia el futuro o hacia el pasado, son iguales a alg\'{u}n exponente de Lyapunov en los subespacios de Oseledecs. Entonces son no nulos. Esto vale para $\mu$-c.t.p. $x \in M$. Entonces $\mu$-c.t.p no tiene exponentes de Lyapunov iguales a cero; es decir, $\mu$ es una medida de probabilidad hiperb\'{o}lica.

{\bf (c)} Sea $\widetilde  f = f|_{\Lambda} : \Lambda \mapsto \Lambda$. Siendo $\Lambda$ un espacio m\'{e}trico compacto y $\widetilde  f$ continua en $\Lambda$, el Teorema \ref{teoremaExistenciaMedErgodicas} asegura que existen medidas de probabilidad $\mu$ soportadas en $\Lambda$ invariantes y erg\'{o}dicas para $\widetilde  f$. Es inmediato chequear que $\mu$, como medida de probabilidad en $M$, es invariante y erg\'{o}dica para $f$. Por la parte b) toda tal medida es hiperb\'{o}lica.
\hfill $\Box$

\begin{remark} \em
 \label{remarkMedidaHiperbolicaErgodica} {\bf Medidas hiperb\'{o}licas erg\'{o}dicas:} \index{medida! hiperb\'{o}lica erg\'{o}dica}
Si una medida   es erg\'{o}dica, entonces los exponentes de Lyapunov son constantes $\mu$-c.t.p. (pues son funciones medibles invariantes c.t.p.). Por el mismo motivo, las dimensiones de los subespacios del splitting en la definici\'{o}n de Lyapunov regularidad son constantes $\mu$-c.t.p. Por lo tanto, para las medidas hiperb\'{o}licas erg\'{o}dicas, la definici\'{o}n del conjunto no uniformemente hiperb\'{o}lico $\Lambda$, tal que $\mu(\Lambda)= 1$, adquiere las particularidades siguientes:

   {\bf (i) } En las desigualdades (\ref{equationChi+}) y (\ref{equationChi-}), las tasas de contracci\'{o}n y dilataci\'{o}n $0 < \lambda < 1 < \sigma$  son constantes independientes de $x \in \Lambda$ (mientras que en general   el coeficiente $C(x)>0$ var\'{\i}a con $x$).

   {\bf (ii) } Las dimensiones de los subespacios estable $S_x$ e inestable $U_x$, son constantes, independientes de $x \in \Lambda$.

 \end{remark}

\subsection{\large Variedades estable e inestable  en la regi\'{o}n de  Pe\-sin}
En esta secci\'{o}n,  $M$ es una variedad compacta y riemanniana y $f \in \mbox{Diff}^{1 + \alpha}(M)$, es decir $f$ es difeomorfismo $C^1$ m\'{a}s H\"{o}lder.
Recordamos que esto significa que   $f \in \mbox{Diff }^1(M)$ y tanto $df_x$ como $df^{-1}_x$ son funciones H\"{o}lder-continuas del punto $x \in M$, i.e. existen constantes $\alpha, K >0$ tales que
$$\|df_x - df_y\| \leq K \, [\mbox{dist}(x,y)]^{\alpha},$$
y an\'{a}logamente para $df^{-1}_x$.

El siguiente teorema, generaliza al caso no uniformemente hiperb\'{o}lico, el Teorema \ref{teoremavariedadesinvariantesAnosov} de existencia de variedades invariantes para conjuntos uniformemente hiperb\'{o}licos (en particular para difeomorfismos de Anosov). Sin embargo, la validez de la siguiente generalizaci\'{o}n, as\'{\i} como la de los resultados que fundamentan la Teor\'{\i}a de Pesin, est\'{a} restringida
a difeomorfismos de clase $C^{1 + \alpha}$.

\newpage

\begin{theorem} {\bf Variedades Estable e Inestable locales (Pesin)}
\label{theoremVarInvariantesRegionPesin} \index{variedad invariante! local}  \index{variedad invariante! inestable}

Sea $f: M \mapsto M$ un difeomorfismo $C^1$ m\'{a}s H\"{o}lder en una variedad riemanniana compacta $M$ tal que la regi\'{o}n de Pesin $P(f)$ es no vac\'{\i}a.  Para  $x \in P(f)$ denotamos $S_x$ y $U_x$ los subespacios estable e inestable, respectivamente, correspondientes a los exponentes de Lyapunov negativos y positivos  de $x \in \Lambda$.  \em (Notar que \  $S_x \oplus U_x = T_xM$). \em

Entonces, para  todo  $x \in P(f)$ existen subvariedades locales $W  _{\mbox{\footnotesize{loc}}}^s(x)$ y $W_{\mbox{\footnotesize{loc}}}  ^u(x)$,   $C^1$-encajadas    en $M$,     a las que   llamamos variedad   estable e inestable  local respectivamente,   tales que:

\em
{\bf (a) } $$T_xW_{\mbox{\footnotesize{loc}}}^s(x) = S_x, \ \ \ T_xW _{\mbox{\footnotesize{loc}}}^u(x) = U_x.$$

{\bf (b) } $$f(W_{\mbox{\footnotesize{loc}}} ^s(x)) \subset W  _{\mbox{\footnotesize{loc}}}^s(f(x)), \ \ \ f (W_{\mbox{\footnotesize{loc}}}^u  (x)) \supset W_{\mbox{\footnotesize{loc}}}^u  (f (x)).$$

{\bf (c) } Para todo $x \in P(f)$ y para todo $y \in M$  se cumple:    
 $$\lim_{n \rightarrow + \infty} \mbox{dist}(f^n(x), f^n(y)) = 0   \  \mbox{  si }  \ y \in  W_{\mbox{\footnotesize{loc}}}  ^s(x),$$
$$\lim_{n \rightarrow + \infty} \mbox{dist}(f^{-n}(x), f^{-n}(y)) = 0   \ \mbox{ si   } \ y \in W_{\mbox{\footnotesize{loc}}}   ^u(x).$$

\end{theorem}

El teorema anterior y su demostraci\'{o}n se encuentran en Pesin \cite{Pesin76}. La demostraci\'{o}n puede encontrarse tambi\'{e}n en \cite[Theorem 4.1.1, pag.81]{BarreiraPesin} o en \cite{HirschPughShub}.

Se observa, en el Teorema \ref{theoremVarInvariantesRegionPesin}, que la hip\'{o}tesis $f$ de clase $C^1$ m\'{a}s H\"{o}lder  es necesaria. En efecto, Pugh en \cite{PughC1masHolder} construy\'{o} un ejemplo $f \in \mbox{Diff }^1(M)$ cuya derivada $df$ es continua pero no es H\"{o}lder continua, con regi\'{o}n de Pesin no vac\'{\i}a, para el que no vale el Teorema \ref{theoremVarInvariantesRegionPesin} de existencia de variedades invariantes.





\section{Atractores topol\'{o}gicos y erg\'{o}dicos} \label{chapterAtractorestopoyergo}

\subsection{Atractores topol\'{o}gicos}

A lo largo de esta secci\'{o}n $f: X \mapsto X$ es una transformaci\'{o}n continua en un espacio m\'{e}trico compacto $X$.

\begin{definition}
{\bf Estabilidad Lyapunov   y orbital.} \index{estabilidad! orbital}
\index{estabilidad! Lyapunov} \index{conjunto! orbitalmente estable} \index{conjunto! Lyapunov estable}

\em
Sea $K \subset X$ un conjunto compacto no vac\'{\i}o invariante por $f$ (es decir $f^{-1}(K) = K$.)

  $K$  se dice \em orbitalmente estable \em (hacia el futuro) si para todo $\epsilon >0$ existe $\delta >0$ tal que   $$\mbox{dist}(p, K) < \delta \ \Rightarrow \ \mbox{dist} (f^n(p), K) < \epsilon \ \ \forall \ n \geq 0.$$

  $K$  se dice \em Lyapunov estable \em (hacia el futuro) si para todo $\epsilon >0$ existe $\delta >0$ tal que   $$\mbox{dist}(p, K) < \delta \ \Rightarrow \ \exists \ q \in K \mbox{ tal que } \mbox{dist} (f^n(p), f^n(q)) < \epsilon \ \ \forall \ n \geq 0.$$

  De las definiciones anteriores se deduce que si $K$ es Lyapunov estable, entonces es orbitalmente estable. Sin embargo el rec\'{\i}proco no es cierto en general.

\end{definition}

\begin{definition} {\bf Atractor topol\'{o}gico I} \em  \label{definitionAtractorTopologico} \index{atractor! topol\'{o}gico}

Un \em atractor topol\'{o}gico \em es un conjunto   $K$ compacto y no vac\'{\i}o tal que

1) $ K = f^{-1}(K) = f(K)$.

2) Existe un   abierto $V \supset K$, llamado \em cuenca local de atracci\'{o}n topol\'{o}gica \em de $K$, tal que  \begin{equation} \label{eqn28}\lim_{n \rightarrow + \infty} \mbox{dist}(f^n(x), K) = 0 \ \ \forall \ x \in V. \end{equation}
\index{cuenca de atracci\'{o}n! topol\'{o}gica}
Consideramos, de   particular importancia, los atractores topol\'{o}gicos \em que sean orbitalmente estables \em (para los cuales daremos una definici\'{o}n equiva\-lente en \ref{definitionAtractorTopologicoII}).

Un atractor topol\'{o}gico $K$ se llama \em minimal \em (como atractor topol\'{o}gico)  si   el \'{u}nico compacto no vac\'{\i}o
   $K' \subset K$   que cumple 1) y 2)  es $K' = K$.
   \index{conjunto! minimal}

\end{definition}

\begin{remark} \em   Muchos autores requieren adem\'{a}s de las condiciones 1) y 2), que el compacto no vac\'{\i}o $K$ sea orbitalmente estable o Lyapunov estable, para llamarse atractor topol\'{o}gico.
La condici\'{o}n   de estabilidad orbital  en la Definici\'{o}n \ref{definitionAtractorTopologico} no es redundante con las otras dos condiciones 1) y 2). En el ejemplo \ref{exampleAsintoticoNoEstable} parte (B),  se muestra un conjunto compacto $  \{p_0\}$ formado por un solo punto fijo $p_0$ que cumple las condiciones 1) y 2) de la definici\'{o}n \ref{definitionAtractorTopologico}, pero que no es  orbitalmente estable.
\end{remark}

{\bf Ejemplo:} En el cap\'{\i}tulo 1 vimos que \index{pozo}
toda \'{o}rbita peri\'{o}dica  hiperb\'{o}lica tipo pozo (de un difeomorfismo $f: M \mapsto M$ en una variedad   $M$), es un atractor topol\'{o}gico.

\begin{exercise}\em
Probar que una \'{o}rbita peri\'{o}dica hiperb\'{o}lica tipo pozo es un atractor topol\'{o}gico Lyapunov estable.
\end{exercise}
\begin{exercise}\em \label{exerciseTransitivoAtrTop} (a) Sea $f: M \mapsto M$ un homeomorfismo transitivo en una variedad compacta $M$. Probar que $M$ es un atractor topol\'{o}gico y que es minimal (como atractor topol\'{o}gico). Deducir que $M$ es el \'{u}nico atractor topol\'{o}gico. Sugerencia: De la transitividad deducir que existe alguna \'{o}rbita futura densa en $M$. Para probar que $M$ es minimal como atractor topol\'{o}gico y el \'{u}nico atractor asumir que   $K \subset M$ es un atractor topol\'{o}gico. Su cuenca local de atracci\'{o}n topol\'{o}gica  $V$ contiene  alg\'{u}n punto $x$ de una    \'{o}rbita futura   densa.  Probar que la \'{o}rbita futura de $x$ es densa. Usando  la definici\'{o}n de atractor, probar que el omega-l\'{\i}mite de $x$ para cualquier $x \in V$ debe estar en $K$. Deducir    que $K= M$.

(b) Probar que el \'{u}nico atractor topol\'{o}gico del autormorfismo lineal $f = \left(
        \begin{array}{cc}
          2 & 1   \\
          1 & 1   \\
        \end{array}
      \right)
 $ en el toro $\mathbb{T}^2$ es todo el toro. \index{automorfismo! lineal del toro}
\end{exercise}

\begin{exercise}\em \label{ejercicioBobo}
Sea $K \subset X$ un atractor topol\'{o}gico. Sea en $K$ la topolog\'{\i}a inducida por su inclusi\'{o}n en el espacio m\'{e}trico $X$. Suponga que existe $x \in K$ tal que la \'{o}rbita futura de $x$ es densa en $K$. Demostrar que $K$ es minimal como atractor topol\'{o}gico. \index{atractor! topol\'{o}gico}
\end{exercise}

\begin{exercise}\em \label{ejercicioAtractorCircunferencia}
Sea $X \subset \mathbb{R}^2 \sim \mathbb{C}$ un disco compacto de centro en el origen y radio $1 < r < 2$. Consideremos en $X$ coordenadas polares $p = \rho e ^{i \varphi} \in X: \ \ \ (\rho, \varphi): \ \ 0 \leq \rho \leq r, \ \ 0 \leq \varphi < 2 \pi$. Sea $f: X \mapsto X$ el homeomorfismo dado por las siguientes ecuaciones
$$f(p) = f(\rho e^{i \varphi})= \rho* e ^{i \varphi*}, \mbox{ donde } $$
$$\rho^* = \frac{  \rho (4 -  \rho)}{3} , \ \ \ \ \ \ \ \varphi^* = \varphi + a$$
donde $a$ es una constante $0 \leq a < 2 \pi$.

(a) Dibujar algunas de las \'{o}rbitas y probar que la circunferencia $K$ de centro en el origen y radio $1$ es un atractor topol\'{o}gico Lyapunov estable. (Sugerencia: ver la demostraci\'{o}n de que $K$ es atractor topol\'{o}gico en el ejemplo \ref{exampleAsintoticoNoEstable}).

(b) Probar que si $a/(2\pi)$ es irracional entonces $K$   es minimal como atractor topol\'{o}gico. (Sugerencia: usar el ejercicio \ref{ejercicioBobo})

(c) Si $a/(2 \pi)$ es racional, ?`es $K$ minimal como atractor topol\'{o}gico?

\end{exercise}
\begin{example} \em  \label{exampleAsintoticoNoEstable}.
\index{conjunto! no orbitalmente estable}

  Sea $X \subset \mathbb{R}^2 \sim \mathbb{C}$ el disco compacto de centro en el origen y radio $1 < r < 2$, como en el ejercicio \ref{ejercicioAtractorCircunferencia}.

{\bf (A) }  Sea, en coordenadas polares, el siguiente homeormorfismo $f$:
$$f(p) = f(\rho e^{i \varphi})= \rho^* e ^{i \varphi^*}, \mbox{ donde } $$
$$\rho^* = \frac{  \rho (4 -  \rho)}{3} , \ \ \ \ \ \ \ \varphi^*  = \varphi + (\rho-1) $$
Tenemos $\rho^* -1 = (3- \rho) (\rho - 1) / 3$. Luego $|\rho^*-1| \leq \lambda |\rho-1| $ si $\rho > 3 \cdot (1 - \lambda)$ donde $0 < \lambda < 1$. Entonces, por inducci\'{o}n en $n$ se deduce que $\lim_{n \rightarrow + \infty} |\rho^{(n)} - 1| = 0 $ donde $\rho^{(n)} $ es la distancia al origen  de $f^n(p)$ para cualquier $p \neq (0,0)$. Deducimos que la distancia de $f^n(p)$ a la circunferencia $K$ de centro $0$ y radio $1$ tiende a cero con $n \rightarrow + \infty$, para cualquier punto inicial $0 \neq p  \in X$. Adem\'{a}s $\rho^*= 1$ si $\rho= 1$, de donde $\rho^{(n)} = 1$ si $\rho ^{(0)} = 1$. Entonces   la circunferencia $K$ de centro en el origen y radio $1$ es invariante por $T$ y por lo probado antes  $\mbox{dist}(f^n(p), K) \rightarrow 0$. Luego $K$ es un atractor topol\'{o}gico de acuerdo con la definici\'{o}n \ref{definitionAtractorTopologico}. Adem\'{a}s $K $ es orbitalmente estable porque  $|\rho^* - 1| \leq |\rho - 1|$. Entonces $|\rho^{(n)} - 1| $ es decreciente con $n$ si $p \neq 0$. Luego, la distancia de $f^n(p)$ a la circunferencia $K$ decrece con $n$. Por lo tanto, es menor que $\epsilon >0$ si inicialmente es menor que $\delta= \epsilon$. Esto prueba la estabilidad orbital de $K$.

\begin{figure}[h]


\begin{center}
\includegraphics[scale=.4]{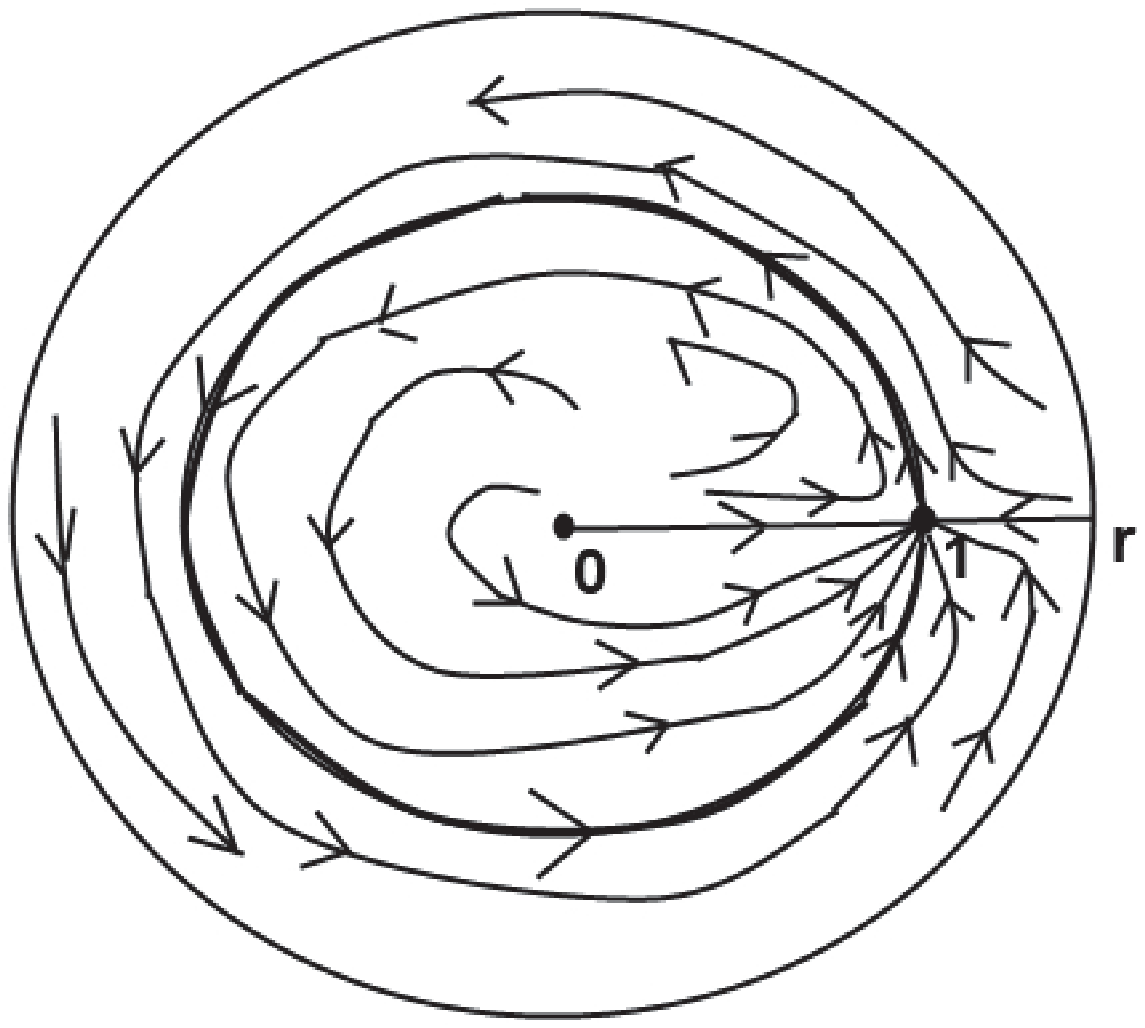}

\caption{Algunas \'{o}rbitas del Ejemplo $g$ en \ref{exampleAsintoticoNoEstable} (B). \label{figuraAtractorNoEstable}  
\index{conjunto! no orbitalmente estable} \index{atractor! topol\'{o}gico! no orbitalmente estable}}

El punto 1 es fijo, es un atractor topol\'{o}gico 

no orbitalmente estable.
\end{center}
\end{figure}

{\bf (B) } Sea ahora $g: X \mapsto X$ dado por las ecuaciones
$$g(p) = f(\rho e^{i \varphi})= \rho^* e ^{i \varphi^*}, \mbox{ donde } $$
$$\rho^* =   \frac{  \rho (4 -  \rho)}{3} , \ \ \ \ \ \ \ \varphi^*  = \varphi \, \Big(2 - \frac{\varphi}{2 \pi} \Big)  \mbox{ para } \varphi \in [0, 2 \pi). $$
En la figura \ref{figuraAtractorNoEstable} se representan algunas de las \'{o}rbitas. Como en la parte (A), se prueba que  la circunferencia $K$ de centro en el origen y radio 1, es un atractor topol\'{o}gico orbitalmente estable. Adem\'{a}s, si graficamos $\varphi^*$ en funci\'{o}n de $\varphi$ para todo $\varphi \in [0, 2 \pi]$, observamos que es creciente, con puntos fijos en $0$ y $2 \pi$, siendo $0$ repulsor (a la derecha de 0) y $2 \pi$ atractor (a la izquierda de $2 \pi$). Entonces, si $\varphi^{(n)}$ es el \'{a}ngulo en coordenadas polares del punto $f^n(p)$, para $p= \rho \, e ^{ i \varphi}$, se tiene que $\lim_{n \rightarrow + \infty} \varphi^{(n)} = 2 \pi \sim 0$.
Luego, en este ejemplo, $1 = 1 e ^{2 \pi \, i}$ es un punto fijo, y  $\lim_{n \rightarrow + \infty} f^n(p) = 1$ para todo $p \neq 0$. Entonces $\{1\} \subset K$ tambi\'{e}n es atractor topol\'{o}gico seg\'{u}n la definici\'{o}n \ref{definitionAtractorTopologico}. Pero $\{1\}$ no es orbitalmente estable: Por un lado $1$ es un punto fijo, por lo tanto, si $\{1\}$ fuera \'{o}rbitalmente estable, toda \'{o}rbita futura con punto inicial $p$ suficientemente cercano a $1$, deber\'{\i}a mantenerse arbitrariamente cerca de $1$.  Consideremos la \'{o}rbita futura $o:= \{f^n(p)\}_{n \geq 0}$ por $0 \neq p = \rho e^{i \varphi} \not \in K$  con $0 <\varphi < \delta$ (para $\delta >0$ suficientemente peque\~{n}o). Aunque $o$ se acerca a $K$ cuando $n $ crece, y adem\'{a}s $\lim_{n} f^n(p) = 1$,  el \'{a}ngulo $\varphi^{(n)}$ (esto es la coordenada angular de $f^n(p)$ en polares) no se mantiene a distancia menor que $\epsilon >0$ de $0$.    Dicho de otra forma,   la \'{o}rbita $\{f^n(p)\}_n$, para algunos $n \geq 1$, se aleja del punto fijo $1$ m\'{a}s que una constante positiva (digamos 1/6, si $\delta >0$ es suficientemente peque\~{n}o), aunque finalmente tiende a $1$ cuando $n \rightarrow + \infty$ (ver Figura \ref{figuraAtractorNoEstable}).
\end{example}

 El siguiente resultado da una caracterizaci\'{o}n de los atractores topol\'{o}gicos orbitalmente estables. Debido a este resultado algunos  autores definen atractor topol\'{o}gico agregando, a las condiciones 1) y 2) de la Definici\'{o}n \ref{definitionAtractorTopologico}, la condici\'{o}n de que sea orbitalmente estable.

\newpage

\begin{proposition}.

 \label{propositionAtractorTopologico} {\bf Caracterizaci\'{o}n de atractores topol\'{o}gicos orbitalm. estables}

Sea $f: X \mapsto X$ un homeomorfismo en un espacio m\'{e}trico compacto $X$.

{\bf (a) } Si $V$ es un abierto no vac\'{\i}o tal que $\overline{f(V) } \subset V$  y si \begin{equation}\label{eqnMaximalInvariante} K:= \bigcap_{n= 0}^{+ \infty} f^n(V),\end{equation} entonces $K$ es compacto no vac\'{\i}o,   es un atractor topol\'{o}gico orbitalmente estable y $V$ es cuenca local de atracci\'{o}n de $K$. \em

{\bf (b) } Rec\'{\i}procamente, si $K$ es un atractor topol\'{o}gico orbitalmente estable, entonces existe un abierto $V$ que es cuenca local de atracci\'{o}n de $K$, tal que $\overline{f(V)} \subset V$  y $K = \bigcap_{n= 0}^{+ \infty} f^n(V).$
\end{proposition}

{\bf Definici\'{o}n de conjunto maximal invariante: } \index{conjunto! maximal invariante} \index{maximal invariante} Si $K$ es compacto invariante no vac\'{\i}o que satisface la igualdad (\ref{eqnMaximalInvariante}) para un abierto $V \supset K$, entonces $K $ se llama \em conjunto maximal invariante (hacia el futuro) \em  de $V$.

\vspace{.2cm}

La Proposici\'{o}n \ref{propositionAtractorTopologico}  justifica la siguiente definici\'{o}n:

\begin{definition} {\bf Atractor topol\'{o}gico II}  \label{definitionAtractorTopologicoII} \index{atractor! topol\'{o}gico}
\index{estabilidad! orbital} \index{conjunto! orbitalmente estable}

\em Sea $f: X \mapsto X$ es un homeomorfismo en un espacio m\'{e}trico compacto $X$.  Un subconjunto no vac\'{\i}o, compacto e invariante
$K \subset X$ se llama \em atractor topol\'{o}gico  orbitalmente estable \em o tambi\'{e}n \em atractor topol\'{o}gico maximal invariante \em  si existe un entorno abierto $V \supset K$, llamado cuenca local de atracci\'{o}n topol\'{o}gica de $K$, tal que \em $$\overline {f(V)} \subset V, \ \ \ \mbox{ y } \ \ \   K := \bigcap_{n= 0}^{+ \infty} f^n(V).$$
\end{definition}

{\em Demostraci\'{o}n: }
{\em de la Proposici\'{o}n } \ref{propositionAtractorTopologico}:

{\bf (a) } Como $\overline {f(V)} \subset V$ y $f$ es un homeomorfismo, tenemos, por inducci\'{o}n: $$f^{n+1}(V) \subset \overline {f^{n+1}(V)} \subset f^n(V) \ \forall \ n \in \mathbb{N}.$$   $$K:= \cap_{n = 0}^{+ \infty} f^n(V) = V \cap (\cap_{n= 1}^{+ \infty}f^n(V)) \subset V \cap (\cap_{n= 1}^{+ \infty}\overline{f^n(V)})= $$ $$= \cap_{n= 1}^{+ \infty}\overline{f^n(V)}  \subset  \cap_{n= 0}^  {+\infty}f^n(V)) = K .$$ Luego, todas las inclusiones anteriores son igualdades y tenemos $$K = \bigcap_{n= 1}^{+ \infty}\overline{f^n(V)}, \ \  \ \ \ \bigcap_{n= 1}^{+ N}  {f^n(V)} =  {f^N(V)}, \ \  \ \ \ \bigcap_{n= 1}^{+ N} \overline{f^n(V)} = \overline{f^N(V)} \ \ \forall \ N \geq 1.$$ Por la propiedad de intersecciones finitas no vac\'{\i}as de compactos, deducimos que $K$ es compacto no vac\'{\i}o. Adem\'{a}s $$f^{-1}(K) = \cap_{n= 1}^{+ \infty} f^{-1}(\overline{f^n(V)}) = \cap_{n= 0}^{+ \infty} {\overline{f^n(V)}} = \overline V \cap K = K, $$
pues, por construcci\'{o}n $K \subset V$. Hemos probado que $K$ es compacto no vac\'{\i}o e invariante con $f$.

Para   $\epsilon >0$, denotamos $B_{\epsilon}(K) := \{x \in X: \mbox{dist}(x, K) < \epsilon\}$.

\vspace{.2cm}

 {\bf Afirmaci\'{o}n A:}
 \em Para todo $\epsilon >0$ existe $N \geq 1$ tal que
$  \overline{f^N(V)} \subset B_{\epsilon}(K).$ \em

\vspace{.2cm}

Por absurdo, si existe $\epsilon >0$ y para todo $n \geq 1$ existe $x_n \in \overline{f^n(V)}$ tal que $\mbox{dist}(x_n, K) \geq \epsilon$, entonces tomando una subsucesi\'{o}n de $\{x_n\}$ convergente a un punto $x$, tenemos $\mbox{dist}(x, K) \geq \epsilon$. Adem\'{a}s $x \in   \overline {f^N(V)}\  \forall \ N \geq 1$ (pues $x_N \in \overline {f^N(V)} = \cap_{n= 1}^N \overline{f^n(V)}$, de donde $x_{m} \in \overline {f^N(V)}$ para todo $m \geq N$). Entonces $x \in \cap_{N \geq 1} \overline{f^N(V)} = K$, lo cual contradice que $\mbox{dist}(x, K) \geq \epsilon$. Hemos probado la afirmaci\'{o}n A.

Para probar que $K$ es un atractor topol\'{o}gico, y que $V$ es cuenca local de atracci\'{o}n de $K$, basta probar que si $y = \lim_{j \rightarrow + \infty} f^{n_j}(x)$ para alg\'{u}n $x \in V$ y para alguna subsucesi\'{o}n $n_j \rightarrow + \infty$, entonces $y \in K$. En efecto para todo $\epsilon >0$, tenemos $\mbox{dist}(y, \overline {f^{n_j}(V)}) < \epsilon$ para todo $j$ suficientemente grande. Luego, usando la afirmaci\'{o}n (A), como ${\overline{f^{n_j}(V)}} \subset B_{\epsilon}(K)$ para todo $n_j$ suficientemente grande, deducimos: $\mbox{dist}(y, K) < 2 \epsilon$. La desigualdad anterior vale para todo $\epsilon >0$. Entonces $\mbox{dist}(y,K) = 0$; es decir $y \in K$, como quer\'{\i}amos demostrar.

Ahora resta probar que $K$ es orbitalmente estable. Dado $\epsilon >0$ sea $N \geq 1$ como en la afirmaci\'{o}n (A). Siendo $f$ un homeomorfismo, $f^N(V)$ es abierto. Adem\'{a}s $K \subset f^N(V)$ porque $K := \bigcap_{N \geq 0}f^N(V)$. Como $K$ es compacto, existe $\delta >0$ tal que $B_{\delta}(K) \subset f^N(V)$. Si $x \in B_{\delta}(K)$ entonces, para todo $m \geq 0$ se cumple $f^m(x) \in f^m(B_{\delta}(K)) \subset f^m (f^N(V))   \subset  \overline{f^{m+ N}(V)} \subset   \overline  {f^N(V)} \subset B_{\epsilon}(K). $ Hemos probado que si $\mbox{dist}(x, K) < \delta$ entonces $\mbox{dist}(f^m(x), K) < \epsilon$ para todo $m \geq 0$. Luego $K$ es orbitalmente estable, terminando la demostraci\'{o}n de la parte (a).

\vspace{.3cm}

{\bf (b) }  Sea $K$ un atractor topol\'{o}gico orbitalmente estable y sea $V'$ un abierto que es   cuenca local de atracci\'{o}n de $K$, seg\'{u}n la condici\'{o}n 2) de la Definici\'{o}n \ref{definitionAtractorTopologico}. Sea $\epsilon >0$ tal que $$ \overline {B_{2\epsilon}(K)} \subset V',$$ donde $  B_{\epsilon}(K) := \{x \in X: \mbox{dist}(x, K) < \epsilon\}$. Por la definici\'{o}n de esta\-bi\-li\-dad orbital de $K$ existe $0 <\delta < \epsilon$ tal que $f^m(B_{2\delta}(K)) \subset B_{\epsilon}(K)$ para todo $m \geq 0$.
Como $\overline{B_{\delta}(K)} \subset B_{2 \delta}(K)$, tenemos
$$f^m (\overline {B_{\delta}(K)}) \subset B_{\epsilon}(K) \ \ \forall \ m \geq 0.$$

\vspace{.2cm}

{\bf Afirmaci\'{o}n (B): } Existe $N \geq 1$ tal que $$f^m(x) \in B_{ \epsilon}(K) \ \forall \ m \geq N, \ \forall \ x \in \overline {B_{2\epsilon}(K)}.$$

En efecto, como $\overline {B_{2\epsilon}(K)} \subset V'$ y $V'$ satisface la condici\'{o}n 2) de la Definici\'{o}n \ref{definitionAtractorTopologico}, para cada $x \in \overline{B_{2\epsilon}(K)}$ existe $N_x \geq 1$ tal que $\mbox{dist} (f^{N_x}(x), K) < \delta$. Como $f$ es continua, existe un entorno abierto $U_x$ de $x$ tal que $\mbox{dist} (f^{N_x}(y), K) < \delta$ para todo $y \in U_x$. Siendo  $\overline {B_{2\epsilon}(K)}$ compacto, exis\-te un subcubrimiento finito   $\{U_{x_1}, U_{x_2}, \ldots, U_{x_h}\}$. Deducimos que para todo $y \in \overline {B_{2\epsilon}(K)}$ existe $N_{x_i} \geq 1$ tal que $\mbox{dist} (f^{N_{x_i}}(y), K) < \delta$. Por la construcci\'{o}n de $\delta$ deducimos que $\mbox{dist} (f^{m}(y), K) < \epsilon$ para todo $m \geq N_{x_i}$. Luego, si $N= \max\{N_{x_i}: 1 \leq i \leq h\}$, se cumple $f^m(y) \in B_{\epsilon}(K)$ para todo $m \geq N$ y para todo $y \in \overline{B_{2\epsilon}(K)}$, terminando la prueba de la Afirmaci\'{o}n (B).

\vspace{.2cm}

Definimos el abierto $V$ de la siguiente manera:
$$V := B_{\epsilon_0}(K) \ \cup \  f(B_{ {\epsilon_1 }}(K)) \ \cup \    f^2(B_{ {\epsilon_2 }}(K)) \   \cup $$ $$ \cup \  f^3(B_{ {\epsilon_3 }}(K)) \  \cup \  \ldots \ \cup \  f^N(B_{ {\epsilon_N }}(K)),$$
  $$\mbox{donde } \ \ \epsilon _i := \epsilon \left( 1 + \frac{i}{2^N} \right) < 2 \epsilon \ \ \ \ \ \forall \ \ 0 \leq i \leq N.$$
  Probemos que $\overline{f(V)} \subset V$. En efecto, por un lado  si $0 \leq i \leq N-1$ entonces $$\overline {f(f^i(B_{\epsilon_i}(K)))} \subset  \overline {f^{i+1}(B_{\epsilon_i}(K))} \subset f^{i+1}(B_{\epsilon_{i+1}}(K)) \subset V.$$ Por otro lado, si $i= N$ entonces $$\overline {f(f^N(B_{\epsilon_N}(K)))} \subset \overline {f (f^N(B_{2\epsilon}(K))} = f^{N+1}(\overline{B_{2 \epsilon}(K)}) \subset B_{\epsilon}(K) \subset V ,$$ debido a la construcci\'{o}n del natural $N$ seg\'{u}n la afirmaci\'{o}n (B).

Ahora probemos que $V$ es cuenca local de atracci\'{o}n de $K$. Sea $x \in V$. Entonces $x \in f^i(B_{2 \epsilon}(K))$ para alg\'{u}n $0 \leq i \leq N$,  es decir $x = f^i(y)$ para alg\'{u}n $y \in B_{2 \epsilon}(K) \subset V'$, donde $V'$ es por hip\'{o}tesis una cuenca local de atracci\'{o}n de $K$. Entonces, por la propiedad 2) de la Definici\'{o}n \ref{definitionAtractorTopologico}, tene\-mos $\lim_{n \rightarrow + \infty}\mbox{dist}(f^n(y), K) = 0$. Luego $\lim_{m \rightarrow + \infty} \mbox{dist}(f^m(x), K) = \lim _{m \rightarrow + \infty}\mbox{dist}(f^{m-i}(y), K) = 0$, para todo $x \in V$. Hemos probado que $V$ satisface la condici\'{o}n 2) de la Definici\'{o}n \ref{definitionAtractorTopologico}. Entonces, por definici\'{o}n $V$ es cuenca local de atracci\'{o}n topol\'{o}gica de $K$.

Para terminar la demostraci\'{o}n solo falta probar que $K = \cap_{n= 0}^{+ \infty} f^n(V)$. Por un lado como $f(K) = K$ (porque por hip\'{o}tesis $K$ es invariante con $f$), y como $K \subset V$ por construcci\'{o}n, entonces $K = f^n(K) \subset f^n(V)$ para todo $n \geq 0$. Luego $K \subset \cap_{n= 0}^{+ \infty} f^n(V)$.

Ahora probemos la inclusi\'{o}n opuesta. Sea $y \in \cap_{n= 0}^{+ \infty} f^n(V)$. Entonces  $$\forall \ n \geq 1 \ \ \ \exists \ \ x_n \in f(V) \mbox{ tal que  } y = f^{n-1}(x_n).$$ Sea $x$ el l\'{\i}mite de una subsucesi\'{o}n convergente de $\{x_n\} \subset f(V)$, es decir
$$x = \lim_{j \rightarrow + \infty} x_{n_j}, \ \ \ n_j \rightarrow + \infty.$$  Tenemos $x \in \overline {f(V)} \subset V$.

Dado $\epsilon' >0$, por la estabilidad orbital de $K$, existe $\delta' >0$ tal que $f^n(B_{\delta'}(K)) \subset B_{\epsilon'}(K)$ para todo $n \geq 0$.
Como $V$ es cuenca local de atracci\'{o}n de $K$, y $x \in V$, deducimos que   existe $N$ tal que   $\mbox{dist} (f^m(x), K) < \delta'   \ \ \forall \ m \geq N,$  en particular $$\mbox{dist}(f^N(x), K) < \delta'.$$
Por la continuidad de $f^N$, existe un entorno abierto $U_x$ de $x$, tal que $$\mbox{dist}(f^N(z), K) < \delta' \ \forall \ z \in U_x.$$ Por construcci\'{o}n del punto $x $, $x_{n_j} \in U_x$ para todo $j$ suficientemente grande. Luego:
 $$\mbox{dist}(f^N(x_{n_j}), K) < \delta' \ \ \forall \ j \mbox{ suficientemente grande.}$$
  Por la construcci\'{o}n del n\'{u}mero $\delta'$ a partir de la estabilidad orbital de $K$, deducimos que $$\mbox{dist}(f^m(x_{n_j}), K) < \epsilon' \ \ \forall \ \ m \geq N.$$ Luego, en particular, para $m= n_j -1$, si elegimos $j$ suficientemente grande tal que
    $$n_j \geq N +1,$$
  tenemos $$\mbox{dist}(f^{n_j -1}(x_{n_j}), K) < \epsilon'.$$
  Por construcci\'{o}n de la sucesi\'{o}n $\{x_n\}$ tenemos $y = f^{n_j-1}(x_{n_j})$. Deducimos que $\mbox{dist}(y, K) < \epsilon '$. Como $\epsilon'>0$ es arbitrario, concluimos que $y \in K$ como quer\'{\i}amos demostrar. \hfill $\Box$

 \subsection{Cuenca global de atracci\'{o}n topol\'{o}gica}
\begin{definition}
{\bf Cuenca global de atracci\'{o}n topol\'{o}gica} \index{cuenca de atracci\'{o}n! topol\'{o}gica} \index{cuenca de atracci\'{o}n! global} \index{$C(K)$ o $ E_K$ cuenca de atracci\'{o}n! topol\'{o}gica del compacto $K$}

\em Sea $f: X \mapsto X$ continua en un espacio m\'{e}trico compacto $X$. Sea $K \subset X$ (no vac\'{\i}o, compacto e invariante) un atractor topol\'{o}gico seg\'{u}n la Definici\'{o}n \ref{definitionAtractorTopologico}. Se llama \em cuenca de atracci\'{o}n topol\'{o}gica (global) de $K$ \em al siguiente conjunto:
\begin{equation} \label{eqn29} C(K) = \{ x \in X: \lim_{n \rightarrow + \infty} \mbox{dist}(f^n(x), K) = 0\}. \end{equation}
Nota: Dado un compacto $K$ no vac\'{\i}o e invariante cualquiera, aunque el conjunto $C(K)$ construido mediante la igualdad (\ref{eqn29}) resulte  no vac\'{\i}o, en general no se llama a este conjunto \lq\lq cuenca de atracci\'{o}n topol\'{o}gica\rq\rq  de $K$,  si $K$ no es un atractor topol\'{o}gico. Es decir, cuando se usa el nombre \lq\lq cuenca de atracci\'{o}n topol\'{o}gica\rq\rq, previamente sabemos que $K$ satisface todas las condiciones   de la Definici\'{o}n \ref{definitionAtractorTopologico}, en particular la condici\'{o}n 2) de existencia de un  \em entorno abierto \em de $K$ contenido en la cuenca global de atracci\'{o}n topol\'{o}gica.
\end{definition}

\begin{proposition}. \label{proposicionCuencaTopologica}

{\bf Propiedades de la cuenca global de atracci\'{o}n topol\'{o}gica.} \index{cuenca de atracci\'{o}n! topol\'{o}gica} \index{cuenca de atracci\'{o}n! global} \index{cuenca de atracci\'{o}n! abierta} \index{cuenca de atracci\'{o}n! invariante}

Sea $f: X \mapsto X$ continua en un espacio m\'{e}trico compacto $X$. Si $K \subset X$ es  un atractor topol\'{o}gico, entonces su cuenca (global) de atracci\'{o}n topol\'{o}gica $C(K)$, tiene las si\-guien\-tes propiedades:

{\bf (a) } $C(K)$ es abierto y no vac\'{\i}o.

{\bf (b) } $C(K)$ es invariante con $f$, es decir $f^{-1}(C(K)) = C(K)$.

{\bf (c) } $K \subset C(K) $.

\em Como consecuencia de (a) y (b) si $X $ es conexo, entonces o bien $K = X$ (en cuyo caso $C(K) = K = X$), o bien la cuenca de atracci\'{o}n $C(K)$ contiene propiamente a $K$ (no coincide $K$ con su cuenca).
\end{proposition}

{\em Demostraci\'{o}n: }

{\bf (b) } Si $x \in C(K)$ entonces, por la igualdad (\ref{eqn29}):

  $\lim_{n \rightarrow + \infty} \mbox{dist}(f^n(x), K) = 0$.
Luego,   $$\lim_{n \rightarrow + \infty} \mbox{dist} (f^{n}(f(x)), K) = \lim_{n \rightarrow + \infty} \mbox{dist} (f^{n-1}(f(x)), K) = $$ $$ = \lim_{n \rightarrow + \infty} \mbox{dist} (f^{n}( x), K)=  0.$$ Esto muestra que $f(x) \in C(K)$ para todo $x \in C(K)$. Hemos probado que $f(C(K)) \subset C(K)$, o lo que es lo mismo,   $C(K) \subset f^{-1}(C(K))$.

Sea ahora $y \in f^{-1}( C(K))$. Entonces $f(y) = x \in C(K)$. Luego, para todo $n \geq 0$ tenemos $f^n(x) = f^{n+ 1}(y)$. Entonces: $$\lim_{n \rightarrow + \infty} \mbox{dist} (f^{n }(y), K) = \lim  _{n \rightarrow + \infty} \mbox{dist} (f^{n +1}(y), K) = $$ $$= \lim_{n \rightarrow + \infty} \mbox{dist} (f^{n }(x), K) = 0.  $$ Se deduce que  $y  \in C(K)$ para todo $y \in f^{-1}(C(K))$, probando que $f^{-1}(C(K)) \subset C(K)$.

{\bf (a) y (c) } Por hip\'{o}tesis $K$ es un atractor topol\'{o}gico. Entonces existe   $V \supset K$ abierto, cuenca local de atracci\'{o}n topol\'{o}gica de $K$, de acuerdo a la condici\'{o}n 2) de la Definici\'{o}n \ref{definitionAtractorTopologico}. De la definici\'{o}n de $C(K)$ en (\ref{eqn29}), junto con la igualdad (\ref{eqn28}), tenemos $V \subset C(K)$, probando que $K \subset C(K)$ y, por lo tanto, $C(K) $ es no vac\'{\i}o. Si $x \in C(K)$, entonces $\lim_n \mbox{dist}(f^n(x), K) = 0$. Como $V $ es un entorno abierto de $K$, deducimos que existe $N \geq 1$ (que depende de $x \in C(K)$), tal que
$$f^n(x) \in V \ \ \ \forall \ \ n \geq N.$$
En particular, la afirmaci\'{o}n de arriba vale para $n= N$. Luego, por la continuidad de $f$ y la apertura de $V$, existe un entorno abierto  $U_x$ de $x$ tal que
$$f^N(x) \in V \ \ \forall \ \ x \in U_x.$$
Como $V \subset C(X)$, deducimos que $f^N(U_x) \subset C(X)$, o dicho de otra forma $U_x \subset f^{-N}(C(X)) = C(X)$. Luego $C(X)$ es abierto.
\hfill $\Box$

\begin{proposition} \label{propositionMedidaInvarianteAtractorTopologico} \index{atractor! topol\'{o}gico} \index{medida! soportada en atractor}

Sea $f: X \mapsto X$ continua en un espacio m\'{e}trico compacto $X$. Sea $K \subset X$ un atractor topol\'{o}gico con cuenca global de atracci\'{o}n topol\'{o}gica $C(K)$. Entonces:

 {\bf (a) } Existen medidas de probabilidad invariantes   para $f|_{C(K)}$.

  {\bf (b) } Existen medidas   de probabilidad erg\'{o}dicas para $f|_{C(K)}$.

  {\bf (c) } Toda medida de probabilidad invariante para $f|_{C(K)}$ est\'{a} soportada en el atractor $K$ (es decir $\mu(K)= 1$).
\end{proposition}

{\em Demostraci\'{o}n: }
{\bf (a) y (b)} Consideremos $f|_K: K \mapsto K$. Es la restricci\'{o}n al atractor $K \subset C(K)$ de $f|_{C(K)}: C(K) \mapsto C(K)$. Como $f$ es continua y $K$ es compacto, entonces, por lo demostrado en el cap\'{\i}tulo 1, existen medidas de probabilidad invariantes para $f|_K$, y tambi\'{e}n  existen medidas de probabilidad  erg\'{o}dicas para $f|_K$. Tomemos una de estas medidas $\nu$ invariante para $f|_K$. Definamos la medida de probabilidad $\mu$ en los borelianos $B$ de $C(K)$ de la siguiente forma: $$\mu(B) := \nu (B \cap K) \ \ \forall \ B \subset C(K) \ \mbox{ boreliano }.$$ (Observar que  $\mu(C(K)) = \nu(C(K) \cap K)= \nu (K) = 1$, porque $K \subset C(K)$.)  Veamos que la probabilidad $\mu$ es invariante para $f|_{C(K)}$. Sea $B \subset C(K)$ un boreliano cualquiera:
$$\mu(f^{-1}(B)) = \nu (f^{-1}(B) \cap K) = \nu(f^{-1}(B) \cap f^{-1}(K)) = $$ $$ = \nu (f^{-1}(B \cap K)) = \nu (B \cap K) = \mu (B),$$
pues $\nu$ es invariante con $f|_K$.
Ahora veamos que si  $\nu$ es erg\'{o}dica para $f|_K$ entonces $\mu$ es erg\'{o}dica para $f|_{C(K)}$. Sea $B \subset C(K)$ boreliano tal que $f^{-1}(B) = B$. Tenemos que $f^{-1}(B \cap K) = f^{-1}(B) \cap f^{-1}(K) = B \cap K$. Luego, como $\nu$ es erg\'{o}dica para $f|_K$, se cumple $\nu(B \cap K) \in \{0,1\}.$. Entonces
$\mu(B) = \nu(B \cap K) \in \{0,1\},  $ lo que prueba la ergodicidad de $\mu$.

{\bf (c) } Sea $\mu$ un medida de probabilidad invariante para $f|_{C(K)}$. Entonces $\mu(C(K)) = 1$. Debemos probar que $\mu(K) = 1$ (lo cual implica que $\mu(C(K) \setminus K) = 0$). Por el Lema de Recurrencia de Poincar\'{e}, $\mu$-c.t.p. es recurrente; es decir $x \in \omega(x)$ para $\mu$-c.t.p. $x \in C(K)$. Para probar que $\mu(K)= 1$ basta probar que $x \in K$ para cualquier $x \in C(K)$ tal que $ x \in \omega(x)$. Para esto, alcanza demostrar que $\omega(x) \subset K$ para todo $x \in C(K)$. En efecto, por la definici\'{o}n de omega-l\'{\i}mite, si $y \in \omega(x)$ entonces $$y = \lim_{j \rightarrow + \infty} f^{n_j}(x) , \ \ \ n_j \rightarrow + \infty. $$
Luego: $$\mbox{dist}(y, K) = \lim_{j \rightarrow + \infty} \mbox{dist}(f^{n_j}(x), K) \leq \limsup_{n \rightarrow + \infty} \ \ \mbox{dist}(f^n(x), K).$$
Como $x \in C(K)$, por la definici\'{o}n de cuenca de atracci\'{o}n topol\'{o}gica dada en (\ref{eqn29}), tenemos: $$\limsup_{n \rightarrow + \infty} \ \ \mbox{dist}(f^n(x), K) = \lim_{n \rightarrow + \infty} \ \ \mbox{dist}(f^n(x), K) = 0.$$
Luego $\mbox{dist}(y,K)= 0$, o lo que es lo mismo $y \in K$ para todo $y \in \omega(x)$, como quer\'{\i}amos demostrar. \hfill $\Box$

\subsection{Atractores hiperb\'{o}licos ca\'{o}ticos}

  \index{atractor! hiperb\'{o}lico} \index{atractor! ca\'{o}tico}  En esta secci\'{o}n consideramos el caso particular en que $X= M$ es una variedad diferenciable compacta y riemanniana, y $f \in  \mbox{Diff }^1(M)$.

Se llaman \em atractores hiperb\'{o}licos, \em a aquellos atractores topol\'{o}gicos $K$ que est\'{a}n soportados   en las variedades inestables de un conjunto compacto no vac\'{\i}o (unif. o no unif.) hiperb\'{o}lico $\Lambda \subset K$. M\'{a}s precisamente: $K$ es un atractor hiperb\'{o}lico si es un atractor topol\'{o}gico y existe un conjunto invariante compacto y (unif. o no unif.) hiperb\'{o}lico $\Lambda \subset K$ tal que $K \subset \bigcup_{p \in \Lambda} W^u(p)$.

El caso no trivial es cuando la dimensi\'{o}n de estas variedades inestables es mayor o igual que 1,   para $\mu$- casi todo punto $x \in K$ para alguna medida invariante $\mu$, soportada en el atractor $K$.
Como son varie\-dades inestables de puntos de un conjunto hiperb\'{o}lico, los exponentes de Lyapunov en las direcciones tangentes a estas variedades (tangentes al atractor) son po\-si\-ti\-vos. Llamaremos a tales atractores topol\'{o}gicos, \em atractores hiperb\'{o}licos ca\'{o}ticos. \em

La b\'{u}squeda de variedades inestables que soporten  el atractor, es una de las motivaciones m\'{a}s relevantes de la teor\'{\i}a de sistemas din\'{a}micos diferenciables. Significa que la din\'{a}mica dentro de un atractor  topol\'{o}gico $K$ (es decir la din\'{a}mica de $f|_K$) puede ser expansiva en el futuro, o en otras palabras, ca\'{o}tica; por tener exponentes de Lyapunov positivos. Dicho de otra forma, veamos el caso en que el atractor $K$ sea orbitalmente estable (\'{o} Lyapunov  estable, respectivamente). La estabilidad orbital (de Lyapunov, resp.)      rige en la  cuenca de atracci\'{o}n $C(K)$. En efecto, por definici\'{o}n, la estabilidad orbital (de Lyapunov resp.) tiene significado no trivial  solo para las \'{o}rbitas de $C(K)$ \em fuera \em de $K$. Sin embargo, no es contradictoria con una din\'{a}mica inestable y expansiva, dada por los exponentes de Lyapunov positivos, \em dentro \em del atractor $K$ (ver el ejemplo de siguiente ejercicio).

El siguiente  ejercicio  muestra  un  ejemplo  en el que la din\'{a}mica dentro del atractor topol\'{o}gico es expansiva o ca\'{o}tica: el atractor est\'{a} formado por las variedades inestables de puntos con exponentes de Lyapunov positivos:

\begin{exercise}\em
{\bf Atractor Solenoide (Smale-Williams)} \index{atractor! solenoide} \index{atractor! de Smale-Williams}
\index{solenoide}

Sea $S \subset \mathbb{R}^3$ el toro s\'{o}lido compacto, definido como la imagen en $\mathbb{R}^3$ de la siguiente parametrizaci\'{o}n con tres par\'{a}metros reales (en coordenadas cil\'{\i}ndricas) $$(r, \varphi, \theta) \in [0, a] \times [0, 2 \pi] \times [0, 2 \pi],$$ donde $0 < a < 1/2$ es una constante real:
  \begin{eqnarray}
       x  &=& (r -1/2) \cos\theta \cos \varphi   \nonumber \\
       y  & = &   r \cos \theta \sen \varphi \nonumber \\
        z &= & r \sen \theta    \nonumber
  \end{eqnarray}

  (a) Interpretar el significado geom\'{e}trico de los par\'{a}metros $(r, \varphi, \theta)$   en $\mathbb{R}^3$ y dibujar el toro s\'{o}lido $S$. (Sugerencia, dibujar primero el disco circular de radio $a$ en el plano $xz$ (es decir el plano $\{\varphi= 0\}$). Luego, observar que $S$ es el s\'{o}lido de revoluci\'{o}n que se obtiene haciendo girar ese disco alrededor del eje de las $z$).

  (b) Sea $f: S \mapsto \mbox{int}(S)$ continua,  tal que lleva cada disco circular $D_{\varphi}$ que se obtiene cortando $S$ con el plano $\varphi$ constante, en un disco circular $f(D_{\varphi})$ de radio $a/4$, contenido en la secci\'{o}n $D_{2 \varphi}$, de modo que   $f|_{D_{\varphi}} =   f_3 \circ f_2 \circ f_1$  donde:

$\bullet$ $f_1$ es una rotaci\'{o}n con eje   $z$ de \'{a}ngulo $\varphi$ (esto implica que   $f_1 (D_{\varphi}) = D_{2 \varphi}$; cuando escribimos $2 \varphi$ nos referimos a este   \'{a}ngulo m\'{o}dulo $2 \pi$).

$\bullet$ $f_2$ es una homotecia de raz\'{o}n $1/4$ que lleva el centro del disco $D_{2 \varphi}$ al punto, en el interior de $D_{2 \varphi}$, con par\'{a}metros $(a/2, 2 \varphi, 0)$

$\bullet$ $f_3$ es una rotaci\'{o}n de \'{a}ngulo $\varphi$ alrededor del eje perpendicular al plano del disco $D_{2 \varphi}$ que pasa por el  centro del c\'{\i}rculo $D_{2 \varphi}$.

  Verificar que $f: S \mapsto f(S)$ es inyectiva, y que $f(S)$ es un compacto contenido en el interior de $S$. Dibujar $f(S)$. Sugerencia: $f(S)$ da dos vueltas alrededor del eje de las $z$. Ver por ejemplo \cite[Figure 1]{Pikovsky_Arnoldcat}.

  (c) Sea $K = \bigcap_{n \geq 0} f^n(S)$ el atractor topol\'{o}gico con cuenca local de atracci\'{o}n $S$.  Sea, para cada $N \geq 0$ el compacto
  $K_N = \bigcap_{n= 0}^{N} f^n(S).$ Sea $A_N = K_N \cap \{(r, \varphi, \theta) \in S: \ \varphi= 0\}.$
  Dibujar $A_1, A_2,  K_1, K_2$, y probar  que $K_{N+1} \subset \mbox{int}(K_N)$ para todo $N \geq 0$. (Sugerencia: inducci\'{o}n en $N$).

  (d) Para cada punto $p \in K$ se define la variedad inestable de $p$ del siguiente modo:
$$W^u(p) = $$ $$=\{q \in S:  \exists \ f^{-n}(q) \in S \ \forall \ n \geq 0 \ \mbox{ y } \lim_{n \rightarrow + \infty} \mbox{dist}(f^{-n}(q), f^{-n}(p))= 0\}.$$
Probar que $W^u(p)$ es una variedad inmersa en $\mathbb{R}^3$ de dimensi\'{o}n 1. Suge\-rencia: Fijar un punto cualquiera  $q \in W^u(p)$. No es restrictivo asumir que la coordenada $\varphi$ de $q$ es cero. Usando las secciones $A_N$ definidas en la parte (c) probar que la componente conexa de la intersecci\'{o}n de $W^u(p)$ con un entorno peque\~{n}o de   punto $q  $ en $\mathbb{R}^3$, es un arco diferenciable.

(e) Probar que el atractor $K$ es:
$$K = \bigcup_{p \in K} W^u(p).$$

(f) Probar que $K$ es un atractor topol\'{o}gico ca\'{o}tico. Esto significa que para alguna medida invariante $\mu$ soportada en $K$, para $\mu$-c.t.p. $p \in K$ y para toda direcci\'{o}n tangente inestable $0 \neq v \in T_p (W^u(p)) = T_p(K)$ el exponente de Lyapunov de $v$ hacia el futuro es positivo. Alcanza con probar que:
$$\liminf_{n \rightarrow + \infty}\frac{\log \|df^n(v)\|}{n} \geq \log 2 >0 \ \ \forall \ p \in K, \ \ \forall \ 0 \neq v \in T_p(W^u(p)).$$
(Sugerencia: probar que $\|df_p(v)\|/\|v\| \geq 2  $ para todo $p \in K$ y para todo $0 \neq v \in T_p(W^u(p))$.)

\end{exercise}

\subsection{Atractores erg\'{o}dicos}

 En esta secci\'{o}n consideraremos $f\colon  M \mapsto M$ continua en una varie\-dad compacta y riemanniana $M$.  Por el teorema erg\'{o}dico de Birkhoff-Khinchin, cualquiera sea la medida invariante $\mu$, para $\mu$-c.t.p. $x \in M$ y para toda funci\'{o}n   $\psi \in L^1(\mu)$ (en particular para toda funci\'{o}n continua $\psi$) existe el promedio temporal asint\'{o}tico $\widetilde \psi $ (en el futuro) definido por: \index{promedio! temporal}
 $$\widetilde \psi (x) = \lim_{n \rightarrow + \infty} \frac{1}{n} \sum_{j= 0}^{n-1} \psi (f^j(x)).$$
 Es decir,   para $\mu$-c.t.p. $x $  fijo en la variedad $M$, est\'{a} definido el funcional lineal
 $$\psi \in C^0(M, \mathbb{R}) \ \mapsto \ \widetilde \psi (x).$$
 Por el teorema de Representaci\'{o}n de Riesz, existe una medida de probabilidad $\mu_x$ tal que
 $$\widetilde \psi (x) = \int \psi \, d \mu_x \ \ \ \mbox{ para } \mu-{\mbox{c.t.p. }} x \in M, \ \ \forall \ \mu \in {\mathcal M}_f,$$
 donde ${\mathcal M}_f$ denota el espacio de todas las medidas de probabilidad en la sigma-\'{a}lgebra de Borel de $M$ que son $f$-invariantes. En los cap\'{\i}tulos anteriores probamos que $\mu_x \in {\mathcal M}_f$. Adem\'{a}s, dotando el espacio ${\mathcal M}$ de probabilidades (no necesariamente $f$-invariantes) de la topolog\'{\i}a d\'{e}bil-estrella, tenemos:
 $$\mu_x = \lim_{n \rightarrow + \infty} \frac{1}{n} \sum_{j= 0}^{n-1} \delta_{f^j(x)} \ \ \ \mbox{ para } \mu-{\mbox{c.t.p. }} x \in M, \ \ \forall \ \mu \in {\mathcal M}_f, $$
 donde el l\'{\i}mite es tomado en ${\mathcal M}$ con la topolog\'{\i}a d\'{e}bil$^*$ y $\delta_y$ denota la probabilidad Delta de Dirac soportada en el punto $y \in M.$

  En las secciones anteriores vimos adem\'{a}s que $\mu_x $ es erg\'{o}dica para $\mu$-c.t.p. $x \in M$, para toda $\mu \in {\mathcal M}_f$. Luego, el promedio temporal asint\'{o}tico  $\widetilde \psi (x)$ coincide con el valor esperado de $\psi$ con respecto a la probabilidad $\mu_x$.

 Por un lado, notamos que la teor\'{\i}a erg\'{o}dica desarrollada hasta ahora es v\'{a}lida para $\mu$-c.t.p. $x \in M$, y en general, no para cualquier punto $x \in M$. Dicho de otra forma, si el criterio de selecci\'{o}n de los puntos iniciales $x \in M$ no es $\mu$.c.t.p. para alguna medida $\mu \in {\mathcal M}_f$, entonces no necesariamente existen   los promedios asint\'{o}ticos de Birkhoff $\widetilde \psi (x)$. Y si existen, en general, no coinciden con el promedio espacial que resulta de integrar $\psi$  respecto a una medida erg\'{o}dica.

 Por otro lado, cuando se tiene un atractor topol\'{o}gico $K$ con cuenca local $V $ (abierta), el criterio de selecci\'{o}n de los estados iniciales, por la propia definici\'{o}n de atractor topol\'{o}gico, reside en tomar los puntos $x \in V$ para alg\'{u}n abierto $V$.  El criterio topol\'{o}gico de relevancia u \lq\lq observabilidad\rq\rq \ de conjuntos de \'{o}rbitas, es que estos conjuntos sean abiertos, o m\'{a}s en general,  con interior no vac\'{\i}o.

 Sin embargo, en la mayor\'{\i}a de los ejemplos de atractores topol\'{o}gicos que vimos en la secci\'{o}n anterior, el atractor $K $ (que est\'{a}   contenido en su cuenca local $V$) tiene interior vac\'{\i}o. Adem\'{a}s vimos (Proposici\'{o}n \ref{propositionMedidaInvarianteAtractorTopologico}), que las medidas   $\mu$ invariantes por $f|_{V}$ (que siempre existen) est\'{a}n soportadas en $K$. Entonces $\mu$-c.t.p. $x $ de la cuenca $V$, est\'{a}  en $K$. Luego, el teorema erg\'{o}dico de Birkhoff, y el teorema de existencia de medidas invariantes y erg\'{o}dicas, \em no aseguran la existencia de los promedios temporales asint\'{o}ticos, \em  para las \'{o}rbitas con punto inicial en  \em un conjunto de estados iniciales relevante u \lq\lq observable\rq\rq, desde el punto de vista topol\'{o}gico. \em

 Notamos que, en general, no hay esperanza que los promedios temporales asint\'{o}ticos existan para las \'{o}rbitas con punto inicial en un conjunto de estados iniciales con interior no vac\'{\i}o (es decir, relevante u \lq\lq observa\-ble\rq\rq, desde el punto de vista topol\'{o}gico). En efecto, si $K$ es un atractor topol\'{o}gico con cuenca de atracci\'{o}n $C(K)$ (abierta),  se puede demostrar que si $f|_{C(K)}$ no es \'{u}nicamente erg\'{o}dica (y en general no lo es), entonces el conjunto de estados iniciales $x \in C(K)$ para los cuales no existe el promedio asint\'{o}tico de Birkhoff, es denso. Luego, el conjunto de estados iniciales para los cuales existe ese promedio temporal asint\'{o}tico, tiene interior vac\'{\i}o.

 Por estos motivos, entre otras razones, se adopta otro criterio de \lq\lq obser\-va\-bi\-lidad\rq\rq \ de las \'{o}rbitas o de selecci\'{o}n de los estados iniciales. Es  un criterio  medible en vez de topol\'{o}gico, pero que considera a la mayor\'{\i}a de los puntos de la cuenca $C(K)$.

 \begin{definition} \label{remarkObservabilidadEstadistica} {\bf Criterio de observabilidad medible } \em \index{observabilidad!}

 Cuando el espacio es una variedad riemanniana $M$, \em el criterio medible de   \lq\lq observabilidad\rq\rq \ de las \'{o}rbitas, es que formen   un conjunto con medida de Lebesgue positiva. \em

 \end{definition}

 Lo anterior justifica las siguientes definiciones:

 \begin{definition} \label{definitionAtractorErgodico}
 {\bf Atractor erg\'{o}dico} \index{atractor! erg\'{o}dico}

 \em Sea $f: M \mapsto M$ continua en una variedad compacta y Riemanniana $M$ de dimensi\'{o}n finita. Denotamos con $m$ a la medida de Lebesgue en $M$. Notamos que $m$ no es necesariamente $f$-invariante.

  Un conjunto compacto no vac\'{\i}o $K$   se llama \em atractor erg\'{o}dico \em si:

   $\bullet$ $f^{-1}(K) = K = f(K)$

   $\bullet$ Existe un abierto $V \supset K$   tal que
   \begin{equation}\label{eqn28z} \lim_{n \rightarrow + \infty} \mbox{dist}(f^n(x), K) = 0 \   \mbox{ para Lebesgue-c.t.p. } \ x \in V.\end{equation}

    $\bullet$
  \begin{equation} \label{eqn30} \exists \ \mu \mbox{ erg\'{o}dica tal que: } \mu(K)= 1 \mbox{ y } $$ $$ \exists \ \widetilde \psi (x) = \lim_{n \rightarrow + \infty} \frac{1}{n} \sum_{j= 0}^{n-1} \psi (f^j(x)) = \int \psi \, d \mu \  \end{equation} $$\mbox{ para Lebesgue-c.t.p. } x \in V \ \    \mbox{ y } \ \ \ \forall \ \psi \in C^0(M, {\mathbb{R}}).$$

  Un abierto $V \supset K$ que satisface las condiciones anteriores,  se llama \em cuenca local de atracci\'{o}n \em del atractor erg\'{o}dico $K$.

  \vspace{.3cm}

  {\bf Nota:} No todos  adoptan la  Definici\'{o}n \ref{definitionAtractorErgodico}. Algunos autores exigen que $K$ sea, por definici\'{o}n, un atractor topol\'{o}gico que satisface adem\'{a}s la condici\'{o}n (\ref{eqn30}), para llamarlo atractor erg\'{o}dico (ver por ejemplo \cite{PughShubErgodicAttractors}).

  {\bf Medida SRB o f\'{\i}sica soportada en un atractor erg\'{o}dico: } \index{medida! SRB} \index{medida! f\'{\i}sica}

  Se observa que cuando existe una medida $\mu$ (erg\'{o}dica;  soportada en el atractor erg\'{o}dico $K$) que satisface la igualdad (\ref{eqn30})   para Lebesgue-c.t.p. $x \in V$, entonces esta medida es \em \'{u}nica. \em
  Tal medida $\mu$    se llama \em medida SRB erg\'{o}dica \em o tambi\'{e}n \em medida f\'{\i}sica erg\'{o}dica, \em del atractor erg\'{o}dico $K$.

  En \ref{definitionMedidaSRB} definiremos medida  SRB o medida  f\'{\i}sica $\mu$, en un contexto m\'{a}s  general, aunque $\mu$ no est\'{e}  soportada  en un atractor  erg\'{o}dico.

\vspace{.3cm}

  {\bf Observaci\'{o}n: } Dado cualquier atractor topol\'{o}gico $K$ con cuenca local $V$, debido a la Proposici\'{o}n \ref{propositionMedidaInvarianteAtractorTopologico} siempre existen medidas invariantes y erg\'{o}dicas para $f|_V$, y est\'{a}n soportadas   en $K$. Si alguna de estas medidas erg\'{o}dicas $\mu$   satisface la condici\'{o}n (\ref{eqn30}) para Lebesgue-casi todo punto  $x \in V$, entonces $K$ es un atractor erg\'{o}dico.

  \end{definition}

  \begin{remark} \em .

   \index{observabilidad}

  A diferencia de los atractores topol\'{o}gicos, para los atractores erg\'{o}dicos $K$  la atracci\'{o}n a $K$ de las \'{o}rbitas en su cuenca $V$  dada por la igualdad (\ref{eqn28z})  es solo para Lebesgue c.t.p. $x \in V$, y no necesariamente para todo punto $x \in V$. Es un \em criterio de observabilidad medible Lebesgue\em-c.t.p. de la cuenca. Entonces un atractor erg\'{o}dico no es necesariamente un atractor topol\'{o}gico. En el Ejercicio \ref{ejercicioEjemploMagnetico} se muestra un ejemplo de atractor erg\'{o}dico que no es topol\'{o}gico.

  Por otra parte, un atractor topol\'{o}gico satisface la igualdad (\ref{eqn28}) para todo $x \in V$. Luego satisface  (\ref{eqn28z}). Pero no necesariamente satisface la condici\'{o}n de existencia de una medida SRB erg\'{o}dica para la cual valga la igualdad (\ref{eqn30}). Entonces un atractor topol\'{o}gico no es necesariamente un atractor erg\'{o}dico. En el ejercicio \ref{ejercicioEjemploFacil} se muestra un ejemplo de atractor topol\'{o}gico que no es erg\'{o}dico.

   \begin{exercise}\em \label{ejercicioEjemploFacil}
   Sea en $Q= [0,1]^2$ la aplicaci\'{o}n $T(x,y) = \big((1/2) x, y\big)$. Probar que el segmento $ K= \{0\} \times [0,1]$ es un atractor topol\'{o}gico pero no es un atractor erg\'{o}dico. Sugerencia: para probar que $K$ no es atractor erg\'{o}dico, demostrar que toda medida $\mu$   erg\'{o}dica es delta de Dirac en un punto fijo y que el conjunto de puntos $x \in Q$ para los cuales vale la igualdad (\ref{eqn28z}) tiene medida de Lebesgue cero.
   \end{exercise}

   \begin{exercise}\em
   \label{ejercicioEjemploMagnetico} Sea en el disco $D= \{z \in \mathbb{C}: |z| \leq 1\}$ la transformaci\'{o}n $f: D \mapsto D$ que deja fijo el origen y tal que para todo $z \neq 0$ expresado en  polares, $f(z)$ est\'{a} dado por la siguiente f\'{o}rmula:
   $$f(\rho e^{i \varphi}) = \widehat \rho e^{i \widehat \varphi},$$
   donde $$\widehat \varphi = \frac{3}{2} \varphi \mbox{ si } 0 \leq \varphi < \pi  (\mbox{m\'{o}d} 2 \pi),$$
   $$\widehat \varphi = \pi + \frac{\varphi}{2} \mbox{ si } \pi \leq \varphi < 2\pi  (\mbox{m\'{o}d} 2 \pi),$$
   $$\widehat \rho = \Big (1- \frac{\widehat \varphi}{2 \pi}   \Big ) \rho.$$
   (a) Bosquejar las \'{o}rbitas (Sugerencia: los puntos $\varphi = 0$ son fijos, y las dem\'{a}s \'{o}rbitas son tales que el argumento tiende a $2 \pi$ por abajo y el m\'{o}dulo tiende a cero).

   (b) Probar que toda \'{o}rbita con punto inicial $z$ que no se encuentre en el segmento $\varphi=0$ es tal que la distancia al origen tiende a cero y la sucesi\'{o}n de promedios de Birkhoff de las funciones continuas $\psi$ tiende a $\int \psi \, d \delta_0$, donde $\delta_0$ es la Delta de Dirac soportada en el origen.

   (c) Concluir que $K= \{0\}$ es un atractor erg\'{o}dico pero no es atractor topol\'{o}gico.

   Nota: En este ejemplo $f$ es discontinuo en el semieje real positivo. Sin embargo, puede construirse un ejemplo continuo con bosquejo similar de \'{o}rbitas.
   \end{exercise}
\end{remark}
\begin{remark} \em  \label{remarkObservabilidadTopyEstad} \index{observabilidad} \index{cuenca de atracci\'{o}n! observabilidad de}
   Por un lado tenemos el criterio de obser\-va\-bilidad de la cuenca local $V$ del atractor. O bien la atracci\'{o}n se produce para todo estado inicial $x$ en el abierto $V$ (criterio de observabilidad topol\'{o}gica) o bien la atracci\'{o}n se produce solo para   Lebesgue c.t.p.   $x \in V$ (criterio de observabilidad Lebesgue-medible). El criterio de observabilidad de la cuenca es entonces el \em criterio con el cual se eligen los estados iniciales \em para observar a d\'{o}nde son atra\'{\i}das las \'{o}rbitas.

   \begin{definition}\label{definitionAtraccionTopol} {\bf Atracci\'{o}n topol\'{o}gica } \em En forma independiente al criterio de  \index{atracci\'{o}n! topol\'{o}gica}  obser\-va\-bilidad de los estados iniciales en la cuenca local,  la atracci\'{o}n en s\'{\i} misma,  definida por la igualdad (\ref{eqn28}) para los atractores topol\'{o}gicos, y por la igualdad (\ref{eqn28z}) para los atractores erg\'{o}dicos  se llama \em atracci\'{o}n topol\'{o}gica. \em
    Esta significa, por definici\'{o}n,     que el l\'{\i}mite de la distancia al atractor $K$   existe y es cero, desde los puntos iniciales elegidos seg\'{u}n el criterio de observabilidad que corresponda.

\end{definition}

\begin{definition}\label{definitionAtraccionEstad} \index{atracci\'{o}n! estad\'{\i}stica}

  {\bf Atracci\'{o}n estad\'{\i}stica. } \em Esta significa, por definici\'{o}n, que los promedios de Birkhoff  de las funciones continuas convergen (o por lo menos, en un contexto m\'{a}s general, tienen subsucesiones convergentes) al   valor  esperado  respecto a alguna  medida  invariante  $\mu$ soportada  en el atractor $K$, desde los puntos iniciales elegidos seg\'{u}n el criterio de observabilidad que corresponda.  En un contexto m\'{a}s general, $\mu$ no es necesariamente erg\'{o}dica.

   Lo usual es que cuando se estudia la atracci\'{o}n estad\'{\i}stica, el criterio de selecci\'{o}n de puntos iniciales sea el de observabilidad medible Lebesgue c.t.p.

   \end{definition}
  \end{remark}

  En los pr\'{o}ximos ejemplos veremos   casos particulares de existencia y de no existencia de atractor erg\'{o}dico.

 \begin{example} \em  \label{ejemploTentMapEnDisco} \index{tent map}
 En el tent map del intervalo, todo el intervalo es un atractor topol\'{o}gico y erg\'{o}dico a la vez, cuya    medida erg\'{o}dica SRB   (o f\'{\i}sica) es la medida de Lebesgue (ver por ejemplo \cite{buzziEnciclopedia}).

 \end{example}

 \begin{example} \em  \index{flujo polo norte-polo sur}
  En el flujo Polo Norte-Polo Sur, el Polo Sur es atractor erg\'{o}dico y topol\'{o}gico a la vez, cuya media erg\'{o}dica SRB (o f\'{\i}sica) es la delta de Dirac soportada en el Polo Sur.
 \end{example}

 \begin{example} \em . \label{ejemploRotacionEsfera}

 {\bf Contraejemplo: Rotaci\'{o}n irracional de la esfera. } \index{rotaci\'{o}n! irracional}

 En este ejemplo  no existen atractores erg\'{o}dicos. Sea $f: S^2 \mapsto S^2$   conti\-nua, definida en la superficie esf\'{e}rica $S^2$ por la siguiente para\-me\-tri\-zaci\'{o}n en coordenadas \lq\lq esf\'{e}ricas\rq\rq:

    $S^2 : \{(x,y,z) \in \mathbb{R}^3 : \ \ x = \cos \varphi \sen \theta, \ \ y= \sen \varphi \sen \theta, \ \ z= \cos \theta, \ \ 0 \leq \varphi \leq 2 \pi, \ \ 0 \leq \theta \leq \pi\}$

    $f$ est\'{a} definida por las siguientes ecuaciones  en coordenadas esf\'{e}ricas:

    \ \ \ \ \ \ \ \ \ \ \ \ \ \ \ \ \ \ $f(\varphi, \theta) = (\widehat \varphi , \widehat \theta )$ donde:
  \begin{eqnarray}
   \widehat \theta & = \theta & \nonumber \\
   \widehat \varphi & = & \varphi + a\nonumber
  \end{eqnarray}
 siendo $a \in (0, 2 \pi)$ una constante tal que $a/2 \pi$ es irracional.

 Es inmediato chequear que los polos $C_0 = \{\theta= 0\}$ y $C_{\pi} = \{\theta = \pi\}$ son puntos fijos por $f$, y que cada circunferencia $C_{\theta}$, con $\theta$ constante diferente de $0$ y de $\pi$, es invariante por $f$. Adem\'{a}s, $f|_{C_{\theta}}$ es una rotaci\'{o}n irracional si $\theta \neq 0, \pi$. Por lo tanto toda   \'{o}rbita  por $f$ es densa  en la secci\'{o}n $C_{\theta}$ donde est\'{a} contenida. Deducimos que no hay \'{o}rbitas densas en $S^2$. Luego, $f: S^2 \mapsto S^2$ no es transitivo.

 Sea $m$   la medida de Lebesgue bidimensional en la esfera $S^2$ (normalizada para que sea una probabilidad, es decir, dividimos la medida de Lebesgue en la esfera, entre el \'{a}rea de toda la esfera, para que  $m(S^2) = 1$). Esta medida es invariante con $f$, pues el Jacobiano $|\mbox{det} df|$ es id\'{e}nticamente igual a 1.  Sin embargo $m$ no es erg\'{o}dica, pues $m$ es positiva sobre abiertos pero $f$ no es transitivo.

 \vspace{.3cm}

 {\bf Afirmaci\'{o}n: } \em  No existen atractores erg\'{o}dicos para la rotaci\'{o}n irracional $f$ en la esfera $S^2$. \em

 {\em Demostraci\'{o}n: }
 Por absurdo, supongamos que existe    un atractor erg\'{o}dico $K$ y llamemos $V$  a su cuenca local de atracci\'{o}n. Entonces $K$ es compacto no vac\'{\i}o e invariante por $f$ y $V$ es abierto que contiene a $K$. Sea $p \in V$. Vimos que $\omega(p) = C_{\theta_p}$ donde $  C_{\theta_p}$ es la secci\'{o}n horizontal de la esfera que contiene al punto $p$. De la igualdad (\ref{eqn28z}) deducimos $\omega(p) \subset K$ para Lebesgue c.t.p. $p \in V$. Entonces $C_{\theta_p} \subset K$, y como $p \in C_{\theta_p}$ deducimos que $p \in K$. Hemos probado que, bajo la hip\'{o}tesis de absurdo, Lebesgue c.t.p. $p \in V$ est\'{a} contenido en $K$. Como $K$ es compacto, tenemos $V \subset K$. Pero por definici\'{o}n de atractor erg\'{o}dico $K \subset V$. Entonces $K= V$ es compacto y abierto a la vez, y es no vac\'{\i}o. Como $S^2$ es conexo concluimos que, si existiera un atractor erg\'{o}dico $K$, \'{e}ste ser\'{\i}a toda la esfera $K= S^2$.

  Por (\ref{eqn30}), el promedio temporal asint\'{o}tico $\widetilde \psi (p)$ deber\'{\i}a ser constante para $m$-c.t.p. $p \in V= S^2$, para cualquier funci\'{o}n continua $\psi: S^2 \mapsto \mathbb{R}$. Como en este ejemplo $m$ es invariante, entonces $m$ ser\'{\i}a una medida erg\'{o}dica, contradiciendo que $m$ es positiva sobre abiertos y $f$ no es transitivo. \hfill $\Box$

  {\bf Nota: } En este ejemplo \ref{ejemploRotacionEsfera} toda la esfera es un atractor topol\'{o}gico (en realidad es el \'{u}nico atractor topol\'{o}gico). Luego, este ejemplo prueba que no todo atractor topol\'{o}gico es un atractor erg\'{o}dico.

 \end{example}

\begin{example} \em  \index{automorfismo! lineal en el toro}
Sea $f =  {\left [
           \begin{array}{cc}
             2 & 1 \\
             1 & 1 \\
           \end{array}
          \right ]  }_{\mbox{m\'{o}d}\displaystyle (1,1)}
$ en el toro $\mathbb{T}^2$. Es erg\'{o}dico respecto a la medida de Lebesgue $m$. Todo el toro $\mathbb{T}^2$ es un atractor erg\'{o}dico y $m$ es la medida erg\'{o}dica SRB o f\'{\i}sica.

\end{example}

 \subsection{Atracci\'{o}n estad\'{\i}stica y medidas SRB o f\'{\i}\-si\-cas}

 Independientemente de si existe o no un atractor erg\'{o}dico, dada un medida de probabilidad $\mu$ (no necesariamente erg\'{o}dica ni invariante) definimos el siguiente conjunto $B(\mu)$ en el espacio $X$ donde act\'{u}a $f$:

 \begin{definition} {\bf Cuenca de atracci\'{o}n estad\'{\i}stica} \label{DefinicionCuencaDeAtraccionEstadistica} \em \index{cuenca de atracci\'{o}n! estad\'{\i}stica} \index{medida! atracci\'{o}n estad\'{\i}stica de} \index{$B(\mu)$ cuenca de atracci\'{o}n! estad\'{\i}stica de la medida $\mu$}

 Sea $f: X \mapsto X$ una transformaci\'{o}n continua en un espacio m\'{e}trico compacto.
 Sea $\mu \in {\mathcal M}$ una probabilidad.

 Se llama \em cuenca de atracci\'{o}n estad\'{\i}stica \em de $\mu$ al siguiente conjunto:
 $$B(\mu): = \{ x \in X: \ \ \lim_{n \rightarrow + \infty} \frac{1}{n} \sum_{j= 0}^{+ \infty} \psi(f^j(x)) = \int \psi \, d \mu \ \ \forall \ \psi \in C^0(X, \mathbb{R}) \}.$$
 Usando la caracterizaci\'{o}n de la topolog\'{\i}a d\'{e}bil$^*$ en el espacio ${\mathcal M}$ obte\-nemos:
 \begin{eqnarray}
 \label{equationB(mu)}
 &B(\mu) = \{ x \in X: \ \ \lim_{n \rightarrow + \infty} \sigma_{n,x}  = \mu\}, & \mbox{ d\'{o}nde } \\ \label{(1)}
 &\sigma_{n,x} := \frac{1}{n} \sum_{j= 0}^{n-1} \delta_{f^j(x)}.&
 \end{eqnarray} \index{$\sigma_{n,x}$ probabilidad emp\'{\i}rica}

 {\bf Observaci\'{o}n:} La cuenca de atracci\'{o}n estad\'{\i}stica de cualquier medida de probabilidad $\mu$ es un conjunto medible. En efecto, tomando una familia numerable $\{\psi_i\}_{i \in \mathbb{N}}$ de funciones reales continuas $\psi: X \mapsto [0,1]$ que sea denso en el espacio $C^0(X, [0,1])$, la igualdad dentro de la definici\'{o}n de la cuenca $B(\mu)$, se verifica para toda $\psi \in C^0(X, \mathbb{R})$ si y solo si se verifica para   $\psi_i$ para todo $i \in {\mathbb{N}}$. Para cada $i$ fijo, la igualdad se satisface para un conjunto medible, pues el l\'{\i}mite puntual de una sucesi\'{o}n de funciones continuas es medible. Luego, $B(\mu)$ es la intersecci\'{o}n numerable de conjuntos medibles; es decir, es medible.
\end{definition}
\begin{definition}
\label{definitionProbaEmpiricas}

 {\bf Probabilidades emp\'{\i}ricas}  \em \index{probabilidad! emp\'{\i}rica}
\index{medida! de probabilidad emp\'{\i}rica}

 \index{$\sigma_{n,x}$ probabilidad emp\'{\i}rica}
 Se observa que $\sigma_{n,x}$, definida por la igualdad (\ref{(1)}), es una medida de probabilidad, para todo $n \geq 1$ y para todo $x \in X$ (en general $\sigma_{n,x}$ no es $f$ invariante). Es la probabilidad promedio concentrada en los puntos de tramos finitos de la \'{o}rbita futura por $x$.

 Las probabilidades $ \sigma_{n,x}$ (para cualquier $n \geq 1$, con $x \in X$ fijo) se llaman \em  probabilidades emp\'{\i}ricas  \em de la \'{o}rbita futura por $x$. Este nombre proviene de que un experimentador  no puede observar los promedios temporales asint\'{o}ticos (en tiempo infinito), sino que observa los promedios hasta tiempo $n$  finito. Estos promedios se pueden calcular como el valor esperado integrando respecto a las probabilidades emp\'{\i}ricas. M\'{a}s precisamente, debido a la igualdad (\ref{(1)})  tenemos:
 $$\frac{1}{n} \sum_{j= 0}^{n-1} \psi (f^j(x)) = \int \psi \, d \sigma_{n,x} \ \ \ \forall \ \ \psi \in C^0(X, \mathbb{R}).$$
\index{promedio! temporal} \index{promedio! de Birkhoff}
 \end{definition}
 \begin{remark} \em
 Si $\mu$ es una probabilidad tal que su cuenca de atracci\'{o}n estad\'{\i}stica $B(\mu)$ es no vac\'{\i}a, entonces $\mu$ es invariante con $T$ (Ejercicio \ref{ejercicioZ0} a)).  Se puede caracterizar a una medida erg\'{o}dica por la siguiente afirmaci\'{o}n (Ejercicio \ref{ejercicioZ0} b)):
 $$\mu \mbox{ es invariante y erg\'{o}dica si y solo si } \mu  (B(\mu)) = 1.$$
 Para cualquier   medida de probabilidad $\mu$ no invariante, y para cualquier medida de probabilidad invariante no erg\'{o}dica, se cumple $\mu (B(\mu)) = 0$   (Ejercicio \ref{ejercicioZ0} c)), aunque $B(\mu)$ puede ser no vac\'{\i}o.
 \end{remark}

 \begin{remark} \em
  \label{remarkBowen} \index{cuenca de atracci\'{o}n! estad\'{\i}stica} \index{medida! atracci\'{o}n estad\'{\i}stica de}
  \em La cuenca de atracci\'{o}n estad\'{\i}stica $B(\mu)$ de una medida invariante $\mu$ no erg\'{o}dica, puede ser no vac\'{\i}a, y adem\'{a}s puede cubrir Lebesgue c.t.p. \em

  En efecto, un ejemplo de homeomorfismo   no diferenciable en un disco bidimensional compacto $D$ (que es adaptaci\'{o}n de un difeomorfismo  en el disco compacto, atribuido a Bowen) exhibe una medida  invariante y \em no erg\'{o}dica \em    $\mu$ tal que $B(\mu) \neq \emptyset$ (ver  \cite[Example 7.2, case B] {CatIlyshenkoAttractors}). M\'{a}s a\'{u}n, en este ejemplo, $B(\mu)$ contiene Lebesgue casi todo punto $x \in D$ a pesar que $\mu(B(\mu)) = 0$, es decir, el soporte de $\mu$ tiene medida de Lebesgue nula. Luego, en este ejemplo,  existe una medida SRB $\mu$ no erg\'{o}dica, la medida de Lebesgue $m$ no es inva\-rian\-te, existe un atractor topol\'{o}gico no erg\'{o}dico cuya cuenca de atracci\'{o}n topol\'{o}gica es abierta y cubre Lebesgue casi todo punto, y no existen atractores erg\'{o}dicos.

 En la versi\'{o}n   $C^{1}$ del ejemplo atribuido a Bowen (ver   \cite{Golenishcheva}),   se cumple  $m(B(\mu)) = 0$ para toda medida de probabilidad $\mu$ invariante soportada en el atractor topol\'{o}gico. Luego, no existen medidas SRB. Adem\'{a}s en este ejemplo, toda   medida invariante es hiperb\'{o}lica (tiene exponentes de Lyapunov no nulos) y tiene exponentes de Lyapunov positivos. Como no existen medidas SRB, no   existen atractores erg\'{o}dicos.  Sin embargo, existe un atractor topol\'{o}gico cuya cuenca es abierta y cubre Lebesgue casi todo punto.
 \end{remark}

 \begin{exercise}\em
 \label{ejercicioZ0}
 Sea $f: X \mapsto X$ continua en un espacio m\'{e}trico $X$. Denotamos ${\mathcal M}$ el espacio de todas las probabilidades de Borel (no necesariamente $f$-invariantes), dotado de la topolog\'{\i}a d\'{e}bil$^*$.
 Sea $\mu \in {\mathcal M}$.

 {\bf a) } Probar que si la cuenca de atracci\'{o}n estad\'{\i}stica $B(\mu)$ es no vac\'{\i}a, entonces $\mu$ es $f$-invariante.

  {\bf b) } Probar que
 $\mu $   es invariante y erg\'{o}dica si y solo si  $ \mu  (B(\mu)) = 1.$

 {\bf c) }  Probar que si $\mu$ no es invariante, o si $\mu$ es invariante pero no erg\'{o}dica, entonces
  $$\mu (B(\mu)) = 0.$$ (Sugerencias: $B(\mu)$ es siempre un conjunto $f$-invariante. En el caso $\mu$ invariante, aplicar el teorema de descomposici\'{o}n erg\'{o}dica).

 \end{exercise}

 De la parte c) del Ejercicio \ref{ejercicioZ0} deducimos:

 \vspace{.2cm} \index{cuenca de atracci\'{o}n! estad\'{\i}stica} \index{medida! atracci\'{o}n estad\'{\i}stica de} \index{medida! erg\'{o}dica}

 \index{equivalencia de definiciones! de ergodicidad}

 \em Una medida de probabilidad $\mu$ es $f$-invariante y erg\'{o}dica si y solo si $$\mu(B(\mu)) = 1,$$
 donde $B(\mu)$ denota la cuenca de atracci\'{o}n estad\'{\i}stica de $\mu$. \em

%

 \begin{definition}
 {\bf Medidas SRB o f\'{\i}sicas }  {\em (Sinai \cite{Sinai_SRB}- Ruelle \cite{Ruelle_SRB}-Bowen \cite{Bowen1971,Bowen-Ruelle_SRB})}  \index{medida! SRB} \index{medida! f\'{\i}sica}  \index{cuenca de atracci\'{o}n! estad\'{\i}stica} \index{medida! atracci\'{o}n estad\'{\i}stica de} \label{definitionMedidaSRB}

 \em
 Sea $f: M \mapsto M$ continua en una variedad compacta y riemanniana $M$. Sea $m$ la medida de Lebesgue en $M$, normalizada para que sea una probabilidad: $m(M) = 1$. (En general, $m$ no es necesariamente $f$-invariante.) Sea $\mu \in {\mathcal M}$.

 Decimos que la medida de probabilidad $\mu$ es \em SRB (Sinai-Ruelle-Bowen) \em o, indistintamente, que es \em f\'{\i}sica, \em
 si $$m(B(\mu)) >0,$$
 donde $B(\mu)$ es la cuenca de atracci\'{o}n estad\'{\i}stica de $\mu$, definida en \ref{DefinicionCuencaDeAtraccionEstadistica}.
 \end{definition}

 {\bf Nota: } Si $\mu$ es SRB, entonces $B(\mu) \neq \emptyset$, y por lo observado en la parte (a) del Ejercicio \ref{ejercicioZ0}, $\mu \in {\mathcal M}_f$ (es decir las medidas SRB son invariantes con $f$). \index{medida! SRB}

 De acuerdo a la Definici\'{o}n \ref{definitionMedidaSRB}, una medida SRB puede no ser erg\'{o}dica. En efecto, en el ejemplo mencionado en \ref{remarkBowen}, que es adaptaci\'{o}n $C^0$ de un ejemplo  atribuido a Bowen,  existe   medida SRB no erg\'{o}dica. \index{medida! SRB no erg\'{o}dica}

 \vspace{.3cm}

 {\bf Sobre la nomenclatura \lq\lq SRB\rq\rq \  y \lq\lq  f\'{\i}sica\rq\rq. } \index{medida! SRB} \index{medida! f\'{\i}sica}
La Definici\'{o}n \ref{definitionMedidaSRB} de medida SRB no es adoptada por todos los autores. En general se utiliza    esta definici\'{o}n que solo requiere $m(B(\mu)) >0$, solo para llamar \em f\'{\i}sica \em  a la medida  $\mu$.
Pero, para una  parte importante de matem\'{a}ticos (por ejemplo \cite{buzziEnciclopedia}, \cite{YoungSurvey}),  medida SRB no es sin\'{o}nimo de medida f\'{\i}sica.  Llaman f\'{\i}sica a cualquier probabilidad $\mu$ cuya cuenca de atracci\'{o}n estad\'{\i}stica $B(\mu)$ tenga medida de Lebesgue positiva. Pero para llamar SRB a $\mu$,    requieren   que la probabilidad $\mu$  sea erg\'{o}dica, tenga exponentes de Lyapunov positivos, y   tenga medidas condicionales inestables absolutamente con\-ti\-nuas (veremos m\'{a}s adelante qu\'{e} significa esta propiedad adicional, al introducir las medidas de Gibbs, mediante las Definiciones \ref{definitionMedidasCondicionadasAC},  \ref{definitionMedidasCondicionasInestables}  y \ref{definitionMedidaGibbs}). Esto es debido a que, entre otros motivos, en el contexto   restringido de los difeomorfismos de clase $C^{1 + \alpha}$ uniformemente
hiperb\'{o}licos, ambas definiciones son equivalentes (ver Teorema \ref{TheoremSRBanosov}).

   Nosotros adoptamos la definici\'{o}n   de, por ejemplo,
\cite[Definition 1.9]{BonattiDiazVianaLibro}  o   \cite[Definition 22]{Jost}.
Aqu\'{\i}, en la Definici\'{o}n \ref{definitionMedidaSRB}, \em no estamos asumiendo ninguna  condici\'{o}n adicional a la fisicalidad de la medida para llamarla SRB o f\'{\i}sica (ni la ergodicidad, ni la positividad de exponentes de Lyapunov, ni la continuidad absoluta condicionada inestable).  En resumen, usamos ambas palabras \lq\lq SRB \ \rq\rq \'{o} \lq\lq f\'{\i}sica\rq\rq, como sin\'{o}nimos. \em

\section{Teor\'{\i}a de Pesin} \index{Pesin! teor\'{\i}a de} \index{teor\'{\i}a de Pesin} \label{chapterAtractores}

\subsection{Desintegraci\'{o}n en medidas condicionales}

En esta secci\'{o}n asumimos que $X$ es un espacio m\'{e}trico compacto, provisto de la sigma-\'{a}lgebra de Borel ${\mathcal B}$, y que $\mu$ es una medida de probabilidad en   $(X, {\mathcal B})$.
  Expondremos el enunciado de un resultado   de la teor\'{\i}a  de la medida (el Teorema de Rohlin), v\'{a}lido   aunque no exista una din\'{a}mica definida en el espacio $X$. En la secci\'{o}n siguiente veremos el uso del Teorema de Rohlin, junto con la Teor\'{\i}a de Pesin, en sistemas din\'{a}micos de clase $C^1$-m\'{a}s H\"{o}lder (cuando $X$ tiene adem\'{a}s, una estructura de variedad).

%
\begin{definition}
\label{definicionParticionMedible}

 {\bf Partici\'{o}n medible } \index{partici\'{o}n! medible}

 \em

 Se llama \em partici\'{o}n    \em    en $X$  
 a una colecci\'{o}n (puede ser finita, infinita numerable o infinita no numerable) $${\mathcal P} := \{W({x})\}_{x \in X}$$ de subconjuntos $  W({x}) \subset X$ (pueden ser medibles o no medibles)    tales que:

  {\bf (a) }   $ x \in W({x})$ para todo $x \in X$ (esto implica que $\bigcup_{x \in X} W(x) = X$)

  {\bf (b) } Los conjuntos $W(x)$ son dos a dos disjuntos; es decir, para toda pareja de puntos $x \neq y$ en $X$, \'{o}   bien $W(x) = W(y)$  \'{o} bien $W_x \cap W_y = \emptyset$.

  \vspace{.3cm}

  La partici\'{o}n ${\mathcal P}$ se dice   \em medible, \em si sus piezas $W(x)$ son todas medibles y ${\mathcal P}$ est\'{a} generada por una colecci\'{o}n numerable de particiones finitas con piezas medibles.  Esto es:

  {\bf (c) } Existe una colecci\'{o}n numerable $\{\mathcal P_n\}_{n \geq 1}$ donde:
  $${\bf  (c1) } \ \ {\mathcal P}_n := \{E_{n,j}\}_{1 \leq j \leq k_n} $$  es una partici\'{o}n finita de $X$ para todo $n \geq 1$, con exactamente $k_n$ piezas $E_{n,j} \subset X$ medi\-bles (disjuntas dos a dos al cambiar $j$ con $n$ fijo, y cuya uni\'{o}n en $j$ es $X$ para todo $n$ fijo).
  $${\bf (c2) } \ \ {\mathcal P}_{n+1}  \prec {\mathcal P}_n \ \ \forall \ n \geq 1.$$ \index{partici\'{o}n! m\'{a}s fina que}
  Esto significa que cada pieza de la partici\'{o}n ${\mathcal P}_{n+1}$ est\'{a} contenida en alguna pieza de la partici\'{o}n ${\mathcal P}_{n}$ (Se dice que ${\mathcal P}_{n+1}$ es \em m\'{a}s fina \em que ${\mathcal P}_{n}$ y se denota ${\mathcal P}_{n+1} \prec {\mathcal P}_n$).
   $${\bf (c3) } \ \ \forall \ x \in X: \ \ W(x) = \bigcap _{n \geq 1}  \  E_{n, j_n(x)} \in {\mathcal P}$$   donde $\ 1 \leq j_n(x) \leq k_n $ es el \'{u}nico \'{\i}ndice, para cada $n$ fijo, tal que $\ x \in E_{n, j_n(x)} \ \ \forall \ n \geq 1.$
  Observar que esta \'{u}ltima condici\'{o}n implica  que para todo $x \in X$ y para todos $n \geq 1$ y $1 \leq j \leq k_n$, o bien $x \not \in E_{n,j}$, o bien $W(x) \subset E_{n,j}$. En otras palabras, cada conjunto  medible  $E_{n,j}$ est\'{a}n formado  por piezas enteras de la partici\'{o}n ${\mathcal P}$. Dicho de otra forma: ${\mathcal P} \prec {\mathcal P}_n$ para todo $n \geq 1$.

  La condici\'{o}n (3) establece una condici\'{o}n m\'{a}s fuerte  que ${\mathcal P} \prec {\mathcal P}_n$ para todo $n \geq 1$:   La \lq\lq intersecci\'{o}n\rq\rq \   decreciente de las particiones ${\mathcal P}_n$ es ${\mathcal P}$. Esto se denota como $${\mathcal P} = \bigvee_{n= 1}^{+ \infty} {\mathcal P}_n,$$
  \index{partici\'{o}n! operaci\'{o}n $\vee$} \index{$\vee$ operaci\'{o}n! entre particiones} \index{$\prec$ relaci\'{o}n! entre particiones}
donde para cualquier pareja ${\mathcal R} , {\mathcal S}$ de particiones finitas se define
$${\mathcal R} \vee {\mathcal S} := \{R \cap S \colon  \ R \in {\mathcal R}, \ S \in {\mathcal S}\}.$$

  La condici\'{o}n (c3) (junto con la propiedad de medibilidad de las piezas de cada partici\'{o}n ${\mathcal P}_n$) implica, en particular, la medibilidad de las piezas de ${\mathcal P}$. Sin embargo el rec\'{\i}proco es falso, como veremos en el ejemplo del Ejercicio \ref{exercise2111particionNoMedible}: existen particiones ${\mathcal P}$ cuyas piezas son todas medibles y que no cumplen la condici\'{o}n (c). Estas particiones, no son particiones medibles, de acuerdo a esta definici\'{o}n, a pesar que sus piezas son todas medibles, porque no est\'{a} generada por ninguna colecci\'{o}n numerable de particiones con piezas medibles.
  \end{definition}

  \begin{exercise}\em
   Probar  que si ${\mathcal P}$ es una partici\'{o}n con piezas medibles, y si ${\mathcal P}$ es finita o infinita numerable, entonces ${\mathcal P}$ es una partici\'{o}n medible.
   \end{exercise}

%
%
%
%

En los siguientes ejercicios veremos ejemplos de particiones medibles y no medibles, con una cantidad no numerable de piezas:

 \begin{exercise}\em \label{ejericioParticionesMedibleshorizontales}
 Sea    $X= [0,1]^2$. Para cada $(x_0, y_0) \in X$ se denota $S(x_0, y_0) := \{(x,y) \in X: x = x_0\} $  al segmento    vertical de $X$ que    se obtiene seccionando $X$ con la recta $x= x_0$ constante.

{\bf (a) }  Probar que   la foliaci\'{o}n ${\mathcal F} := \{S(x,y)\}_{(x,y) \in X}$    es una partici\'{o}n medible de $X$. (Sugerencia: considerar la colecci\'{o}n numerable  $\{E_{n,i}\}_{n \geq 1, 1 \leq i \leq n}$ de borelianos $E_{n,i} := [(i-1)/n, i/n)  \times [0,1] $   si $i < n$, \  $E_{n,n} := [(n-1)/n, 1] \times [0,1]$.)

 {\bf (b) } En el intervalo $[0,1]$ sea $K = \cap_{n= 0}^{+ \infty} K_n$ el conjunto de Cantor de los \lq\lq tercios mitad\rq\rq, i.e. $K_0 = [0,1]$ y $K_{n+1}$ se obtiene de $K_n$    retirando de cada intervalo cerrado $I$ que forma a $K_n$ el subintervalo abierto  con punto medio en el punto medio de $I$ y con longitud $(1/3)\mbox{long(I)}$. Sea en $X$, la  partici\'{o}n $\{\mathcal P\} = \{P(x,y)\}_{(x,y) \in X}$ definida por $P(x_0,y_0) = S(x_0, y_0)$ si $x_0 \in K$, y $P(x_0, y_0) := \{(x,y) \in X: \ x \not \in K\}$ si $x_0 \not \in K$. Probar que ${\mathcal P}$ es una partici\'{o}n medible en $X$.

 {\bf (c) } Sea  $\xi: X \mapsto Y$ un homeomorfismo. Sea en $Y$ la foliaci\'{o}n

 ${  \mathcal G } := \{\xi(S(\xi^{-1}(x,y)))\}_{(x,y) \in Y}$. (Se dice que $\xi^{-1}$ es una trivializaci\'{o}n $C^0$ de la foliaci\'{o}n ${\mathcal G}$). Probar que ${\mathcal G}$ es una partici\'{o}n medible.

 {\bf (d) } Sean $X$ e $Y$ espacios m\'{e}tricos compactos, y sea en $X \times Y$ la partici\'{o}n ${\mathcal F} = \{S(x,y)\}_{(x,y) \in X \times Y}$ por secciones verticales $S(x_0,y_0):= \{(x,y) \in X \times Y: \ x= x_0\}.$ Probar que ${\mathcal F}$ es una partici\'{o}n medible.

\end{exercise}

\begin{exercise}\em
Sean $(X, {\mathcal A})$ e $(Y, {\mathcal B})$ dos espacios medibles. Sea $\xi: X \mapsto Y $ una transformaci\'{o}n bimedible. Sea en $X$ una partici\'{o}n

${\mathcal P} = \{P(x)\}_{x \in X}$ cualquiera (no necesariamente medible). Sea   en $Y$  la partici\'{o}n   ${\mathcal Q} := \{\xi(S(\xi^{-1}(y)))\}_{y \in Y} $. Probar que ${\mathcal P}$ es partici\'{o}n medible si y solo si ${\mathcal Q}$ lo es.
 \end{exercise}

 \begin{exercise}\em \label{exercise2111particionNoMedible} \index{partici\'{o}n! medible} \index{automorfismo! lineal en el toro}

 Sea $f: \mbox{Diff }^{\infty}(\mathbb{T}^2)$ en el toro $\mathbb{T}^2$ el autormorfismo lineal $$f = \left(
                             \begin{array}{cc}
                               2 & 1 \\
                               1 & 1 \\
                             \end{array}
                           \right)
   (\mbox{ m\'{o}d. }\mathbb{Z}^2) .$$   Sea $W^u(p)$ la variedad inestable (global) por cada punto $p \in \mathbb{T}^2$, i.e.: $$W^u(p) = \{q \in \mathbb{T}^2: \lim_{n \rightarrow + \infty} \mbox{dist}(f^{-n}(p), f^{-n}(q))= 0\}.$$

     (a) Probar que $W^u(p)$ es un conjunto medible para todo $p \in \mathbb{T}^2$.

     (b) Probar que la partici\'{o}n ${\mathcal P}:= \{W^u(p)\}_{p \in {\mathbb{T}^2}}$ no es medible.

   Sugerencia: Considerar la medida de Lebesgue $m$ en el toro. Chequear que $m(W^u(p)) = 0$ para todo $p$. Sea $E \subset \mathbb{T}^2$   medible con $m(E) >0$ y tal que si $p \in E$ entonces $W(p) \subset E$. Probar que  $m(E) = 1 $ (Recordar que $\mathbb{T}^2 = \mathbb{R}^2/\mathbb{Z}^2$. Tomar un levantado $\widehat E$ a $\mathbb{R}^2$ del conjunto $E$: si $\widehat p \in \widehat E$ entonces toda la recta con direcci\'{o}n inestable en $\mathbb{R}^2$ que pasa por $\widehat p$ est\'{a} contenida en $\widehat E$. Mirar la  intersecci\'{o}n de $\widehat E$ con cada  recta  horizontal de altura entera en $\mathbb{R}^2$. Proyectar todas esas intersecciones, siguiendo las verticales en $\mathbb{R}^2$,  sobre una sola recta horizontal. Observar que esa proyecci\'{o}n  (m\'{o}d. $1$) en $[0,1] = \mathbb{R}/\mathbb{Z}$, es invariante con la rotaci\'{o}n irracional y tiene medida de Lebesgue positiva en el intervalo.) De la propiedad $m(E) = 0$ \'{o} $m(E) = 1$, deducir que para $m$-c.t.p. $p \in \mathbb{T}^2$ no se puede verificar la  condici\'{o}n (c3)   de la Definici\'{o}n \ref{definicionParticionMedible}) de partici\'{o}n medible.

   \end{exercise}

\begin{exercise}\em  \label{exerciseMedibleSiiMuMedible}

(a) Probar que si ${\mathcal P}$ y ${\mathcal Q}$ son dos particiones medibles, entonces ${\mathcal P} \vee {\mathcal Q}$ es una partici\'{o}n medible. (Sugerencia, construir ${\mathcal P}_n \vee {\mathcal Q}_n$ donde $ {\mathcal P}_n $ y ${\mathcal Q}_n$ son las particiones finitas que satisfacen, para ${\mathcal P}$ y ${\mathcal Q}$ respectivamente, la condici\'{o}n (c) de la Definici\'{o}n \ref{definicionParticionMedible}.)

(b) Probar que si ${\mathcal Q}$ es medible, entonces ${\mathcal P}$
es medible si y solo si ${\mathcal P} \vee {\mathcal Q}$ es medible.

\end{exercise}

%
%

 \vspace{.3cm}

 {\bf Espacio cociente por una partici\'{o}n medible.} \index{espacio! cociente por partici\'{o}n}

 Dada una partici\'{o}n medible ${\mathcal P}= \{W(x)\}_{x \in X}$ de $(X, {\mathcal B})$ consideramos el conjunto cociente $X /\sim$ por la relaci\'{o}n de equivalencia $x \sim y$ si y solo si $W(x) = W(y)$. Denotando $[x]$ a la clase de equivalencia que contiene a $x$, tenemos
 $$X/\sim  = \{[x]\}_{x \in X} \ \mbox{ donde }  \ [x]  = W(x) $$
 Dotamos el conjunto  $X/\sim$ de la estructura medible cociente, definiendo la sigma-\'{a}lgebra:
 $$(X/\sim, {\mathcal A}/\sim): \ \ \ \widehat A \in {\mathcal A}/\sim \ \ \ \mbox{ si y solo si }  $$ $$A :=   \{x \in X: W(x) \in \widehat A\}   \in {\mathcal A}.  $$

 Para cada $W(x) \in X / \sim$, es decir, en cada pieza   $W(x)$ de la partici\'{o}n medible ${\mathcal P}$,  definimos la sigma-\'{a}lgebra ${\mathcal A}_{W(x)}$  que resulta de restringir la sigma-\'{a}lgebra ${\mathcal A}$ a    $W (x)$. Esto es:
 $$ (W(x), {\mathcal A}_{W(x)}): \   \forall \  A  \subset W(x), \ \ A  \in {\mathcal A}_{W(x)} \mbox{ si y solo si } A \in {\mathcal A}.$$

\begin{theorem}
{\bf   Desintegraci\'{o}n Medible  (Teorema de Rohlin)} \index{teorema! Rohlin} \index{teorema! de descomposici\'{o}n medible}  \index{teorema! de existencia de! medidas condicionales} \label{theoremRohlin}

 Sea $X$ un espacio m\'{e}trico compacto y sea ${\mathcal A}$ la sigma-\'{a}lgebra de Borel.
Sea  $\{W(x)\}_{x \in X}$   una partici\'{o}n medible.
Sea $\mu$ una medida de probabilidad  en $(X, {\mathcal A})$.
  Entonces:

  {\bf  (i) } Existe una medida de probabilidad $\widehat \mu $ en el espacio medible cociente $X /\sim$, tal que para $\widehat \mu$-c.t.p. $W(x)   \in X/\sim$,   existe   una medida  de probabilidad $\mu^{W(x)}$ en el espacio medible  $(W(x), {\mathcal A}_{W(x)})$, tal  que:

  $\bullet $ Para todo $A \in {\mathcal A}$:
 \begin{equation}\label{eqnRohlin1}\int_X  \chi_A \, d \mu = \int _{[x] \in X/\sim} \Big(\int_ {y \in [x]= W(x)} \chi_A (y) \, d \mu^{W(x)} \Big) \, d \widehat \mu. \end{equation}

 $\bullet $ M\'{a}s en general, para toda $\psi \in L^1(\mu)$:
 \begin{equation}\label{eqnRohlin2}\int_X  \psi \, d \mu = \int _{[x] \in X/\sim} \, d \widehat \mu \ \Big(\int_ {y \in [x]= W(x)} \psi (y) \, d \mu^{W(x)} \Big).  \end{equation}

 {\bf (ii) } La  medida de probabilidad $\widehat \mu$, y  para $\widehat \mu$-casi todas las piezas $W(x) \in X /\sim$,  las  medidas  de probabilidad $\mu^{W^u(x)}$  que verifican las igualdades \em (\ref{eqnRohlin1}) \em y \em (\ref{eqnRohlin2}), \em  son \'{u}nicas.

\end{theorem}

La   prueba original del Teorema de Rohlin se encuentra en \cite{Rohlin}, o tambi\'{e}n en \cite{Rohlin1}. La demostraci\'{o}n se encuentra reformulada tambi\'{e}n, por ejemplo, en \cite{Ledrappier} o en \cite{VianaRohlin}.

El teorema de Rohlin     generaliza a un contexto medible, y para cualquier medida de probabilidad, el teorema de Fubini que vale para la medida de Lebesgue en un rect\'{a}ngulo de $\mathbb{R}^a \times \mathbb{R}^b$. En efecto, por el Teorema de Fubini en un rect\'{a}ngulo, la integral con respecto a la medida de Lebesgue  se obtiene como la integral doble en una secci\'{o}n    \lq\lq horizontal\rq\rq \ del rect\'{a}ngulo, de las integrales sobre las secciones \lq\lq verticales,\rq\rq \  con respecto a las medidas de Lebesgue a lo largo de dichas secciones,   respectivamente.

 \begin{definition}
 {\bf Medidas condicionales } \em   \index{medida! condicional} \label{definitionMedidasCondicionas}

 En las hip\'{o}tesis y conclusiones del Teorema \ref{theoremRohlin} de Rohlin, las medidas $\mu^{W }$ se llaman \em medidas condicionales de $\mu$ a lo largo de las piezas $W $ de la  partici\'{o}n ${\mathcal P}$. \em

 Tambi\'{e}n las llamamos, en breve   \em medidas ${\mathcal P}$-condicionales  de $\mu$, \em o si la partici\'{o}n est\'{a} clara en el contexto, simplemente \em medidas condicionales \em de $\mu$.

 Observar que cada  medida  condicionada   $\mu^{W }$ est\'{a}  soportada  en una pieza $W $, es una probabilidad  (es decir $\mu^{W}(W ) = 1$). Adem\'{a}s, $\mu^{W}$ est\'{a}  definida  para $\widehat \mu$-casi todas las piezas $W$ de la partici\'{o}n ${\mathcal P}$, y no necesariamente para todas las piezas.

\vspace{.3cm}

 {\bf Continuidad absoluta de las medidas condicionales.} \index{continuidad absoluta! de medidas condicionales} \index{medida! condicional} \index{medida! absolutamente continua}

 Cuando el espacio m\'{e}trico $X $ tiene una estructura de variedad compacta y riemanniana $M$ (es decir $X= M$), y si las piezas $W $ de la partici\'{o}n medible ${\mathcal P}$ son subvariedades de $M$, se considera, como medida  privi\-legiada  a lo largo de cada una de estas subvariedades, la medida  de Lebesgue en $W $, que denotamos $m^{W}$.

 En un contexto general, las medidas  condicionales $\mu^{W}$ de $\mu$ a lo largo de las piezas $W $ de la partici\'{o}n ${\mathcal P}$, podr\'{\i}an no tener   relaci\'{o}n con las medidas de Lebesgue $m^{W}$ a lo largo de estas piezas. Sin embargo, dentro de la Teor\'{\i}a de Pesin que veremos en la pr\'{o}xima secci\'{o}n, tienen especial importancia las medidas $\mu$ que cumplen la siguiente definici\'{o}n:
\end{definition}
\begin{definition}
\label{definitionMedidasCondicionadasAC}
\em
 Se dice que $\mu$ \em tiene medidas ${\mathcal P}$-condicionales absolutamente continuas \em cuando para $\widehat \mu$- casi toda pieza $W $ de la partici\'{o}n ${\mathcal P}$, se cumple:
 $$\mu^{W } \ll m^ {W }, \ \ \mbox{ i.e. } \mu^{W }(A) = 0 \mbox{ si } m^{W } (A) = 0, $$
 donde $A \subset W$ es medible, $W$ es una subvariedad $u$-dimensional de $M$, y $m^W$ es la medida de Lebesgue $u$-dimensional a lo largo de la subvariedad $W$. \index{continuidad absoluta! de medidas condicionales}
\index{medida! condicional inestable}
\index{medida! absolutamente continua}
 \end{definition}


\subsection{Medidas de Gibbs}

 \begin{definition}
 \label{definitionMedidasCondicionasInestables} \
\index{medida! condicional inestable}

 {\bf Medidas condicionales inestables } \em

 Sea $M$  una variedad compacta y riemanniana y sea $f: \in {\mbox{Diff }}^1(M)$ tal que el siguiente  conjunto
 \begin{equation} \label{eqnvariedadinestableglobal} W^u(x) := \{y \in M: \lim_{n \rightarrow + \infty} \mbox{dist} (f^{-n}(x), f^{-n}(y)) = 0\}\end{equation} es, por hip\'{o}tesis, una subvariedad  $C^1$ inmersa  en $M$   para $\mu$-casi todo punto $x \in M$, para cierta medida de probabilidad $\mu$  invariante con $f$. Bajo esta hip\'{o}tesis $W^u(x)$ se llama \em subvariedad  inestable global \em por el punto $x$. Debido a la continuidad de $f$, es inmediato chequear que $f^{-1}(W^u(x)) = W^u(f^{-1}(x))$ para todo $x \in M$.

 Sea $B_{\delta}$ una bola (compacta) de radio suficientemente peque\~{n}o con $\mu(B_{\delta}) >0$, y sea en ella la partici\'{o}n ${\mathcal F}_{\delta}^u = \{W^u_{\delta}(x)\}_{x \in M}$ definida como sigue:

  $\bullet$ $W^u_{\delta}(x)= c.c._x(W^u(x) \cap B_{\delta}(x)) $ es la componente conexa que contiene al punto $x$ de la intersecci\'{o}n $W^u(x) \cap B_{\delta}(x)$, para todo punto $x$ tal que $W^u(x)$ es una subvariedad $C^1$-inmersa en $M$ (es decir para $\mu$-c.t.p $x \in M$).  Por hip\'{o}tesis, esta subvariedad $W^u_{\delta}(x)$ est\'{a} $C^1$-encajada en $M$ para $\mu$-c.t.p. $x \in B_{\delta}$. Se llama \em subvariedad inestable local \em por el punto $x$.

  $\bullet$ Abusando de la notaci\'{o}n, para los restantes puntos  $x \in B_{\delta}$, denotamos $W^u_{\delta}(x) := \{y \in B_{\delta}\colon $ no existe subvariedad inestable por el punto $y \}$. Por construcci\'{o}n, esta \'{u}nica pieza de la partici\'{o}n ${\mathcal F}^u$ (que no es necesariamente una variedad) tiene $\mu$-medida cero.

   Rescalando $\mu$ para que $\mu(B_{\delta}) = 1$, tenemos lo siguiente: \index{medida! condicional inestable}
   \index{medida! condicional!}

     \vspace{.5cm}

   Si la partici\'{o}n ${\mathcal F}^u = \{W_{\delta}^u(p)\}_{p \in B_{\delta}}$ de subvariedades inestables locales as\'{\i} cons\-truida, fuera una partici\'{o}n $\mu$-medible,  llamamos \em medidas condicionales  inestables \em  de $\mu$, a las medidas condicionales $\mu^{W^u_{\delta}(p)}$  de $\mu$ a lo largo de las piezas $W^u_{\delta}(p)$ de esa partici\'{o}n, es decir a lo largo de las variedades inestables locales.

   Por simplicidad en la notaci\'{o}n escribiremos $\mu^u = \mu^{W^u_{\delta}}$ a  las medidas condicionales inestables, y $m^u = m^{W^u_{\delta}} $ a   las medidas de Lebesgue a lo largo de las subvariedades inestables $W^u_{\delta}$.
   \end{definition}

   \begin{definition}
   \index{continuidad absoluta! de medidas condicionales}
\index{medida! condicional inestable}
\index{medida! absolutamente continua}
   \label{definitionmedidascondicionalesinestablesAC} {\bf Medidas condicionales inestables absolutamente continuas.} \em
\index{continuidad absoluta! de medidas condicionales} \index{medida! condicional} \index{medida! absolutamente continua}
   En el contexto de la Definici\'{o}n \ref{definitionMedidasCondicionasInestables}, una medida $\mu$ se dice que tiene \em medidas condicionales inestables absolutamente continuas \em cuando
   $$\mu^u \ll m^u $$ para $\widehat \mu$- casi toda variedad inestable local $W^u_{\delta}$ de la partici\'{o}n ${\mathcal F}^u$, donde $\widehat \mu$ y $ \mu^u $ son las medidas del Teorema \ref{theoremRohlin} de desintegraci\'{o}n de Rohlin en la partici\'{o}n de variedades inestables locales.

    Se recuerda que, por definici\'{o}n, dadas dos medidas $\nu_1$ y $\nu_2$, se dice que \em $\nu_1$ es absolutamente con respecto de $\nu_2$, \em y se denota $\nu_1 \ll \nu_2$ cuando para todo boreliano $A$ se cumple $$\nu_2(A)= 0 \ \Rightarrow \ \nu_1(A) = 0.$$
   Dos medidas finitas $\nu_1$ y $\nu_2$ cumplen $\nu_1 \ll \nu_2$ si y solo si existe una funci\'{o}n $h \in L^1(\nu_2)$, llamada \em derivada de Radon-Nikodym \em de $\nu_1$ con respecto de $\nu_2$, tal que para todo conjunto medible $A$ se cumple
   $$\nu_1(A) = \int _A \ h \, d \nu_2.$$
   (Ver por ejemplo \cite[page 85]{Folland} o
   \cite[page 113]{Rudin}.)

   Luego, si $\mu$ tiene medidas condicionales inestables absolutamente continuas, entonces para $\mu$-c.t.p. $x \in M$ existe una funci\'{o}n $h_x \in L^1(m^u_x)$, donde $m^u_x$ es la medida de Lebesgue a lo largo de la variedad inestable local $W^u_{\delta}(x)$, tal que la medida condicional inestable $\mu^u_x$ a lo largo de $W^u_{\delta}(x) $ cumple:
   $$\mu^u_x\big(A \cap W^u_{\delta} (x)\big) = \int_{y \in W^u_{\delta}(x)} \chi_A(y) \, h_x(y) \, d m_x^u(y) $$
   para todo boreliano $A \subset M$.

 \end{definition}
 %

 \begin{definition}
 {\bf Medidas de Gibbs} \label{definitionMedidaGibbs} \em \index{continuidad absoluta! de medidas condicionales}
\index{medida! condicional inestable}
\index{medida! absolutamente continua}
\index{medida! de Gibbs}

 Sea $M$ una variedad compacta y riemanniana y sea $f \in \mbox{Diff }^1(M)$. Una medida de probabilidad $\mu$   $f$-invariante  se dice que es \em   medida de Gibbs \em cuando cumple:

 {\bf (i) } Para $\mu$-c.t.p. $x \in M$ el    conjunto definido por la igualdad (\ref{eqnvariedadinestableglobal}) es una subvariedad  $C^1$-inmersa  en $M$ (subvariedad  inestable   global  por el punto $x$).

 {\bf (ii) } Para $\delta  >0$ suficientemente peque\~{n}o, si $B_{\delta} \subset M$ es una bola compacta de radio   radio $\delta$ tal que $\mu(B_{\delta})>0$, entonces la siguiente familia de variedades inestables locales, es una partici\'{o}n $\mu$-medible $${\mathcal F}^u := \{W^u_{\delta}(x)\}_{x \in B_{\delta}}, \ \ \mbox{ donde } \ W^u_{\delta}(x) := c.c._x (W^u (x) \cap B_{\delta}).$$
 ($c.c._x$ denota la componente conexa que contiene al punto $x$).

 {\bf (iii) } $\mu$  tiene medidas condicionales inestables absolutamente continuas, de acuerdo  a la definici\'{o}n dada en el \'{u}ltimo p\'{a}rrafo de \ref{definitionMedidasCondicionasInestables}
 \end{definition}

Veamos ahora que la  propiedad de tener medidas condicionales absolutamente continuas a lo largo de las piezas de una partici\'{o}n, se transmite de una medida invariante $\mu$ a sus componentes erg\'{o}dicas, y rec\'{\i}procamente. Por lo tanto una medida invariante $\mu$ es de Gibbs, si y solo si sus componentes erg\'{o}dicas son medidas de Gibbs.

\begin{corollary}
{\bf del Teorema de Rohlin} \label{corolarioRohlin1} \index{teorema! Rohlin} \index{continuidad absoluta! de medidas condicionales}
\index{medida! condicional inestable}
\index{medida! absolutamente continua} \index{componentes erg\'{o}dicas}
\index{medida! de Gibbs}

Sea $f\colon X \mapsto X$  continua en el espacio m\'{e}trico compacto, sea $\mu$ una medida de probabilidad $f$-invariante y sea ${\mathcal P} = \{W(x)\}_{x \in X}$ una partici\'{o}n  medible y $f$-invariante \em (i.e. $f^{-1}(W(x)) =  W(f^{-1}(x))$ para todo $x \in X$). \em

Entonces:

{\bf  (a) } Las medidas  condicionales $\mu^{W}$ de $\mu$ a lo largo de las piezas $W$ de ${\mathcal P}$, y la medida inducida $\widehat \mu$ por $\mu$ en el espacio cociente $X / \sim$ de la partici\'{o}n, son medidas invariantes con $f$.

{\bf  (b) }   $\mu_x^{W(x)} = \mu^{W(x)} \ \mbox{ para } \ \mu-\mbox{c.t.p. } x \in M,$
donde $W(x)$ denota la pieza de la partici\'{o}n ${\mathcal P}$ que contiene al punto $x$, y $\mu_x$ denota la componente erg\'{o}dica de la medida $\mu$ a la que pertenece el punto $x$ seg\'{u}n el Teorema \em \ref{theoremDescoErgodicaEspaciosMetricos} (es decir,  $\lim_{n \rightarrow + \infty} (1/n) \sum_{j= 0}^{n-1} \delta_{f^j(x)} = \mu_x$ en la topolog\'{\i}a d\'{e}bil estrella).

{\bf (c)} Si adem\'{a}s   $X = M$ es una variedad compacta y riemanniana,  \em    $f \in \mbox{Diff }^1(M)$,  \em y    ${\mathcal P} = \{W_{\delta} (x)\}_{x \in M}$ es una partici\'{o}n  medible de la bola $B_{\delta}(x)$ formada por subvariedades $W_{\delta}(x)$ \ $C^1$ inmersas en $M $ para $\mu$-c.t.p. $x \in X$, entonces:

 \ \ $\bullet$ Las medidas  condicionales $\mu^{W}$ de $\mu$ son absolutamente continuas a lo largo de las subvariedades $W \in {\mathcal P}$  si y solo si  para $\mu$-c.t.p. $x \in X$ las medidas condicionales $\mu_x^W$ de las componentes erg\'{o}dicas $\mu_x$ de $\mu$,   son absolutamente continuas. \index{componentes erg\'{o}dicas}

 \ \ $\bullet$ $\mu$ es una medida de Gibbs    si y solo si las componentes erg\'{o}dicas $\mu_x$ de $\mu$ son medidas de Gibbs   para $\mu$-c.t.p. $x \in M$.
\end{corollary}

 {\em Demostraci\'{o}n: }
  {\bf (a) } Por el Teorema \ref{theoremRohlin} de Rohlin y la $f$-invariancia de $\mu$ tenemos, para todo conjunto medible $A$:
  \begin{equation} \label{eqn-10}\mu(A) =  \int _{[x] \in X/ \sim} \, d \widehat \mu\int _{x \in [x] } \chi_A(x) \, d\mu^{W(x)}  = $$ $$ =\int_{[x] \in X/ \sim} \Big( \mu^{W(x)}(W(x) \cap A \Big) \, d \widehat \mu = $$
$$= \mu(f^{-1}(A)) = \int _{[x] \in X/ \sim} \, d \widehat \mu\int _{x \in [x] } \chi_{f^{-1}(A)} (x) \, d\mu^{W(x)}  =$$ $$ = \int_{[x] \in X/ \sim} \Big( \mu^{W(x)}(W(x) \cap f^{-1}(A)) \Big) \, d \widehat \mu.\end{equation}
Las medidas de probabilidad $\mu^{W(x)}$ est\'{a}n definidas para $\widehat \mu$-casi toda pieza de la partici\'{o}n ${\mathcal P}$. Para las otras piezas tomamos por convenci\'{o}n cualquier medida de probabilidad soportada en ellas. De esta forma tenemos definidas $\mu^{W(x)}$ para todo $x \in M$.

Debido a la $f$-invariancia de las piezas $W \in {\mathcal P}$, se cumple:
$$W(x) \cap f^{-1}(A) = f^{-1}(W(f(x)) \cap  A)$$
Definimos la siguiente medida $(\mu^{W(y)})^*$ a lo largo de la pieza $W(y)$ para   todo $y \in M$:
$$(\mu^{W(y)})^*(B):= \mu^{W(f^{-1}(y))}(f^{-1}(B)) \ \ \forall \ B \in {\mathcal A}, $$ donde
${\mathcal A}$ denota la sigma-\'{a}lgebra de Borel en $X$. Entonces:
$$\mu^{W(x)}(W(x) \cap f^{-1}(A)) = (\mu^{W(f(x))})^*(W(f(x)) \cap A),$$
y sustituyendo en (\ref{eqn-10}) obtenemos
\begin{equation} \label{eqn-10b}\mu(A) =    \int_{[x] \in X/ \sim} \Big( (\mu^{W(f(x))})^*(W(f(x)) \cap A)  \Big) \, d \widehat \mu =$$ $$ =\int _{[y] \in X /\sim} \Big( (\mu^{W(y)})^*(W(y) \cap A)  \Big) \, d (\widehat \mu)^*,\end{equation}
donde $(\widehat \mu)^*$ es la medida de probabilidad en $X / \sim$ definida por
$$\int  \varphi  \, d \widehat \mu ^* := \int \varphi \circ f \, d \widehat \mu  \ \ \forall \ \ \varphi \in L^1(\widehat \mu).$$ (En la \'{u}ltima igualdad $f: X \sim \mapsto X \sim $ denota la aplicaci\'{o}n que lleva la pieza $W \in {\mathcal P}$ en  la pieza $f^{-1}(W) $).
De (\ref{eqn-10b}) deducimos
$$\mu(A) =  \int _{[y] \in X /\sim} \Big( (\mu^{W(y)})^*(W(y) \cap A)  \Big) \, d (\widehat \mu)^* = $$ $$ =   \int _{[y] \in X /\sim} d \widehat \mu ^* \int_{y \in [y]} \chi_A(y) \, d (\mu^{W(y)})^* \ \ \forall \ A \in {\mathcal A}.$$
 Luego hemos encontrado otra desintegraci\'{o}n de Rohlin de la medida $\mu$ con respecto a la partici\'{o}n ${\mathcal P}$. Por la unicidad de las medidas de probabilidad $\widehat \mu$ y $\mu^{W}$ de la desintegraci\'{o}n de Rohlin, se cumple
 $$(\widehat\mu)^* = \widehat \mu, \ \ \ \mu^{W(x)} = (\mu^{W(f(x)})^*,$$
 demostrando la $f$-invariancia de $\widehat \mu$ y de $\mu^W$.

 {\bf (b) } Para todo conjunto $A \in {\mathcal A}$, por el Teorema de Descomposici\'{o}n Erg\'{o}dica (Teorema \ref{theoremDescoErgodicaEspaciosMetricos}) tenemos: $$\mu(A) = \int \mu_x(A) \, d \mu.   $$
 Por el Teorema de Descomposici\'{o}n de Rohlin, aplicado a cada medida erg\'{o}dica $\mu_x$, se cumple:
 \begin{equation}\label{eqn-11}\mu(A) = \int \mu_x (A) \, d \mu = \int_{x \in X} d \mu \int_{[y] \in X / \sim} d \widehat \mu_x \int_{y \in [y]} \chi_A \, d \mu_x^{W(y)} =$$
 $$=\int_{ [y]  \in   X/ \sim }    d   \widehat \nu \int_{y \in [y]} \chi_A \, d \mu_x^{W(y)},\end{equation}
 donde la medida de probabilidad $ \widehat  \nu $ en el espacio cociente $   X / \sim, {\mathcal A}  $ est\'{a} definida por:
 $$\int_{X/ \sim} \varphi \, d\widehat \nu  := \int _{x \in X} \Big (\int _{[y] \in X / \sim} \varphi([y]) \, d \widehat \mu_x \Big ) \, \ d \mu \ \ \ \forall \ \ \varphi \in L^1 (\widehat \mu).$$
 Luego, la igualdad (\ref{eqn-11}) es una descomposici\'{o}n de Rohlin de la medida $\mu$ con respecto a la partici\'{o}n ${\mathcal P}$. Por la unicidad de las probabilidades de la descomposici\'{o}n de Rohlin, concluimos que $\widehat \mu = \widehat \nu$ y      $\mu^W_x = \mu^W$.

 {\bf (c) } Como consecuencia de la parte (b) $\mu^W \ll m^W$ si y solo si $\mu_x ^W \ll m ^W$.
 En el caso particular en que  la partici\'{o}n medible ${\mathcal P}$ es la partici\'{o}n de una bola $B \subset M$ en variedades inestables locales, obtenemos que $\mu$ es medida de Gibbs si y solo si $\mu_x$ lo es para $\mu$-c.t.p. $x \in M$.
 \hfill $\Box$

\vspace{.3cm}

 \subsection{Relaci\'{o}n entre medidas de Gibbs y SRB} \index{Pesin! teor\'{\i}a de} \index{teor\'{\i}a de Pesin}

En esta secci\'{o}n asumimos que $M$ es una variedad compacta y riemanniana y que $f \in \mbox{Diff }^{1 + \alpha}(M)$. Algunos de los teoremas que veremos en esta secci\'{o}n se obtienen de resultados de la Teor\'{\i}a de Pesin,    bajo la hip\'{o}tesis de que $f$ es de clase  $C^{1  + \alpha}$. Pero son falsos  si $f$ es solo de clase $C^1$. En particular las condiciones de continuidad absoluta de las medidas condicionales inestables (existencia de medida erg\'{o}dica de Gibbs que es tambi\'{e}n SRB) no rige  para todo $f \in \mbox{Diff }^1(M)$.

M\'{a}s adelante, en las   secciones posteriores de este cap\'{\i}tulo, veremos algunas formas de reformular los resultados de esta secci\'{o}n, generalizando algunas definiciones (mediante la introducci\'{o}n de las medidas SRB-like) y modificando adecuadamente los enunciados,  para que sean aplicables a todo difeomorfismo de clase $C^1$ en la variedad $M$ (y m\'{a}s a\'{u}n, algunos de ellos son aplicables tambi\'{e}n a transformaciones continuas $f: M \mapsto M$).

 \begin{theorem}
 \label{teoremaAtractoresErgodicos} \label{theoremGibbs->SRB} {\bf     de Gibbs implica SRB o f\'{\i}sica} \index{medida! SRB} \index{medida! de Gibbs} \index{medida! f\'{\i}sica} \index{teorema! de medida de Gibbs} \index{componentes erg\'{o}dicas}

Sea \em $f \in \mbox{Diff }^{1 + \alpha}(M)$ \em y sea $\mu $ una medida invariante hiperb\'{o}lica.

  Si $\mu$ es una medida de Gibbs, entonces las componentes erg\'{o}dicas de $\mu$ son medidas SRB o f\'{\i}sicas.


 \end{theorem}

 En el p\'{a}rrafo \ref{proofTheoremGibbs->SRB} veremos   la demostraci\'{o}n de este Teorema, reduci\'{e}ndola a     resultados de la Teor\'{\i}a de Pesin. Tambi\'{e}n se puede encontrar la demostraci\'{o}n del teorema \ref{teoremaAtractoresErgodicos} en \cite{PughShubErgodicAttractors} o en \cite[Proposition 11.24]{BonattiDiazVianaLibro}.

 Es cierto tambi\'{e}n el rec\'{\i}proco del Teorema \ref{theoremGibbs->SRB}, bajo hip\'{o}tesis  de $C^{1 + \alpha}$-hiperbolicidad uniforme: si las componentes erg\'{o}dicas de una medida invariante $\mu$ son medidas SRB o f\'{\i}sicas, entonces $\mu$ es una medida de Gibbs. M\'{a}s precisamente:

 \begin{theorem} \label{theoremSRB->Gibbs}  {\bf SRB erg\'{o}dica implica de Gibbs} \index{medida! SRB erg\'{o}dica} \index{medida! SRB} \index{medida! de Gibbs}   \index{teorema! de medida SRB erg\'{o}dica}

 Sea \em $f \in \mbox{Diff }^{1 + \alpha}(M)$, \em sea $\Lambda \subset M$ un atractor erg\'{o}dico   uniformemente hiperb\'{o}lico, y sea $\mu $ la medida  erg\'{o}dica SRB o f\'{\i}sica   soportada en $\Lambda$ \em (de acuerdo a la Definici\'{o}n \ref{definitionAtractorErgodico} de atractor erg\'{o}dico).

 \em Entonces   $\mu$ es medida  de Gibbs.
 \end{theorem}

 El Teorema \ref{theoremSRB->Gibbs} es consecuencia   del siguiente Teorema de Pesin-Sinai, que demuestra   la existencia de medida SRB erg\'{o}dica, bajo hip\'{o}tesis de hiperbolicidad uniforme en el contexto $C^{1 + \alpha}$

 \begin{theorem} \label{theoremPesinSinai}
 {\bf [Pesin-Sinai, \cite{Pesin-Sinai}]}  \index{teorema! Pesin-Sinai}

 Sea $f \in \mbox{Diff }^{1 + \alpha}(M)$. Sea $\Lambda $ un atractor topol\'{o}gico uniformemente hi\-per\-b\'{o}\-lico. Entonces:

  Existen medidas SRB-erg\'{o}dicas soportadas en $\Lambda$. Toda  medida SRB erg\'{o}dica soportada en $\Lambda$ es de Gibbs. Rec\'{\i}procamente, toda medida erg\'{o}dica de Gibbs soportada en $\Lambda $ es SRB. \em (cf. Teorema \ref{theoremGibbs->SRB}).
 \end{theorem}

 En el caso particular que $f$ sea Anosov, demostraremos este resultado  m\'{a}s adelante en esta secci\'{o}n (Teorema \ref{TheoremSRBanosov}).

\vspace{.3cm}

 {\bf Generalizaci\'{o}n del Teorema \ref{theoremPesinSinai} para difeomorfismos $C^{1 + \alpha}$ no uniformemente hiperb\'{o}licos}

 La   implicaci\'{o}n SRB-erg\'{o}dica $\Rightarrow$ Gibbs, como en el Teorema \ref{theoremPesinSinai} de Pesin-Sinai, rige a\'{u}n en hip\'{o}tesis m\'{a}s gene\-rales que la hiperbolicidad uniforme que nosotros enunciamos en ese Teorema,   asumiendo siempre que $f \in \mbox{Diff }^{1 + \alpha}(M)$. En efecto, el atractor topol\'{o}gico $\Lambda$ puede no ser   hiperb\'{o}lico, sino que alcanza que sea   \em parcialmente hiperb\'{o}lico  \em (ver Definici\'{o}n, por ejemplo  en \cite[Definition B.3, page 289]{BonattiDiazVianaLibro}). Puede ser no uniformemente hiperb\'{o}lico con singularidades, como el atractor de Lorenz \cite{Pesin1992}, por ejemplo.

  Para los difeomorfismos parcialmente hiperb\'{o}licos,  las medidas \index{hiperbolicidad! parcial} \index{difeomorfismos! parcialmente hiperb\'{o}licos} invarian\-tes no son necesariamente hiperb\'{o}licas, sino que puede existir un subespacio de Oseledets con exponente de Lyapunov igual a cero. Sin embargo, puede existir una separaci\'{o}n acotada uniformemente lejos de cero, entre los exponentes de Lyapunov positivos y los no positivos.

 La demostraci\'{o}n   de Pesin-Sinai del Teorema \ref{theoremPesinSinai} se encuentra en \cite{Pesin-Sinai} para difeomorfismos $C^{1 + \alpha}$-parcialmente hiperb\'{o}licos; en particular para los que son uniformemente hiperb\'{o}licos \'{o} Anosov. Tambi\'{e}n puede encontrarse la demostraci\'{o}n general para difeomorfismos $C^{1 + \alpha}$    par\-cial\-mente hi\-per\-b\'{o}\-li\-cos,  en \cite[Theorem 11.16] {BonattiDiazVianaLibro}.

Observemos que la hiperbolicidad parcial no implica que la medidas invariantes sean hiperb\'{o}licas, y por lo tanto, aunque existan medidas de Gibbs erg\'{o}dicas, no podremos aplicar el Teorema \ref{theoremGibbs->SRB} para deducir que estas son SRB.  En general, el problema de existencia de medidas SRB para difeomorfismos parcialmente hiperb\'{o}licos, est\'{a} esencialmente abierto. Recientemente, en  \cite{VianaYang} se demuestra que los
  difeomorfismos parcialmente hiperb\'{o}licos con dimensi\'{o}n central 1 (1 es la dimensi\'{o}n del subespacio de Oseledets con exponente de Lyapunov igual a cero), $C^{1 + \alpha}$-gen\'{e}ricamente existe   una, y a lo sumo  una cantidad finita,  de medidas erg\'{o}dicas   SRB cuyas cuencas de atracci\'{o}n estad\'{\i}stica cubren Lebesgue c.t.p.

Los difeomorfismos hiperb\'{o}licos son un caso particular de los llamados difeomorfismos con \em splitting dominado. \em  \index{splitting! dominado} La definici\'{o}n, los enunciados y las demostraciones de propiedades din\'{a}micas topol\'{o}gicas de los difeomorfismos con splitting dominado se encuentran en \cite{PujalsSamba}. El problema de existencia de medidas SRB para los difeomorfismos con splitting dominado, salvo en casos particulares, est\'{a} esencialmente abierto.

 Otro caso en el que la existencia de medidas SRB se ha estudiado, es el de    \index{difeomorfismos! derivados de Anosov}
 los llamados difeomorfismos  \em derivados de Anosov. \em La demostraci\'{o}n de la misma tesis del Teorema \ref{theoremPesinSinai} para difeomorfismos derivados de Anosov transitivos con exponente de Lyapunov positivo y de clase $C^{1 + \alpha}$, y adem\'{a}s  la unicidad de su medida SRB, fueron dadas por M.F. Carvalho     en \cite{Carvalho} (o tambi\'{e}n, m\'{a}s detalladamente explicadas, en   \cite{CarvalhoTesis}).

 Tambi\'{e}n existen otras generalizaciones del Teorema \ref{theoremPesinSinai} para ciertas clases de difeomorfismos no uniformemente   ni parcialmente hiperb\'{o}licos ni derivados de Anosov  y que no tienen splitting dominado (por ejemplo en \cite{Hu}, \cite{HeberTesis} y \cite{CatEnr2001}). Estas clases de difeomorfismos tienen $C^{r}$ regularidad para valores de $r \geq 2$ suficientemente grande y se componen de ciertos difeomorfismos $f_1$ llamados \em casi-Anosov, \em \index{difeomorfismos! casi-Anosov} que est\'{a}n en el borde de los de Anosov en el espacio ${\mbox{Diff }^{r}(M)}$ para cierto $r > 1$ suficientemente grande. Por definici\'{o}n, un difeomorfismo $f_1$ es casi Anosov, o \em almost Anosov, \em    si $f_1$ se obtiene por medio de una isotop\'{\i}a (es decir, por medio de una deformaci\'{o}n continua $f_t \in \mbox{Diff }^r(M)$ para $t \in [0,1] \subset \mathbb{R}$),    tal que $f_t$ es de Anosov   transitivo  para todo $0 \leq t < 1$. La isotop\'{\i}a en un casi-Anosov debe cumplir la siguiente condici\'{o}n:   el splitting $T_{x}M = S^t_{x } \oplus U^t_{x}$ es un splitting uniformemente hiperb\'{o}lico de $f_t$ para todo $t \in [0, 1)$ fijo, y para todo $x \in M$. Adem\'{a}s, por hip\'{o}tesis, existe una   \'{o}rbita $o(x_0)$ tal que para todo $x \not \in \overline{o (x_0)} $ existe un splitting $S_x^1 \oplus U_x^1$ invariante por $f_1$, obtenido como el l\'{\i}mite cuando $t \rightarrow 1$ del splitting hiperb\'{o}lico de $S_x^{t} \oplus U_x^t$. Finalmente, el difeomorfismo $f_1$, o bien posee tambi\'{e}n un splitting (no hiperb\'{o}lico) definido como $S^{1} \oplus U^{1}$, donde $S^{1}= \lim_{t \rightarrow 1^-} {S^1}, \ U^{1} = \lim_{t \rightarrow 1^-} U^t $ (cuando existen dichos l\'{\i}mites en la \'{o}rbita de $x_0$ y son mutualmente transversales);  o bien existen dichos l\'{\i}mites pero son tangentes;  o bien no existen.

 La existencia de una \'{u}nica medida de Gibbs erg\'{o}dica que es SRB, en ejemplos del primer caso de difeomorfismos $C^{2}$ casi-Anosov, es demostrada en  \cite{Carvalho}, donde $x_0$ es un punto fijo  por $f_t$ para todo $0 \leq t \leq 1$, y los exponentes de Lyapunov en $U_{x_0}^1$ (para $t= 1$) son estrictamente positivos.  En cambio, en el ejemplo de casi-Anosov tambi\'{e}n del primer caso, estudiado en \cite{HuYoung},  se prueba que no existe ninguna medida de probabilidad de Gibbs, si se toma como foliaci\'{o}n inestable, la formada por las $C^1$-subvariedades que son tangentes a $E^1_x$ en todo punto $x \in M$.  En este ejemplo, el difeomorfismo $f_1$, de clase $C^2$ casi-Anosov,   se construye tomando $x_0$ en un punto fijo, pero tal que los exponentes de Lyapunov en $U_{x_0}^1$ son nulos. Sin embargo, a\'{u}n en este ejemplo,  en \cite{HuYoung} se demuestra que existe una \'{u}nica medida SRB erg\'{o}dica y que su cuenca de atracci\'{o}n estad\'{\i}stica cubre Lebesgue casi todo punto.

 En el segundo caso de difeomorfismos casi-Anosov, en el toro $\mathbb{T}^2$,    \cite{HeberTesis}    prueba la existencia de la medida de Gibbs erg\'{o}dica que es SRB,    cuando $o(x_0)$ es una \'{o}rbita no peri\'{o}dica (m\'{a}s precisamente una \'{o}rbita de tangencia heterocl\'{\i}nica). Tambi\'{e}n en el segundo caso  (cuando   $S_{x_0}^1$ y $U_{x_0}^1$ son tangentes) y en $\mathbb{T}^2$, la exis\-tencia y unicidad de una medida SRB erg\'{o}dica con cuenca de atracci\'{o}n que cubre Lebesgue-c.t.p.,  es   demostrada en \cite{CatEnr2001}, cuando $x_0$ es un punto fijo (no hiperb\'{o}lico) para $f_1 $ de clase $C^3$, y adem\'{a}s una de las derivadas segundas parciales de $f_1$ en $x_0 $ se anula. Este caso incluye los ejemplos   de difeomorfismos casi-Anosov en el toro bidimensional introducidos por Lewowicz en \cite{Lewowicz-Ejemplo}, en los que las direcciones estable  e inestable en un punto fijo $x_0$ para $t < 1$, se hacen tangentes cuando el par\'{a}metro $t$ de la isotop\'{\i}a alcanza el valor 1 (sin saber a priori si los exponentes de Lyapunov en casi todo punto diferente de $x_0$, son o no son, diferentes de cero).

 El tercer caso de difeomorfismos casi-Anosov, en que los l\'{\i}mites $S_{x_0}^1$ y $U_{x_0}^1$ no existen, es estudiado en \cite{Hu} cuando   $x_0$ es un punto fijo, $f$ es de clase $C^r$ para $r \geq 2$ suficientemente grande, y se cumplen ciertas hip\'{o}tesis sobre las derivadas de orden mayor que 1  en $x_0$.

 Estos resultados parciales, fueron obtenidos en subclases muy parti\-culares de difeomorfismos casi-Anosov. La dificultad grande para demos\-trar   la existencia de medidas SRB, cuando no se tienen a priori   hip\'{o}tesis de existencia de una constante uniforme  que separe  los exponentes de Lyapunov positivos de los no positivos,   provoca  que el problema general de existencia de medidas SRB para los casi-Anosov, permanezca a\'{u}n mayor\-mente abierto.

 \vspace{.3cm}

 Ahora re-enunciamos el Teorema \ref{theoremPesinSinai}   en el caso particular de

 $f \in \mbox{Diff }^{1+ \alpha}(M)$  difeomorfismo de Anosov. M\'{a}s precisamente:


 \begin{theorem}.
{\bf (Sinai)} \label{TheoremSRBanosov} \index{medida! SRB erg\'{o}dica} \index{medida! SRB} \index{medida! de Gibbs}   \index{teorema! de medida SRB erg\'{o}dica}
 \index{teorema! Pesin-Sinai} \index{teorema! de existencia de! medidas SRB} \index{teorema! de existencia de! medidas de Gibbs}

 Sea $f \in \mbox{Diff }^{1 + \alpha}(M)$ un difeomorfismo de Anosov.  Entonces:

 {\bf (a) } Existen   medidas    erg\'{o}dicas SRB (o f\'{\i}sicas).

 {\bf (b) } Toda medida SRB (o f\'{\i}sica)  es de Gibbs   erg\'{o}dica, y rec\'{\i}procamente

  {\bf (c) } La uni\'{o}n de las cuencas de atracci\'{o}n estad\'{\i}stica de las medidas SRB   cubren Lebesgue c.t.p. $x \in M$.

  {\bf (d) } Si adem\'{a}s $f$ es topol\'{o}gicamente transitivo, entonces existe una \'{u}nica medida SRB, es erg\'{o}dica y de Gibbs, y su cuenca de atracci\'{o}n estad\'{\i}stica cubre Lebesgue c.t.p. $x \in M$.

 \end{theorem}

Demostraremos   el Teorema \ref{TheoremSRBanosov} m\'{a}s adelante en esta secci\'{o}n, en el p\'{a}rrafo \ref{proofTheoremSRBanosov}.

\begin{corollary}
\label{corolarioAnosovTransitivomedidaLebesgue}
\label{corollarySRBanosov}
Sea $f \in {\mbox{Diff }^{1 + \alpha}(M)}$ de Anosov transitivo que preserva la medida de Lebesgue $m$. Entonces, $m$ es erg\'{o}dica,   es la \'{u}nica medida SRB de $f$ y es de Gibbs.
\end{corollary}
{\em Demostraci\'{o}n: }
Por la parte (d) del Teorema \ref{TheoremSRBanosov}, Lebesgue c.t.p. $x \in M$ cumple
$$\lim_{n \rightarrow + \infty} \frac{1}{n} \sum_{j= 0}^{n-1} \delta_{f^j(x)} = \mu,$$
donde $\mu$ es la \'{u}nica medida SRB de $f$, es erg\'{o}dica y de Gibbs.

Por el Teorema \ref{theoremDescoErgodicaEspaciosMetricos}, como la medida de Lebesgue $m$ es invariante, entonces para $m$-c.t.p. $x \in M$
$$\lim_{n \rightarrow + \infty} \frac{1}{n} \sum_{j= 0}^{n-1} \delta_{f^j(x)} = m_x,$$
donde $m_x$ es la componente erg\'{o}dica de $m$ a la que pertenece $f$.

Cuando existe, el l\'{\i}mite de una sucesi\'{o}n de medidas de probabilidad en la topolog\'{\i}a d\'{e}bil estrella, es \'{u}nico (esta es una propiedad de todo espacio m\'{e}trico). Concluimos que $m_x = \mu$ para $m$-c.t.p. $x \in M$. Dicho de otra forma, la descomposici\'{o}n erg\'{o}dica de $m$ tiene una \'{u}nica componente erg\'{o}dica que es $\mu$. Entonces $m = \mu$, como quer\'{\i}amos demostrar.
\hfill $\Box$

\begin{remark} \em
   \label{remarkSRBanosovDensidad} \index{medida! densidad de} \index{medida! SRB} \index{medida! de Gibbs}
 En la demostraci\'{o}n del Teorema \ref{TheoremSRBanosov}, probaremos adem\'{a}s los siguientes resultados:

  \em Las medidas condicionales   inestables $\mu_x^u$ de cualquier medida $\mu$ que sea SRB  \em (para un difeomorfismos de Anosov de clase $C^{1+ \alpha}$) \em no solo son   absolutamente continuas respecto a las  medidas de Lebesgue inestables $m_x^u$ \em (a lo largo de las respectivas variedades inestables locales $W^u_{\delta}(x)$ para $\mu$-c.t.p. $x \in M$), \em sino que son adem\'{a}s equivalentes a estas; es decir: \em
$$\mu_x^u \ll m^u_x , \ \ \mu^u_x \ll m^u_x \ \ \ \mu-\mbox{c.t.p.} \ x \in M.$$ \em \index{densidad} \index{derivada de Radon-Nykodim}
Adem\'{a}s,   la derivada de Radon-Nikodym  $d \mu^u_x/d m^u_x$ (es decir la densidad de las medidas condicionales inestables) es una funci\'{o}n continua $h_x(y)$ y positiva, y est\'{a} definida por:
\em $$\frac{d\mu_x^u}{dm_x^u}(y) = h_x(y) := \prod_{j= 0}^{+ \infty} \frac{J^u(f^j(x))}{J^u(f^j(y))} \in C_0(M, \mathbb{R}^+) \ \ \ \mu-\mbox{c.t.p.} \ x \in M,$$
donde $J^u(x) := \big|\mbox{det} (df_x|_{E^u_x}\big |$ se llama   \em Jacobiano inestable de $f$ en el punto $x$. \index{jacobiano inestable}
\end{remark}

Cuando intentamos generalizar el Corolario \ref{corollarySRBanosov}
a difeomorfismos no uniformemente hiperb\'{o}licos que preservan la medida de Lebesgue, el m\'{e}todo de demostraci\'{o}n que usamos para los Anosov transitivos no funciona si uno no sabe a priori, o demuestra primero, la existencia de una medida SRB:

\begin{conjecture} \em
{\bf (Viana \cite{VianaSurvey})} \index{conjetura! Viana}
\em
Sea $f$ un difeomorfismo que preserva la medida de Lebesgue $m$. Si $m$ es una medida hiperb\'{o}lica \em (cf. Definici\'{o}n \ref{definitionMedidaHiperbolica}) \em entonces existe alguna medida SRB.
\end{conjecture}

La demostraci\'{o}n de la conjetura de Viana, o el hallazgo de un contraejemplo, es un problema abierto.

\vspace{.3cm}

\subsection{Sobre la F\'{o}rmula de Pesin para la entrop\'{\i}a}

Como consecuencia de los teoremas \ref{theoremGibbs->SRB} y \ref{theoremSRB->Gibbs}, la b\'{u}squeda de medidas SRB o f\'{\i}sicas en el caso de difeomorfismos   de clase $C^1$ m\'{a}s H\"{o}lder, se centra en la b\'{u}squeda de medidas de Gibbs, o sea, de medidas invariantes cuyas medidas condicionales inestables sean absolutamente continuas (ver Definiciones \ref{definitionMedidasCondicionasInestables} y \ref{definitionMedidaGibbs}). Por ese motivo, la caracterizaci\'{o}n de las medidas invariantes que son medidas de Gibbs, adquiere en el contexto $C^{1+ \alpha}$ especial relevancia. Una tal caracterizaci\'{o}n est\'{a} dada por la  igualdad del siguiente Teorema \ref{theoremFormulaPesin}, llamada \em F\'{o}rmula de Pesin para la entrop\'{\i}a. \em

\vspace{.2cm}

{\bf Entrop\'{\i}a m\'{e}trica} Para poder enunciar la F\'{o}rmula de Pesin,  \index{entrop\'{\i}a m\'{e}trica} introducimos brevemente el concepto de \em entrop\'{\i}a m\'{e}trica $h_{\mu}(f)$ \em de una transformaci\'{o}n continua $f: M \mapsto M$ con respecto a una medida de probabilidad $\mu$  invariante con $f$. La definici\'{o}n de entrop\'{\i}a m\'{e}trica, y el estudio de las primeras propiedades que dan su forma de c\'{a}lculo para difeomorfismos de Anosov, son debidos a Kolmogorov    \cite{Kolmogorov} y a Sinai \cite{Sinai-MetricEntropy}. La definici\'{o}n precisa y sus propiedades,  puede encontrarse adem\'{a}s, por ejemplo en \cite[\S 4.4]{Walters},   \cite[Cap\'{\i}tulo 4]{Mane},
\cite[pag. 168-170]{Katok-Hasselblatt},   \cite{Sinai-MetricEntropy2}, \cite[pag. 55-76]{SinaiBook1994}, \cite[\S 4.1-4.2]{Jost},  \cite[Chapter 3]{Keller}, o \cite{KingEnciclopedia}, entre muchos otros textos que tratan matem\'{a}ticamente el concepto de la entrop\'{\i}a m\'{e}trica  para los sistemas din\'{a}micos.

Para comprender los enunciados siguientes, adm\'{\i}tase que tenemos defi\-ni\-do un n\'{u}mero real no negativo $h_{\mu}(f) \geq 0$, llamado \em entrop\'{\i}a m\'{e}trica de $f$ \em con respecto de la medida de probabilidad $f$-invariante $\mu$,   que depende solo de $f$  y de $\mu$ y   tal que   $h_{\mu}(f)$ es invariante por isomorfismos de espacios de medida. El n\'{u}mero $h_{\mu}(f)$, por la forma en que se define, mide con respecto a la probabilidad $\mu$, la tasa exponencial maximal en que \em los iterados futuros de $f$ desordenan  los pedazos de cualquier partici\'{o}n finita ${\mathcal P}$ del espacio $M$, \em ponderados con la probabilidad $\mu$.

M\'{a}s precisamente,  consideremos la partici\'{o}n ${\mathcal P}_n := \bigvee_{j= 0}^{n-1} f^{-j}({\mathcal P})$, definida por la siguiente condici\'{o}n: $x, y \in A \in {\mathcal P}_n$  si y solo si para todo $0 \leq j \leq n-1$ existe $P_j \in {\mathcal P}$ tal que $x,y \in P_j$. Cada pedazo  diferente  de esta partici\'{o}n ${\mathcal P}_n$, es el conjunto de puntos $x$ que mutuamente se acompa\~{n}an dentro de un mismo pedazo de ${\mathcal P}$, al ser iterados hacia el futuro.  Cuando  mayor es la cantidad de pedazos diferentes de la partici\'{o}n ${\mathcal P}_n$ al crecer $n$, m\'{a}s desordena el iterado $f^n$ a las \'{o}rbitas en el espacio, en relaci\'{o}n a la partici\'{o}n inicial ${\mathcal P}$.  \index{caos}

Se define
$$h_{\mu}(f, {\mathcal P}) := \limsup_{n \rightarrow + \infty} \frac{1}{n} \sum_{A \in {\mathcal P}_n} \mu(A) \log \mu(A).$$
   $\sum_{A \in {\mathcal P}_n} \mu(A) \log \mu(A)$ es un promedio ponderado de logaritmos. Luego, se  interpreta como un exponente (ponderado). Entonces, al dividirlo entre $n$,   da una \em tasa o coeficiente  de  crecimiento \em exponencial con $n$. Por este motivo, ese cociente  se interpreta como la tasa o velocidad de crecimiento exponencial con $n$ (hasta el instante $n$), del \lq\lq desorden espacial\rq\rq \ que producen los iterados de $f$, ponderado con la medida de probabilidad $\mu$. Por lo tanto, intuitivamente hablando, $h_{\mu}(f, \mathcal P)$ es la tasa exponencial, $\mu$-ponderada y asint\'{o}tica,  en que los iterados de $f$ (mejor dicho $\mu$-casi toda \'{o}rbita  de $f)$ desordenan a la partici\'{o}n inicial   ${\mathcal P}$ con la que mir\'{a}bamos, como referencia, la distribuci\'{o}n de puntos   iniciales.

Finalmente, se define la entrop\'{\i}a m\'{e}trica
$$h_{\mu}(f) : = \sup_{{\mathcal P} \in \mathbb{P}} h_{\mu}(f, {\mathcal P}),$$
donde ${\mathbb{P}}$ denota el conjunto de todas las particiones finitas de $M$ en piezas medibles. \index{entrop\'{\i}a m\'{e}trica}

Luego $h_{\mu}(f) \geq 0$ puede interpretarse como la tasa exponencial asint\'{o}tica maximal de crecimiento del desorden espacial de $f^n$ al hacer $n \rightarrow + \infty$, ponderada con la probabilidad $\mu$. Si   $h_{\mu}(f) > 0$,   el sistema  restringido al soporte de $\mu$  se llama \em ca\'{o}tico \em  en sentido medible.
Cuanto mayor es $h_{\mu}(f)$,  m\'{a}s r\'{a}pidamente se desordena el espacio al iterar $f$.

Luego,   al observador que quiera cuantificar el caos, interesan, si existen, aquellas medidas de probabilidad invariantes $\mu$ que maximicen la entrop\'{\i}a m\'{e}trica $h_{\mu}(f)$ en relaci\'{o}n    al valor esperado de cierta funci\'{o}n real, llamada \em potencial  \em (el cual, gruesamente hablando,   cuantifica el significado    de \lq\lq optimizar la observaci\'{o}n del caos\rq\rq). Estas medidas que maximizan la diferencia entre la entrop\'{\i}a m\'{e}trica y el valor esperado de una funci\'{o}n potencial, se llaman \em estados de equilibrio. \em  \index{estados de equilibrio} Su estudio se constituye en la sub-teor\'{\i}a, dentro de la teor\'{\i}a erg\'{o}dica de los sistemas din\'{a}micos, llamada \em formalismo termodin\'{a}mico   \em (ver por ejemplo \cite[Chapter 4]{Keller}).


\begin{theorem} \label{theoremFormulaPesin}
{\bf [Pesin] F\'{o}rmula de Pesin para la entrop\'{\i}a} \index{entrop\'{\i}a! f\'{o}rmula de Pesin} \index{f\'{o}rmula de Pesin} \index{Pesin! f\'{o}rmula de}

Sea $f \in \mbox{Diff }^{1 + \alpha}(M)$ que preserva una medida  $\mu$ hiperb\'{o}lica y de Gibbs. Entonces \em \index{medida! de Gibbs} \index{medida! hiperb\'{o}lica}
\begin{equation} \label{eqnformuladePesin} h_{\mu}(f) = \int \ \sum_{i=1}^{\mbox{\ \footnotesize dim}(M)} \chi_i^+ (x) \, d \mu,\end{equation}
\em donde $h_{\mu}(f)$ es la entrop\'{\i}a m\'{e}trica de $f$ con respecto de $\mu$, y $$ \chi_i^+(x)  := \max\{0, \chi_i(x)\},$$ siendo  \em $\{\chi_i\}_{1 \leq i \leq \mbox{\footnotesize dim}(M)}$ \em los exponentes de Lyapunov de $f$ en el punto $x$, repetido cada uno de ellos tantas veces como la dimensi\'{o}n del espacio de Oseledets $E^i_x$ correspondiente.  \index{exponentes de Lyapunov! positivos} \em

Se recuerda que, por el Teorema de Oseledets,   $\mu$-c.t.p. es regular. Por lo tanto existe la funci\'{o}n real no negativa $\sum_{i=1}^{\mbox{\footnotesize dim}(M)} \chi_i^+ (x)$ y es medible. Si ninguno de los exponentes de Lyapunov $\chi_i(x)$ positivo, resulta  $\sum_{i= 1}^{\mbox{\footnotesize dim}(M)} \chi^+_i(x) = 0$.
\end{theorem}

La demostraci\'{o}n  de la F\'{o}rmula de Pesin se encuentra en \cite{PesinLyapunovExponents}. Una prueba diferente que prescinde de la  propiedad  de continuidad absoluta de la foliaci\'{o}n estable se encuentra en \cite{ManePesinFormula}. La demostraci\'{o}n tambi\'{e}n se encuentra, por ejemplo, en \cite[Theorem 5.4.5]{BarreiraPesin}. Finalmente, en \cite{Tahzibi} se gene\-ra\-liza el Teorema \ref{theoremFormulaPesin} para   difeomorfismos $C^1$-gen\'{e}ricos.

 Si una medida $\mu$ satisface la f\'{o}rmula de Pesin, entonces ella mide \'{o}ptimamente el  caos medible del sistema  en relaci\'{o}n al valor esperado de la suma de los exponentes de Lyapunov positivos. En efecto, la igualdad (\ref{eqnformuladePesin}) da el m\'{a}ximo posible de la entrop\'{\i}a m\'{e}trica $h_{\mu}(f)$ en relaci\'{o}n a ese valor esperado, ya que para toda aplicaci\'{o}n de clase $C^1$ rige la siguiente cota superior de $h_{\mu}(f)$: \index{caos}


 \begin{theorem}
 {\bf      Desigualdad de Margulis-Ruelle \cite{Margulis-Inequality}, \cite{Ruelle_Inequality}} \index{desigualdad Margulis-Ruelle} \index{entrop\'{\i}a! desigualdad Margulis-Ruelle}
 \index{teorema! Margulis-Ruelle}
 \index{exponentes de Lyapunov! positivos}

 Sea $f : M \mapsto M$ de clase $C^1$. Sea $\mu$ una medida de probabilidad inva\-rian\-te por $f$. Entonces
   \em
\begin{equation} \label{eqnDesigualdadDeRuelle} h_{\mu}(f) \leq \int \ \sum_{i=1}^{\mbox{\ \footnotesize dim}(M)} \chi_i^+ (x) \, d \mu,\end{equation}
\em donde $h_{\mu}(f)$ es la entrop\'{\i}a m\'{e}trica de $f$ con respecto de $\mu$,  y $$ \chi_i^+(x)  := \max\{0, \chi_i(x)\},$$ donde  \em $\{\chi_i\}_{1 \leq i \leq \mbox{\footnotesize dim}(M)}$ \em son los exponentes de Lyapunov de $f$ en el punto $x$, repetido cada uno de ellos tantas veces como la dimensi\'{o}n del espacio de Oseledets $E^i_x$ correspondiente.  \em
 \end{theorem}

 La demostraci\'{o}n de la Desigualdad (\ref{eqnDesigualdadDeRuelle}) de Margulis-Ruelle se puede encontrar en \cite{Ruelle_Inequality}, o tambi\'{e}n en \cite[Theorem 5.4.1]{BarreiraPesin}.

 \vspace{.3cm}

La f\'{o}rmula de Pesin, caracteriza a las medidas de Gibbs, debido al siguien\-te rec\'{\i}proco del Teorema \ref{theoremFormulaPesin}:

\begin{theorem}
\label{theoremLedrappier-Young} \index{teorema! Ledrappier-Young}
\index{entrop\'{\i}a} \index{medida! de Gibbs} \index{exponentes de Lyapunov! positivos} \index{continuidad absoluta! de medidas condicionales}
\index{medida! condicional inestable}

{\bf [Ledrappier-Young, \cite{Ledrappier-Young}]}

Sea $f \in \mbox{Diff }^2(M)$ y sea $\mu$ una medida invariante tal que \em $$\sum_{i=1}^{\mbox{\footnotesize dim}(M)} \chi_i^+(x) >0 \ \ \mu-\mbox{c.t.p. } x \in M.$$   \em Si se verifica la F\'{o}rmula de Pesin \em (\ref{eqnformuladePesin}) \em de la entrop\'{\i}a m\'{e}trica $h_{\mu}(f)$, entonces la medida $\mu$ es de Gibbs, es decir, $\mu$ tiene medidas condicionales inestables absolutamente continuas.
\end{theorem}

La demostraci\'{o}n del Teorema \ref{theoremLedrappier-Young} se encuentra en \cite{Ledrappier-Young}, y tambi\'{e}n en \cite{Ledrappier} con la hip\'{o}tesis adicional de $\mu$ hiperb\'{o}lica.

\begin{corollary} \index{atractor! topol\'{o}gico} \index{atractor! hiperb\'{o}lico} \index{hiperbolicidad! uniforme} \index{medida! SRB}
\index{medida! de Gibbs} \index{Pesin! f\'{o}rmula de} \index{f\'{o}rmula de Pesin} \index{entrop\'{\i}a! f\'{o}rmula de Pesin}  \index{teorema! Ledrappier-Young}
Sea $\Lambda$ un atractor topol\'{o}gico uniformemente hiperb\'{o}lico de $f \in {\mbox{Diff }^2(M)}$. Entonces las siguientes afirmaciones son equivalentes:

{\bf (i)} $\mu$ es una medida erg\'{o}dica SRB (o f\'{\i}sica).

{\bf (ii)} $\mu$ es una medida erg\'{o}dica de Gibbs.

{\bf (iii)} $\mu$ es una medida erg\'{o}dica que satisface la F\'{o}rmula \em (\ref{eqnformuladePesin}) \em de Pesin para la entrop\'{\i}a.
\end{corollary}
{\em Demostraci\'{o}n: }
Se obtiene inmediatamente de reunir los Teoremas \ref{theoremPesinSinai}, \ref{theoremFormulaPesin} y \ref{theoremLedrappier-Young}.
\hfill $\Box$

\subsection{Demostraci\'{o}n del Teorema  \ref{theoremGibbs->SRB}}

Para poder demostrar los Teoremas \ref{theoremGibbs->SRB} y \ref{TheoremSRBanosov}, necesitamos introducir algunos resultados   relevantes de la llamada \em Teor\'{\i}a de Pesin: \em Esta teor\'{\i}a estudia el comportamiento de  $\mu$-casi toda \'{o}rbita, para $f \in {\mbox{Diff }^{1 + \alpha}(M)}$ donde $\mu$ es una medida de probabilidad $f$-invariante e \em hiperb\'{o}lica. \em Es decir, la regi\'{o}n de Pesin  $\Sigma$ (i.e. el conjunto de los puntos regulares cuyos exponentes de Lyapunov son todos diferentes de cero) tiene $\mu$-probabilidad igual a 1.

\begin{definition}
{\bf Holonom\'{\i}a} \em  \label{definitionHolonomia} \index{holonom\'{\i}a} \index{foliaci\'{o}n! invariante} \index{foliaci\'{o}n! holonom\'{\i}a de}

Sea $\mu$   una medida de probabilidad invariante hiperb\'{o}lica de

$f \in \mbox{Diff }^{1 + \alpha}(M)$. Es decir, la regi\'{o}n de Pesin $\Sigma$ cumple $\mu(\Sigma) = 1$. Tomemos   $x_0 \in \Sigma$ tal que para todo $\delta >0$ la bola $B_{\delta}(x_0) \subset M$   de centro $x_0$ y radio $\delta >0$ cumple $\mu(B_{\delta}(x_0))>0$. Por el Teorema \ref{theoremVarInvariantesRegionPesin},   existen para   $\mu$-c.t.p. $x \in B_{\delta}(x_0)$   las variedades estables e inestables locales por el punto $x$, que denotamos $W^{s}_{\delta}(x), \ W^u_{\delta}(x) \subset B_{\delta}(x_0)$ (Si fuera necesario tomamos la componente conexa que contiene al punto $x$ de la intersecci\'{o}n de la variedad estable o inestable local por el punto $x$  con la bola $B_{\delta}(x_0)$).

 Definimos la holonom\'{\i}a $$h_s: B^*_{\delta}(x_0)   \mapsto W^u_{\delta}(x_0),$$ donde $$B^*_{\delta}(x_0) = \{ y \in B_{\delta}(x_0) \cap \Sigma\colon \ \#(W^s_{\delta}(y) \cap W^u_{\delta}(x_0)) = 1$$ como la transformaci\'{o}n que a cada punto $y \in B^*{\delta}(x_0)$ hace corresponder el \'{u}nico punto $h_s(y) \in   W^s_{\delta}(y) \cap W^u_{\delta}(x_0)$.
 Debido al Teorema \ref{theoremVarInvariantesRegionPesin}, si $\delta>0$ es suficientemente peque\~{n}o, entonces $h_s$ existe, pues el conjunto $B^*_{\delta}(x_0)$ contiene por lo menos a la variedad estable local $W^s_{\delta}(x_0)$ que interseca transversalmente a $W^u_{\delta}(x_0)$ en el punto $x_0$.

 La transformaci\'{o}n $h_s$ se llama \em holonom\'{\i}a a lo largo de las variedades estables locales \em en la bola $B_{\delta}(x_0)$ sobre la variedad inestable local del punto $x_0$, o en breve, \em holonom\'{\i}a local estable.  \em

Decimos que \em la   holonom\'{\i}a estable \index{holonom\'{\i}a} \index{continuidad absoluta! de holonom\'{\i}a} \index{foliaci\'{o}n! absolutamente continua}   \index{foliaci\'{o}n! invariante} \index{foliaci\'{o}n! holonom\'{\i}a de} \index{continuidad absoluta! de foliaci\'{o}n} $h_s$ es absolutamente continua \em si  para $\mu$-c.t.p. $x_0 \in \Sigma$ y para todo boreliano $A \subset B_{\delta}(x_0)$ se cumple:
\begin{equation} \label{eqnholonomiaAC} m(h_s^{-1}(A) = 0 \ \Leftrightarrow \ m^u(A \cap W_{\delta}^u(x_0)) = 0, \end{equation}
donde $m$ es la medida de Lebesgue en la variedad $M$ y $m^u$  es la medida de Lebesgue a lo largo de la subvariedad inestable local $W_{\delta}^u(x_0)$.

An\'{a}logamente se define \em holonom\'{\i}a inestable \em y \em continuidad absoluta de la holonom\'{\i}a inestable, \em intercambiando entre s\'{\i}   los roles de las variedades estables e inestables locales en la definici\'{o}n anterior.

\end{definition}


\begin{theorem} \label{theoremTeoriaPesin}
{\bf [Pesin] \cite{Pesin76}}

{\bf (Teorema fundamental de la Teor\'{\i}a de Pesin)} \index{teor\'{\i}a de Pesin} \index{Pesin! teor\'{\i}a de} \index{holonom\'{\i}a}
\index{foliaci\'{o}n! invariante}
\index{foliaci\'{o}n! holonom\'{\i}a de}
\index{continuidad absoluta! de holonom\'{\i}a}
\index{foliaci\'{o}n! absolutamente continua}
\index{continuidad absoluta! de foliaci\'{o}n}

Sea \em $f \in \mbox{Diff }^{1 + \alpha}(M)$ \em y sea $\mu$ una medida erg\'{o}dica e hiperb\'{o}lica. Entonces existe $\delta >0$ suficientemente peque\~{n}o, tal que la holonom\'{\i}a estable y la holonom\'{\i}a inestable en las bolas $B_{\delta}(x_0)$ para $\mu$-c.t.p. $x_0 \in M$ son absolutamente continuas.
\end{theorem}

La prueba del Teorema \ref{theoremTeoriaPesin} se encuentra    en \cite{Pesin76}, y tambi\'{e}n, por ejemplo, en [Theorem 4.3.1]\cite{BarreiraPesin}.

Observamos que el Teorema \ref{theoremTeoriaPesin} es falso en la topolog\'{\i}a $C^1$. En efecto, en \cite{Bowen_C1horseshoe} y \cite{RobinsonYoungContrajemploCAdeFoliacion} (y tambi\'{e}n en \cite{SchmittGora} para endomorfismos)  se construyen   ejemplos de atractores hiperb\'{o}licos para los cuales la holonom\'{\i}a a lo largo de las variedades estables locales no es absolutamente continua.

Por ese motivo la demostraci\'{o}n   del Teorema \ref{theoremGibbs->SRB} no funciona para difeomorfismos de clase $C^1$ que no sean $C^{1 + \alpha}$. En particular, la demostraci\'{o}n que daremos del Teorema \ref{TheoremSRBanosov} no funciona para difeomorfismos de clase $C^1$ que no sean de clase $C^{1 + \alpha}$. M\'{a}s adelante expondremos algunos resultados utilizando las   medidas \lq\lq SRB-like\rq\rq (ver Definici\'{o}n \ref{definitionMedidaSRBlike}), que no son necesariamente   medidas de Gibbs,  pero que existen para difeomorfismos y endomorfismos sin m\'{a}s regularidad que  $C^1$ manteniendo   propiedades de atracci\'{o}n estad\'{\i}stica similares a las medidas SRB.

\vspace{.3cm}

Ahora, veremos que la demostraci\'{o}n del Teorema \ref{theoremGibbs->SRB},   que relaciona  para difeomorfismos de clase $C^{1 + \alpha}$  las medidas de Gibbs hiperb\'{o}licas y erg\'{o}dicas con las medidas SRB o f\'{\i}sicas, se reduce al Teorema fundamental de la Teor\'{\i}a de Pesin que establece la continuidad absoluta de la holonom\'{\i}a estable local.


 \begin{nada} \em  \label{proofTheoremGibbs->SRB}
 {\bf Demostraci\'{o}n del Teorema \ref{teoremaAtractoresErgodicos}}
 \end{nada}

{\em Demostraci\'{o}n: }
Por la parte (c) del Corolario \ref{corolarioRohlin1} si $\mu$ es medida de Gibbs, entonces $\mu$-casi toda componente erg\'{o}dica $\mu_x$  de $\mu$ es medida de Gibbs. Adem\'{a}s como $\mu$ es hiperb\'{o}lica por hip\'{o}tesis, entonces $\mu(\Sigma)= 1$, donde $\Sigma$  es la regi\'{o}n de Pesin. Luego, por el Teorema \ref{theoremDescoErgodicaEspaciosMetricos} de descomposici\'{o}n erg\'{o}dica $\mu$-casi toda componente erg\'{o}dica   $\mu_x$ de $\mu$ cumple $\mu_x(\Sigma) = 1$; es decir $\mu_x$ es hiperb\'{o}lica y de Gibbs, adem\'{a}s de ser erg\'{o}dica. Tomemos una tal $\mu_x$ y renombremosla como $\mu$.

Para probar el Teorema \ref{theoremGibbs->SRB}, aplicando la Definici\'{o}n \ref{definitionMedidaSRB} de medida SRB o f\'{\i}sica, debemos probar que la cuenca $B$ de atracci\'{o}n estad\'{\i}stica de $\mu$, definida por:
$$B= \{x \in M: \ \lim_{n \rightarrow + \infty} \frac{1}{n} \sum_{j= 0}^{n-1} \delta_{f^j(x)} = \mu\}$$
tiene medida de Lebesgue $m(B)$ positiva.

Como $\mu$ es erg\'{o}dica, $\mu(B) = 1$.   Por el Teorema \ref{theoremRohlin} de Descomposici\'{o}n de Rohlin en una bola $B_{\delta}(x_0)$ con $\mu$-medida positiva, se cumple:
$$\mu^u(B \cap W^u_{\delta}(x_0)) = 1$$
para $ \mu$-c.t.p. $x_0$, donde $\mu^u$ es la medida condicionada inestable de $\mu$.
Como $\mu$ es medida de Gibbs, por la Definici\'{o}n \ref{definitionMedidaGibbs}, se cumple $\mu^u \ll m^u$, para $\mu$-c.t.p. $x_0$, donde $m^u$ es la medida de Lebesgue a lo largo de la variedad $W^u_{\delta}(x_0)$. Luego deducimos
$$m^u(B \cap W^u_{\delta}(x_0))>0$$
Aplicando el Teorema \ref{theoremTeoriaPesin} de la Teor\'{\i}a de Pesin, obtenemos:
$$m(h_s^{-1}(B \cap W^u_{\delta}(x_0)) >0.$$
Entonces, para terminar de demostrar que $\mu$ es una medida SRB, es decir para terminar de probar que $m(B) >0$, basta demostrar ahora que $h_s^{-1}(B \cap W^u_{\delta}(x_0)) \subset B$. En efecto, sea $y \in h_s^{-1}(B \cap W^u_{\delta}(x_0))$. Probemos que $y \in B$. Por definici\'{o}n de la holonom\'{\i}a estable $h_s$, se cumple $ y \in B^*_{\delta}$ tal que $h_s (y) = z := W^s_{\delta}(y) \cap W^u_{\delta}(x_0) \in B$. Entonces, como $z \in B \cap W^s_{\delta}(y)$, obtenemos $$\lim_{n \rightarrow + \infty} \mbox{dist}(f^n(z), f^n(y)) = 0, \ \ \ \lim_{n \rightarrow + \infty} \frac{1}{n} \sum_{j= 0}^{n-1} \delta_{f^j(z)} = \mu.$$
Aplicando el resultado probado en el Ejercicio \ref{exercise4}, resulta $$\lim_{n \rightarrow + \infty} \frac{1}{n} \sum_{j= 0}^{n-1} \delta_{f_j(y)} = \mu,$$
es decir $y \in B$, como quer\'{\i}amos probar.
\hfill $\Box$

\subsection{Lema de Distorsi\'{o}n Acotada}

Al final de este  cap\'{\i}tulo  probaremos el Teorema \ref{TheoremSRBanosov} que establece, en el contexto $C^{1 + \alpha}$-Anosov, la equivalencia entre las medidas SRB erg\'{o}dicas y las medidas de Gibbs erg\'{o}dicas, y la existencia de estas.

La demostraci\'{o}n del Teorema \ref{TheoremSRBanosov} est\'{a} fuertemente basada en el si\-guien\-te Lema \ref{lemmaDistoAcotada} de Distorsi\'{o}n Acotada.     La prueba de este lema se basa en la hip\'{o}tesis $f \in \mbox{Diff }^{1 + \alpha}(M)$.

\begin{lemma} \index{lema de distorsi\'{o}n acotada} \index{distorsi\'{o}n acotada} \index{jacobiano inestable}
\label{lemmaDistoAcotada} {\bf Lema de Distorsi\'{o}n Acotada}
Sea $f \in {\mbox{Diff }^{1 + \alpha}}(M)$ de Anosov. Sea $TM = E^u  \oplus E^s$ el splitting uniformemente hiperb\'{o}lico del fibrado tangente, y sea $W^u_{\delta}(x)$ la variedad inestable local por el punto $x$ para $\delta >0$ constante, suficientemente peque\~{n}o. Denotamos con $$J^u(x) = \big|\det \, df_x|_{E^u_x}  \big| >0$$ al Jacobiano inestable en el punto $x$. Sean $x, y \in M $     dos puntos tales que $y \in W_{\delta}^u(x)$ y sean, para todo $n \geq 1$, las funciones
$$h_n(x, y) := \frac{\prod_{j= 1}^{n}  J^u(f^{-j}(x))}{\prod_{j= 1}^{n}  J^u(f^{-j}(y))} \in (0, + \infty)$$
$$h(x,y) = := \frac{\prod_{j= 1}^{+ \infty}  J^u(f^{-j}(x))}{\prod_{j= 1}^{+ \infty}  J^u(f^{-j}(y))} \in [0, + \infty],$$
definidas solo para las parejas $(x,y) $ de puntos en el conjunto:
$$ H_{\delta}:= \{(x,y) \in M \mbox{ tales que } y \in W^u_{\delta}(x)\}.$$
Entonces:

{\bf (i) } Existe una constante real $K >0$ tal que
  $$   \ \frac{1}{K} < h(x,y) <  {K} \ \ \forall \ (x,y) \in H_{\delta} $$

{\bf (ii) } La funci\'{o}n $h: H_{\delta} \rightarrow \mathbb{R}^+  $    es continua.

{\bf (iii) } Para todo $\epsilon >0$ existe $N \geq 1$ \em (independiente de la pareja de puntos $(x,y) \in H_{\delta}$) \em  tal  que:
$$
e^{- \epsilon} < \frac{h_n(x,y)}{h(x,y)} < e^{\epsilon} \ \ \forall \ (x,y) \in H_{\delta}, \ \ \forall \ n \geq N.$$

\end{lemma}

{\em Demostraci\'{o}n: }
{\bf (i) } Definamos $\log :[0, + \infty] \mapsto [- \infty, + \infty]$ acordando que $\log 0 = -\infty $ y $\log (+ \infty) = + \infty$. Para demostrar la afirmaci\'{o}n (i) basta probar que existe una constante real $c >0$ tal que $|\log h(x,y)| \leq c$. Por construcci\'{o}n de la funci\'{o}n $h: H_{\delta} \mapsto [0, + \infty] $ tenemos:
\begin{equation} \label{eqn-5}|\log h(x,y)|  = \big|\sum_{j= 1}^{+ \infty} \big(\log J^u(f^{-j}(x)) - \log J^u(f^{-j}(y)) \big)  \big| \leq $$ $$  \sum_{j= 1}^{+ \infty} \big|\log J^u(f^{-j}(x)) - \log J^u(f^{-j}(y)) \big|. \end{equation}
  $J^u: M \mapsto \mathbb{R}^+$ es continuo (porque $df_x$ es continuo pues $f$ es de clase $C^1$ y   $E^u_x$ depende continuamente de $x$). Entonces la funci\'{o}n $J^u$ est\'{a} uniformemente acotada superiormente e inferiormente por constantes reales positivas. Luego,  por el Teorema del valor medio del c\'{a}lculo diferencial aplicado a la funci\'{o}n real   $\log t$ de variable real positiva $t$, existe una constante $c_1 >0$ tal que $$|\log J^u (z_1) - \log J^u(z_2)| \leq c_1 |J^u(z_1) - J^u(z_2)|$$ para toda pareja $(z_1, z_2)$ de puntos en la variedad $M$.
Sustituyendo en (\ref{eqn-5}) resulta:
\begin{equation} \label{eqn-6}|\log h(x,y)|    \leq  c_1 \sum_{j= 1}^{+ \infty} \big|  J^u(f^{-j}(x)) -  J^u(f^{-j}(y)) \big|. \end{equation}
Siendo $f$ de clase $C^{1 + \alpha}$, las variedades inestables son de clase $C^{1 + \alpha}$ (cf. la \'{u}ltima parte del enunciado del Teorema \ref{teoremavariedadesinvariantesAnosov}). Entonces el subespacio tangente $T_y W^u_{\delta}(x) = E^u_y$ depende $\alpha$-H\"{o}lder continuamente del punto $y$ cuando $y$ var\'{\i}a a lo largo de la subvariedad $W^u_{\delta}(x)$. Es decir, para cada $x$ existe una constante $c_2(x) >0$ tal que
\begin{equation}
\label{eqn33}
\mbox{dist}(E^u_y , E^u_z) \leq c_2(x) \mbox{dist}(y,z)^{\alpha} \ \ \forall \ \ y,z \in W_{\delta}^u(x).\end{equation} Afirmamos que para $\delta >0$ fijo suficientemente peque\~{n}o, existe una constante $c_2$ uniforme tal que $$c_2(x) \leq c_2 \ \ \forall \ x \in M.$$ En efecto, fijando $x_1 \in M$ y una constante $C(x_1)> c_2(x_1)$, la desigualdad (\ref{eqn33}) se satisface para cualquier punto $\widehat x$ en un entorno suficientemente peque\~{n}o de $x_1$,   sustituyendo $c_2( x) $ por $C(x_1)$. Luego, cubriendo la variedad compacta $M$ con una cantidad finita de tales entornos, centrados en puntos $x_1, \ldots, x_m$ respectivamente, y tomando $c_2= \max_{i= 1}^{m} C(x_i)$, deducimos que se satisface la desigualdad (\ref{eqn33}) para todo $x \in M$, con la constante $c_2$  en lugar de $c_2(x)$.

Por otra parte, como $f$ es de clase $C^{1 + \alpha}$, su derivada $df_x$ depende $\alpha$-H\"{o}lder continuamente del punto $x$. Es decir, existe una constante $c_3 >0$ tal que
$$\|df_y - df_z\| \leq c_3 \mbox{dist}(y,z) ^{\alpha} \ \ \forall \ y,z \in M.$$ Por otra parte existe una constante $c_4>0$ tal que $$|\det A |\leq c_4\|A\|$$ para toda aplicaci\'{o}n lineal $A$ de un espacio vectorial de dimensi\'{o}n igual a $\mbox{dim}(E^u)$ en otro de la misma dimensi\'{o}n.
Reuniendo las cuatro \'{u}ltimas desigualdades que involucran las constantes $ c_2(x), c_2, c_3$ y $c_4$, obtenemos
$$ |J^u(y) - J^u(z) |  =\big| \, |\det(df_y|_{E^u_y} | - |\det(df_z|_{E^u_z} | \, \big| \leq $$ $$\big| \, |\det(df_y|_{E^u_y} | - |\det(df_z|_{E^u_y} | \, \big| + \big| \, |\det(df_z|_{E^u_y} | - |\det(df_z|_{E^u_z} | \leq $$ $$ \big| \leq c_4 \|df_y - df_z\| + c_4 (\max_{z \in M}  \|df_z\|) \mbox{dist} (E^u_y , E^u_z) \leq $$ $$\big(c_4 c_3 + c_4 (\max_{z \in M} \|df_z\|) c_2 \big) \, \cdot \, \mbox{dist}(y,z)^{\alpha}  \ \ \forall \ (y,z) \in W^u_{\delta}(x), \ \ \forall \ x \in M.$$
Sustituyendo en la desigualdad (\ref{eqn-6}), obtenemos una constante $c_5 >0$ tal que
\begin{equation}
\label{eqn-7}
|\log h(x,y)| \leq c_5 \, \sum_{j= 1}^{+ \infty} \mbox{dist}(f^{-j}(x), f^{-j}(y))^{\alpha} \ \ \ \forall \ (x,y) \in H_{\delta}
\end{equation}
Nota: Para obtener la desigualdad (\ref{eqn-7}) observamos que no necesariamente  se cumple $f^{-j}W^u_{\delta}(x) \subset W^u_{\delta}(f^j(x))$, y por lo tanto, aunque $(x,y) \in H_{\delta}$ no necesariamente son aplicables directamente  todas las desigualdades anteriores para la pareja de puntos $(f^{-j}(x), f^{-j}(y))$ para cualquier $j \geq 1$. Sin embargo, fijado $\delta  >0$ suficientemente peque\~{n}o, aplicamos la parte B) del Teorema \ref{teoremavariedadesinvariantesAnosov}, y observamos que hay   convergencia uniforme a cero de las distancias a lo largo de las variedades inestables locales al iterar hacia el pasado, pues los vectores tangentes $u \in E^u_x$ se contraen uniformemente hacia el pasado seg\'{u}n la Definici\'{o}n \ref{definicionAnosov} de difeomorfismo de Anosov. Luego,   existe $N \geq 1$ uniforme tal que $f^{-n}W_{\delta}^u(x) \subset W_{\delta}^u (f^{-n}(x))$ para todo $n \geq N$ para todo $x \in M$. Entonces, basta elegir $0 <\delta'< \delta $ tal que $f^{-j}W_{\delta'} ^u (x) \subset W_{\delta}^u (f^{-j}(x))$ para $ j \in \{0, \ldots, N-1\}$, para que esa misma inclusi\'{o}n valga tambi\'{e}n para todo $j \geq 0$. Finalmente renombramos $\delta'$ como $\delta$ para deducir la desigualdad (\ref{eqn-7}).

Por la Definici\'{o}n \ref{definicionAnosov} tenemos:
\begin{equation}
\label{eqnDistanciasInestablesHiperbolicidad}
\mbox{dist}^u(f^{-j}(x), f^{-j}(y)) \leq C \sigma^{-j} \mbox{dist}^u(x,y),\end{equation}
donde $C >0$ y $\sigma > 1$ son las constantes dadas en la Definici\'{o}n \ref{definicionAnosov}, y $\mbox{dist}^u(x,y)$ es la distancia entre los puntos $x$ e $y$ a lo largo de la variedad inestable local $W^u_{\delta}(x)$. Adem\'{a}s la distancia $\mbox{dist}$ en la variedad ambiente $M$ entre dos puntos que est\'{a}n en la misma variedad inestable local, es siempre el m\'{\i}nimo de las longitudes de curvas que unen esos puntos (con la m\'{e}trica riemanniana dada en $M$). Por lo tanto es menor o igual que la distancia $\mbox{dist}^u$ a lo largo de la variedad inestable local. Sustituyendo en la desigualdad (\ref{eqn-7}) resulta:
\begin{equation}
\label{eqn-8bb}
|\log h(x,y)| \leq c_5 \, C^{\alpha} \, \sum_{j= 1}^{+ \infty} \sigma^{-\alpha j} \, (\mbox{dist}^u(x, y))^{\alpha} \leq c_6 (\mbox{dist}^u(x,y))^{\alpha} $$ $$\ \ \ \forall \ (x,y) \in H_{\delta},
\end{equation}
donde $c_6 = c_5 \, C^{\alpha} \,  {1}/(1- \sigma^{-\alpha})$ (observar que $0 <\sigma^{-\alpha} < 1$ porque $\sigma > 1$ y $\alpha >0$).

Como $W^u_{\delta}(x)$ es una variedad $C_1$-encajada en $M$ que depende continuamente de $x$, existe $$\mbox{Diam}(W^u_{\delta}(x)) := \sup_{y \in W^u_{\delta}(x)} \mbox{dist}^u(x,y) \in \mathbb{R}^+$$ y es una funci\'{o}n continua real de $x \in M$. Luego, est\'{a} acotada superiormente por   una constante uniforme $c_7 >0$. Sustituyendo en (\ref{eqn-8bb}), concluimos:
\begin{equation}
\label{eqn-8}
|\log h(x,y)| \leq   c_6 \, (\mbox{dist}^u(x,y))^{\alpha}  \leq c_6 \, c_7^{\alpha} = c\ \ \ \forall \ (x,y) \in H_{\delta},
\end{equation}
terminando la demostraci\'{o}n de la parte (i) del Lema \ref{lemmaDistoAcotada}.

\vspace{.2cm}
{\bf (iii) } Fijemos $n \geq 1$. Por la definici\'{o}n de las funciones $h$  y $h_n$  tenemos
$$\frac{h_n(x,y)}{h(x,y)} = := \frac{\prod_{j= n +1}^{+ \infty}  J^u(f^{-j}(y))}{\prod_{j= n+1}^{+ \infty}  J^u(f^{-j}(x))}$$
Luego, aplicando la f\'{o}rmula (\ref{eqn-8bb}) a los puntos $f^{-n-1}(x)$ y $f^{-n-1}(y)$, en lugar de $x$ e $y$ respectivamente, obtenemos
$$\Big|\log \Big(\frac{h_n(x,y)}{h(x,y)}   \Big) \Big| \leq c_6 \big(\mbox{dist}^u(f^{-n-1}(x), f^{-n-1}(y))\Big)^{\alpha}.$$
Usando ahora la desigualdad (\ref{eqnDistanciasInestablesHiperbolicidad}), deducimos:
$$\Big|\log \Big(\frac{h_n(x,y)}{h(x,y)}   \Big) \Big| =   \leq c_6 C^{\alpha} \sigma^{-(n+1) \alpha} \big(\mbox{dist}^u(x , y)\big)^{\alpha} \leq c_6 \, C \, c_7^{\alpha} \, \sigma^{-(n+1) \alpha}.$$
Como $\sigma >1$ y $\alpha >0$, el t\'{e}rmino de la derecha tiende a cero cuando $n \rightarrow + \infty$. Luego, existe $N \geq 1$, independiente de $(x,y)$ tal que
$$\Big|\log  \, \frac{h_n(x,y)}{h(x,y)}     \Big| < \epsilon \ \ \forall \ n \geq N,$$
terminando de probar la afirmaci\'{o}n (iii) del Lema \ref{lemmaDistoAcotada}.

\vspace{.2cm}
{\bf (ii) }  Debido a la afirmaci\'{o}n (ii), la sucesi\'{o}n de funciones $h_n(x,y)$ converge uniformemente a $h(x,y)$ para todo $(x,y) \in H_{\delta}$. Adem\'{a}s $h_n(x,y)$ es continua en $H_{\delta}$, porque el Jacobiano  inestable  $J^u(x) = \big|\det df^n_x|_{E^u_x}\big|$ es una funci\'{o}n continua    de $x$, pues $df_x$ es continuo y $E_x^u$ tambi\'{e}n. El l\'{\i}mite uniforme de una sucesi\'{o}n de funciones continuas es continuo. Concluimos que $h$ es una funci\'{o}n continua,
finalizando la prueba del Lema \ref{lemmaDistoAcotada}.
\hfill $\Box$

\subsection{Demostraci\'{o}n del Teorema \ref{TheoremSRBanosov}}

 En la demostraci\'{o}n del Teorema \ref{TheoremSRBanosov}, usaremos, adem\'{a}s del Lema de Distorsi\'{o}n Acotada,  el siguiente resultado, v\'{a}lido para cualquier difeomorfismo de Anosov (no necesariamente de clase $C^{1 + \alpha}$), y generalizable tambi\'{e}n (reformulando adecuadamente el enunciado) para los entornos locales de atractores topol\'{o}gicos uniformemente hiperb\'{o}licos.

 \begin{theorem}
 \label{theoremProductoLocal} {\bf Estructura de producto local} \index{producto local} \index{teorema! de producto local} \index{estructura de producto local}
 Sea $f \in \mbox{Diff }^1(M)$ Anosov. Existen  constantes $ 0 < \delta' < \delta $   suficientemente peque\~{n}as, tales que   $f$ tiene estructura de producto local en entornos de radio $\delta'$.  \em

 La estructura de producto local significa, por definici\'{o}n, que para toda pareja de puntos $(x,y) $ tales que $\mbox{dist}(x,y) < \delta'$, existen, son \'{u}nicos y dependen continuamente de $(x,y)$ los puntos $[x,y]$ y $[y,x]$ definidos por:
 $$[x,y] = W^u_{ \delta}(x) \cap W^s_{ \delta}(y), \ \ [y,x] = W^u_{ \delta}(y) \cap W^s_{ \delta}(x),$$
 donde $W_{ \delta}^{u,s}(x)$ denota la componente conexa de $W^{u,s}(x) \cap B_{ \delta}(x)$ que contiene al punto $x$.
 \end{theorem}
 Una prueba del Teorema  \ref{theoremProductoLocal}   puede encontrarse en \cite[Theorem 3.12]{Bowen} o en
 \cite[Proposition 6.4.21]{Katok-Hasselblatt}. El Teorema \ref{theoremProductoLocal} es  consecuencia del  llamado \lq\lq shadowing lemma\rq\rq (lema de sombreado).

 \vspace{.3cm}

 Ahora estamos en condiciones de demostrar el Teorema \ref{TheoremSRBanosov} que establece la existencia, y mutua equivalencia, de las medidas erg\'{o}dicas de Gibbs y de las medidas SRB, para los difeomorfismos de Anosov de clase $C^{1 + \alpha}$.


\begin{nada} \em 
\label{proofTheoremSRBanosov} {\bf Demostraci\'{o}n del Teorema \ref{TheoremSRBanosov}}
\end{nada}
 Esta prueba consta de tres partes:

 \vspace{.3cm}

\index{medida! de Gibbs} \index{medida! condicional inestable} \index{derivada de Radon-Nykodim} \index{densidad} \index{medida! densidad de} \index{medida! equivalencia de}
 {\bf Afirmaci\'{o}n I:} \em Existen medidas de probabilidad $\mu$ de Gibbs   tales que:

  $\bullet$ Las probabilidades condicionales inestables $\mu^u $ son equivalentes a las medidas de Lebesgue $m^u$ a lo largo de las respectivas variedades inestables locales $W^u_{\delta}(x)$ para $\mu$-c.t.p. $x \in M$. \em
 (Nota: La equivalencia entre $\mu^u$ y $m^u$ se define como $\mu^u \ll m^u$ y $m^u \ll \mu^u$). \em

 $\bullet$ La derivada de Radon-Nikodym $d\mu^u/dm^u$ \em (la densidad de las medidas condicionales inestables $\mu^u$) \em es una funci\'{o}n real continua  y estrictamente positiva. \em

  \vspace{.3cm}

  Una vez probada la Afirmaci\'{o}n I, como consecuencia, debido al Teorema \ref{theoremDescoErgodicaEspaciosMetricos} de Descomposici\'{o}n en componentes erg\'{o}dicas, y por las partes (b) y (c) del Corolario \ref{corolarioRohlin1} del Teorema de Rohlin, concluimos que existen medidas de probabilidad erg\'{o}dicas   de Gibbs $\mu$ tales que $\mu^u $ es equivalente a $m^u$ para $\mu$-c.t.p. Aplicando entonces el Teorema \ref{theoremGibbs->SRB} ya demostrado, toda medida de probabilidad erg\'{o}dica y de Gibbs, es SRB. Deducimos que existen medidas   SRB erg\'{o}dicas, es decir, concluimos la parte (a) del Teorema \ref{TheoremSRBanosov}. Adem\'{a}s la Afirmaci\'{o}n I implica que existen medidas SRB erg\'{o}dicas, que son de Gibbs,   cuyas probabilidades condicionales inestables son equivalentes a la  medida de Lebesgue a lo largo de las variedades inestables respectivas para $\mu$-c.t.p., y cuyas densidades son funciones reales continuas y estrictamente positivas.

  \vspace{.3cm}
    {\bf Afirmaci\'{o}n II:} \em  Sea $\mu$ una medida de Gibbs erg\'{o}dica $\mu$ tal que la derivada de Radon-Nikodym \em (densidad  $d\mu^u/dm^u$ de las medidas condicionales inestables $\mu^u$) \em  es una funci\'{o}n continua \em (no necesariamente estrictamente positiva). \em Entonces la cuenca de atracci\'{o}n estad\'{\i}stica de $\mu$ cubre Lebesgue c.t.p. de un entorno abierto de su soporte. \em
    \index{medida! de Gibbs erg\'{o}dica}
    \index{derivada de Radon-Nykodim} \index{densidad}
    \index{medida! densidad de}
    \index{continuidad! de la funci\'{o}n densidad}
    \index{cuenca de atracci\'{o}n! estad\'{\i}stica}

   \vspace{.3cm}

    De esta propiedad, una vez demostrada, deducimos que dos medidas de Gibbs erg\'{o}dicas diferentes cuyas densidades inestables sean continuas (en particular dos medidas de Gibbs erg\'{o}dicas que satisfagan la Afirmaci\'{o}n I), deben te\-ner soportes disjuntos   a distancia positiva. En efecto, por la Definici\'{o}n \ref{DefinicionCuencaDeAtraccionEstadistica}, las cuencas de atracci\'{o}n estad\'{\i}stica de medidas diferentes deben ser disjuntas. Debido a la Afirmaci\'{o}n II los soportes de esas dos medidas est\'{a}n mutuamente aislados. De aqu\'{\i} deducimos (usando que $M$ es un espacio m\'{e}trico compacto y por lo tanto existe base numerable de abiertos) que las medidas de Gibbs erg\'{o}dicas que satisfacen la Afirmaci\'{o}n I (y por lo tanto tambi\'{e}n la II) son a lo sumo una cantidad numerable.

    \vspace{.3cm}

     {\bf Afirmaci\'{o}n III:} \em Lebesgue casi todo punto de la variedad $M$ pertenece a   la uni\'{o}n de las cuencas de atracci\'{o}n estad\'{\i}stica de las medidas de Gibbs erg\'{o}dicas cuyas medidas condicionales inestables   cumplen las dos propiedades de la Afirmaci\'{o}n \em I.
\index{cuenca de atracci\'{o}n! estad\'{\i}stica}
\index{medida! de Gibbs erg\'{o}dica}
 \index{medida! condicional inestable}

    \vspace{.3cm}

     Supongamos probada la  afirmaci\'{o}n   III. Tomando en particular una medida SRB o f\'{\i}sica cualquiera $\nu$,  con cuenca de atracci\'{o}n estad\'{\i}stica $A$, deducimos que $A$ est\'{a} contenido en la uni\'{o}n de las cuencas de atracci\'{o}n estad\'{\i}stica (que son disjuntas dos a dos) de una colecci\'{o}n finita o infinita numerable de medidas de Gibbs erg\'{o}dicas que satisfacen la Afirmaci\'{o}n I. Como las cuencas de atracci\'{o}n estad\'{\i}stica de medidas diferentes son disjuntas, deducimos que toda medida SRB o f\'{\i}sica es de Gibbs erg\'{o}dica y satisface la Afirmaci\'{o}n I, y por lo tanto tambi\'{e}n la II. Esto, junto con el Teorema \ref{theoremGibbs->SRB}, prueba la  parte (b) del Teorema \ref{TheoremSRBanosov}.

     Adem\'{a}s, la afirmaci\'{o}n III implica inmediatamente la parte (c) del Teorema \ref{TheoremSRBanosov}, y m\'{a}s a\'{u}n, toda medida SRB no solo es de Gibbs erg\'{o}dica, sino que adem\'{a}s satisface las afirmaciones I y II.

     Supongamos ahora adem\'{a}s que $f$ es topol\'{o}gicamente transitivo. Si  hubiera dos medidas SRB (erg\'{o}dicas) diferentes $\mu$ y $\nu$, entonces sus cuencas de atracci\'{o}n estad\'{\i}stica $B(\mu)$ y $B(\nu)$ ser\'{\i}an disjuntas. Como $\mu$ y $\nu$ satisfacen la Afirmaci\'{o}n II, existen abiertos disjuntos $U$ y $V$ tales que Lebesgue c.t.p. de $U$ pertenece a $B(\mu)$ y Lebesgue c.t.p. de $V$ pertenece a $B(\nu)$. Como $B(\mu)$ es $f$-invariante, entonces $m(  f^n(U) \setminus B(\mu))= 0$ para todo $n \geq 1$. Como $f$ es topol\'{o}gicamente transitivo existe $n \geq 1$ tal que el abierto $f^n(U) \cap V \neq \emptyset$. Todo abierto (en particular $f^n(U) \cap V$) tiene medida de Lebesgue positiva, y Lebesgue c.t.p. de $f^n(U) \cap V$ pertenece a $ B(\mu)$ (por pertenecer a $f^n(U)$), y pertenece tambi\'{e}n a $B(\nu)$ (por pertenecer a $V$). Luego las cuencas de atracci\'{o}n estad\'{\i}stica $B(\mu)$ y $B(\nu)$ no son disjuntas, de donde $\mu = \nu$. Esto prueba la unicidad de la medida SRB, es decir, la parte (d) del Teorema \ref{TheoremSRBanosov}.

     En resumen, para demostrar completamente el Teorema \ref{TheoremSRBanosov}, basta probar las Afirmaciones I, II y III.

  \vspace{.3cm}


  {\bf Primera parte de la demostraci\'{o}n del Teorema \ref{TheoremSRBanosov}}

  (Prueba de la Afirmaci\'{o}n I de \S\ref{proofTheoremSRBanosov})

 {\em Demostraci\'{o}n: }.

\vspace{.2cm}

{\bf Paso 1:  Construcci\'{o}n de candidata  a medida de Gibbs}

Elijamos   una variedad inestable local cualquiera $W_0$ y    consideremos   la medida $m^u_0$ de Lebesgue a lo largo de $W_0$. Tenemos $m_0^u(W_0) >0$. Definimos la siguiente medida de probabilidad $\mu_0$ en la variedad ambiente $M$:
$$\mu_0(A) := \frac{m^u(A \cap W_0)}{m^u(W_0)} \ \ \ \forall \ A \in {\mathcal A},$$
donde ${\mathcal A}$ es la sigma-\'{a}lgebra de Borel.
Para todo $n \geq 1$ definimos las siguientes medidas de probabilidad $\nu_n$ y $\mu_n$:
\begin{equation} \label{eqn nu_n-mu_n}\nu_n (A):= \mu_0 (f^{-n}(A)), \ \ \ \mu_n(A) := \frac{1}{n} \sum_{j= 0}^{n-1} \nu_j(A) \ \ \ \forall \ A \in {\mathcal A}.\end{equation}
Finalmente tomamos una subsucesi\'{o}n   $\{\mu_{n_i}\}_{i \in \mathbb{N}}$ convergente a   $\mu$ en el espacio ${\mathcal M}$ de las medidas de probabilidad en $(M, \mathcal A)$ con la topolog\'{\i}a d\'{e}bil$^*$. Probaremos que la medida $\mu$ satisface la Afirmaci\'{o}n I.

\vspace{.2cm}


{\bf Paso 2: Descomposici\'{o}n de Rohlin de las medidas $\nu_n$. } \index{descomposici\'{o}n de Rohlin}

Consideremos un punto $x_0$ en el soporte de $\mu$ tal que existe una bola $B= B_{\delta}(x_0) \subset M$ de radio $\delta >0$ suficientemente peque\~{n}o, tal que $\mu(B) >0$. Entonces, por la definici\'{o}n de la topolog\'{\i}a d\'{e}bil estre\-lla (tomando por ejemplo una funci\'{o}n continua no negativa que vale 1 en $x_0$ y est\'{a} soportada en $B$), existe $N \geq 1$ tal que $\mu_N(B) >0$, de donde deducimos que, para todo $n \geq N$ existe $j \in \{0, \ldots, n-1\}$ tal que $\nu_j(B) >0$, y por lo tanto \begin{equation} \label{eqn-4} \mu_n(B) >0 \ \ \forall \ n \geq N.\end{equation} Esto implica, por la definici\'{o}n de la medida $\nu_j$ que $f^j(W_0) \cap B \neq \emptyset$.
Fijemos $n \geq 1$ tal que
$$f^n(W_0) \cap B \neq \emptyset.$$
Por construcci\'{o}n, para todo boreliano $A \subset M$, se cumple:
$$\nu_n(A \cap B) = \frac{1}{m_0^u(W_0)} \int  \chi_{f^{-n}(A)}(y) \chi_{W_0 \cap f^{-n}(B)}(y) \, d m_0^u(y).$$
Haciendo el cambio de variable $x= f^n(y)$ y denotando $J^u(y)$ al Jacobiano inestable en el punto $y$ definido por
$$J^u (y) := \big|\det df_y|_{E^u(y)} \big|,$$ donde $E_y^u = T_yW^u(y)$, obtenemos:
$$\nu_n(A \cap B) = \frac{1}{m_0^u(W_0)} \int  \chi_{A}(x) \chi_{f^n(W_0) \cap B}\,(x)  \, \big |\det d(f^{-n})|_{E^u_x}(x)\big|\, \, d m^u(x)$$
donde $m^u(x)$ denota la medida de Lebesgue a lo largo de la variedad inestable local $W^u_{\delta}(x)$ por el punto $x$.
Tenemos $$\big |\det d(f^{-n})|_{E^u_x}(x)\big| = \frac{1}{\big |\det d(f^{n})|_{E^u_{y}}(y)\big|} = \frac{1}{\prod_{i= 0}^{n} J^u(f^i(y)) } = $$ $$=\frac{1}{\prod_{h= 1}^{n} J^u(f^{-h}(x)) },$$
de donde:
$$\nu_n(A \cap B) = \frac{1}{m_0^u(W_0)} \int  \chi_{A}(x) \chi_{f^n(W_0) \cap B}\,(x)   \,\frac{1}{\prod_{h= 1}^{n} J^u(f^{-h}(x)) } \, d m^u(x).$$
Denotamos $k_n$ el n\'{u}mero de componentes conexas de $f^n(W_0) \cap B$ (por convenci\'{o}n, $k_n = 0$ si el conjunto es vac\'{\i}o y $\sum_{i= 1} ^0 \cdot = 0$). Para cada $n \geq 1$ fijo, si $k_n \geq 1$  denotamos $\{W_{i, n}\}_{1 \leq i \leq k_n}$ al conjunto de componentes conexas de $f^n(W_0) \cap B$. Observamos que para todo $x \in W_{i,n}$ la medida $m^u(x)$ es la   de Lebesgue a lo largo de la subvariedad $W_{i,n}$. Entonces:
\begin{equation}\label{eqn-3}\nu_n(A \cap B) = \frac{1}{m_0^u(W_0)} \sum_{i= 1}^{k_n}\int_{W_{i,n}}  \chi_{A}(x)     \,\frac{1}{\prod_{h= 1}^{n} J^u(f^{-h}(x)) } \, d m^u(x).\end{equation}
Fijamos $n \geq 1$ tal que $k_n \geq 1$, es decir $f^n(W_0) \cap B \neq \emptyset$. Luego, por cons\-truc\-ci\'{o}n de la medida $\nu_n$ mediante la igualdad (\ref{eqn nu_n-mu_n}), tenemos $\nu_n(B) >0$. Tomemos una  subvariedad cualquiera $ V$,   encajada en $B_{\delta}$, de dimensi\'{o}n complementaria a la inestable, que interseque en un solo punto cada una, a todas las variedades estables locales inestables $W^u{\delta}(z) $ (para todo $z \in B$).   Por el Teorema \ref{theoremProductoLocal} de estructura de producto local, tal variedad topol\'{o}gicamente transversal $V$ existe, ya que puede tomarse igual a una variedad estable local.

Para cada $i \in \{1, \ldots, k_n\}$ consideremos el punto   $x_{i,n} =  W_{i,n} \cap V$. Construimos, para cada $x \in \bigcup_{i= 1}^{k_n} W_{i, n}$, el siguiente valor real \em positivo \em    $h_n(x)$ de la que llamaremos funci\'{o}n $h_n$:
\begin{equation} \label{eqnhsubn} h_{n}(x) := \frac{\prod_{h= 1}^{n} J^u(f^{-h}(x_{i,n}))}{\prod_{h= 1}^{n} J^u(f^{-h}(x))},\end{equation}
donde $i \in \{1, \ldots, n\}$ es el \'{u}nico \'{\i}ndice tal que $x \in W_{i,n}$.

Consideramos la partici\'{o}n ${\mathcal P}_n$ de la bola $B$ cuyas piezas son $\{W_{i,n}\}_{1 \leq i \leq k_n}$ y adem\'{a}s el complemento en $B$ de $\cup_{i= 1}^{k_n} W_{i,n}$. Denotamos con
$\rho_n$ la si\-guien\-te medida (finita) en el espacio medible cociente formado por  las piezas de la partici\'{o}n ${\mathcal P}_n$:
\begin{equation} \label{eqnrhosubn}\rho_n = \frac{1}{m_0^u(W_0)}   \sum_{i= 1}^{k_n} \frac{\int_{W_{i,n}} h_n \, d m^u} {\prod_{h= 1}^{n}   J^u \circ f^{-h} (x_{i,n})} \,  \delta_{x_{i,n}}, \end{equation} donde $\delta_{x_{i,n}}$ denota la Delta de Dirac soportada en la pieza $W_{i,n}$ (repre\-sentada por el punto $x_{i,n}$). Sustituyendo (\ref{eqnhsubn}) y (\ref{eqnrhosubn}) en la igualdad (\ref{eqn-3}), obtenemos para todo boreliano $A$ la siguiente descomposici\'{o}n de la medida de probabilidad $\nu_n$ restringida a la bola $B$:
\begin{equation}
\label{eqnDescoRohlin nu_n}
\frac{\nu_n(A \cap B)}{\nu_n(B)} =  \frac{1}{\nu_n(B)} \int d \rho_n   \int_{W_{i,n}} \chi_{A}(x)   \,\frac{h_n(x)}{\int_{W_{i,n}} h_n \, d m^u} \, d m^u(x). \end{equation}
La igualdad anterior  da la descomposici\'{o}n de Rohlin con respecto de la partici\'{o}n ${\mathcal P}_n$, de la medida de probabilidad $\nu_n/\nu_n(B)$ en la bola $B$. En efecto, la medida de probabilidad $\nu^u_{n} $ condicionada a lo largo de las piezas de la  partici\'{o}n (i.e. a lo largo de las variedades inestables locales $W_{i,n}$) est\'{a} dada por la integral de la derecha en la igualdad (\ref{eqnDescoRohlin nu_n}).   Esta medida condicionada $\nu^u_n$ cumple $$d\nu^u_{n} = \frac{h_n}{\int_{W_{i,n}} h_n \, dm^u } \, d m^u,$$ y por lo tanto
$$\nu^u_{n}\ll m^u.$$ Decimos entonces que las medidas condicionales inestables de $\nu_n$ son absolutamente continuas (a\'{u}n cuando $\nu_n$ no es necesariamente invariante con $f$). La integral de la izquierda en (\ref{eqnDescoRohlin nu_n}) da la medida de probabilidad $\widehat \nu_n$ en el espacio cociente de la bola $B$ con respecto a la partici\'{o}n ${\mathcal P}_n$. Esta medida cociente $\widehat \nu_n$ es
$$\widehat \nu_n = \frac{1}{\nu_n(B)}\,  {\rho_n}. $$
El objetivo en los siguientes pasos de la demostraci\'{o}n es usar la descomposici\'{o}n de Rohlin de las medidas $\nu_j$ dada por la igualdad (\ref{eqnDescoRohlin nu_n}) para todo $j \in \{0, \ldots, n-1\}$, para encontrar la descomposici\'{o}n de Rohlin de la medida $\mu_n = \sum_{j= 0}^{n-1}\mu_n$. Luego, pasando al l\'{\i}mite para una subsucesi\'{o}n $\{n_i\}_{i \in \mathbb{N}}$, encontraremos la descomposici\'{o}n de Rohlin de la medida $\mu$ construida en el Paso 1, para demostrar que es medida de Gibbs y satisface la Afirmaci\'{o}n I. El punto dif\'{\i}cil en este argumento es justificar rigurosamente el pasaje al l\'{\i}mite de los promedios de las medidas $\widehat \nu_n$ y de las medidas condicionales $\nu_n^u$. Para ello, necesitamos  $\epsilon$-aproximar la descomposici\'{o}n de Rohlin dada por la igualdad (\ref{eqnDescoRohlin nu_n}),   para  todo $\epsilon >0$
arbitrariamente peque\~{n}o, y para $n$ suficientemente grande.

\vspace{.2cm}

{\bf Paso 3: $\epsilon $- aproximaci\'{o}n de la descomposici\'{o}n de Rohlin de $\nu_n$} \index{descomposici\'{o}n de Rohlin}

Fijemos, como en el paso 2,  $n \geq 1$ tal que $f^n(W_0) \cap B \neq \emptyset$. Recordamos que $\{W_{i,n}\}_{1 \leq i \leq k_n}$ denota a las componentes conexas de $f^n(W_0) \cap B$, que son por lo tanto, variedades inestables locales contenidas en la bola $B$. Recordemos que para obtener la descomposici\'{o}n de Rohlin de $\nu_n$ dada por la igualdad (\ref{eqnDescoRohlin nu_n}), hemos elegido y dejado fijo un punto y uno solo $x_{i,n} \in W_{i,n} $. Definimos ahora, para todo $x \in  \bigcup_{i= 1}^{k_n} W_{i,n}$ el siguiente valor $h(x) \in [0, + \infty]$ de la que llamaremos funci\'{o}n $h$  restringida a $\bigcup_{i= 1}^{k_n} W_{i,n}$:
\begin{equation} \label{eqnh } h (x) := \frac{\prod_{h= 1}^{+ \infty} J^u(f^{-h}(x_{i,n}))}{\prod_{h= 1}^{+ \infty} J^u(f^{-h}(x))},\end{equation}
donde $i$ es el \'{u}nico \'{\i}ndice en $\{1, \ldots, k_n\}$ tal que $x \in W_{i,n}$
Ahora aplicamos el   Lema \ref{lemmaDistoAcotada} de Distorsi\'{o}n Acotada  que establece que:
 Existe una constante real $K >0$ tal que \begin{equation} \label{eqn(i)}   \ \frac{1}{K} < h(x) < \frac{1}{K} \ \ \forall \ x \in H:= \bigcup_{n \geq 1}^{+ \infty} \bigcup_{i= 1}^{k_n} W_{i,n}\end{equation}
 \begin{equation}
\label{eqn(ii)}  \ h: H \rightarrow \mathbb{R}^+ \mbox{ es continua } \end{equation}
   Para todo $\epsilon >0$ existe $N \geq 1$ (independiente de $x \in H$) tal que:
\begin{equation}
\label{eqn(iii)}
e^{- \epsilon} < \frac{h_n(x)}{h(x)} < e^{\epsilon} \ \ \forall \ x \in \bigcup_{i=1}^{k_n} W_{i,n} \ \ \forall \ n \geq N.\end{equation}

Aplicando la afirmaci\'{o}n  (\ref{eqn(i)})  del Lema de Distorsi\'{o}n Acotada, cons\-truimos, para cada $n \geq 1$ fijo tal que $f^n(W_0) \cap B \neq \emptyset$, la siguiente medida finita $I_n$:
\begin{equation} \label{eqnDescoRohlinI_n} I_n(A) := \int d \rho_n   \int_{W_{i,n}} \chi_{A}(x)   \,\frac{h (x)}{\int_{W_{i,n}} h  \, d m^u} \, d m^u(x),  \ \ \forall \ A \in {\mathcal{A}} \end{equation}
donde $\rho_n$ es la medida definida en la igualdad (\ref{eqnrhosubn}).
Comparando las igualdades (\ref{eqnDescoRohlin nu_n}) y (\ref{eqnDescoRohlinI_n}), obtenemos, para todo boreliano $A$   la siguiente igualdad:
$$        \nu_n (A \cap B)  =   \int d \rho_n   \int_{W_{i,n}} \chi_A    \,\frac{h_n  }{\int_{W_{i,n}} h_n  \, d m^u} \, d m^u = $$ $$=\int d \rho_n   \int_{W_{i,n}} \chi_A   \,\frac{h_n}{h}\frac{h }{\int_{W_{i,n}} h (h_n/h)  \, d m^u} \, d m^u  $$
Usando la propiedad (\ref{eqn(iii)}) del Lema de Distorsi\'{o}n Acotada, deducimos, para todo $n \geq N$ tal que $f^n(W_0) \cap B \neq \emptyset$ y para todo Boreliano $A$:
 $$e^{-2\epsilon} I_n(A) \leq  {\nu_n(A \cap B)}   \leq e^{2\epsilon} I_n(A)$$
 En conclusi\'{o}n:
 \begin{equation}
 \label{eqnDescoRohlinAproxnu_n}
 e^{-2 \epsilon} \int d \rho_n   \int_{W_{i,n}} \chi_{A}(x)   \,\frac{h (x)}{\int_{W_{i,n}} h  \, d m^u} \, d m^u(x) \leq \nu_n(A \cap B) \leq $$ $$e^{ 2 \epsilon} \int d \rho_n   \int_{W_{i,n}} \chi_{A}(x)   \,\frac{h (x)}{\int_{W_{i,n}} h  \, d m^u} \, d m^u(x)   \ \ \forall \ A \in {\mathcal{A}}
 \end{equation}

 \vspace{.2cm}

 {\bf Paso 4: Descomposici\'{o}n de Rohlin de la medida $\mu$} \index{descomposici\'{o}n de Rohlin}

 Sea la medida de probabilidad $\mu$ construida en el Paso 1:
 \begin{equation} \label{eqnmu}\mu= \lim_{i \rightarrow + \infty} \mu_{n_i},\end{equation}
 donde el l\'{\i}mite es en la topolog\'{\i}a d\'{e}bil estrella, y $\mu_n$ est\'{a} definida para todo $n \geq 1$ por la igualdad (\ref{eqn nu_n-mu_n}), como el promedio aritm\'{e}tico de las medidas $\nu_j$ para $j \in \{0, \ldots, n-1\}$. Consideremos $n \geq N$ fijo y la partici\'{o}n medible ${\mathcal Q}$ de la bola $B$ formada por las variedades inestables locales. Denotemos $B / \sim$ el espacio medible cociente de $B$ con respecto a esa partici\'{o}n. Usando las desigualdades (\ref{eqnDescoRohlinAproxnu_n}) obtenemos, para todo boreliano $A \subset M$ la  siguiente $\epsilon$-aproximaci\'{o}n de $\mu_n(A)$
 \begin{equation}
 \label{eqnDescoRohlinAproxmu_n}
 e^{-2 \epsilon} \int_{B/\sim} d \Big(\frac{1}{n}\sum_{j= 0}^{n-1}\rho_j \Big) \, \int_{W} \chi_{A}    \,\frac{h  }{\int_{W} h  \, d m^u} \, d m^u(x) \leq \mu_n(A \cap B) \leq $$ $$e^{ 2 \epsilon} \int d \frac{1}{n}\Big(\sum_{j= 0}^{n-1}\rho_j \Big) \,  \int_{W} \chi_{A}    \,\frac{h  }{\int_{W} h  \, d m^u} \, d m^u(x)  \ \ \forall \ A \in {\mathcal{A}}, \ \ \ \forall \ n \geq N.
 \end{equation}
donde, por convenci\'{o}n, $\rho_j$ est\'{a} definida por la igualdad (\ref{eqnrhosubn}) si

$f^j(W_0) \cap B \neq \emptyset$, y $\rho_j := 0$ en caso contrario.

Denotamos \begin{equation} \label{eqnrho*_n} \rho^*_n := \frac{1}{n} \sum_{j= 0}^{n-1} \rho_j.\end{equation}
Por construcci\'{o}n $\rho^*_n$ es una medida finita en el espacio cociente $B/ \sim$ de la bola $B$ con respecto a la partici\'{o}n ${\mathcal Q}$ en subvariedades locales inestables.
Por la igualdad (\ref{eqn-4}) $\mu_n(B)>0$ para todo $n$ suficientemente grande (digamos $n \geq N$). Entonces, debido a las desigualdades (\ref{eqnDescoRohlinAproxmu_n}),   $\rho^*_n$ es una medida no nula para todo $n \geq N$. Probemos que la sucesi\'{o}n de medidas $\{\rho^*_n\}_{n \geq N}$ est\'{a} uniformemente acotada superiormente por $e^{\epsilon}$. En efecto, por un lado la medida $\mu_n$ construida por la igualdad (\ref{eqn nu_n-mu_n}) es una medida de probabilidad en toda la variedad $M$, y por otro lado, de la desigualdad a la izquierda en (\ref{eqnDescoRohlinAproxmu_n}), obtenemos
$$1 \geq \mu(B) \geq e^{-2 \epsilon} \int _{B/ \sim} d \rho_n^* = \rho^*(B/\sim).$$
El espacio de todas las medidas finitas en el espacio medible $B/\sim$ que est\'{a}n uniformemente acotadas por $e^{2 \epsilon}$, provisto de la topolog\'{\i}a d\'{e}bil estrella, de secuencialmente compacto. Luego, existe una subsucesi\'{o}n $\{n_{i_k}\}_{k \in \mathbb{N}}$ de $\{n_i\}_{i \in \mathbb{N}}$ (que por simplicidad seguimos denotando como $n_i$), tal que $\{\rho^*_{n_i}\}_{i}$ es convergente a una medida finita $\rho^*$.
En resumen, hemos demostrado la existencia de una medida finita $\rho^*$ en el espacio medible cociente $B/\sim$, tal que:
\begin{equation} \label{eqnrho*}\rho^* = \lim_{i \rightarrow + \infty} \rho^*_{n_i} = \lim_{u \rightarrow + \infty} \frac{1}{n_i} \sum_{j= 0}^{n_i -1} \rho_j,\end{equation}
donde el l\'{\i}mite de medidas es en la topolog\'{\i}a d\'{e}bil estrella.

Consideremos una funci\'{o}n real continua no negativa $\psi: M \mapsto [0,1]$.  Las desigualdades (\ref{eqnDescoRohlinAproxmu_n}) y la igualdad (\ref{eqnrho*_n}) implican
\begin{equation}
 \label{eqnDescoRohlinAproxmu_n_2}
 e^{-2 \epsilon} \int_{B/\sim}  \, d \rho^*_{n_i} \, \int_{W} \psi    \,\frac{h  }{\int_{W} h  \, d m^u} \, d m^u  \leq \int_B \psi \, d\mu_{n_i}  \leq $$ $$e^{2 \epsilon} \int_{B/\sim}   \, d \rho^*_{n_i} \, \int_{W} \psi    \,\frac{h  }{\int_{W} h  \, d m^u} \, d m^u  \ \ \forall \ n_i \geq N.
 \end{equation}
 Aplicando la propiedad (\ref{eqn(ii)}) del Lema de Distorsi\'{o}n Acotada, deducimos que el producto $\psi \cdot h$ es una funci\'{o}n continua en $B$. Luego, como tambi\'{e}n $m^u(W^u_{\delta}(x))$ depende continuamente de $x$ (porque el subespacio tangente $E^u_x $  depende continuamente de $x$), deducimos que la integral a la derecha en (\ref{eqnDescoRohlinAproxmu_n_2}) es un n\'{u}mero real positivo que depende continuamente de la variedad inestable $W \in B/\sim$. Entonces, por la igualdad (\ref{eqnrho*}) y por la definici\'{o}n de la topolog\'{\i}a d\'{e}bil estrella, deducimos que:
 $$\lim_{i \rightarrow + \infty} \int_{B/\sim}  \, d \rho^*_{n_i} \, \int_{W} \psi    \,\frac{h  }{\int_{W} h  \, d m^u} \, d m^u = $$ $$= \int_{B/\sim}  \, d \rho^*  \, \int_{W} \psi    \,\frac{h  }{\int_{W} h  \, d m^u} \, d m^u. $$
  Tomando l\'{\i}mite en las desigualdades (\ref{eqnDescoRohlinAproxmu_n_2}) y recordando (\ref{eqnmu}), obte\-ne\-mos:
\begin{equation}
 \label{eqnDescoRohlinAproxmu_n_3}
 e^{-2 \epsilon} \int_{B/\sim}   \, d \rho^*  \, \int_{W} \psi    \,\frac{h  }{\int_{W} h  \, d m^u} \, d m^u  \leq \int_B \psi \, d\mu   \leq $$ $$e^{2 \epsilon} \int_{B/\sim}   \, d \rho^*  \, \int_{W} \psi    \,\frac{h  }{\int_{W} h  \, d m^u} \, d m^u
 \end{equation}
 Por linealidad, la integral de la derecha en (\ref{eqnDescoRohlinAproxmu_n_3}) define un funcional lineal acotado y positivo en   el espacio $C^0(M, \mathbb{R})$ de todas las funciones continuas $\psi: M \mapsto \mathbb{R}$. Luego, por el Teorema de Representaci\'{o}n de Riesz, define una medida finita $\widetilde \mu $.
 Como las desigualdades (\ref{eqnDescoRohlinAproxmu_n_3})  valen para toda $\psi \in C^0(M, [0,1]$, y la funci\'{o}n caracter\'{\i}stica $\chi_A$ de cualquier boreliano $A$ puede aproximarse en $L^1(\mu)$ y en $L^1(\widetilde \mu)$ por una funci\'{o}n continua $\psi \in C^0(M, [0,1])$,  obtenemos, para todo boreliano $A$ las siguientes desigualdades:
 \begin{equation}
 \label{eqnDescoRohlinAproxmu}
 e^{-2 \epsilon} \int_{B/\sim}    d \rho^*  \, \int_{W} \chi_A   \,\frac{h  }{\int_{W} h  \, d m^u} \, d m^u  \leq \mu(A \cap B) \leq $$ $$\int_{B/\sim}   \, d \rho^*  \, \int_{W} \chi_A    \,\frac{h  }{\int_{W} h  \, d m^u} \, d m^u.
 \end{equation}
 Por construcci\'{o}n, las medidas involucradas en las desigualdades (\ref{eqnDescoRohlinAproxmu}) son independientes del n\'{u}mero $\epsilon>0$, pero valen para todo $\epsilon >0$ suficientemente peque\~{n}o. Luego, haciendo $\epsilon \rightarrow 0^+$, deducimos:
\begin{equation}
 \label{eqnDescoRohlinmu}
  \mu(A \cap B) = \int_{B/\sim}   \, d \rho^*  \, \int_{W} \chi_A    \,\frac{h  }{\int_{W} h  \, d m^u} \, d m^u \ \ \forall \ A \in {\mathcal A}.
 \end{equation}
 Dividiendo la igualdad (\ref{eqnDescoRohlinmu}) entre $\mu(B) >0$, hemos encontrado   la descomposici\'{o}n de Rohlin de la medida de probabilidad $\mu$ restringida a la bola $B$, con respecto a la partici\'{o}n ${\mathcal Q}$ de $B$ en subvariedades inestables locales. En efecto, la medida $\widehat \mu$ en el espacio cociente $B/\sim$ con respecto a esa partici\'{o}n es $\rho^*/\mu(B)$ (es inmediato chequear que esta es una medida de probabilidad en $B/\sim$); y las medidas condicionales inestables $\mu^u$ de $\mu$ est\'{a}n definidas por las integrales de la derecha en la igualdad (\ref{eqnDescoRohlinmu}). Luego, $$\mu^u \ll m^u$$ para $\widehat \mu$-casi toda variedad inestable $W \in B/\sim$. Adem\'{a}s, la derivada de Radon-Nikodym de $\mu^u$ con respecto a $m^u$ (la densidad) es:
 $$\frac{d\mu^u}{dm^u}(x) = h. $$
 Usando las propiedades (\ref{eqn(i)}) del Lema de Distorsi\'{o}n Acotada, $1/h$ es una funci\'{o}n continua y acotada superiormente. Entonces
 $$\frac{dm^u}{d\mu^u}(x) = \frac{1}{h}, $$
 lo cual implica que $$m^u \ll \mu^u,$$
 terminando de demostrar que $\mu$ es una medida de Gibbs,   que sus medidas condicionales inestables $\mu^u$ son equivalentes a las medidas de Lebesgue $m^u$ a lo largo de las variedades inestables, para $\mu$-c.t.p., y que su funci\'{o}n densidad inestable $h$ es continua y estrictamente positiva. Esto completa la prueba de la Afirmaci\'{o}n I del par\'{a}grafo   \S\ref{proofTheoremSRBanosov}. \hfill $\Box$

 {\bf Nota: } Adem\'{a}s, usando la igualdad (\ref{eqnh }) hemos probado que la densidad de las medidas condicionales inestables de $\mu$ a lo largo de la variedad inestable local $W^u_{\delta}(x_i)$, es
 $$h (x) = \frac{\prod_{h= 0}^{+ \infty} J^u(x_i)}{\prod_{h= 0}^{+ \infty} J^u(x)},$$
 demostrando tambi\'{e}n lo afirmado en la Observaci\'{o}n \ref{remarkSRBanosovDensidad}.

 \vspace{.3cm}
\newpage

{\bf Segunda parte de la demostraci\'{o}n del Teorema \ref{TheoremSRBanosov}}

(Prueba de la Afirmaci\'{o}n II de   \S\ref{proofTheoremSRBanosov})

\begin{exercise}\em
Probar la Afirmaci\'{o}n II de \S\ref{proofTheoremSRBanosov}.

Sugerencia: Sea dada $\mu$, una medida de probabilidad de Gibbs erg\'{o}dica con densidad inestable continua. Hay que probar   que existe un conjunto medible con $\mu$-medida 1 y un abierto que lo contiene, tal que Lebesgue c.t.p. de ese abierto est\'{a} en la cuenca de atracci\'{o}n estad\'{\i}stica $B$ de   $\mu$.  Para $\mu$-c.t.p. $x_0 \in M$,   toda bola $B_{\delta}(x_0)$ tiene $\mu$-medida positiva. Probemos que Lebesgue casi todo punto de $B_{\delta}(x_0)$, para $\delta >0$ suficientemente peque\~{n}o que depende de $x_0$, pertenece a la cuenca $B$. Consideremos la densidad $h(x_0)$ de $\mu^u_{x_0}$ respecto la medida de Lebesgue inestable $m^u_{x_0}$. Cualquier funci\'{o}n densidad es no negativa. Por definici\'{o}n de la funci\'{o}n densidad inestable $ h= d \mu^u/dm^u$, y por el teorema de descomposici\'{o}n de Rohlin,  $h_{x_0} >0$ para $\mu$-c.t.p. $x_0 \in M$. Como por hip\'{o}tesis $h$ es continua,   entonces $h_{x} >0$ para todo $x \in W^u_{\delta}(x_0)$, si tomamos $\delta >0$ suficientemente peque\~{n}o (dependiendo del punto $x_0$).  Ahora podemos aplicar los mismos argumentos de la demostraci\'{o}n del Teorema \ref{theoremGibbs->SRB}, usando la condici\'{o}n (\ref{eqnholonomiaAC}) de continuidad absoluta de la holonom\'{\i}a $h_s$ a lo largo de la foliaci\'{o}n estable, establecida en el Teorema \ref{theoremTeoriaPesin} de la Teor\'{\i}a de Pesin. En extenso, por la construcci\'{o}n sugerida antes, la  medida condicionada inestable $\mu^u$ a lo largo de la variedad $W^u_{\delta}(x_0)$  tiene densidad estrictamente positiva.   Entonces   es equivalente a la medida de Lebesgue $m^u$ en esa variedad inestable local.   Usar esta propiedad, y la ergodicidad de $\mu$ para probar que    $m^u$ c.t.p. de $W^u_{\delta}(x_0)$ pertenece a $B$  (i.e. $\mu$-c.t.p. pertenece a la cuenca $B$, lo cual implica que $\mu^u$ c.t.p. de $W^u_{\delta}(x_0)$ pertenece a $B$ debido al Teorema \ref{theoremRohlin}).  Usar    (\ref{eqnholonomiaAC}) en ambas direcciones, para deducir que la medida de Lebesgue $m$ de $h_s{-1}(W^u_{\delta}(x_0) \setminus B) \subset M$ es cero. Finalmente aplicar el Teorema \ref{theoremProductoLocal} de estructura de producto local del difeomorfismo de Anosov $f$, para deducir que $h_s^{-1}(W^u_{\delta}(x_0)) = B_{\delta}(x_0)$, y el mismo argumento que en la prueba del Teorema \ref{theoremGibbs->SRB} para demostrar que $h_s^{-1}(B \cap W^u_{\delta}(x_0)) \subset B$ y concluir que $m$-c.t.p. de $B_{\delta}(x_0)$ pertenece a $B$. Concluimos que para $\mu$-c.t.p. $x_0 \in M$ existe $\delta >0$ (que puede depender de $x_0$) tal que $m$-c.t.p. de $B_{\delta}(x_0)$ pertenece a la cuenca de atracci\'{o}n estad\'{\i}stica $B$ de la medida de Gibbs erg\'{o}dica $\mu$. Tomando la uni\'{o}n de esas bolas abiertas, se obtiene un abierto que contiene al soporte de $\mu$ y tal que $m$-c.t.p. de ese abierto est\'{a} contenido en la cuenca $B$.
\end{exercise}

 {\bf Tercera y \'{u}ltima parte de la demostraci\'{o}n del Teorema \ref{TheoremSRBanosov}}

 (Prueba de la Afirmaci\'{o}n III de \S\ref{proofTheoremSRBanosov})

{\em Demostraci\'{o}n: } Debido a la Afirmaci\'{o}n I ya probada, existen medidas   de Gibbs $\mu$  cuyas medidas condicionales inestables $\mu^u$  satisfacen   las dos propiedades en dicha afirmaci\'{o}n.  Aplicando el Teorema \ref{theoremGibbs->SRB}, la parte (b) del Corolario \ref{corolarioRohlin1}, y la Definici\'{o}n \ref{definitionMedidaSRB}, deducimos que la cuenca de atracci\'{o}n estad\'{\i}stica   de las componentes erg\'{o}dicas de $\mu$ tienen medida de Lebesgue positiva, y estas componentes erg\'{o}dicas son medidas de Gibbs erg\'{o}dicas cuyas  condicionales inestables   satisfacen   las dos propiedades de la Afirmaci\'{o}n I. Por lo tanto, el   conjunto $B$ formado por todos los puntos que pertenecen a las cuencas de atracci\'{o}n estad\'{\i}stica de medidas de Gibbs erg\'{o}dicas que satisfacen esas dos propiedades, cumple $$m(B) >0.$$ Nota:   $B$ es medible: En efecto, , las medidas de Gibbs erg\'{o}dicas son a lo sumo, una cantidad numerable, porque  por el Teorema \ref{theoremGibbs->SRB} y la Definici\'{o}n \ref{definitionMedidaSRB}, sus cuencas de atracci\'{o}n estad\'{\i}stica son disjuntas dos a dos y tienen medida de Lebesgue positiva. La cuenca de atracci\'{o}n estad\'{\i}stica de cualquier medida de probabilidad, por su construcci\'{o}n  dada en la Definici\'{o}n \ref{DefinicionCuencaDeAtraccionEstadistica}, es un conjunto medible. Luego, la uni\'{o}n    numerable $B$ de conjuntos medibles, es medible.

Queremos probar que $m(B)= 1$, es decir $m(C)= 0$ donde $$C = M \setminus B.$$
Como la cuenca de atracci\'{o}n estad\'{\i}stica de cualquier medida de probabilidad es invariante con $f$, tenemos:
$$f(B)= B, \ \ \ \ \ f(C)= C.$$
Tomemos $\delta >0$ suficientemente peque\~{n}o. Como $m(B) >0$, deducimos que:

 $$  \mbox{\em Existe una bola $B_{\delta}(x_0)$ tal que $m(B \cap B_{\delta}(x_0))>0$.  \em } $$

 Afirmamos que
\begin{equation}
\label{eqnAprobarC}
\mbox{A probar: } \  \ \forall \ x_0  \in M, \ \ m(B \cap B_{\delta}(x_0))>0 \   \Rightarrow \   m(C \cap B_{\delta}(x_0) = 0.  \end{equation}

Primero veamos que una vez probada la afirmaci\'{o}n (\ref{eqnAprobarC}), se deduce la afirmaci\'{o}n (III). En efecto, si la variedad $M$ fuera   conexa, entonces la afirmaci\'{o}n (III) implica que las bolas abiertas de radio $\delta$ se clasifican en dos clases disjuntas: la de aquellas bolas tales que Lebesgue-c.t.p. de ella pertenece a $B$, y la de aquellas bolas tales que Lebesgue c.t.p. pertenece al complemento de $B$, o sea a $C$. Como la variedad es conexa, una de las dos clases es vac\'{\i}a. Como existe un bola $B_{\delta}(x_0)$ en la primera clase, entonces la segunda clase es vac\'{\i}a. Esto prueba que para toda bola de radio $\delta$, $m(C \cap B_{\delta})= 0$. Luego $0 = m(C) = m(M\setminus B), $ de donde deducimos que $m(B)= 0$, concluyendo la Afirmaci\'{o}n III, cuando la variedad $M$ es conexa. Y si la variedad $M$ no es conexa, como es compacta por hip\'{o}tesis, tiene una cantidad finita $k \geq 2$ de componentes conexas $M_1, \ldots, M_k$. Siendo $f$ un difeomorfismo de clase $C^{1 + \alpha}$, en particular, es continuo. Entonces existe un $p \geq 1$ tal que $f^p: M_i \mapsto M_i$ para todo $1 \leq i \leq k$, y $f^p|_{M_i} \in \mbox{Diff}^{1 + \alpha} (M_i)$. Como $f$ es Anosov, $f^p|_{M_i}$ es Anosov, donde $M_i$ es una variedad compacta y conexa. Por lo probado antes $m(C \cap M_i) = 0$ para todo $1 \leq i \leq k$. Luego $m(C) = 0$, terminando de probar la Afirmaci\'{o}n III de \S\ref{proofTheoremSRBanosov},  a partir de (\ref{eqnAprobarC}).

Ahora solo resta demostrar (\ref{eqnAprobarC}). Supongamos por absurdo que $$m(B \cap B_{\delta}(x_0)) >0, \ \ m(C \cap B_{\delta}(x_0)) >0.$$ Denotemos
$$W_0= W_{\delta}^u(x_0)$$ a la variedad inestable local por $x_0$ en la bola $B_{\delta}(x_0)$ (es decir, $W^u_{\delta}(x_0)$ es la componente conexa que contiene al punto $x_0$ de la intersecci\'{o}n $W^u(x_0) \cap B_{\delta}(x_0)$).

Debido a la propiedad (\ref{eqnholonomiaAC}), establecida en el Teorema \ref{theoremTeoriaPesin} de continuidad absoluta de la holonom\'{\i}a estable, tenemos
$$m^u(W_0 \cap B) >0, \ \ \ m^u(W_0 \cap C) >0$$
donde $m^u$ es la medida de Lebesgue inestable a lo largo de la variedad inestable local $W_0$.
Construyamos, como en el paso 1 de la demostraci\'{o}n de la Afirmaci\'{o}n I en \S\ref{proofTheoremSRBanosov}, las siguientes medidas de probabilidad en $M$. Para cualquier boreliano $A \subset M$, definimos:
$$\mu_0(A):= \frac{m^u(W_0 \cap A)}{m^u(W_0)} = \lambda \mu_{0,B} (A) + (1-\lambda) \mu_{0,C} (A),$$
donde $$\lambda:= \frac{m^u(B \cap W_0)}{m^u(W_0)}, \ \ 1 - \lambda:= \frac{m^u(C \cap W_0)}{m^u(W_0)},$$
$$\mu_{0,B} (A):= \frac{m^u(W_0 \cap A \cap B)}{m^u(W_0 \cap B)}, \ \ \ \frac{m^u(W_0 \cap A \cap C)}{m^u(W_0 \cap C)}.$$
$$\nu_n(A) := \mu_0(f^{-n}(A)) = \lambda \nu_{n,B}(A) + (1- \lambda) \nu_{n,C}(A),$$
donde $$\nu_{n,B}(A) := \mu_{0,B} (f^{-n}(A)), \ \ \nu_{n,C}(A) := \mu_{0,C} (f^{-n}(A)), $$
$$\mu_n := \frac{1}{n} \sum_{j= 0}^{n-1}\nu_{j} = \lambda \mu_{n, B} + (1- \lambda) \mu_{n,C},$$
donde
$$\mu_{n,B} := \frac{1}{n} \sum_{j= 0}^{n-1}\nu_{j, B}, \ \ \ \mu_{n,C} := \frac{1}{n} \sum_{j= 0}^{n-1}\nu_{j, C}.$$
Tomamos una subsucesi\'{o}n $\{n_i\}_{i \geq 1}$ tal que las sucesiones de medidas de probabilidad $\{\mu_{n_i}\}_{i \geq 1},$ \ $\{\mu_{n_i, B}\}_{i \geq 1}$  y  $\{\mu_{n_i, C}\}_{i \geq 1}$ son convergentes en la topolog\'{\i}a d\'{e}bil estrella. Es decir, existen medidas de probabilidad $\mu$, $\mu_B$ y $\mu_C$ tales que:
$$ \mu = \lim_{i \rightarrow + \infty}   \mu_{n_i}, \ \  \mu_B = \lim_{i \rightarrow + \infty}   \mu_{n_i, B}, \ \  \mu_C = \lim_{i \rightarrow + \infty}   \mu_{n_i,C}. $$
Entonces:
$$\mu = \lim_{i \rightarrow + \infty}   \mu_{n_i} = \lim_{i \rightarrow + \infty} \lambda \mu_{n_i, B} + (1 - \lambda) \mu_{n_i, C} =$$ $$= \lambda \lim_{i \rightarrow + \infty} \mu_{n_i, B} + (1 - \lambda) \lim_{i \rightarrow + \infty}\mu_{n_i, C} = \lambda   \mu_{ B} + (1 - \lambda)  \mu_{ C}.$$
 Luego $$\mu_C \ll \mu.$$ Por el teorema de Descomposici\'{o}n Erg\'{o}dica, como $\mu_C \ll \mu$, deducimos que:

 \em Las componentes erg\'{o}dicas ${\mu_C}_x$ de $\mu_C$ coinciden con las componentes erg\'{o}dicas $\mu_x$  de $\mu$, para $\mu_C$-c.t.p. $x \in M$. \em

  Por lo demostrado en la Afirmaci\'{o}n I de \S\ref{proofTheoremSRBanosov}, la medida de probabilidad $\mu$ es de Gibbs y sus medidas condicionales inestables satisfacen las dos propiedades de la Afirmaci\'{o}n I. Luego, usando las partes (b) y (d)  del Corolario \ref{corolarioRohlin1}, obtenemos:

  \em Las componentes erg\'{o}dicas $\mu_x$ de $\mu$ son medidas erg\'{o}dicas de Gibbs y sus medidas condicionales inestables $\mu_x^u$ satisfacen las dos propiedades de la Afirmaci\'{o}n \em I, \em para $\mu$-c.t.p. $x \in M$. \em

  De las dos enunciados probados arriba, concluimos que:

   \em Las componentes erg\'{o}dicas de  $\mu_C$ con medidas erg\'{o}dicas de Gibbs cuyas condicionales inestables satisfacen las dos propiedades de la Afirmaci\'{o}n \em I, \em para $\mu_C$-c.t.p. $x \in M$. \em

   Ahora aplicamos la Afirmaci\'{o}n II ya probada, para deducir el siguiente resultado:

  {\bf Enunciado (a) ya probado: } \em Lebesgue c.t.p. de un abierto $V$ que contiene al soporte de $\mu_C$, pertenece al conjunto $B= M \setminus C$. \em

  Pero, por construcci\'{o}n de la medida $\mu_C = \lim_i \mu_{n_i,C}   $, si tomamos una funci\'{o}n real continua no negativa $\psi: M \mapsto [0,1] $  soportada en $V$ y tal que $\psi_V >0$, se cumple:
  $$0 < \int \chi_V \psi \, d \mu_C = \int \psi \, d \mu_C =  \lim_{i \rightarrow + \infty} \int \psi \, d \mu_{n_i, C} = $$ $$= \lim_{i \rightarrow + \infty} \frac{1}{n_i} \sum_{j= 0}^{n_i - 1} \int \psi \, d \nu_{j, C}. $$
  Entonces, existe $j \geq 1$ tal que $$0 < \int \psi \, d \nu_{j, C} = \int \psi \circ f^{-j} \, d \mu_{0,C} = \int  (\chi_V  \cdot  \psi) \circ f^{-j} \, d \mu_{0,C}.$$
  (La \'{u}ltima igualdad proviene de que $\psi$ est\'{a} soportada en $V$, es decir $ \chi_V \cdot \psi = \psi$).
  Por la construcci\'{o}n de $\mu_{0, C}$, y debido a la positividad establecida en la \'{u}ltima desigualdad, deducimos que
  $$m^u(f^{-j}(V) \cap C \cap W_0) >0.$$
  Como $f$ es un difeomorfismo y $f^j(C)= C$, deducimos
  \begin{equation} \label{eqn-21}m^u(V \cap C \cap f^j(W_0))  >0. \end{equation}
  Ahora  tomamos una bola abierta $B_{\delta}(x_1) \subset V$, de radio $\delta >0$ suficientemente peque\~{n}o, centrada en un punto $x_1 \in V \cap C \cap f^j(W_0))$ tal que \begin{equation}
  \label{eqn-22}
  m^u(  C \cap   W^u_{\delta}(x_1)) >0,\end{equation} donde $W^u_{\delta}(x_1)$ es la componente conexa que contiene a $x_1$ de la intersecci\'{o}n $B_{\delta}(x_1) \cap f^j(W^u(x_0))   $, quien a su vez contiene a la componente conexa por el punto $x_1$ de $B_{\delta}(x_1) \cap f^j(W_0) $, cuya intersecci\'{o}n con $C$, debido a (\ref{eqn-21}),    tiene medida de Lebesgue $m^u$ positiva.

  Aplicamos la propiedad (\ref{eqnholonomiaAC}) de continuidad absoluta de la holonom\'{\i}a estable  $h_s:  B_{\delta}(x_1) \mapsto W_{\delta}^u(x_1)$, establecida en el Teorema \ref{theoremTeoriaPesin}. De (\ref{eqn-22}) deducimos
  $$m(h_s^{-1}( C \cap W^u_{\delta}(x_1))) >0. $$
  Como $h_s^{-1}(C) \subset C$ y $h_s^{-1}(W^u_{\delta}(x_1)) \subset B_{\delta}(x_1) \subset V$, deducimos que $m (C \cap V) >0$, lo cual contradice el enunciado (a) probado antes.
Esta contradicci\'{o}n termina de demostrar (\ref{eqnAprobarC}) por absurdo, y por lo tanto, con ella termina  la prueba de la Afirmaci\'{o}n III de \S\ref{proofTheoremSRBanosov}  y  del Teorema
 \ref{TheoremSRBanosov}.
\hfill $\Box$

\section{Atractores  estad\'{\i}sticos y medidas SRB-like}

A lo largo de este cap\'{\i}tulo $f: M \mapsto M$ es una transformaci\'{o}n continua en una variedad compacta riemanniana $M$ de dimensi\'{o}n finita. Denotamos con $m$ a la medida de Lebesgue en $M$, re-escalada para que sea una medida de probabilidad: $m(M)= 1$ (i.e. si $0 <m(M) \neq 1$, sustituimos $m$ por la probabilidad $m/m(M)$).

Definiremos primero atractor de Milnor. Este no es en general un atractor estad\'{\i}stico. Sin embargo, la demostraci\'{o}n de su existencia es pr\'{a}cticamente la misma que la demostraci\'{o}n de existencia de atractores estad\'{\i}sticos, que veremos en la secci\'{o}n siguiente.

\subsection{Atractores de Milnor}

\begin{definition}
\label{definitionAtractorMilnor} \em  \index{atractor! de Milnor}
{\bf Atractor de Milnor \cite{Milnor}} Se llama \em atractor de Milnor \em a un conjunto compacto no vac\'{\i}o $K \subset M$,  invariante por $f$ (i.e. $f^{-1}(K) = K$), tal que $$m(E_K) = 1,$$ donde el conjunto $E_K \subset M$, llamado \em cuenca de atracci\'{o}n (topol\'{o}gica) \em de $K$, est\'{a} definido por
\begin{equation} \label{eqnCuencaAtracMilnor}E_K  := \{x \in M: \ \ \lim_{n \rightarrow + \infty} \mbox{dist}(f^n(x), K) = 0\}.  \end{equation} \index{cuenca de atracci\'{o}n! topol\'{o}gica} \index{cuenca de atracci\'{o}n! de Milnor} \index{$C(K)$ o $ E_K$ cuenca de atracci\'{o}n! topol\'{o}gica del compacto $K$}
\end{definition}
\begin{exercise}\em \label{ejercicio7}
Probar que para cualquier conjunto compacto $K$ no vac\'{\i}o, el conjunto $E_K$ definido por la igualdad (\ref{eqnCuencaAtracMilnor}) es medible. Sugerencia: La funci\'{o}n $d: M \mapsto \mathbb{R}$ definida por $d(x) = \mbox{dist}(x,K)$ es medible, y $f: M \mapsto M$ es medible porque es continua. El l\'{\i}mite superior de una sucesi\'{o}n de funciones medibles es medible.
\end{exercise}

\begin{exercise}\em
Sean $K \subset K'$ dos atractores de Milnor. Probar que $E_K \subset E_{K'}$.
\end{exercise}

En la Definici\'{o}n \ref{definitionAtractorMilnor} observamos que para un atractor de Milnor, el criterio de observabilidad de su cuenca es el criterio Lebesgue-medible (cf. Observaci\'{o}n  \ref{remarkObservabilidadTopyEstad}), mientras que el tipo de atracci\'{o}n es topol\'{o}gica (cf. Definiciones \ref{definitionAtraccionTopol}
y \ref{definitionAtraccionEstad}).

Es inmediato chequear que todo atractor topol\'{o}gico es de Milnor. Sin embargo, no todo atractor de Milnor es topol\'{o}gico, como veremos en los Ejemplos \ref{ejemploMilnorNoTopologico}  y \ref{ejemploMilnorNoTopologico1}.

Tambi\'{e}n es inmediato chequear que todo atractor erg\'{o}dico es de Milnor. Pero no todo atractor de Milnor es erg\'{o}dico, como muestra el Ejemplo  \ref{ejemploRotacionEsfera}, en que toda la esfera $S^2$ es un atractor topol\'{o}gico, y por lo tanto un atractor de Milnor, pero no existen atractores erg\'{o}dicos.

\begin{definition} \index{atractor! de Milnor! $\alpha$-obs minimal} \index{conjunto! minimal $\alpha$-obs} \index{atractor! $\alpha$-observable}
{\bf $\alpha$-obs. minimalidad de un atractor de Milnor   \cite{CatIlyshenkoAttractors}} \em Sea dado un n\'{u}mero real $0 <\alpha \leq 1$. Un atractor de Milnor $K$ se dice    \em   $\alpha$-observable \em (escribimos \lq\lq \em $K$ es $\alpha$-obs\rq\rq.\em) si  su cuenca  $E_K$ de atracci\'{o}n (topol\'{o}gica) cumple
$$m(E_K) \geq \alpha.$$

Un atractor $K$ de Milnor $\alpha$-obs. se dice   \em $\alpha$-obs. minimal \em si no contiene subconjuntos compactos propios no vac\'{\i}os que sean atractores de Milnor $\alpha$-obs. para el mismo valor de $\alpha$.

En particular, cuando $\alpha= 1$, tenemos definidos los atractores de Milnor $1$-observables y $1$-observables minimales.

Se observa que todo atractor de Milnor $1$-observable es $\alpha$-observable para cualquier $0 < \alpha \leq 1$. Pero un atractor de Milnor $1$-obs. minimal no tiene por qu\'{e} ser $\alpha$-obs. minimal para todo $0 < \alpha < 1$ (ver Ejemplos \ref{ejemploMilnorNoTopologico} y \ref{ejemploMilnorNoTopologico1}).
\end{definition}


\begin{theorem} {\bf Existencia de atractores de Milnor} \index{teorema! de existencia de! atractor de Milnor} \index{atractor! de Milnor! $\alpha$-obs minimal}
\index{conjunto! minimal $\alpha$-obs} \label{TeoExistenciaAtrMilnor}

Sea $f: M \mapsto M$ continua en una variedad compacta y riemanniana $M$, de dimensi\'{o}n finita. Sea $0 < \alpha \leq 1$ dado. Entonces existen atractores de Milnor $\alpha$-obs. minimales para $f$. Adem\'{a}s, si $\alpha= 1$, el atractor de Milnor $1$-obs. minimal es \'{u}nico.
\end{theorem}

La prueba del Teorema \ref{TeoExistenciaAtrMilnor}, en una versi\'{o}n m\'{a}s restringida que   enuncia solo la existencia y unicidad del   atractor de Milnor $1$-observable minimal, fue dada primeramente en (\cite{Milnor}). En el ap\'{e}ndice de \cite{CatIlyshenkoAttractors} se define $\alpha$-obs. minimalidad de atractores de Milnor para $0 < \alpha \leq 1$, y se observa que la demostraci\'{o}n de existencia de atractor  de Milnor dada en (\cite{Milnor}) puede aplicarse, con una inmediata adaptaci\'{o}n, para probar la existencia de atractores de Milnor $\alpha$-obs. minimales, para cualquier $0 < \alpha \leq 1$.

{\em Demostraci\'{o}n: } {\em del Teorema }\ref{TeoExistenciaAtrMilnor}: Sea ${\aleph}_{\alpha}$ la familia de los atractores de Milnor $\alpha$-obs. (no necesariamente minimales). Esta familia no es vac\'{\i}a pues, trivialmente, $M$ es un atractor de Milnor $1$-obs.  En $\aleph_{\alpha}$ consideramos la relaci\'{o}n de orden parcial $K_1 \subset K_2$. Entonces, las cuencas de atracci\'{o}n topol\'{o}gica $E_{K_1}$ y $E_{K_2}$ cumplen $$E_{K_1} \subset E_{K_2}, \ \ \alpha \leq m(E_{K_1}) \leq m(E_{K_2}).$$
Sea en $\aleph_{\alpha}$ una cadena $\{K_i\}_{i \in I}$ (no necesariamente numerable). Es decir, $\{K_i\}_{i \in I}$ es un subconjunto totalmente ordenado de $\aleph_{\alpha}$, con la relaci\'{o}n de orden $\subset$.

Probemos que:

 {\bf Afirmaci\'{o}n (i)} (A probar) \em Existe en $\aleph_{\alpha}$ un elemento $K$ minimal de la cadena $\{K_i\}_{i \in I}$. \em Es decir, probemos que existe $K \in \aleph_{\alpha}, \ K \subset K_i$ para todo $i \in I$.

 En efecto, el conjunto $K= \bigcap_{i \in I} K_i$ es compacto no vac\'{\i}o, porque cualquier subcolecci\'{o}n finita $K_{1} \supset K_2 \supset \ldots \supset K_l$ de la cadena dada $\{K_i\}_{i \in I}$, tiene intersecci\'{o}n $K_1$ que es un compacto no vac\'{\i}o. Para probar que $K \in  {\aleph}_{\alpha}$, debemos probar ahora que $m(E_K) \geq \alpha$. Sea $j \in \mathbb{N}^+$ y sea $V_j \supset K$ el abierto formado por todos los puntos de $M$ que distan de $K$ menos que $1/j$. Afirmamos que existe \begin{equation}
 \label{eqn-25}    K_{i_j} \subset V_j \mbox{ para alg\'{u}n } i_j \in I. \end{equation} Por absurdo, supongamos que para cierto $j \in \mathbb{N}^+$ fijo, para todo $i \in I$, el compacto $K_i  \setminus V_j)$ es no vac\'{\i}o. Luego, como $\{K_i\}_{i \in \aleph}$ es totalmente ordenado con la relaci\'{o}n de orden $\subset$, obtenemos que la familia de compactos $\{K_i  \setminus V_j \}_{i \in I}$ es totalmente ordenada. Argumentando como hicimos m\'{a}s arriba, pero con esta nueva familia totalmente ordenada de compactos, en vez de con la familia $\{K_i\}_{i \in I}$, deducimos que el siguiente   compacto es no vac\'{\i}o $$\bigcap_{i \in I} (K_i  \setminus V_j) = \Big(\bigcap _{i \in I} K_i   \Big)  \setminus V_j  = K \setminus V_j,$$
contradiciendo que $K \subset V_j$. Hemos probado la afirmaci\'{o}n (\ref{eqn-25}).

 Como todo punto $x \in E_{K_{i_j}}$ cumple $\lim_{n} \mbox{dist}(f^n(x), K_{i_j}) = 0$, entonces  $f^{n}(x) \in V_j$ para todo $n$ suficientemente grande (que depende del punto $x$). Este argumento vale para cualquier $j \in \mathbb{N}^+$. Concluimos  que todo punto en $\bigcap _{j \in {\mathbb{N}^+}} E_{K_{i_j}}$ pertenece a $E_K$. Rec\'{\i}procamente, todo punto de $E_{K}$, por definici\'{o}n de la cuenca de atracci\'{o}n topol\'{o}gica, est\'{a} contenido en $E_{K_i}$ para todo $i \in I$ (porque $K \subset  {K_i}$). En particular, esta afirmaci\'{o}n se satisface para $i_j$, para todo $j \in \mathbb{N}^+$. Luego:
$$E_K := \bigcap _{j \in {\mathcal N}^+} E_{K_{i_j}}.$$
 Como la colecci\'{o}n numerable $E_{K_{i_j}}$ est\'{a} totalmente ordenada,  obtenemos $$m(E_K) = m \big(\bigcap_{j \in {\mathbb{N}^+}} E_{K_{i_j}}\big) = \lim_{j \rightarrow + \infty} m(E_{K_{i_j}}) \geq \alpha,$$
 terminando de demostrar la afirmaci\'{o}n (i).

 De la afirmaci\'{o}n (i) se deduce que para toda cadena en $\aleph_{\alpha}$ existe alg\'{u}n elemento $K \in \aleph_{\alpha}$ minimal de la cadena. Aplicando el Lema de Zorn, existen en $\aleph_{\alpha}$ elementos minimales de $\aleph_{\alpha}$. Es decir, existe $K \in \aleph_{\alpha}$ que no contiene subconjuntos propios que pertenezcan a $\aleph_{\alpha}$. Esto es, existe $K$ atractor de Milnor $\alpha$-obs. minimal.

 Ahora probemos la unicidad  del atractor de Milnor $1$-obs. minimal. Si existieran dos atractores de Milnor $K_1$ y $K_2$ que fueran $1$-obs. minimales, entonces la intersecci\'{o}n $E$ de sus cuencas de atracci\'{o}n topol\'{o}gica $E:= E_{K_1} \cap E_{K_2} $ cumple $$m(E)= 1,$$
 porque $m(E_{K_1})= m(E_{K_2}) = 1$. Todo punto $x \in E$ verifica, por la definici\'{o}n de la cuenca $E_{K_1}$ y $E_{K_2}$, la siguiente propiedad:

  \em Para todo $\epsilon >0$ existe $N \geq 1$ tal que \em
  $$\mbox{dist}(f^n(x), K_1), \ \ \mbox{dist}(f_n(x), K_2) < \epsilon \ \ \forall \ n \geq N$$
  Luego $\mbox{dist}(K_1, K_2) < 2 \epsilon \ \ \forall \ \epsilon >0$, de donde $K:= K_1 \bigcap K_2 \neq \emptyset$. Es est\'{a}ndar chequear que, siendo $K_1$ y $K_2$ compactos no vac\'{\i}os con intersecci\'{o}n no vac\'{\i}a, si un punto $x$ cumple
  $$\lim_{n \rightarrow + \infty} \mbox{dist}(f^n(x), K_1)=\lim_{n \rightarrow + \infty} \mbox{dist}(f^n(x), K_1)= 0, $$
  entonces $$\lim_{n \rightarrow + \infty} \mbox{dist}(f^n(x), K_1 \cap K_2)=0.$$
  (Chequear esta \'{u}ltima afirmaci\'{o}n en la parte (a) del Ejercicio \ref{exercise5}.)
  Luego $$E= E_{K_1} \bigcap E_{K_2} \subset E_{K_1 \cap K_2} $$ y como $m(E)= 1$, deducimos que $m(E_{K_1 \cap K_2})= 1$. Entonces $K_1 \cap K_2$ es un atractor de Milnor $1$-obs. Como $K_1$ y $K_2$ eran atractores de Milnor $1$-obs. minimales, concluimos que $K_1 \cap K_2 = K_1 = K_2$, terminando de demostrar la unicidad del atractor de Milnor 1-obs. minimal, y el Teorema \ref{TeoExistenciaAtrMilnor}.
\hfill $\Box$
\begin{exercise}\em \label{exercise5}
{\bf (a)} Demostrar que si $K_1$ y $K_2$ son compactos no vac\'{\i}os, y si existe una sucesi\'{o}n de puntos $x_n \in M$ tal que
$$\lim_{n \rightarrow + \infty} \mbox{dist}(x_n, K_1)= \lim_{n \rightarrow + \infty} \mbox{dist}(x_n, K_2)= 0,  $$
entonces $K_1 \cap K_2 \neq \emptyset$ y $$\lim_{n \rightarrow + \infty} \mbox{dist}(x_n, K_1 \cap K_2) = 0.$$

{\bf (b)} Demostrar que si $K_1$ y $K_2$ son dos atractores de Milnor tales $m(E_{K_1} \cap E_{K_2}) >0$, entonces $K_1 \cap K_2$ es no vac\'{\i}o, y es un atractor de Milnor cuya cuenca de atracci\'{o}n topol\'{o}gica es $$E_{K_1 \cap K_2} = E_{K_1} \cap E_{K_2}.$$
\end{exercise}

\begin{example} \em
\label{ejemploMilnorNoTopologico} Sea $X=  [0, 4 \pi] \times [0,1] $ el rect\'{a}ngulo compacto de ancho $4 \pi$ y altura 1, con un v\'{e}rtice en el origen, cuyo interior est\'{a} contenido en el primer cuadrante y tiene lados paralelos a los ejes. Sea en el intervalo $[0,4 \pi]$ la transformaci\'{o}n
$$f_1(x) = 1+ x - \cos x  \ \ \forall \ x \in [0, 4 \pi] $$ y sea en $X$ la transformaci\'{o}n
$$f(x,y) = (f_1(x), y/2) \ \ \forall \ (x,y) \in [0, 4 \pi] \times [0,1]. $$
Sea $m$ la medida de Lebesgue en $X$, re-escalada para que $$m(X)= 1.$$
Es est\'{a}ndar chequear (ver Ejercicio \ref{exercise6}), que
$$K_1 := \{(2 \pi, 0)\} $$ es un atractor de Milnor con cuenca de atracci\'{o}n topol\'{o}gica $$E_{K_1} = (0, 2 \pi] \times [0,1],$$
y   por lo tanto $K_1$ es $1/2$-obs. minimal como atractor de Milnor.
$$K_2 := \{(4 \pi, 0)\}$$ es otro atractor de Milnor con cuenca de atracci\'{o}n topol\'{o}gica $$E_{K_2} = (2 \pi, 4 \pi] \times [0,1],$$ y por lo tanto $K_2$ es tambi\'{e}n $1/2$-obs. minimal como atractor de Milnor.

  $K_1$ no es atractor topol\'{o}gico, porque todo entorno de $K_1$ contiene puntos de la cuenca de atracci\'{o}n topol\'{o}gica de $K_2$.

 $K_2$  es atractor topol\'{o}gico.

 $K_1 \cup K_2$ es el \'{u}nico atractor de Milnor $1$-obs. minimal, y tambi\'{e}n es el \'{u}nico atractor $\alpha$-obs, si $1/2 < \alpha \leq 1$.

 $K_1$ y $K_2$ son los \'{u}nicos atractores de Milnor $\alpha$-obs. si $0 < \alpha \leq 1/2$ (y como contienen un solo punto cada uno, son   $\alpha$-obs. minimales).

 $K_1 \cup K_2$ es atractor topol\'{o}gico, pero no es minimal como atractor topol\'{o}gico, pues contiene a $K_2$ que tambi\'{e}n es atractor topol\'{o}gico.
\end{example}
\begin{exercise}\em
\label{exercise6} Probar todas las afirmaciones del Ejemplo \ref{ejemploMilnorNoTopologico}.
\end{exercise}
\begin{example} \em
\label{ejemploMilnorNoTopologico1}
Sea en el cuadrado $X = [0,1] ^2$ la transformaci\'{o}n
$$f(x,y) = (x, (1/2) y) \ \ \forall \ (x,y) \in [0,1]^2.$$
Es inmediato chequear que todos los puntos del segmento $[0,1] \times\{0\}$ son puntos fijos y que $\lim_{n \rightarrow + \infty}f^n(x,y)  = (x,0) \ \ \forall \ (x,y) \in [0,1]^2$. Entonces, para todo $0 < \alpha \leq 1$, cualquier segmento compacto $I \times \{0\}$ tal que la longitud del intervalo compacto $I \subset [0,1]$ sea exactamente $\alpha$, es un atractor de Milnor $\alpha$-observable minimal. Sea $K = (I_1 \cup I_2 \cup \ldots \cup I_k) \times \{0\}$ donde $I_k \subset [0,1]$ es un intervalo compacto con interior no vac\'{\i}o, la colecci\'{o}n de los intervalos $I_k$ es disjunta dos a dos  y tal que la suma de las longitudes de los $I_k$ es exactamente $\alpha$. Entonces, $K$ tambi\'{e}n es un atractor de Milnor $\alpha$-obs. minimal. Adem\'{a}s, si $0 <\alpha < 1$ y si $I \subset [0,1]$ es un conjunto de Cantor con medida de Lebesgue  igual a $\alpha $ (tales conjuntos de Cantor siempre existen), entonces $I \times \{0\}$ es tambi\'{e}n un atractor de Milnor $\alpha$-obs. minimal. Ninguno de los atractores $\alpha$-obs. minimales   construidos anteriormente, si $0 <\alpha < 1$, es atractor topol\'{o}gico, pues las cuencas de atracci\'{o}n topol\'{o}gica no contienen ning\'{u}n entorno del atractor. El \'{u}nico atractor $1$-observable minimal es el intervalo $[0,1] \times \{0\}$, que s\'{\i} es atractor topol\'{o}gico.
\end{example}

\subsection{Atractores estad\'{\i}sticos o de Ilyashenko}

En esta secci\'{o}n, al igual que en la anterior, $f: M \mapsto M$ es una transformaci\'{o}n continua en una variedad compacta riemanniana $M$ de dimensi\'{o}n finita. Denotamos con $m$ a la medida de Lebesgue en $M$, re-escalada para que sea una medida de probabilidad: $m(M)= 1$, i.e. si $0 <m(M) \neq 1$, sustituimos $m$ por la probabilidad $m/m(M)$.

\begin{definition}
\label{definitionAtractorIlyashenko} \em
{\bf Atractor estad\'{\i}stico o de Ilyashenko \cite{Ilyashenko}, \cite{Ilyashenkogorodetski}} \index{atractor! de Ilyashenko}
\index{cuenca de atracci\'{o}n! estad\'{\i}stica}
\index{cuenca de atracci\'{o}n! de Ilyashenko} \index{$A_K$ cuenca de atracci\'{o}n! estad\'{\i}stica del compacto $K$}
\index{atractor! estad\'{\i}stico}
Se llama \em atractor   estad\'{\i}stico o de Ilyashenko \em a un conjunto compacto no vac\'{\i}o $K \subset M$,  invariante por $f$ (i.e. $f^{-1}(K) = K$), tal que $$m(A_K) = 1,$$ donde el conjunto $A_K \subset M$, llamado \em cuenca de atracci\'{o}n estad\'{\i}stica \em de $K$, est\'{a} definido por
\begin{equation} \label{eqnCuencaAtracIlyashenko}A_K  := \{x \in M: \ \ \lim_{n \rightarrow + \infty} \frac{1}{j}\sum_{j= 0}^{n-1}\mbox{dist}(f^j(x), K) = 0\}. \end{equation}
\end{definition}
\begin{exercise}\em
Probar que para cualquier conjunto compacto $K$ no vac\'{\i}o, el conjunto $A_K$ definido por la igualdad (\ref{eqnCuencaAtracMilnor}) es medible. Suge\-rencia:  La misma que para el Ejercicio \ref{ejercicio7}.
\end{exercise}

\begin{exercise}\em
Sean $K \subset K'$ dos atractores estad\'{\i}sticos. Probar que $A_K \subset A_{K'}$.
\end{exercise}

Es est\'{a}ndar chequear que todo atractor de Milnor (y en particular todo atractor topol\'{o}gico y todo atractor erg\'{o}dico)   es estad\'{\i}stico o de Ilyashenko (Ejercicio \ref{ejercicio8}). Sin embargo, no todo atractor estad\'{\i}stico   es de Milnor, como veremos en el Ejemplo  \ref{ejemploIlyashenkoNoTopologico}.  Adem\'{a}s, no todo    atractor estad\'{\i}stico es erg\'{o}dico, como muestran los Ejemplos  \ref{ejemploRotacionEsfera} y \ref{ejemploMilnorNoTopologico1}.

\begin{exercise}\em \label{ejercicio8} {\bf (a)} Sea $\{a_n\}_{n \in \mathbb{N}}$ una sucesi\'{o}n de reales que converge a $a \in \mathbb{R}$. Probar que la sucesi\'{o}n de promedios $\frac{1}{n} \sum_{j= 0}^{n-1} a_j$ converge a $a$.

{\bf (b)} Probar que todo atractor de Milnor $K$ es atractor estad\'{\i}stico o de Ilyashenko. Sugerencia: probar, usando la parte (a),  que la cuenca de atracci\'{o}n topol\'{o}gica de $K$ est\'{a} contenida en la cuenca de atracci\'{o}n  estad\'{\i}stica de $K$.

{\bf (c) } Probar que en los Ejemplos \ref{ejemploRotacionEsfera} y \ref{ejemploMilnorNoTopologico1}, los atractores de Milnor $K$ (que por la Parte (b) son tambi\'{e}n atractores estad\'{\i}sticos) no son atractores erg\'{o}dicos. Sugerencia: el \'{u}nico atractor de Milnor $K$ que cumple la condici\'{o}n (\ref{eqn28z}) de la Definici\'{o}n \ref{definitionAtractorErgodico} de atractor erg\'{o}dico, no cumple la condici\'{o}n (\ref{eqn30}) de esa definici\'{o}n.

{\bf (d)} Probar que en el Ejemplo \ref{ejemploMilnorNoTopologico},   los atractores de Milnor $K$ (que por la parte (b) tambi\'{e}n son atractores  estad\'{\i}sticos) no son atractores erg\'{o}dicos, aunque para dos de ellos se cumple la condici\'{o}n (\ref{eqn30}) de la Definici\'{o}n \ref{definitionAtractorErgodico} de atractor erg\'{o}dico.

\end{exercise}

\subsection{Atracci\'{o}n estad\'{\i}stica de un compacto}

En la Definici\'{o}n \ref{definitionAtractorIlyashenko} observamos que para un atractor   estad\'{\i}stico, el criterio de observabilidad de su cuenca es el criterio Lebesgue-medible (cf. Observaci\'{o}n  \ref{remarkObservabilidadTopyEstad}), mientras que la atracci\'{o}n es estad\'{\i}stica (cf. Definiciones \ref{definitionAtraccionTopol}
y \ref{definitionAtraccionEstad}). En efecto:

\begin{proposition}
\label{propositionAtraccionEstadistica1} {\bf Caracterizaci\'{o}n de la atracci\'{o}n estad\'{\i}stica} \index{cuenca de atracci\'{o}n! estad\'{\i}stica}
  \index{atractor! estad\'{\i}stico}

Sea $K$   un atractor estad\'{\i}stico o de Ilyashenko, y sea $A_K$   su cuenca de atracci\'{o}n estad\'{\i}stica. Entonces, para todo $x \in M$ se cumple
$x \in A_K$ si y solo si
toda medida de probabilidad $\mu$ que sea l\'{\i}mite de alguna subsucesi\'{o}n convergente en la topolog\'{\i}a d\'{e}bil estrella de las probabilidades emp\'{\i}ricas $\sigma_n:= \frac{1}{n} \sum_{j= 0}^{n-1} \delta_{f^j(x)}$ \em (cf. Definici\'{o}n \ref{definitionProbaEmpiricas}),  \em cumple
$$\mu(K) = 1.$$
\end{proposition}

 {\em Demostraci\'{o}n: }

{\em Demostraci\'{o}n del \lq\lq solo si\rq\rq: } Sea $x \in A_K$ y sea $$\mu = \lim_{i \rightarrow + \infty} \frac{1}{n} \sum_{j= 0}^{n-1} \delta_{f_j(x)},$$
para una cierta subsucesi\'{o}n $\{n_i\}_{i \in \mathbb{N}}$ de los naturales  tal que ese l\'{\i}mite $\mu$ existe en la topolog\'{\i}a d\'{e}bil estrella del espacio de probabilidades. Debemos probar que $\mu(K) = 1$.

Por la definici\'{o}n de la topolog\'{\i}a d\'{e}bil estrella, para toda funci\'{o}n continua $\psi: M \mapsto \mathbb{R}$ se cumple
$$\lim_{i \rightarrow + \infty} \frac{1}{n_i} \sum_{j= 0}^{n_i-1} \psi \circ f^j(x) = \int \psi \, d \mu.$$ En particular, la igualdad anterior se cumple para la funci\'{o}n continua $$\psi(x) := \mbox{dist}(x, K).$$ Por la Definici\'{o}n \ref{definitionAtractorIlyashenko}, como $x \in A_K$, existe $n \geq 0$ tal que
 $$\frac{1}{n} \sum_{j= 0}^{n-1} \mbox{dist}(f^j(x), K) < \epsilon$$
 En particular, para $n= n_i$ para todo $i$ suficientemente grande, tenemos: $$\frac{1}{n_i} \sum_{j= 0}^{n_i-1} \psi(f^j(x)) \geq \epsilon $$
 De las igualdades anteriores, deducimos que $\int \psi \, d \mu = 0$, es decir
 $$\int \mbox{dist}(x, K) \, d  \mu = 0$$
 y como $\mbox{dist}( \cdot, K)$  es una funci\'{o}n no negativa, deducimos que $\mbox{dist}(x,K) = 0$ para $\mu$-c.t.p. $x \in M$.
 Como $K$ es compacto, deducimos que $x \in K$ para $\mu$-c.t.p. $x \in M$, luego
 $$\mu(K) = \int \chi_K(x) \, d\mu = \int 1 \, d \mu = 1,$$
 como quer\'{\i}amos probar.

 {\em Demostraci\'{o}n del \lq\lq si\rq\rq: } Sea dado $x \in M$ tal que toda medida $\mu$ de probabilidad que sea el l\'{\i}mite d\'{e}bil estrella de una subsucesi\'{o}n convergente de las probabilidades emp\'{\i}ricas $\sigma_n(x)$, est\'{a} soportada en $K$. Tenemos que probar que  $\lim_{n \rightarrow + \infty} \frac{1}{n} \sum_{j= 0}^{n-1} \mbox{dist}(f^j(x), K) = 0$. Sea una subsucesi\'{o}n $n_i$ tal que existe el l\'{\i}mite $$\lim _{i \rightarrow + \infty} \frac{1}{n_i} \sum_{j= 0}^{n_i-1} \mbox{dist}(f^j(x), K)   = a.$$  No es restrictivo suponer que la subsucesi\'{o}n $\sigma_{n_i}(x)$ es convergente (en caso contrario, reemplazamos $\{n_i\}$ por una subsucesi\'{o}n adecuada de ella). Entonces, llamando $\mu$ al l\'{\i}mite de $\sigma_{n_i}$, e integrando   la funci\'{o}n continua $\psi = \mbox{dist}(\cdot, K)$ obtenemos:
 \begin{equation}
 \label{eqn-23}
   a = \lim_{i \rightarrow + \infty} \frac{1}{n_i} \sum_{j= 0}^{n_i-1} \mbox{dist}(f^j(x), K) = \lim_{i \rightarrow + \infty} \int \psi(y) \, d (\sigma_{n_i}(x))(y) = $$ $$= \int \psi \, d \mu = \int \mbox{dist}(y,K)\, d \mu(y)  .\end{equation}

 Pero la \'{u}ltima integral    en (\ref{eqn-23}) es igual a cero, pues  de  la hip\'{o}tesis $\mu(K)= 1$ deducimos que $y \in K$ para $\mu$-c.t.p.  Luego $\mbox{dist}(y,K) = 0$ para $\mu$-c.t.p. $y \in M$. Entonces, la igualdad (\ref{eqn-23}) implica que

 $\lim_n\frac{1}{n} \mbox{dist}\sum_{j= 0}^{n-1}(f^j(x), K)= 0 $ (pues el l\'{\i}mite de cualquier subsucesi\'{o}n convergente   es $a= 0$). Entonces $x \in A_K$,  terminando  de demostrar la Proposici\'{o}n \ref{propositionAtraccionEstadistica1}.
 \hfill $\Box$

\begin{definition}
\label{definitionSoporteDemu} \em Sea $\mu$ una medida de probabilidad en los borelianos de $M$. Se llama \em soporte compacto de $\mu$ \em al m\'{\i}nimo compacto $K \subset M$ tal que $\mu(K) = 1$. El soporte compacto existe, debido al Lema de Zorn y a la propiedad de que cualquier familia de compactos tiene intersecci\'{o}n   no vac\'{\i}a si las intersecciones de todas las subfamilias finitas son no vac\'{\i}as. El soporte compacto de $\mu$ es \'{u}nico, pues si existieran dos diferentes $K_1$ y $ K_2$, entonces $\mu(K_1 \cap K_2)= 1$ y  ni $K_1 $ ni $ K_2$ tendr\'{\i}an la propiedad de minimalidad.
\end{definition}


\begin{proposition}
\label{propositionSRB->estadistico} {\bf Medidas SRB y atractores estad\'{\i}sticos} \index{medida! SRB} \index{atractor! estad\'{\i}stico}

Sea $f: M \mapsto M$ continua en la variedad  compacta $M$. Si existe una medida $\mu$ que es SRB para $f$ \em (cf. Definici\'{o}n \ref{definitionMedidaSRB}), \em entonces el soporte compacto $K$ de $\mu$ es un atractor estad\'{\i}stico o de Ilyashenko.
\end{proposition}

\begin{remark} \em
\label{remarkAtracIlyashenkoNo->SRB} El rec\'{\i}proco de la Proposici\'{o}n \ref{propositionSRB->estadistico} es falso. En efecto, las medidas SRB no siempre existen (cf. Ejemplos \ref{ejemploRotacionEsfera} y \ref{ejemploMilnorNoTopologico1}), pero los atractores de Ilyashenko siempre existen, como demostraremos en el Teorema \ref{TeoExistenciaAtrIlyashenko}.

A pesar de que el rec\'{\i}proco de la Proposici\'{o}n \ref{propositionSRB->estadistico} es falso, se puede gene\-ralizar el enunciado de esta Proposici\'{o}n de forma que su rec\'{\i}proco es verdadero. M\'{a}s precisamente, un compacto $K$ no vac\'{\i}o e invariante es atractor estad\'{\i}stico o de Ilyashenko, si y solo s\'{\i}, es el m\'{\i}nimo soporte compacto de una colecci\'{o}n adecuada de medidas invariantes, que llamamos SRB-like, y que describen en forma \'{o}ptima la estad\'{\i}stica de las \'{o}rbitas de un conjunto de medida de Lebesgue positivo. En efecto, en la pr\'{o}xima secci\'{o}n definiremos las medidas de probabilidad   \lq\lq SRB-like\rq\rq \ o \lq\lq pseudo-f\'{\i}sicas\rq\rq, que incluyen a las medidas SRB o f\'{\i}sicas (cuando \'{e}stas existen), pero tambi\'{e}n pueden incluir a otras medidas de probabilidad invariantes que no son necesariamente SRB. El nuevo enunciado generalizado de la Proposici\'{o}n \ref{propositionSRB->estadistico} utilizando las medidas SRB-like, en lugar de las medidas SRB, y su rec\'{\i}proco, ser\'{a} establecido y demostrado en el Teorema \ref{theoremSRB-like<->Ilyashenko}.
\end{remark}

 {\em Demostraci\'{o}n: }{\em de la Proposici\'{o}n } \ref{propositionSRB->estadistico}

 Para probar que $K$ es atractor estad\'{\i}stico, por la Definici\'{o}n \ref{definitionAtractorIlyashenko}, hay que probar que $m(A_K) >0$, donde $m$ denota la medida de Lebesgue y $A_K$ denota la cuenca de atracci\'{o}n estad\'{\i}stica de $K$ definida en (\ref{eqnCuencaAtracIlyashenko}). Como $\mu$ es medida SRB, por la Definici\'{o}n \ref{definitionMedidaSRB}, la cuenca de atracci\'{o}n estad\'{\i}stica $B(\mu)$ de $\mu$, dada en la Definici\'{o}n (\ref{DefinicionCuencaDeAtraccionEstadistica}), cumple $m(B(\mu)) >0$. Luego, basta probar que $B(\mu) \subset A_K$. Sea $x \in B(\mu)$. Entonces, por la Definici\'{o}n \ref{DefinicionCuencaDeAtraccionEstadistica}, para cualquier funci\'{o}n continua $\psi: M \mapsto \mathbb{R}$ se cumple:
 $$\lim_{n \rightarrow + \infty} \frac{1}{n} \sum_{j= 0}^{n-1} \psi \circ f^j(x) = \int \psi \, d\mu. $$
 En particular si tomamos $\psi(x) = \mbox{dist}(x, K)$ tenemos:
 $$\lim_{n \rightarrow + \infty} \frac{1}{n} \sum_{j= 0}^{n-1} \mbox{dist}(f^j(x), K)  = \int \mbox{dist}(y, K)\, d \mu$$
 Para concluir que $x \in A_K$ basta probar que la integral en la igualdad anterior es cero, y para ello basta chequear que $\mbox{dist}(y, K)= 0$ para $\mu$-c.t.p. $y \in M$. En efecto, $\mu$-c.t.p. $y \in M$ pertenece a $K$ porque $\mu(K)= 1$. Luego, hemos probado que $x \in A_K$ para todo $x \in B(\mu)$, y por lo tanto $m(A_K) \geq m(B(\mu)) >0$, y $K$ es un atractor de Ilyashenko.
 \hfill $\Box$

 \begin{definition}
 {\bf Tiempo medio de estad\'{\i}a} \label{definitionFrecuencia} \em
 \index{tiempo medio de estad\'{\i}a}
 Sea   $V \subset M$ un conjunto medible no vac\'{\i}o. Sea $x \in M$ cualquiera. Llamamos \em frecuencia     de visitas a  $V $ \em de la \'{o}rbita por $x$ hasta tiempo $n$, al siguiente n\'{u}mero $ \sigma_{n, x}(V)$ (cf. (\ref{(1)}) en la Definici\'{o}n \ref{definitionProbaEmpiricas} de las probabilidades emp\'{\i}ricas $\sigma_{n,x}$):
 \begin{equation}\sigma_{n,x}(V) := \int \chi_V \, d \sigma_{n,x} =    \frac{1}{n} \sum_{j= 0}^{n-1} \chi_{V}(f^j(x)) =$$ $$= \frac{\#\{0 \leq j \leq n-1\colon f^j(x) \in V \}}{n},\label{eqn-24}\end{equation}
 donde $\chi_V$ denota la funci\'{o}n caracter\'{\i}stica de $V$ y $\#A$ denota cantidad de elementos de un conjunto finito $A$.

  El \'{u}ltimo t\'{e}rmino de la igualdad (\ref{eqn-24}) indica que el tiempo medio de estad\'{\i}a en $V$ de la \'{o}rbita por $x$ hasta tiempo $n$ es la cantidad relativa de iterados del punto $x$ que caen dentro del conjunto $V$.
  Por el teorema erg\'{o}dico de Birkhoff, si $\mu$ es una medida invariante, entonces para $\mu$-c.t.p. $x \in M$ el tiempo medio de estad\'{\i}a en $V$ de la \'{o}rbita por $x$ hasta tiempo $n$, converge cuando $n \rightarrow + \infty$, a $\widetilde \chi_V$, cuyo valor esperado $\int \widetilde \chi_V \, d \mu$ es igual a $\mu(V)$. Adem\'{a}s si la medida $\mu$ es erg\'{o}dica, el tiempo medio de estad\'{\i}a en $V$ tiende a $\mu(V)$ para $\mu$-c.t.p. $x \in M$.

  Sin embargo, el enunciado del Teorema de Birkhoff y las propiedades de las medidas invariantes erg\'{o}dicas, no nos ser\'{a}n \'{u}tiles a los prop\'{o}sitos en esta secci\'{o}n. Nos interesa   considerar el tiempo medio de estad\'{\i}a en $V$ para Lebesgue  c.t.p., aunque la medida de Lebesgue no sea invariante con $f$, y aunque los puntos considerados $x$ no pertenezcan al soporte de ninguna medida invariante $\mu$.
 \end{definition}

 \begin{proposition}.
 \label{propositionAtraccionEstadisticaFrec}

 {\bf Caracterizaci\'{o}n de atracci\'{o}n estad\'{\i}stica por el tiempo medio de estad\'{\i}a en un entorno} \index{tiempo medio de estad\'{\i}a} \index{atractor! estad\'{\i}stico}
\index{cuenca de atracci\'{o}n! estad\'{\i}stica}

 Sea $K$ un compacto no vac\'{\i}o invariante con $f$. Para todo $\epsilon >0$ denotamos con $V (\epsilon)$ al conjunto de puntos en la variedad $M$ que distan de $K$ menos que $\epsilon$.

 Entonces, la cuenca de atracci\'{o}n estad\'{\i}stica $A_K$ del conjunto $K$\em (definida en la igualdad (\ref{eqnCuencaAtracIlyashenko}) de la definici\'{o}n \ref{definitionAtractorIlyashenko}) \em est\'{a} caracterizada por la siguiente igualdad:
 \begin{equation} \label{eqnCuencaAtracIlyashenkoFrec}A_K = \bigcap_{\epsilon>0} \{x \in M: \ \lim_{n \rightarrow + \infty} \sigma_{n,x} (V(\epsilon)) \ = \ 1\}.\end{equation}
 \end{proposition}
 {\bf Notas: } La igualdad (\ref{eqnCuencaAtracIlyashenkoFrec}) implica que $x \in A_K$ si y solo si, cuando $n$ es suficientemente grande, la probabilidad  (la frecuencia relativa) del suceso de encontrar al iterado $f^j(x)$, para alg\'{u}n $0 \leq j \leq n-1$, en un entorno $V(\epsilon)$ dado fijo  arbitrariamente peque\~{n}o del atractor de Ilyashenko $K$, es 1. Sin embargo, la atracci\'{o}n de la \'{o}rbita por $x$ al conjunto $K$, no es necesariamente atracci\'{o}n topol\'{o}gica. Para ser atracci\'{o}n topol\'{o}gica, debe cumplirse con total certeza (no solo con probabilidad cercana a 1, ni tampoco solo con probabilidad 1)  el suceso de encontrar  al iterado $f^j(x)$ en el entorno $   V(\epsilon)$, para   todo $j$ suficientemente grande.

  La igualdad (\ref{eqnCuencaAtracIlyashenkoFrec}) da por lo tanto, la siguiente interpretaci\'{o}n intuitiva del significado de la atracci\'{o}n estad\'{\i}stica a un atractor de Ilyashenko $K$. Las \'{o}rbitas en su cuenca de atracci\'{o}n estad\'{\i}stica $A_K$ se acercan al conjunto $K$ tanto como se desee cuando el iterado $n$ tiende a infinito, pero se admite que existan iterados excepcionales, con una frecuencia relativa despreciable para $n$ arbitrariamente grande, en que la \'{o}rbita \lq\lq se toma una excursi\'{o}n relativamente muy breve (de vacaciones)\rq\rq \ alej\'{a}ndose del atractor $K$ durante la excursi\'{o}n  (ver Ejemplo \ref{ejemploHuYoung}).

 {\em Demostraci\'{o}n: }
 {\em de la Proposici\'{o}n } \ref{propositionAtraccionEstadisticaFrec}: 
{\bf (i)} Probemos que si $x \in A_K$ entonces

$\lim_{n \rightarrow + \infty} \sigma_{n,x}(V(\epsilon))$ para todo $\epsilon >0$. Fijemos $\epsilon >0$. Por simplicidad, escribimos $V$ en lugar de $ V(\epsilon)$, ya que $\epsilon $ est\'{a} fijo.   Tenemos:
$$   {\epsilon \cdot \chi_{M \setminus V}(y) }  \leq  \mbox{dist}(y, K) \ \ \forall \ y \in M,$$
pues, cuando $y \in V$ el t\'{e}rmino de la izquierda es cero, y cuando $y \not \in V$, el t\'{e}rmino de la izquierda es $\epsilon \leq \mbox{dist}(y, K)$. Entonces:
$$0 \leq \epsilon \cdot \frac{1}{n}\sum_{j= 0}^{n-1}(1-  \chi_V(f^j(x)) \leq \frac{1}{n} \sum_{j= 0}^{n-1}\mbox{dist}(f^j(x), K).$$
Como $x \in A_K$ por hip\'{o}tesis, el l\'{\i}mite cuando $n \rightarrow + \infty$ del t\'{e}rmino a la derecha en la desigualdad anterior, es cero. Concluimos que
$$0= \epsilon \cdot \lim_{n \rightarrow + \infty}\frac{1}{n}\sum_{j= 0}^{n-1}(1-  \chi_V(f^j(x))  = \epsilon \cdot \lim_{n \rightarrow + \infty}  \Big(1 - \frac{1}{n} \sum_{j= 0}^{n-1} \chi_V(f^j(x)) \Big)  = $$ $$ = \epsilon \cdot \lim_{n \rightarrow + \infty}   \sigma_{n,x}(V).$$
Como $\epsilon >0$, concluimos que $\lim_{n \rightarrow + \infty} \sigma_{n,x}(V) = 0$, como quer\'{\i}amos demostrar.

{\bf (ii)} Probemos que si  $\lim_{n \rightarrow + \infty} \sigma_{n,x}(V(\epsilon))$ para todo $\epsilon >0$ entonces $x \in A_K$. Fijemos $\epsilon >0$ y por simplicidad denotemos $V= V(\epsilon)$.
Sea $$D:= \mbox{diam}(M) = \max\{\mbox{dist}(x,y)\colon x,y \in M\} >0.$$  Luego:
$$\mbox{dist}(y, K) \leq D \cdot \chi_{M \setminus V}(y) \ \ \forall \ y \not \in V, $$ $$\mbox{dist}(y,K) < \epsilon \cdot \chi_V(y) \ \ \forall \ y \in V. $$
Entonces $$\mbox{dist} (y,K) \leq \epsilon \cdot \chi_V(y) + D \cdot \chi_{M \setminus V}(y) \ \ \forall \ y \in M.$$
Consideramos $x$ que satisface las hip\'{o}tesis. De la \'{u}ltima desigualdad obtenemos:
$$0 \leq \frac{1}{n}  \sum_{j= 0}^{n-1} \mbox{dist}(f^j(x), K) \leq \epsilon \cdot \frac{1}{n} \sum_{j= 0}^{n-1}\chi_{V}(f^j(x)) + D \cdot \frac{1}{n} \sum_{j= 0}^{n-1} \chi_{M \setminus V}(f^j(x)) = $$ $$= \epsilon \cdot \sigma_{n,x}(V) + D \cdot (1- \sigma_{n,x}(V)). $$
Por hip\'{o}tesis $\lim_{n \rightarrow + \infty} \sigma_{n,x}(V) = 1$. Como los n\'{u}meros reales positivos $\epsilon$ y $D$ est\'{a}n fijos  (son independientes de $n$), deducimos que
$$  \lim_{n \rightarrow + \infty} \frac{1}{n} \sum_{j= 0}^{n-1} \mbox{dist}(f^j(x), K) = 0, $$
La \'{u}ltima igualdad significa, debido a la condici\'{o}n (\ref{eqnCuencaAtracIlyashenko}) de la Definici\'{o}n \ref{definitionAtractorIlyashenko}, que $x \in A_K$  como quer\'{\i}amos demostrar.
 \hfill $\Box$


\subsection{Existencia de atractores de Ilyashenko}

\begin{definition}.

{\bf $\alpha$-obs. minimalidad de un atractor estad\'{\i}stico   \cite{CatIlyshenkoAttractors}} \label{definitionAtractorIlyashenkoMinimal} \index{atractor! de Ilyashenko! $\alpha$-obs minimal} \index{atractor! $\alpha$-observable}
\index{atractor! estad\'{\i}stico! $\alpha$-obs minimal}
\index{conjunto! minimal $\alpha$-obs}

\em Sea dado un n\'{u}mero real $0 <\alpha \leq 1$. Un atractor estad\'{\i}stico o de Ilyashenko $K$ se dice    \em   $\alpha$-observable \em (escribimos \lq\lq \em $K$ es $\alpha$-obs\rq\rq.\em) si  su cuenca  $A_K$ de atracci\'{o}n estad\'{\i}stica cumple
$$m(A_K) \geq \alpha.$$

Un atractor $K$ estad\'{\i}stico o de Ilyashenko $\alpha$-obs. se dice   \em $\alpha$-obs. minimal \em si no contiene subconjuntos compactos propios no vac\'{\i}os que sean atractores estad\'{\i}sticos o de Ilyashenko $\alpha$-obs. para el mismo valor de $\alpha$.

En particular, cuando $\alpha= 1$, tenemos definidos los atractores estad\'{\i}sticos o de Ilyashenko $1$-observables y $1$-observables minimales.

Se observa que todo atractor estad\'{\i}stico o de Ilyashenko $1$-observable es $\alpha$-observable para cualquier $0 < \alpha \leq 1$. Pero un atractor estad\'{\i}stico o de Ilyashenko $1$-obs. minimal no tiene por qu\'{e} ser $\alpha$-obs. minimal para todo $0 < \alpha < 1$.
\end{definition}

\begin{remark} \em {\bf Sobre conjuntos minimales} \index{conjunto! minimal}

Recordemos la caracterizaci\'{o}n de los conjuntos minimales desde el punto de vista topol\'{o}gico, dada en la Definici\'{o}n \ref{definicionMinimalTopologico}: un conjunto   $K$ compacto, no vac\'{\i}o y $f$-invariante es minimal para $f$, desde el punto de vista topol\'{o}gico,  si   no contiene subconjuntos propios compactos no vac\'{\i}os que sean invariantes por $f$ hacia el futuro.

Observamos que no todo compacto $K$ minimal para $f$ desde el punto de vista topol\'{o}gico  es   un atractor de Ilyshenko $\alpha$-obs. minimal para alg\'{u}n $0 \leq \alpha \leq 1$ (ver Ejemplo \ref{ejemploHuYoung} o tambi\'{e}n \cite{klepstynMinStatAttr}).

Rec\'{\i}procamente, no todo atractor de Ilyshenko $\alpha$-obs. minimal es necesariamente un  compacto minimal para $f$ desde el punto de vista topol\'{o}gico (ver Ejemplo \ref{ejemploAnosovlineal}).

En \cite{Ilyashenkogorodetski} y \cite{klepstynMinStatAttr} se estudian relaciones entre minimalidad para $f$ desde el punto de vista topol\'{o}gico y los atractores estad\'{\i}sticos o de Ilyashenko de $f$.
\end{remark}

\begin{theorem} {\bf Existencia de atractores de Ilyashenko} \label{TeoExistenciaAtrIlyashenko} \index{atractor! de Ilyashenko}
\index{atractor! estad\'{\i}stico}
\index{atractor! estad\'{\i}stico! $\alpha$-obs minimal}
\index{teorema! de existencia de! atractor estad\'{\i}stico}
\index{teorema! de existencia de! atractor de Ilyashenko}

Sea $f: M \mapsto M$ continua en una variedad compacta y riemanniana $M$, de dimensi\'{o}n finita. Sea $0 < \alpha \leq 1$ dado. Entonces existen atractores estad\'{\i}sticos o de Ilyashenko $\alpha$-obs. minimales para $f$. Adem\'{a}s, si $\alpha= 1$, el atractor  de Ilyashenko $1$-obs. minimal es \'{u}nico.
\end{theorem}

 La prueba del Teorema \ref{TeoExistenciaAtrIlyashenko} sigue los mismos argumentos de la prueba del Teorema \ref{TeoExistenciaAtrMilnor} de existencia de atractores de Milnor, con leves adaptaciones.

{\em Demostraci\'{o}n: }  Sea ${\aleph}_{\alpha}$ la familia de los atractores de Ilyashenko $\alpha$-obs. (no necesariamente minimales). Esta familia no es vac\'{\i}a pues, trivialmente, $M$ es un atractor de Ilyashenko $1$-obs, y por lo tanto es $\alpha$-obs. para cualquier $0 < \alpha \leq 1$.  En $\aleph_{\alpha}$ consideramos la relaci\'{o}n de orden parcial $K_1 \subset K_2$. De la condici\'{o}n (\ref{eqnCuencaAtracIlyashenko}), teniendo en cuenta que $\mbox{dist}(y, K_2) \leq \mbox{dist}(y, K_1)$ para todo $y \in M$, las cuencas de atracci\'{o}n estad\'{\i}stica $A_{K_1}$ y $A_{K_2}$ cumplen $$A_{K_1} \subset A_{K_2}, \ \ \alpha \leq m(A_{K_1}) \leq m(A_{K_2}).$$
Sea en $\aleph_{\alpha}$ una cadena $\{K_i\}_{i \in I}$ (no necesariamente numerable). Es decir, $\{K_i\}_{i \in I}$ es un subconjunto totalmente ordenado de $\aleph_{\alpha}$, con la relaci\'{o}n de orden $\subset$.

Probemos que:

 {\bf Afirmaci\'{o}n (i)} (A probar) \em Existe en $\aleph_{\alpha}$ un elemento $K$ minimal de la cadena $\{K_i\}_{i \in I}$. \em Es decir, probemos que existe $K \in \aleph_{\alpha}, \ K \subset K_i$ para todo $i \in I$.

 En efecto, el conjunto $K= \bigcap_{i \in I} K_i$ es compacto no vac\'{\i}o, porque cualquier subcolecci\'{o}n finita $K_{1} \supset K_2 \supset \ldots \supset K_l$ de la cadena dada $\{K_i\}_{i \in I}$, tiene intersecci\'{o}n $K_1$ que es un compacto no vac\'{\i}o. Para probar que $K \in  {\aleph}_{\alpha}$, debemos probar ahora que $m(A_K) \geq \alpha$. Sea $j \in \mathbb{N}^+$ y sea $V_j \supset K$ el abierto formado por todos los puntos de $M$ que distan de $K$ menos que $1/j$. En la prueba del Teorema \ref{TeoExistenciaAtrMilnor} demostramos la afirmaci\'{o}n (\ref{eqn-25}):   $$K_{i_j} \subset V_j \mbox{ para alg\'{u}n } i_j \in I.$$ Como todo punto $x \in A_{K_{i_j}}$ se cumple $\lim_{n} \sigma_{n,x}(V_j) = 0$ (cf. Proposici\'{o}n \ref{propositionAtraccionEstadisticaFrec}), y este argumento vale para todo $j \in \mathbb{N}^+$, entonces    todo punto en $\bigcap _{j \in {\mathbb{N}^+}} A_{K_{i_j}}$ pertenece a $A_K$. Rec\'{\i}procamente, todo punto de $A_{K}$  est\'{a} contenido en $A_{K_i}$ para todo $i \in I$ (porque $K \subset  {K_i}$). En particular esta afirmaci\'{o}n se satisface para $i_j$, para todo $j \in \mathbb{N}^+$. Luego:
$$A_K := \bigcap _{j \in {\mathcal N}^+} A_{K_{i_j}}.$$
 Como la colecci\'{o}n numerable $A_{K_{i_j}}$ est\'{a} totalmente ordenada,  obtenemos $$m(A_K) = m \big(\bigcap_{j \in {\mathbb{N}^+}} A_{K_{i_j}}\big) = \lim_{j \rightarrow + \infty} m(A_{K_{i_j}}) \geq \alpha,$$
 terminando de demostrar la afirmaci\'{o}n (i).

 De la afirmaci\'{o}n (i) se deduce que para toda cadena en $\aleph_{\alpha}$ existe alg\'{u}n elemento $K \in \aleph_{\alpha}$ minimal de la cadena. Debido al   Lema de Zorn, existen en $\aleph_{\alpha}$ elementos minimales de $\aleph_{\alpha}$. Es decir, existe $K \in \aleph_{\alpha}$ que no contiene subconjuntos propios que pertenezcan a $\aleph_{\alpha}$. Esto es, existe $K$ atractor de Ilyashenko $\alpha$-obs. minimal.

 Ahora probemos la unicidad  del atractor de Ilyashenko $1$-obs. minimal. Si existieran dos atractores de Ilyashenko $K_1$ y $K_2$ que fueran $1$-obs. minimales, entonces la intersecci\'{o}n $A$ de sus cuencas de atracci\'{o}n estad\'{\i}stica $A:= A_{K_1} \cap A_{K_2} $ cumple $$m(A)= 1,$$
 porque $m(A_{K_1})= m(A_{K_2}) = 1$. Todo punto $x \in A$ verifica, por la Proposici\'{o}n \ref{propositionAtraccionEstadisticaFrec} que caracteriza  las cuencas $A_{K_1}$ y $A_{K_2}$, la siguiente propiedad:

  \em Para todo $0 <\epsilon <1/2$ existe $N \geq 1$ tal que \em
  $$\sigma_{n,x}(V_{\epsilon}(K_1)), \ \ \sigma_{n,x}(V_{\epsilon}(K_1)) > 1- (\epsilon/2) \ \ \forall \ n \geq N,$$
  donde $V_{\epsilon}(K_i)$ denota el conjunto de puntos de la variedad $M$ que distan del compacto $K_i$ menos que $\epsilon$, y $\sigma_{n,x}$ denota la probabilidad emp\'{\i}rica definida en la igualdad (\ref{(1)}), Definici\'{o}n \ref{definitionProbaEmpiricas}.

  Luego existen m\'{a}s de $N: = n \cdot (1- (\epsilon/2))\geq (3/4) n$ iterados $f^j(x)$ con tiempos $j \in \{0, 1, \ldots, n-1\}$ tales que $f^j(x) \in V_{\epsilon}(K_1)$. En efecto, si la cantidad de   iterados $f^j(x)$ con tiempos $j \in \{0, 1, \ldots, n-1\}$ fuera  menor o igual que $N $, entonces tendr\'{\i}amos $\sigma_{n,x}(V_{\epsilon}(K_1) \leq N/n  = 1- (\epsilon/2) $. An\'{a}logamente, existen m\'{a}s de   $N$ iterados $f^h(x)$ con tiempos $h \in \{0, 1, \ldots, n-1\}$ tales que $f^h(x) \in V_{\epsilon}(K_1)$. Como la cantidad de iterados posibles con tiempos en $\{0,1, \ldots, n-1\}$ es a lo sumo $n$, y como $2N \geq (3/2) n > n$, deben existir algunos iterados comunes $f^j(x)= f^h(x) \in V_{\epsilon}(K_1) \cap V_{\epsilon}(K_2)$. Tomemos alguno de esos iterados comunes $f^j(x)$. Tenemos
  $$\mbox{dist}(K_1, K_2) \leq \mbox{dist}(f^j(x), K_1) + \mbox{dist}(f^j(x), K_2) < 2 \epsilon. $$ Entonces $ \mbox{dist}(K_1, K_2) < 2 \epsilon \ \ \forall \ 0 <\epsilon < 1/2$, de donde $$K:= K_1 \bigcap K_2 \neq \emptyset.$$ Aunque no es inmediato, es est\'{a}ndar chequear que, siendo $K_1$ y $K_2$ compactos no vac\'{\i}os con intersecci\'{o}n no vac\'{\i}a, si un punto $x$ cumple
  $$\lim_{n \rightarrow + \infty} \sigma_{n,x} V_{\epsilon}(K_1)=\lim_{n \rightarrow + \infty} \sigma_{n,x} V_{\epsilon}(K_2)= 0 \ \ \ \forall \ \epsilon >0, $$
  entonces $$\lim_{n \rightarrow + \infty} \sigma_{n,x} V_{\epsilon}(K_1 \cap K_2)=0 \ \ \forall \ \epsilon >0.$$
  (Chequear esta \'{u}ltima afirmaci\'{o}n en la parte (a) del Ejercicio \ref{exercise5Ilyashenko}.)
  Luego $$A:= A_{K_1} \bigcap A_{K_2} \subset A_{K_1 \cap K_2} $$ y como $m(A)= 1$, deducimos que $m(A_{K_1 \cap K_2})= 1$. Entonces $K_1 \cap K_2$ es un atractor de Ilyashenko $1$-obs. Como $K_1$ y $K_2$ eran atractores de Ilyashenko $1$-obs. minimales, concluimos que $K_1 \cap K_2 = K_1 = K_2$, terminando de demostrar la unicidad del atractor de Ilyashenko 1-obs. minimal, y el Teorema \ref{TeoExistenciaAtrIlyashenko}.
\hfill $\Box$
\begin{exercise}\em \label{exercise5Ilyashenko} Para un conjunto compacto no vac\'{\i}o $K \subset M$, y para   $\epsilon >0$, denotamos $V_{\epsilon}(K) $ al conjunto de puntos de $M$ que distan de $K$ menos que $\epsilon$. Denotamos $\sigma_{n,x}$ las probabilidades emp\'{\i}ricas de la \'{o}rbita con estado inicial $x$ hasta tiempo $n$, seg\'{u}n la igualdad (\ref{(1)}), dada por la Definici\'{o}n \ref{definitionProbaEmpiricas}.

{\bf (a)} Demostrar que si $K_1$ y $K_2$ son compactos no vac\'{\i}os, y si existe un punto  $x  \in M$ tal que
$$\lim_{n \rightarrow + \infty} \sigma_{n,x} (V_{\epsilon} (K_1))= \lim_{n \rightarrow + \infty} \sigma_{n,x} (V_{\epsilon} (K_1))= 0  \ \ \forall \ \epsilon >0  $$
entonces $K_1 \cap K_2 \neq \emptyset$ y $$\lim_{n \rightarrow + \infty} \sigma_{n,x} (V_{\epsilon} (K_1 \cap K_2)) = 0 \ \ \forall \ \epsilon >0$$

{\bf (b)} Demostrar que si $K_1$ y $K_2$ son dos atractores de Ilyashenko tales $m(A_{K_1} \cap A_{K_2}) >0$, entonces $K_1 \cap K_2$ es no vac\'{\i}o, y es un atractor de Ilyashenko cuya cuenca de atracci\'{o}n estad\'{\i}stica es $$A_{K_1 \cap K_2} = A_{K_1} \cap A_{K_2}.$$
\end{exercise}

\subsection{Ejemplos de atractores estad\'{\i}sticos}

\begin{example} \em   \index{automorfismo! lineal en el toro}
\label{ejemploAnosovlineal} Sea $f_0 \in \mbox{Diff }^{\infty}(\mathbb{T}^2)$ el mapa \lq\lq Arnold's cat \rq\rq \ \ en el toro bidimensional dado en el ejemplo de la Secci\'{o}n \ref{section2111}, como el automorfismo lineal hiperb\'{o}lico con matriz asociada $$ \Big(
                                                    \begin{array}{cc}
                                                      2 & 1 \\
                                                      1 &1 \\
                                                    \end{array}
                                                  \Big)  \Big|_{\mbox{\footnotesize m\'{o}d.} (1,1)}.
 $$  Recordando lo expuesto en la Secci\'{o}n \ref{section2111}: $f_0$ es un difeomorfismo de Anosov lineal, tiene a $(0,0)$ como punto fijo y preserva la medida de Lebesgue $m$ (re-escalada para que $m(\mathbb{T}^2)= 1$).

 Probemos que el \'{u}nico   atractor estad\'{\i}stico   (en particular el \'{u}nico    atractor estad\'{\i}stico $\alpha$-obs. minimal para cualquier $0 <\alpha \leq 1$), es $K= \mathbb{T}^2$. En efecto, en el Corolario \ref{corollarySRBanosov} probamos que $m$ es erg\'{o}dica para $f$. Entonces,     $m$-c.t.p. $x \in M$ est\'{a} en la cuenca de atracci\'{o}n estad\'{\i}stica $B(m)$ de $m$. (Recordar que una medida de probabilidad invariante $\mu$ es erg\'{o}dica si y solo si $\mu(B(\mu)) = 1$).

 Sea $K$ un atractor estad\'{\i}stico. Por definici\'{o}n, su cuenca de atracci\'{o}n estad\'{\i}stica $A_K$ tiene medida de Lebesgue positiva. Luego $$x \in B(m) \ \mbox{ para } m-\mbox{c.t.p. } x \in A_K.$$

 Tomemos y dejemos fijo $x \in B(m) \cap A_K$. Como  $x \in B(m)$ tenemos:
 $$\lim_{n \rightarrow + \infty} \frac{1}{n} \sum_{j= 0}^{n-1} \delta_{f^j(x)} = m,$$
 en la topolog\'{\i}a d\'{e}bil$^*$ del espacio ${\mathcal M}$ de las probabilidades. Entonces para toda funci\'{o}n continua $\psi$
 se cumple
 $$\lim_{n \rightarrow + \infty} \frac{1}{n} \sum_{j= 0}^{n-1} \psi \circ f^j(x) = \int \psi \, d m. $$
 En particular si elegimos la funci\'{o}n continua $\psi$ definida por
 $$\psi(x) = \mbox{dist}(x, K) \ \ \ \forall \ x \in M,$$
 se cumple
 $$\lim_{n \rightarrow + \infty} \frac{1}{n} \sum_{j= 0}^{n-1}\mbox{dist}(f^j(x), K) = \int \psi \, d m$$
 Como $x \in A_K$, por la Definici\'{o}n \ref{definitionAtractorIlyashenko}, el l\'{\i}mite de la izquierda es cero. Entonces deducimos que
 $$\int \psi \, d m = 0.$$
 Pero $\psi \geq 0$. Entonces su integral da cero si y solo si $\psi = 0$ para $m$-c.t.p. Una funci\'{o}n continua que es cero para Lebesgue casi todo punto, es id\'{e}nticamente nula. Deducimos que $$0 =\psi(x)=\mbox{dist}(x, K) \ \ \ \forall \ x \in \mathbb{T}^2.$$
 Como $K$ es compacto, $\mbox{dist}(x, K)= 0$ si y solo si $x \in K$. Concluimos que $x \in K$ para todo $x \in \mathbb{T}^2$, es decir $K= \mathbb{T}^2$, probando que el \'{u}nico atractor de Ilyashenko en este ejemplo es todo el toro.

 Entonces  el conjunto $K_0:=\{(0,0)\}$ no es atractor de Ilyashenko. Pero $K_0$ es minimal para $f$  desde el punto de vista topol\'{o}gico porque $(0,0)$ es un punto fijo por $f$ (i.e. $K_0$ es compacto no vac\'{\i}o,  $f$-invariante, y no contiene subconjuntos propios con esas tres propiedades). Luego, en este ejemplo, hay un minimal (desde el punto de vista topol\'{o}gico)   que no es  atractor de Ilyashenko. M\'{a}s a\'{u}n, ning\'{u}n minimal desde el punto de vista topol\'{o}gico puede ser atractor de Ilyashenko. En efecto,  m\'{a}s arriba probamos  que el \'{u}nico atractor estad\'{\i}stico es todo el toro. Pero  todo el toro no es minimal  desde el punto de vista topol\'{o}gico, pues contiene propiamente a $\{(0,0)\}$ que es $f$-invariante.

\end{example}

\begin{remark} \em  \index{probabilidad! emp\'{\i}rica} \index{medida! de probabilidad emp\'{\i}rica}
La demostraci\'{o}n que hicimos en el Ejemplo \ref{ejemploAnosovlineal} de existencia de atractores de Ilyashenko $\alpha$-obs. minimales que no son minimales desde el punto de vista topol\'{o}gico, se bas\'{o} en el uso del Corolario \ref{corollarySRBanosov}, que a su vez se obtiene de la Teor\'{\i}a de Pesin  y de la Teor\'{\i}a de medidas SRB para difeomorfismos de clase $C^{1+ \alpha}$. M\'{a}s precisamente, el argumento que utilizamos en el Ejemplo \ref{ejemploAnosovlineal}, pas\'{o} por deducir  la convergencia de la sucesi\'{o}n de probabilidades emp\'{\i}ricas  $$\sigma_{n,x} =  \frac{1}{n} \sum_{j= 0}^{n-1} \delta_{f^j(x)}$$ para $m$-c.t.p. $x $ en la cuenca de atracci\'{o}n estad\'{\i}stica $A_K$ del atractor de Ilyashenko $K$. Sin embargo, en general (cuando no sean aplicables los argumentos de la Teor\'{\i}a de Pesin), \em no es necesario que la sucesi\'{o}n $\{\sigma_{n,x}\}_{n \geq 1}$ sea convergente, \em para un conjunto de medida de Lebesgue positiva contenido en la cuenca $A_K$ de atracci\'{o}n estad\'{\i}stica de $K$  (ver  por ejemplo \cite{Golenishcheva}).

En la pr\'{o}xima secci\'{o}n, tendremos en cuenta la posible no convergencia de la sucesi\'{o}n de probabilidades emp\'{\i}ricas  para definir las medidas SRB-like (Definici\'{o}n \ref{definitionMedidaSRBlike}) Finalmente, caracterizaremos a   los atractores de Ilyashenko $1$-obs. minimales, estudiando las medidas SRB-like. Otras relaciones  entre la eventual convergencia o no convergencia de la sucesi\'{o}n de probabilidades emp\'{\i}ricas, y los atractores estad\'{\i}sticos o de Ilyashenko se encuentran por ejemplo en \cite{AshwinStatAttr}.

\end{remark}

\begin{example} \em
\label{ejemploIlyashenkoNoTopologico} \label{ejemploHuYoung} {\bf Hu-Young \cite{HuYoung}}

 Sea $f_0 \in \mbox{Diff }^{\infty}(\mathbb{T}^2)$ el Anosov lineal en el toro bidimensional dado en el ejemplo de la Secci\'{o}n \ref{section2111},   con matriz asociada $$ \Big(
                                                    \begin{array}{cc}
                                                      2 & 1 \\
                                                      1 &1 \\
                                                    \end{array}
                                                  \Big)  \Big|_{\mbox{\footnotesize m\'{o}d.} (1,1)}.
 $$

 En \cite{HuYoung} se construye una isotop\'{\i}a $\{f_t\}_{t \in [0,1]}$ de difeomorfismos $f_t \in \mbox{Diff }^{2}(\mathbb{T}^2)$ que transforma continuamente en el espacio $\mbox{Diff }^{2}(\mathbb{T}^2)$ el difeomorfismo de Anosov lineal $f_0$ en $f_1$ de modo que:

  $\bullet$ $f_t$ es un difeomorfismo   de Anosov para todo $0 \leq t < 1$.

  $\bullet$ $f_t(0,0) = (0,0) $ para todo $0 \leq t \leq 1$.

  $\bullet$ $f_t$ es conjugado a $f_0$ para todo $0 \leq t \leq 1$, es decir existe un homeomorfismo $h_t: \mathbb{T}^2 \mapsto \mathbb{T}^2$ tal que $$f_t \circ h_t = h_t \circ f_0.$$

  $\bullet$ Para $t= 1$ la derivada $df_1(0,0)$ en el punto fijo   $(0,0)$ tiene matriz asociada diagonalizable, con un valor propio igual a 1, y el otro positivo, pero menor estrictamente que 1. Esto implica que el punto fijo $(0,0)$ no es hiperb\'{o}lico para $f_1$ (para el par\'{a}metro $t= 1$). Pierde hiperbolicidad en la direcci\'{o}n que era expansora  para valores del par\'{a}metro $t < 1$, pero no la pierde en la direcci\'{o}n contractiva.

    Adem\'{a}s, en la construcci\'{o}n de \cite{HuYoung} se imponen otras condiciones a $f_1$ que permiten deducir la llamada quasi-hiperbolicidad  en $\mathbb{T}^2\setminus \{(0,0)\}$.

  En este ejemplo construido en \cite{HuYoung}, los autores muestran que la medida $\delta_{(0,0)}$ concentrada en el origen (que es $f_1$-invariante porque el $(0,0)$ es punto fijo) es    SRB o f\'{\i}sica para $f_1$ (de acuerdo a nuestra Definici\'{o}n \ref{definitionMedidaSRB}, pero no de acuerdo con lo que los autores en \cite{HuYoung} llaman medida SRB). Adem\'{a}s muestran  que la cuenca de atracci\'{o}n estad\'{\i}stica $B(\delta_{(0,0)})$ cubre Lebesgue c.t.p. del toro $\mathbb{T}^2$. Entonces,  $\delta_{(0,0)}$ es la \'{u}nica medida SRB o f\'{\i}sica, ya que ninguna otra medida puede tener cuenca de atracci\'{o}n estad\'{\i}stica con medida de Lebesgue positiva (recordar que las cuencas de atracci\'{o}n estad\'{\i}stica de medidas diferentes, por la Definici\'{o}n \ref{DefinicionCuencaDeAtraccionEstadistica}, son disjuntas).

  Aplicando la Proposici\'{o}n \ref{propositionSRB->estadistico}, deducimos que \em $  \{(0,0)\}$ es el \'{u}nico atractor de Ilyashenko de $f_1$. \em  Por lo tanto $  \{(0,0)\}$ es el \'{u}nico atractor de Ilyashenko $\alpha$-obs. minimal para cualquier $0 < \alpha \leq 1$.

  Como  $f_1$ es conjugado a $f_0$ y $f_0$ es transitivo, entonces $f_1$ es transitivo. Aplicando lo probado en el Ejercicio \ref{exerciseTransitivoAtrTop}, el \'{u}nico atractor topol\'{o}gico para $f_1$ es ${\mathbb{T}^2}$. Luego $  \{(0,0)\}$ es un atractor estad\'{\i}stico o de Ilyashenko, que \em no es atractor topol\'{o}gico.\em 

  Por la Proposici\'{o}n \ref{propositionAtraccionEstadisticaFrec}, la frecuencia de visita a un entorno arbitrariamente peque\~{n}o del origen, para Lebesgue-casi toda \'{o}rbita, tiende a 1, cuando el n\'{u}mero $n$ de iterados   tiende a infinito. Sin embargo,  el origen no es un atractor topol\'{o}gico. Es decir, no es un pozo. Por el contrario, la din\'{a}mica en un entorno del origen es localmente conjugado  a la de un punto fijo hiperb\'{o}lico tipo silla,   pues el origen es una silla hiperb\'{o}lica para $t= 0$ y $f_1$ es conjugado con $f_0$. Entonces a pesar de que la frecuencia de estad\'{\i}a en un entorno del origen, para Lebesgue casi toda \'{o}rbita, tiende a 1, todo punto que no est\'{e} en la variedad estable local del origen, termina saliendo de ese entorno, para hacer excursiones breves alejado de \'{e}l.

  Concluimos: Si el experimentador de la din\'{a}mica tiene como objetivo observar la estad\'{\i}stica de Lebesgue casi toda \'{o}rbita (es decir los promedios temporales de las distancias al atractor), entonces observar\'{a} que el sistema din\'{a}mico por iterados de $f_1$ se comporta con un punto fijo atractor, que atrae Lebesgue casi toda \'{o}rbita, como si fuera un pozo.

  En cambio si el experimentador de la din\'{a}mica tiene como objetivo observar la topolog\'{\i}a din\'{a}mica de Lebesgue casi toda \'{o}rbita (es decir los conjuntos $\omega$-l\'{\i}mite a donde estas \'{o}rbitas tienden al iterar hacia el futuro), entonces observar\'{a} que el sistema din\'{a}mico por iterados de $f_1$ se comporta como el automorfismo lineal hiperb\'{o}lico en el toro, en que todo el toro es el \'{u}nico atractor transitivo de Lebesegue casi toda \'{o}rbita.

  En el primer caso, el experimentador estad\'{\i}stico, no calificar\'{a} este sistema $f_1$ como ca\'{o}tico, pues es altamente previsible, desde el punto de vista estad\'{\i}stico (i.e. de los promedios temporales) para Lebesgue-casi toda \'{o}rbita. En el segundo caso, el experimentador topol\'{o}gico, lo calificar\'{a} como ca\'{o}tico, pues le resultar\'{a} imposible predecir en qu\'{e} abierto del espacio estar\'{a}   el iterado $n$-\'{e}simo para Lebesgue-casi toda \'{o}rbita.
\end{example}


\subsection{Medidas SRB-like o pseudo-f\'{\i}sicas}

En esta secci\'{o}n, como en las dos anteriores, $M$ es una variedad compacta y riemanniana, de dimensi\'{o}n finita, y $f: M \mapsto M$ es continua (no necesariamente invertible).

\begin{definition}
\label{definitionpomegalimit} {\bf El conjunto $p\omega$ l\'{\i}mite de probabilidades} \index{$p\omega(x)$ p-omega l\'{\i}mite en ${\mathcal M}$} \index{probabilidad! emp\'{\i}rica}

\em Sea $x \in M$, sea ${\mathcal M}$ el espacio de medidas de probabilidad de Borel en $M$ con la topolog\'{\i}a d\'{e}bil estrella, y sea $\{\sigma_{n,x}\}_{n \geq 1} \subset {\mathcal M}$ la sucesi\'{o}n de probabilidades emp\'{\i}ricas de la \'{o}rbita futura de $x$ hasta tiempo $n$, definida en (\ref{(1)}) y Definici\'{o}n \ref{definitionProbaEmpiricas}.
Como ${\mathcal M}$ es secuencialmente compacto, existen subsucesiones convergentes de $\{\sigma_{n,x}\}_{n \geq 1}$.

Llamamos \em $p$-omega-l\'{\i}mite \em de la \'{o}rbita de $x$, u \em omega-l\'{\i}mite en el espacio de probabilidades \em de la \'{o}rbita de $x$ al conjunto $p\omega(x)$ formado por los l\'{\i}mites de todas las subsucesiones convergentes de $\{\sigma_{n,x}\}_{n \geq 1}$. Es decir
\begin{equation}
\label{eqnpomegalimit}
p\omega(x):= \{\mu \in {\mathcal M}\colon \ \ \exists n_i \rightarrow + \infty \mbox{ tal que } \lim_{i \rightarrow + \infty} \sigma_{n_i,x} = \mu    \},
\end{equation}
donde el l\'{\i}mite   a la derecha es en la topolog\'{\i}a d\'{e}bil$^*$ del espacio ${\mathcal M}$ de probabilidades.

Es est\'{a}ndar chequear que, para todo $x \in M$, $p\omega(x) \subset {\mathcal M}$ es un conjunto no vac\'{\i}o y d\'{e}bil$^*$-cerrado (y por lo tanto d\'{e}bil$^*$-compacto, pues ${\mathcal M}$ es d\'{e}bil$^*$-compacto)
\end{definition}

\begin{exercise}\em
Probar las dos \'{u}ltimas afirmaciones  de la Definici\'{o}n \ref{definitionpomegalimit}.
\end{exercise}

Recordemos la Definici\'{o}n \ref{DefinicionCuencaDeAtraccionEstadistica} de cuenca de atracci\'{o}n estad\'{\i}stica $B(\mu)$ de una medida de probabilidad $\mu$, y la Definici\'{o}n \ref{definitionMedidaSRB} de medida SRB o f\'{\i}sica. Ahora generalizaremos esas dos definiciones, agregando una $\epsilon$- aproximaci\'{o}n. Para ello elegimos y dejamos fija una m\'{e}trica $\mbox{dist}^*$ en el espacio ${\mathcal M}$ de probabilidades, que induzca la topolog\'{\i}a d\'{e}bil$^*$ (cf. Teorema \ref{teoremaCompacidadEspacioProbabilidades}).

\begin{definition}
\label{DefinicionCuencaDeAtraccionEpsilonEstadistica}{\bf Cuenca de atracci\'{o}n estad\'{\i}stica $\epsilon$-d\'{e}bil} \em \index{cuenca de atracci\'{o}n! estad\'{\i}stica $\epsilon$-d\'{e}bil} \index{$A_{\epsilon}(\mu)$ cuenca de atracci\'{o}n! estad\'{\i}stica  $\epsilon$-d\'{e}bil! de la medida $\mu$}


Dada una medida de probabilidad $\mu$ y dado $\epsilon >0$, llamamos \em cuenca de atracci\'{o}n estad\'{\i}stica $\epsilon$-d\'{e}bil \em al  conjunto $A_{\epsilon}(\mu)$ definido por:
\begin{equation}
\label{equationAepsilon(mu)}
A_{\epsilon}(\mu):= \{x \in M: \ \mbox{dist} (p \omega(x), \mu) < \epsilon\}. \end{equation}
Comparemos esta Definici\'{o}n \index{cuenca de atracci\'{o}n! estad\'{\i}stica} \ref{DefinicionCuencaDeAtraccionEpsilonEstadistica} con la Definici\'{o}n \ref{DefinicionCuencaDeAtraccionEstadistica} de la cuenca de atracci\'{o}n estad\'{\i}stica (fuerte) $B_{\mu}$ de una medida $\mu$. En efecto, combinando las igualdades (\ref{eqnpomegalimit})  y (\ref{equationB(mu)}), obtenemos:
$$B(\mu):= \big\{x \in M: \ \ p\omega (x) = \{\mu\}\big\}.$$
Luego \begin{equation}
\label{equationB(mu)subsetAepsilon(mu)}
B(\mu) \subset A_{\epsilon}(\mu) \ \ \ \forall \ \ \epsilon \geq 0.\end{equation}
En el segundo t\'{e}rmino de la Igualdad (\ref{equationAepsilon(mu)}), observamos que la distancia entre el conjunto compacto $p\omega(x)$ y el punto $\mu \in {\mathcal M}$ es menor que $\epsilon$. Pero esto no implica que todo el conjunto $p \omega(x)$ (cuando contiene m\'{a}s de un punto) deba estar contenido en la bola de centro $\mu$ y radio $\epsilon$. Por lo tanto, aunque $\bigcap_{\epsilon >0} A_{\epsilon}(\mu)$ pueda ser no vac\'{\i}o, este conjunto contiene a, \em pero no necesariamente coincide con \em $B(\mu)$, quien a\'{u}n, puede   ser vac\'{\i}o.

\end{definition}


\begin{definition}
\label{definitionMedidaSRBlike}{\bf Medida SRB-like o pseudo-f\'{\i}sica} \em \index{medida! SRB-like} \index{medida! pseudo-f\'{\i}sica}

Llamaremos a una medida de probabilidad $\mu$  \em SRB-like o f\'{\i}sica \em si {\bf para todo } $\epsilon >0$ su cuenca de atracci\'{o}n estad\'{\i}stica $\epsilon$-d\'{e}bil $A_{\epsilon}(\mu)$ tiene medida de Lebesgue positiva. En breve:
$$m(A_{\epsilon}(\mu)) >0 \ \ \forall \ \epsilon >0,$$
donde $m$ denota la medida de Lebesgue en la variedad $M$.
\end{definition}

Comparando la Definici\'{o}n \ref{definitionMedidaSRBlike} con la Definici\'{o}n \ref{definitionMedidaSRB} de medida SRB, observamos que debido a la inclusi\'{o}n (\ref{equationB(mu)subsetAepsilon(mu)}):

\begin{center}
\em Toda medida SRB es SRB-like. \em
\end{center}

Sin embargo el rec\'{\i}proco es falso. En efecto,   toda $f$ continua tiene medidas SRB-like (cf. Teorema \ref{theoremExistenciaSRB-like} que demostraremos m\'{a}s adelante en esta secci\'{o}n). Sin embargo existen ejemplos para los que no existen medidas SRB (cf. Ejemplos \ref{ejemploRotacionEsfera} y \ref{ejemploMilnorNoTopologico1}, y en el caso hiperb\'{o}lico, el mapa $C^1$ en un disco, atribuido a Bowen y estudiado en \cite{Golenishcheva}).

\begin{exercise}\em
\label{ejercicioSRB-likeInvariante}.

{\bf (a)} Probar que toda medida SRB-like es invariante por $f$.

{\bf (b)} Probar que la definici\'{o}n de  medida SRB-like no depende de la m\'{e}trica elegida en el espacio ${\mathcal M}$ de probabilidades con la topolog\'{\i}a d\'{e}bil$^*$.
\end{exercise} \index{${\mathcal O}_f $ conjunto de medidas  SRB-like! (o pseudo-f\'{\i}sicas u observables) para $f$}

{\bf Notaci\'{o}n:} Denotamos con ${\mathcal O}_f$ al conjunto de todas las medidas SRB-like para $f$. Esta notaci\'{o}n proviene de \cite{CatEnr2011}, que introduce la definici\'{o}n de medidas SRB-like, llam\'{a}ndolas tambi\'{e}n \em medidas observables. \em

\vspace{.2cm}

{\bf Observaciones: }  ${\mathcal O}_f $ est\'{a} contenido en el conjunto ${\mathcal M}_f$ de medidas de probabilidad $f$-invariantes, pero usualmente difiere mucho de ${\mathcal M}_f$, como veremos en el Ejemplo  \ref{ejemploExpansorasEnElCirculo} $C^1$ gen\'{e}rico. Sin embargo, en el caso $C^0$, y para los que llamamos \em endomorfismos expansores de Misiurewicz en el c\'{\i}rculo $S^1$ \em  (para los que no existen medidas SRB), el conjunto de medidas SRB-like coincide con el conjunto ${\mathcal M}_f$ de todas las medidas invariantes (ver \cite{Misiurewicz}, o tambi\'{e}n el Corolario \ref{corolarioSRB-likeMisiurevicz} y el Ejemplo \ref{ejemploExpansorasC0EnElCirculo} m\'{a}s adelante en esta secci\'{o}n).

\vspace{.2cm}

Algunas propiedades   que distinguen a las medidas SRB-like son: \index{medida! SRB-like}

$\bullet$ Existen medidas SRB-like, sin necesidad de agregar hip\'{o}tesis adicionales a la continuidad de $f$ (Teorema \ref{theoremExistenciaSRB-like}).

$\bullet$ El conjunto de medidas SRB-like describen en forma \'{o}ptima (con un m\'{\i}nimo posible de medidas invariantes) la estad\'{\i}stica en el futuro de Lebesgue casi toda \'{o}rbita (Teorema \ref{theoremOptimalidadSRB-like}).

$\bullet$ El  m\'{\i}nimo soporte  compacto com\'{u}n de   medidas SRB-like (cf. Definiciones \ref{definitionSoporteDemu} y \ref{definitionSoporteDemu+}) caracteriza a   los atractores estad\'{\i}sticos o de Ilyashenko  $\alpha$-obs. minimales.  (Teorema \ref{theoremSRB-like<->Ilyashenko}). \index{atractor! estad\'{\i}stico} \index{atractor! de Ilyashenko}

$\bullet$ Bajo hip\'{o}tesis adicionales de $C^1$ hiperbolicidad, las medidas SRB-like, satisfacen la F\'{o}rmula (\ref{eqnformuladePesin}) de Pesin de la Entrop\'{\i}a (cf. Ejemplos \ref{ejemploExpansorasEnElCirculo} y \ref{ejemploSRB-likeAnosovC1}), aunque en el contexto general de regularidad $C^1$ no necesariamente tienen propiedades de continuidad absoluta respecto a Lebesgue.

 \begin{theorem}
 \label{theoremExistenciaSRB-like} {\bf Existencia de medidas SRB-like \cite{CatEnr2011}} \index{teorema! de existencia de! medidas SRB-like} \index{medida! SRB-like} \index{medida! pseudo-f\'{\i}sica} \index{teorema! compacidad de ${\mathcal O}_f$}

  Sea $f: M \mapsto M$ continua. Entonces:

  {\bf (a) } Existen medidas   de probabilidad SRB-like para $f$.

  {\bf (b) } El conjunto ${\mathcal O}_f$ de las medidas SRB-like es d\'{e}bil$^*$-compacto en el espacio de las probabilidades invariantes por $f$.

 \end{theorem}

  Extraemos la prueba de \cite{CatEnr2011}:
 {\em Demostraci\'{o}n: }
 {\bf (a) } Supongamos por absurdo que ${\mathcal O}_f $ es vac\'{\i}o. Entonces   toda probabilidad $\mu $   es no SRB-like. Es decir, $\mu$ no satisface la   Definici\'{o}n \ref{definitionMedidaSRBlike}. Recordamos la Definici\'{o}n \ref{definitionpomegalimit} del conjunto $p\omega(x)$ ($p$-omega-l\'{\i}mite de la \'{o}rbita por cada punto $x \in M$). Luego,   toda $\mu \in {\mathcal M}$ est\'{a} contenida en un entorno abierto ${B} (\mu) \subset {\mathcal M}$ (con la topolog\'{\i}a d\'{e}bil$^*$ del espacio de probabilidades ${\mathcal M}$) tal que
 $$m\Big(\{x \in M \colon \ p\omega(x) \cap {B}(\mu) \neq \emptyset  \}\Big) = 0$$
 Como ${\mathcal M}$ es d\'{e}bil$^*$ compacto, existe un subcubrimento finito de $M$ $$\{{B}_1, {B}_2, \ldots, {B}_k\}$$
 tal que ${B}_i = {B}(\mu_i)$ para alguna medida de probabilidad $\mu_i$. Luego:
 $$m\Big(\{x \in M \colon \ p\omega(x) \bigcap \big(\bigcup_{i= 1}^k{B}_i\big) \neq \emptyset \} \Big) =$$ $$= m\Big(\bigcup_{i= 1}^k \{x \in M \colon \ p\omega(x) \cap {B}_i \neq \emptyset \} \Big) \leq $$ $$\leq \sum_{i= 1}^k m\Big(\{x \in M \colon \ p\omega(x) \cap {B}_i \neq \emptyset \} \Big) = 0.$$
 Como $\bigcup_{i= 1}^k {B}_i = {\mathcal M}$, deducimos que
 $$m\Big(\{x \in M \colon \ p\omega(x) \cap  {\mathcal M}  \neq \emptyset \} \Big) = 0.$$
 Como   $p \omega(x) \subset {\mathcal M}$ para todo $x \in M$, deducimos que para Lebesgue c.t.p. $x \in M$, $p \omega(x) = \emptyset$, lo cual es una contradicci\'{o}n porque toda sucesi\'{o}n de probabilidades tiene alguna subsucesi\'{o}n convergente.

  {\bf (b) } Para probar que ${\mathcal O}_f$ es d\'{e}bil$^*$ compacto, basta probar que es d\'{e}bil$^*$ cerrado (pues ${\mathcal O}_f \subset{\mathcal M}$ y ${\mathcal M}$ es un espacio metrizable d\'{e}bil$^*$ compacto). Sea $\mu_n \in {\mathcal O}_f$ convergente a $\mu$. Debemos probar que $\mu \in {\mathcal O}_f$.

  Sea $\epsilon >0$ arbitrario. Sea $B_{\epsilon}(\mu)$ la bola de centro $\mu$ y radio $\epsilon$. Sea $n$ tal que $\mu_n \in B_{\epsilon}(\mu)$ y sea $\epsilon _n >0$ tal que la bola $B_{\epsilon_n}(\mu_n)$ de centro $\mu_n$ y radio $\epsilon _n$ satisface $$B_{\epsilon_n}(\mu_n) \subset B_{\epsilon}(\mu).$$
  Como $\mu_n \in {\mathcal O}_f$, por la Definici\'{o}n \ref{definitionMedidaSRBlike} de medida SRB-like, se cumple
  $$m\Big(\{x \in M: \ p \omega(x) \cap B_{\epsilon_n}(\mu_n) \neq    \emptyset\}   \Big) >0.$$
  Como $$ \{x \in M: \ p \omega(x) \cap B_{\epsilon_n}(\mu_n) \neq    \emptyset\} \ \subset \  \{x \in M: \ p \omega(x) \cap B_{\epsilon}(\mu) \neq    \emptyset\},$$
  deducimos que
  $$m\Big(\{x \in M: \ p \omega(x) \cap B_{\epsilon}(\mu) \neq    \emptyset\}   \Big) >0.$$
  Siendo $\epsilon >0$, esto prueba que $\mu $ es SRB-like, como quer\'{\i}amos demostrar.
  \hfill $\Box$
 Para enunciar el siguiente teorema, recordamos la Definici\'{o}n \ref{definitionpomegalimit} del conjunto $p\omega(x)$ ($p$-omega-l\'{\i}mite de la \'{o}rbita por el punto $x$). De la igualdad (\ref{eqnpomegalimit}), observamos que $p \omega(x)$ es el m\'{\i}nimo conjunto de probabilidades que describen completamente la estad\'{\i}stica (es decir los promedios temporales asint\'{o}ticos) de la \'{o}rbita por el punto $x$.

 \begin{theorem}
 \label{theoremOptimalidadSRB-like}.  \index{teorema! optimalidad estad\'{\i}stica} \index{medida! SRB-like}

 {\bf Optimalidad estad\'{\i}stica del conjunto de medidas SRB-like \cite{CatEnr2011}}

 Para toda $f: M \mapsto M$ continua
 el conjunto ${\mathcal O}_f$ de las medidas SRB-like para $f$ es el m\'{\i}nimo conjunto d\'{e}bil$^*$-compacto ${\mathcal K}$ del espacio de probabilidades   tal que \em
 $$p \omega (x) \subset {\mathcal K} \ \ \mbox{para Lebesgue-c.t.p. } x \in M.$$
 \end{theorem}
 {\em Demostraci\'{o}n: }
 Por la parte (b) del Teorema \ref{theoremExistenciaSRB-like}, el conjunto ${\mathcal O}_f$ de todas las medidas SRB-like para $f$ es no vac\'{\i}o y d\'{e}bil$^*$-compacto.

 Primero probemos que $p\omega(x) \subset {\mathcal O}_f$ para Lebesgue c.t.p. $x     \in M$. Para   $\epsilon >0$ arbitrario, denotamos con $B_{\epsilon}({\mathcal O}_f)$ al conjunto abierto de todas las probabilidades $\nu$ que distan del compacto no vac\'{\i}o ${\mathcal O}_f$ menos que $\epsilon$. Consideremos el complemento $${\mathcal C}:= {\mathcal M} \setminus B_{\epsilon}({\mathcal O}_f).$$ El conjunto ${\mathcal C}$ es compacto  porque es el complemento del abierto $B_{\epsilon}({\mathcal O}_f)$  en el espacio compacto ${\mathcal M}$. Por construcci\'{o}n ${\mathcal C} \cap {\mathcal O}_f = \emptyset$. Luego, toda $\nu \in {\mathcal C}$ es no SRB-like. Entonces existe un entorno abierto $B(\nu)$ de $\nu$ tal que:
 $$m \Big(\{x \in M: \ p \omega(x) \cap B(\nu) \neq \emptyset   \}\Big ) = 0.$$
 Al igual que al final de la demostraci\'{o}n de la parte (a) del Teorema \ref{theoremExistenciaSRB-like}, pero escribiendo ${\mathcal C}$ en el rol de ${\mathcal M}$, deducimos que existe un cubrimiento finito $\{B_1, B_2, \ldots, B_k\}$ de ${\mathcal C}$ tal que
 $$m \Big(\{x \in M: \ p \omega(x) \bigcap \big( \bigcup_{i= 1}^k B_i \big) \neq \emptyset\}     \Big)= 0$$
 Como $\bigcup_{i= 1} ^k B_i \supset {\mathcal C}$, deducimos
 $$m \Big(\{x \in M: \ p \omega(x) \cap {\mathcal C} \neq \emptyset\}     \Big)= 0.$$
 Dicho de otra forma, para Lebesgue-casi todo punto $x \in M$, el conjunto $p \omega (x)$ est\'{a} contenido en   ${\mathcal M} \setminus {\mathcal C} = B_{\epsilon}({\mathcal O}_f)$. Luego, tomando $\epsilon = 1/h$ para $h \in \mathbb{N}^+$, hemos probado que
 $$\forall \ h \in \mathbb{N}^+: \ p \omega(x) \subset B_{1/h}({\mathcal O}_f) \ \ m-\mbox{c.t.p. } x \in M. $$
 Como la uni\'{o}n numerable de conjuntos con medida de Lebesgue nula, tiene medida de Lebesgue nula, la intersecci\'{o}n numerable de conjuntos con medida de Lebesgue total, tiene medida de Lebesgue total. Concluimos que:
 $$p \omega(x) \subset \bigcap_{h= 1}^{+\infty} B_{1/h}({\mathcal O}_f) = {\mathcal O}_f \ \ m-\mbox{c.t.p. } x \in M,$$
 como quer\'{\i}amos demostrar.

 Segundo, probemos que ${\mathcal O}_f$ es minimal en el conjunto de compactos no vac\'{\i}os ${\mathcal K} \subset {\mathcal M} $ que tienen la  propiedad $p\omega(x) \subset{\mathcal K}$ para Lebesgue-casi todo punto $x \in M$. Tomemos cualquier compacto no vac\'{\i}o ${\mathcal K} \subset{\mathcal O}_f$, tal que ${\mathcal K} \neq {\mathcal O}_f$. Probemos que tal ${\mathcal K}$ no tiene la propiedad mencionada. Es decir, probemos que   $m \Big( p \omega(x) \cap ({\mathcal M}\setminus {\mathcal K})\Big) >0$.

 En efecto, como ${\mathcal K} \subset \neq{\mathcal O}_f$, existe una medida $\mu \in {\mathcal O}_f \setminus {\mathcal K}$. Sea  $\epsilon >0$ tal que la bola $B_{\epsilon}(\mu)$ es disjunta con el compacto ${\mathcal K}$.   Como $\mu $ es SRB-like, por las Definiciones \ref{definitionMedidaSRBlike} y \ref{DefinicionCuencaDeAtraccionEpsilonEstadistica}, se cumple:
 \begin{equation}
 \label{eqn-26}
 m\Big(\{x \in M\colon \ p \omega(x) \cap B_{\epsilon}(\mu) \neq \emptyset \} \Big) >0.\end{equation}
 Siendo $B_{\epsilon}(\mu) \cap {\mathcal K} = \emptyset$, deducimos que $B_{\epsilon}(\mu) \subset ({\mathcal M} \setminus {\mathcal K})$. Entonces, sustituyendo en (\ref{eqn-26}) concluimos:
 $$m\Big(\{x \in M\colon \ p \omega(x) \cap ({\mathcal M} \setminus {\mathcal K})   \neq \emptyset \} \Big) >0,$$ como quer\'{\i}amos demostrar.
 \hfill $\Box$

 \begin{corollary} \label{corolarioSRB-likefinito}
  Sea $f: M \mapsto M$ continua. \index{medida! SRB-like} \index{cuenca de atracci\'{o}n! estad\'{\i}stica} \index{teorema! unicidad de medida SRB-like} \index{medida! SRB}

  {\bf (a) } Si la medida SRB-like $\mu$ es \'{u}nica, entonces $\mu$ es SRB y su cuenca de atracci\'{o}n estad\'{\i}stica (fuerte) $B(\mu)$ cubre $M$ Lebesgue c.t.p.

  {\bf (b) } Rec\'{\i}procamente, si existe una medida SRB $\mu$ cuya cuenca de atracci\'{o}n estad\'{\i}stica (fuerte) $B(\mu)$ cubre $M$ Lebesgue c.t.p., entonces $\mu$ es la \'{u}nica medida SRB-like.

  {\bf (c) } Si el conjunto de medidas SRB-like es finito, entonces todas las medidas SRB-like son SRB y la uni\'{o}n de las cuencas de atracci\'{o}n estad\'{\i}stica de las medidas SRB cubre $M$ Lebesgue c.t.p.

  {\bf (d) } Rec\'{\i}procamente, si existe una cantidad finita de medidas SRB tales que la uni\'{o}n de sus cuencas de atracci\'{o}n estad\'{\i}stica cubre $M$ Lebesgue c.t.p., entonces estas son las \'{u}nicas medidas SRB-like, y por lo tanto el conjunto de medidas SRB-like es finito.
 \end{corollary}

  \begin{exercise}\em
  Probar el Corolario \ref{corolarioSRB-likefinito}. Sugerencia: Basta probar (c) y (d), pues estas implican (a) y (b). Probar primero que una medida SRB-like es aislada en el conjunto de las medidas SRB-like si y solo si es SRB. Combinar las Definiciones \ref{definitionMedidaSRB} y \ref{definitionMedidaSRBlike} de medidas SRB y SRB-like respectivamente, junto con las Definiciones \ref{DefinicionCuencaDeAtraccionEstadistica} y \ref{DefinicionCuencaDeAtraccionEpsilonEstadistica} de las cuencas de atracci\'{o}n estad\'{\i}stica  fuerte  y $\epsilon$-d\'{e}bil, respectivamente.
  \end{exercise}

 \begin{conjecture} \em  \index{conjetura! Palis}
 {\bf  Palis \cite{PalisConjecture}} \em Para $r \geq 1$ suficientemente grande $C^r$-gen\'{e}ricamente  en ${\mbox{Diff }^r(M)}$ existe una cantidad finita de medidas SRB tales que la uni\'{o}n de sus cuencas de atracci\'{o}n estad\'{\i}stica cubre $M$ Lebesgue c.t.p. \em

 \vspace{.2cm}

 Del Corolario \ref{corolarioSRB-likefinito} deducimos el siguiente enunciado equivalente de la Conjetura de Palis:

 \vspace{.2cm}

\index{conjetura! Palis}
 \em $C^r$-gen\'{e}ricamente  para $f \in {\mbox{Diff }^r(M)}$ el conjunto d\'{e}bil$^*$ no vac\'{\i}o ${\mathcal O}_f$ de medidas de probabalidad SRB-like para $f$, carece de puntos de acumulaci\'{o}n.
 \end{conjecture}

 \begin{corollary}  {\bf (del Teorema \ref{theoremOptimalidadSRB-like})}
 \label{corolarioSRB-likeMisiurevicz}

  Sea $f: M \mapsto M$ continua, no \'{u}nicamente erg\'{o}dica. Si para Lebesgue c.t.p. $x \in M$ el conjunto de probabilidades $p\omega(x)$ \em (cf. Definici\'{o}n \ref{definitionpomegalimit}) \em  coincide con el conjunto de todas las medidas inva\-rian\-tes, entonces no existen medidas SRB y
 el conjunto de las medidas SRB-like coincide con el conjunto de todas las medidas in\-va\-rian\-tes.
 \end{corollary}
 {\bf Nota:}  Existen  transformaciones $f$ que cumplen las hip\'{o}tesis del Corolario \ref{corolarioSRB-likeMisiurevicz}. En efecto, el mapa  $C^0$ expansor  en el c\'{\i}rculo construido por Misiurevicz en \index{mapa expansor} \cite{Misiurewicz}, y los mapas $C^0$-expansores gen\'{e}ricos encontrados recientemente por Abdenur y Andersson en \cite{AbdenurAndersson}, satisfacen las hip\'{o}tesis de este Corolario. Ver tambi\'{e}n el Ejemplo \ref{ejemploExpansorasC0EnElCirculo} m\'{a}s adelante en esta secci\'{o}n.

 {\em Demostraci\'{o}n: }
 {\em del Corolario } \ref{corolarioSRB-likeMisiurevicz}:

 Por hip\'{o}tesis, el m\'{\i}nimo conjunto compacto ${\mathcal K}$ de medidas de probabilidad tal que $p \omega(x) \subset {\mathcal K}$ es el conjunto ${\mathcal M}_f$ de todas las medidas $f$-invariantes. Por el Teorema \ref{theoremOptimalidadSRB-like}, ${\mathcal K}$ es el conjunto ${\mathcal O}_f$ de   las medidas SRB-like. Luego ${\mathcal O}_f = {\mathcal M}_f$ como quer\'{\i}amos demostrar.
 \hfill $\Box$

\subsection{Relaci\'{o}n entre atractor estad\'{\i}stico y medidas SRB-like}

 \begin{definition} \index{medida! soporte compacto de}
 \label{definitionSoporteDemu+} \em Sea $ {\mathcal K} $ un conjunto no vac\'{\i}o de medidas de probabilidad. Se llama \em soporte compacto de ${\mathcal K}$ \em al m\'{\i}nimo compacto $K \subset M$ tal que $$\mu(K) = 1 \ \ \forall \ \mu \in {\mathcal K}.$$
 El m\'{\i}nimo soporte compacto existe como resultado de aplicar el lema de Zorn a la familia de compactos no vac\'{\i}os con $\mu$-medida igual a 1 para toda $\mu \in {\mathcal K}$, y de usar tambi\'{e}n la propiedad de que es no vac\'{\i}a la intersecci\'{o}n de una familia de compactos tal que toda subfamilia finita tiene intersecci\'{o}n no vac\'{\i}o. El soporte compacto de ${\mathcal K}$ es \'{u}nico, debido a su propiedad de minimalidad y a que la intersecci\'{o}n de dos compactos con $\mu$-medida igual a 1 es un compacto con $\mu$-medida igual a 1. \index{soporte compacto}
 \end{definition}

 \begin{theorem}.\index{teorema! de caracterizaci\'{o}n de! atractor estad\'{\i}stico}
\label{theoremSRB-like<->Ilyashenko}

{\bf Medidas SRB-like y atractores de Ilyashenko} \index{atractor! estad\'{\i}stico} \index{atractor! de Ilyashenko} \index{medida! SRB-like} \index{soporte compacto}

Para toda $f\colon M \mapsto M$ continua, el atractor $1$-obs. minimal de Ilyashenko $K$ \em (cf. Theorem \ref{TeoExistenciaAtrIlyashenko})\ \em es el soporte compacto del conjunto ${\mathcal O}_f$ de medidas SRB-like para $f$.
 \end{theorem}

 {\bf Nota:} El Teorema \ref{theoremSRB-like<->Ilyashenko} puede generalizarse, adaptando el enunciado adecuadamente para caracterizar todo atractor de Ilyashenko $\alpha$-obs. minimal (para cualquier $0 < \alpha \leq 1$) como el soporte compacto de un subconjunto adecuado de medidas SRB-like para $f$ (ver \cite{CatIlyshenkoAttractors}).

 {\em Demostraci\'{o}n: } {\em del Teorema } \ref{theoremSRB-like<->Ilyashenko}:
 Sea $K_1$ el atractor de Ilyashenko 1-obs. minimal. Por el Teorema \ref{TeoExistenciaAtrIlyashenko}, el compacto $K_1 \neq \emptyset$ existe y es \'{u}nico. Sea $K_2$ el soporte compacto del conjunto ${\mathcal O}_f$ de medidas SRB-like para $f$, seg\'{u}n la Definici\'{o}n \ref{definitionSoporteDemu+}. De acuerdo a lo observado al final de dicha definici\'{o}n, el compacto $K_2 \neq \emptyset$ existe y es \'{u}nico.

 {\em Demostraci\'{o}n de $K_1 \subset K_2$:} Como $K_1$ es 1-obs. minimal, por la Definici\'{o}n \ref{definitionAtractorIlyashenkoMinimal} basta probar que \begin{equation}
 \label{eqn-29}
  \mbox{{A probar: }} \ \lim_{n\rightarrow + \infty} \frac{1}{n} \sum_{j= 0}^{n-1}\mbox{dist}(f^j(x), K_2) = 0 \mbox{ para } m-\mbox{c.t.p. } x \in M.\end{equation}
 Para todo $x \in M$ consideremos la sucesi\'{o}n $\{\sigma_{n,x}\}_{n \geq 1}$ de probabilidades emp\'{\i}ricas de la \'{o}rbita por  $x$, seg\'{u}n la Definici\'{o}n \ref{definitionProbaEmpiricas}. Aplicando el Teorema \ref{theoremOptimalidadSRB-like}, tenemos para $m$-c.t.p. $x \in M$ la siguiente propiedad:
 \begin{equation}
  \label{eqn-27}
  \lim_{i \rightarrow + \infty}\sigma_{n_i, x} \in {\mathcal O}_f   \end{equation} para toda subsucesi\'{o}n $ n_i \rightarrow + \infty $ tal que   existe ese l\'{\i}mite (donde dicho l\'{\i}mite se toma en la topolog\'{\i}a d\'{e}bil$^*$ del espacio  de probabilidades).

 Consideremos la funci\'{o}n continua $\psi$ definida por \begin{equation}\label{eqndefpsi}\psi(x) := \mbox{dist}(x,K_2) \ \ \forall \ x \in M.\end{equation} Integrando $\psi$ en la inclusi\'{o}n (\ref{eqn-27}), y teniendo en cuenta la definici\'{o}n de la topolog\'{\i}a d\'{e}bil$^*$, deducimos que para $m$-c.t.p. $x \in M$, y para toda subsucesi\'{o}n convergente $\{\sigma_{n_i,x}\}_{i \in \mathbb{N}} $ de probabilidades emp\'{\i}ricas, existe $\mu \in {\mathcal O}_f$ tal que
 \begin{equation}
 \label{eqn-28}
 \int \psi\, d \mu = \lim_{i \rightarrow + \infty} \int \psi \, d (\sigma_{n_i, x} ) = \lim_{i \rightarrow  + \infty} \frac{1}{n_i}\sum_{j= 0}^{n_i-1}  \mbox{dist}(f^j(x), K_2).\end{equation}
 Pero $\mu(K_2) = 1$, porque por hip\'{o}tesis, $K_2$ es   soporte compacto com\'{u}n de todas las medidas de probabilidad en ${\mathcal O}_f$. Luego $\psi(y) = \mbox{dist}(y, K_2) = 0$ para $\mu$-c.t.p. $y \in M$, de donde $\int \psi \, d \mu = 0$. Sustituyendo en la igualdad (\ref{eqn-28}) obtenemos
 \begin{equation}
 \label{eqn-28z}
 \lim_{i \rightarrow  + \infty} \frac{1}{n_i}\sum_{j= 0}^{n_i-1}  \mbox{dist}(f^j(x), K_2) = 0 \end{equation}
 para $m$-c.t.p. $x \in M$, para toda subsucesi\'{o}n $n_i \rightarrow + \infty$ tal que $\{\sigma_{n_i, x}\}_{i \in \mathbb{N}}$ sea convergente.

 Fijado un tal punto $x$, sea dada una sucesi\'{o}n cualquiera $n_i \rightarrow + \infty$ tal que la subsucesi\'{o}n $\{d_{n_i}\}_{i \in \mathbb{N}}$ de promedios de distancias a $K_2$,  dada por $$d_{n_i} :=  \frac{1}{n_i}\sum_{j= 0}^{n_i-1}  \mbox{dist}(f^j(x), K_2), $$  es convergente. Por la compacidad   del espacio de probabilidades, siempre existe una subsucesi\'{o}n de esta subsucesi\'{o}n $\{d_{n_i}\}_i$ (para \'{\i}ndices $i= i_h$) tal que $\{\sigma_{n_{i_h}, x}\}_{h }$ es tambi\'{e}n convergente. Entonces,   la igualdad (\ref{eqn-28z}) aplicada a $n_{i_h}$ en lugar de $n_i$, implica   (para la subsucesi\'{o}n de \'{\i}ndices $\{n_{i}\}_i$), que el l\'{\i}mite del promedio $d_{n_{i}}$ de las distancias a $K_2$ es cero. Concluimos que toda subsucesi\'{o}n convergente de $\{d_n\}_{n \geq 1}$ converge a cero. Luego, la afirmaci\'{o}n (\ref{eqn-29}) est\'{a} probada.

\vspace{.2cm}

{\em Demostraci\'{o}n de $K_2 \subset K_1$:} Por hip\'{o}tesis $K_2$ es el m\'{\i}nimo compacto tal que $\mu(K_2) = 1$ para toda $\mu \in {\mathcal O}_f$. Luego, basta demostrar que \begin{equation}
\label{eqnAprobar-29}
\mbox{A probar: } \mu(K_1) = 1  \ \ \forall \ \mu \in {\mathcal O}_f.\end{equation}
Sean dados $\mu \in {\mathcal O}_f$ y $\epsilon >0$ arbitrario. Sea $\varphi$ la funci\'{o}n continua no negativa definida por \begin{equation}
\label{equation}
\varphi(y) := \mbox{dist}(y, K_1).\end{equation}   Por la definici\'{o}n de la topolog\'{\i}a d\'{e}bil$^*$, con la m\'{e}trica  $\mbox{dist}^*$ que se haya elegido en el espacio ${\mathcal M}$ de las probabilidades, existe $\delta >0$ tal que
\begin{equation}
\label{eqn-30}
\mbox{dist}^*(\nu, \mu)< \delta \ \Rightarrow \ \Big| \int \varphi \, d\mu - \int \varphi \, d \nu \Big| < \epsilon.\end{equation}
  Por la Definici\'{o}n \ref{definitionMedidaSRBlike} de medida SRB-like,
la medida de Lebesgue $m$ de la cuenca $A_{\delta}(\mu)$ de atracci\'{o}n estad\'{\i}stica $\delta$-d\'{e}bil de $\mu$, es positiva.
Como $K_1$ es un atractor de Ilyashenko $1$-obs., su cuenca de atracci\'{o}n estad\'{\i}stica $A(K_1)$ tiene medida de Lebesgue total. Luego deducimos que $$m\big(A(K_1) \cap A_{\delta}(\mu)\big) >0.$$ Tomemos un punto $x \in A(K_1) \cap A_{\delta}(\mu)$. Por la igualdad (\ref{eqnCuencaAtracIlyashenko}) que define  $A(K_1)$ tenemos: \begin{equation}
  \label{eqn-31}
  0 =\lim_{n \rightarrow + \infty} \frac{1}{n} \sum_{j= 0}^{n-1} \mbox{dist} (f^j(x), K_1) = \lim _{n \rightarrow + \infty} \int \varphi \, d (\sigma_{n,x}),\end{equation}
donde    $\sigma_{n,x}$ es la probabilidad emp\'{\i}rica dada definida por (\ref{(1)}). Por la Definici\'{o}n \ref{DefinicionCuencaDeAtraccionEpsilonEstadistica} de $A_{\delta}(\mu)$,   existe una subsucesi\'{o}n $n_i \rightarrow + \infty$ tal que $\{\sigma_{n_i}\}_{i \in \mathbb{N}}$ converge a una medida $\nu \in p \omega(x)$ tal que $\mbox{dist}^*(\mu, \nu) < \delta$. Luego, combinando la afirmaci\'{o}n (\ref{eqn-30}) con   la igualdad (\ref{eqn-31}), deducimos:
$$0 = \lim_{i \rightarrow + \infty} \int \varphi \, d \sigma_{n_i,x} = \int \varphi \, d \nu ,  \ \ \ \ \mbox{dist}^* (\nu, \mu) < \delta, \ \  \Rightarrow $$ $$  \int \varphi \, d \nu= 0, \ \ \   \Big|\int \varphi \, d \mu - \int \varphi \, d \nu \Big|< \epsilon \ \ \ \Rightarrow $$ $$  \Big| \int \varphi \, d \mu \Big|< \epsilon. $$
Como esta desigualdad vale para todo $\epsilon >0$, deducimos que
$0= \int \varphi \, d \mu.  $
Siendo $\varphi \geq 0$, concluimos que $\varphi(y) =\mbox{dist}(y, K_1) = 0$ para $\mu$-c.t.p. $y \in M$. Como $K_1$ es compacto, la distancia de un punto $y$ a $K_1$ es cero si y solo si $y \in K_1$. Hemos probado que $y \in K_1$ para $\mu$-c.t.p. $y \in M$, o dicho de otra forma, $\mu(K_1)= 1$ como quer\'{\i}amos demostrar.
 \hfill $\Box$

\subsection{Ejemplos de medidas SRB-like}
\begin{example} \em
\label{ejemploExpansorasEnElCirculo} {\bf Endormofirsmos $C^1$-expansores en $S^1$} \index{mapa expansor} \index{endomorfismo expansor}

 Sea $f: S^1 \mapsto S^1$ en el c\'{\i}rculo $S^1$ un endomorfismo de clase $C^1$ expansor, esto es, existe una constante $\sigma > 1$ tal que $|f'(z)| \geq \sigma > 1 \ \ \forall \ z \in S^1$. Por ejemplo en el c\'{\i}rculo $$S^1:= \{z \in \mathbb{C}: \ |z|= 1\}$$ la transformaci\'{o}n $f(z)= z^2$ es expansora. Un resultado cl\'{a}sico de din\'{a}mica topol\'{o}gica  (ver por ejemplo
 \cite[Theorem 2.4.6]{Katok-Hasselblatt}), establece que todo endormofismo expansor en el c\'{\i}rculo es sobreyectivo,   el n\'{u}mero de preim\'{a}genes de cualquier punto $x \in S^1$ es constante igual a $k \geq 2$  (llamado \em grado \em de $f$) y $f$ es topol\'{o}gicamente conjugado a $g_k$ definido por $g_k(z)= z^k$ para todo $z \in S^1= \{|z| = 1\}$.

  El endomorfismo expansor, puede mirarse como un endomorfismo (no invertible) uniformemente hiperb\'{o}lico (cf. Definition \ref{definicionAnosov}), en que el subespacio inestable $U_x$ es todo el espacio tangente $T_xS^1$, y el subespacio estable $E_x = \{0\}$. La variedad inestable $W^u (x_0)$ por un punto cualquiera $x_0 \in S^1$ se define por $$W^u (x_0) := \Big \{\ y_0 \in S^1: \ \ \exists \ x_{-n}, \ y_{-n} \in S^1 \mbox{ tales que }  $$ $$ f(x_{-n}) = x_{-n + 1}, \    f(y_{-n}) = y_{-n + 1} \ \forall \ n \geq 0, \ $$ $$\ \lim_{-n \rightarrow - \infty} \mbox{dist}(y_{-n}, x_{-n}) = 0 \ \Big \}.$$
  La variedad inestable de cualquier punto coincide con todo $S^1$:
  $$W^u(x_0) = S^1 \ \ \forall \ x_0 \in S^1.$$

  Por esta raz\'{o}n, para los endomorfismos expansores en el c\'{\i}rculo, la  medida  condicional  inestable  de $\mu$ (cf. Definici\'{o}n \ref{definitionmedidascondicionalesinestablesAC}) es la misma $\mu$. Entonces,  definimos:

   \vspace{.2cm}

   {\bf Definici\'{o}n (medida de Gibbs):} Si $f$ es un endomorfismo $C^1$ expansor en el c\'{\i}rculo $S^1$, decimos que una medida $f$-invariante $\mu$ \em es de Gibbs \em si  $\mu \ll m$, donde $m$ es la medida de Lebesgue en $S^1$. \index{medida! de Gibbs} \index{continuidad absoluta}

   \vspace{.2cm}

  Si $f$ es   de clase $C^{1 + \alpha}$, se tiene el siguiente resultado cl\'{a}sico, que da una versi\'{o}n del Teorema \ref{TheoremSRBanosov}  aplicable a los endomorfismos expansores del c\'{\i}rculo (en vez de a los difeomorfismos de Anosov en variedades de dimensi\'{o}n mayor que 1).

 \vspace{.3cm}

 {\bf Teorema (Ruelle):} \index{teorema! Ruelle} \index{mapa expansor} \index{endormorfismo expansor} \index{medida! de Gibbs} \index{medida! SRB} \index{medida! equivalencia de}

 \em Sea $f: S^1 \mapsto S^1$ de clase $C^{1 + \alpha}$ expansor en el c\'{\i}rculo $S^1$. Entonces:

  {\em (a)} Existe una \'{u}nica medida   de Gibbs $\mu$ \em (i.e. $\mu$ absolutamente continua respecto a la medida de Lebesgue) \em y esta medida $\mu$ es erg\'{o}dica.

  {\em (b)}  Toda medida SRB es de Gibbs y rec\'{\i}procamente. \em (Por lo tanto existe una \'{u}nica medida SRB y es erg\'{o}dica). \em

  {\em (c)} La medida de Gibbs $\mu$ es equivalente a la medida de Lebesgue $m$ \em (i.e. $\mu \ll m$ y $m \ll \mu$). \em

  {\em (d)} La cuenca de atracci\'{o}n estad\'{\i}stica de $\mu$ cubre Lebesgue c.t.p. $z \in S^1$. \em \index{cuenca de atracci\'{o}n! estad\'{\i}stica}

  \vspace{.3cm}

  La prueba de este Teorema  de Ruelle para mapas $C^{1+ \alpha}$ expansores en $S^1$, puede encontrarse por ejemplo en
  \cite[Theorem 5.1.16]{Katok-Hasselblatt}.

  Nota: Un resultado m\'{a}s general, que establece la existencia de medidas de Gibbs  para mapas $C^{1 }$ a trozos en el c\'{\i}rculo, que sean expansores en cada trozo, y que cumplan hip\'{o}tesis de variaci\'{o}n acotada, fue demostrado recientemente por Liverani en \cite{Liverani}.

  \vspace{.3cm}

  Como consecuencia  de la parte (d) del Teorema de Ruelle, en el caso $C^{1+ \alpha}$ expansor, no existen otras medidas SRB-like que no sean la medida SRB. En efecto, otra medida  $\nu \neq \mu$ tendr\'{\i}a una cuenca de atracci\'{o}n estad\'{\i}stica $\epsilon$-d\'{e}bil $A_{\epsilon}(\nu)$ que debe ser disjunta con la cuenca de atracci\'{o}n estad\'{\i}stica (fuerte) $B(\mu)$ de $\mu$ (porque $\nu \neq \mu$). Pero como $m(M \setminus B(\mu)) = 0$, entonces $m(A_{\epsilon}(\nu))= 0$ y $\nu$ no puede ser   SRB-like.
  \index{teorema! unicidad de medida SRB}

  Ahora veremos que la relaci\'{o}n entre medidas de Gibbs y la F\'{o}rmula (\ref{eqnformuladePesin}) de Pesin para la entrop\'{\i}a (que vale para difeomorfismos de clase $C^{1 + \alpha}$ seg\'{u}n vimos en el Teorema \ref{theoremFormulaPesin} de la secci\'{o}n anterior), tambi\'{e}n se generaliza para endormorfismos expansores de clase $C^{1 + \alpha}$ en el c\'{\i}rculo.

   \vspace{.3cm}


  {\bf Teorema }

  {\bf  F\'{o}rmula de Pesin para la Entrop\'{\i}a. Endomorfismos expansores}
  \index{Pesin! f\'{o}rmula de} \index{f\'{o}rmula de Pesin} \index{entrop\'{\i}a! f\'{o}rmula de Pesin} \index{mapa expansor} \index{endomorfismo expansor}
  \em
  Sea $f: S^1 \mapsto S^1$   expansor de clase $C^{2}$. Entonces, la \'{u}nica medida SRB $\mu$ para $f$ \em (que por el teorema de Ruelle es de Gibbs y equivalente a la medida de Lebesgue), \em satisface la F\'{o}rmula \em (\ref{eqnformuladePesin}) \em de Pesin para   la Entrop\'{\i}a, y es la \'{u}nica medida que satisface tal f\'{o}rmula. \em

  \vspace{.3cm}

  Una prueba de este Teorema, fue dada en \cite{PesinLyapunovExponents} por Pesin. Una versi\'{o}n para endomorfismos de clase $C^2$ se encuentra en \cite{QianZhu-TeoLedrYoungParaEndomorfismosC2} o en \cite{QianZhu-Libro}. M\'{a}s a\'{u}n, en \cite{QianZhu-TeoLedrYoungParaEndomorfismosC2} se prueba la versi\'{o}n  del Teorema de Ledrappier-Young \cite{Ledrappier-Young} para endomorfismos, que establece  la equivalencia, para todo endomorfismo  de clase $C^2$, entre las medidas de Gibbs y las medidas que satisfacen la F\'{o}rmula de Pesin para la Entrop\'{\i}a (cf.  Teorema \ref{theoremLedrappier-Young}).

  \vspace{.3cm}

  En el caso que el endormorfismo $f$ sea de clase $C^1$ pero no $C^{1+ \alpha}$, las afirmaciones del Teorema de Ruelle fallan. En efecto, en \cite{AvilaBochi} se prueba que los endormorfismos expansores $C^1$-gen\'{e}ricos no poseen medidas invariantes absolutamente continuas respecto de $m$ (no poseen medidas de Gibbs). M\'{a}s en general, en \cite{SchmittGora} se prueba (con un ejemplo expl\'{\i}cito) que, aunque $f$ sea discontinua, si $f$ es $C^1$ expansora a trozos, entonces la existencia de   medidas invariantes absolutamente continuas respecto de $m$ no es necesaria.

  Sin embargo, en \cite{CampbellQuas}  Campbell y Quas probaron que $C^1$-gen\'{e}ricamente, un endomorfismo  expansor  en el c\'{\i}rculo (que es $C^1$ pero no $C^{1 + \alpha}$), posee  una \'{u}nica medida SRB $\mu$  y que esta medida  es erg\'{o}dica,  satisface la F\'{o}rmula de Pesin para la Entrop\'{\i}a (\ref{eqnformuladePesin})  y tiene cuenca de atracci\'{o}n estad\'{\i}stica que cubre Lebesgue c.t.p.  Pero en vez de ser $\mu$ absolutamente continua respecto a la medida de Lebesgue $m$, se cumple todo lo contrario:
  $$\mu \perp m,$$
  es decir, $\mu$ es mutuamente singular respecto a la medida de Lebesgue $m$.

  Del resultado  de Campbell y Quas  deducimos que $C^1$-gen\'{e}ricamente para los endomorfismos expansores en el c\'{\i}rculo, existe una \'{u}nica medida SRB-like que es la medida SRB $\mu$ encontrada por Campbell y Quas. En efecto, argumentamos igual que antes, como lo hicimos a partir del Teorema de Ruelle para los expansores de clase $C^{1+ \alpha}$. Como la cuenca de atracci\'{o}n estad\'{\i}stica $B(\mu)$ de la medida SRB $\mu$ cubre Lebesgue c.t.p., entonces ninguna otra medida $\nu \neq \mu$ puede ser SRB-like.

 \vspace{.2cm}

 Sin necesidad de asumir hip\'{o}tesis de $C^1$-genericidad, todo endormorfismo expansor de clase $C^1$ en el c\'{\i}rculo tiene medidas SRB-like,   todas sus medidas SRB-like satisfacen la F\'{o}rmula de Pesin para la Entrop\'{\i}a (\ref{eqnformuladePesin}), y la uni\'{o}n de sus cuencas de atracci\'{o}n estad\'{\i}stica $\epsilon$-d\'{e}bil  cubre Lebesgue c.t.p. para todo $\epsilon >0$. Estos resultados fueron probados en \cite{CatEnrPortugaliae}.
Sin embargo, la unicidad de la medida SRB-like (que es cierta para los expansores de clase $C^{1+ \alpha} $ debido al Teorema de Ruelle), es falsa en general para   los expansores de clase $C^1$ que no son $C^{1 + \alpha}$ ni son $C^1$-gen\'{e}ricos. En efecto, en \cite{Quas} Quas construy\'{o} un expansor de clase $C^1$ que exhibe m\'{a}s de una medida SRB-like.
\end{example}

\begin{example} \em
\label{ejemploExpansorasC0EnElCirculo} {\bf Endormofirsmos $C^0$-expansores en $S^1$} \index{mapa expansor} \index{endomorfismo expansor}
\end{example}

Un mapa continuo $f: S^1 \mapsto S^1$ en el c\'{\i}rculo $S^1$ se llama \em endomorfismo $C^0$-expansor \em (cf.  \cite[Definition 2.4.1]{Katok-Hasselblatt}), si existen constantes $\delta >0$ y $\sigma > 1$ tales que
$$\mbox{dist}(f(x), f(y)) \geq \sigma \mbox{dist}(x,y) \ \ \forall \ x,y  \in S^1 \mbox{ tales  que } \mbox{dist}(x,y) < \delta.$$

En \cite{Misiurewicz}, Misiurewicz construy\'{o} un endormofismo   $C^0$-expansor $f$ en $S^1$, que satisface las hip\'{o}tesis del Corolario \ref{corolarioSRB-likeMisiurevicz}. Por lo tanto, en ese ejemplo  no existen medidas SRB  y toda medida $f$-invariante es SRB-like. Existe entonces una cantidad infinita no numerable de medidas SRB-like. Llamaremos a los endomorfismos que tienen esta propiedad \em endomorfismos de Misiurewicz. \em En \cite{AbdenurAndersson} se prueba que los endomorfismos de Misiurewicz son $C^0$-gen\'{e}ricos en el espacio de los endomorfismos $C^0$-expansores del c\'{\i}rculo $S^1$.

\begin{example} \em 
{\bf Anosov $C^1$} \label{ejemploSRB-likeAnosovC1} \index{difeomorfismos! de Anosov} \index{automorfismo! lineal en el toro} \index{transitividad} \index{transformaci\'{o}n transitiva}

  Sea $f$   un difeomorfismo de Anosov transitivo  en una variedad compacta y Riemanniana $M$ (cf. Definition \ref{definicionAnosov}).

  Primero repasemos el caso en que $f$ es adem\'{a}s de clase $C^{1 + \alpha}$ o en particular si es de clase $C^2$:

  En el Teorema \ref{TheoremSRBanosov} de  Sinai \cite{Sinai_SRB}, vimos que si   $f $ es un difeomorfismo de Anosov   de clase $C^{1 + \alpha}$, entonces existe una \'{u}nica medida de probabilidad SRB $\mu$, es de Gibbs erg\'{o}dica, y su cuenca de atracci\'{o}n estad\'{\i}stica cubre Lebesgue c.t.p. \index{medida! de Gibbs erg\'{o}dica} \index{cuenca de atracci\'{o}n! estad\'{\i}stica} \index{medida! SRB-like}
Deducimos que tal medida $\mu$ es la \'{u}nica medida SRB-like, argumentando de igual forma que en el ejemplo \ref{ejemploExpansorasEnElCirculo}.

\vspace{.2cm}

Adem\'{a}s, debido al Teorema \ref{theoremGibbs->SRB}, la \'{u}nica medida SRB-like de un difeomorfismo de Anosov transitivo  de clase $C^{1 + \alpha}$, satisface la F\'{o}rmula (\ref{eqnformuladePesin}) de Pesin para la Entrop\'{\i}a.

\vspace{.3cm}

Si el difeomorfismo $f$ de Anosov transitivo es de clase $C^2$, entonces su \'{u}nica medida     SRB $\mu$, es tambi\'{e}n la \'{u}nica medida de probabilidad que satisface la F\'{o}rmula de Pesin para la Entrop\'{\i}a. En efecto, por el Teorema \ref{theoremLedrappier-Young} de Ledrappier-Young \cite{Ledrappier-Young}, toda medida $\nu$ que satisfaga la F\'{o}rmula de Pesin para la Entrop\'{\i}a es de Gibbs. Y por el Teorema \ref{theoremGibbs->SRB} las componentes erg\'{o}dicas $\nu_x$ de toda medida de Gibbs es SRB, y por lo tanto es SRB-like. Luego, como existe una \'{u}nica medida SRB-like $\mu$, $\nu$ tiene una \'{u}nica componente erg\'{o}dica que es $\mu$, y por lo tanto $\nu$ es erg\'{o}dica y coincide con $\mu$.

\vspace{.2cm}

Ahora pasemos al caso de difeomorfismo de Anosov transitivo $f$ de clase $C^1$ pero no $C^{1 + \alpha}$. En este caso las demostraciones conocidas de los teoremas mencionados en el repaso anterior, no son aplicables, porque utilizan la Teor\'{\i}a de Pesin. Por ejemplo,  las pruebas utilizan el Teorema \ref{theoremTeoriaPesin} que establece la continuidad absoluta de la holonom\'{\i}a de la foliaci\'{o}n estable. Este resultado ser\'{\i}a falso si $f$ no fuera de clase $C^{1 + \alpha}$ (ver   \cite{RobinsonYoungContrajemploCAdeFoliacion}, o tambi\'{e}n \cite{Bowen_C1horseshoe}).

Sin embargo, en los \'{u}ltimos   a\~{n}os se han obtenido algunos resultados parciales, aplicables a los difeomorfismos de Anosov de clase $C^{1}$, que generalizan el Teorema \ref{TheoremSRBanosov}:

En \cite{QiuHyperbolicExisteUnicaSRBPesinFormula} Qiu y Zhu demostraron que $C^1$-gen\'{e}ricamente   los difeomorfismos de Anosov transitivos tienen una \'{u}nica medida SRB, esta medida es erg\'{o}dica,  satisface la F\'{o}rmula de Pesin para la Entrop\'{\i}a,   es la \'{u}nica medida que satisface tal f\'{o}rmula, y su cuenca de atracci\'{o}n estad\'{\i}stica cubre Lebesgue c.t.p. de la variedad. Por lo tanto $\mu$ es la \'{u}nica medida SRB-like. La prueba de Qiu y Zhu no pasa (a diferencia de la prueba del Teorema \ref{TheoremSRBanosov}) por la construcci\'{o}n de medidas de Gibbs. M\'{a}s a\'{u}n (tanto como la autora de este libro conoce) no se sabe si la \'{u}nica medida SRB $\mu$ de un difeomorfismo de Anosov transitivo y $C^1$-gen\'{e}rico, es medida de Gibbs. Pero se sabe que tal $\mu$ no es absolutamente continua respecto a la medida de Lebesgue $m$ en toda la variedad, pues \'{A}vila y Bochi \cite{AvilaBochi2006} probaron que $C^1$ gen\'{e}ricamente no existen medidas invariantes absolutamente continuas respecto de $m$.

\vspace{.2cm}

Cuando   la medida SRB-like $\mu$ es \'{u}nica  deducimos, como consecuencia del Teorema \ref{theoremOptimalidadSRB-like}, que $\mu$ es SRB y que su cuenca de atracci\'{o}n estad\'{\i}stica $B(\mu)$ cubre Lebesgue c.t.p. En este caso, por ejemplo para los difeomorfismos de Anosov transitivos $C^1$-gen\'{e}ricos del Teorema de Qiu y Zhu \cite{QiuHyperbolicExisteUnicaSRBPesinFormula}, es v\'{a}lido un   teoremas erg\'{o}dico, probado  recientemente por Kleptsyn y Ryzhov en \cite{kelpstynRyzov2012}, que estima  la velocidad de convergencia de las probabilidades emp\'{\i}ricas $\sigma_{n,x}$ (en la topolog\'{\i}a d\'{e}bil$^*$) a la medida SRB-$\mu$, para Lebesgue c.t.p.  $x \in M$.

\vspace{.2cm}

En el caso de difeomorfismo de Anosov $C^1$ no gen\'{e}rico y no $C^{1 + \alpha}$, cuando $f$ preserva una medida invariante $\mu$ absolutamente continua respecto a la medida de Lebesgue, Sun y Tian probaron en \cite{SunTianDominatedSplitting} que esta medida $\mu$ satisface la F\'{o}rmula de Pesin para la entrop\'{\i}a. Esto implica que, bajo las hip\'{o}tesis adicional de existencia de medida invariante $\mu$ equivalente a la medida de Lebesgue $ m$, toda medida SRB-like satisface tal f\'{o}rmula. En efecto, por el Teorema \ref{theoremDescoErgodicaEspaciosMetricos} de Descomposici\'{o}n Erg\'{o}dica,  $\mu$-c.t.p. $x$ pertenece a  la cuenca de atracci\'{o}n estad\'{\i}stica $B(\mu_x)$ de una componente erg\'{o}dica $\mu_x$ de $\mu$. Como $\mu$ es equivalente a la medida de Lebesgue, deducimos que Lebesgue c.t.p. $x$ pertenece a $B(\mu_x)$. Luego, por el Teorema \ref{theoremOptimalidadSRB-like}, estas componentes erg\'{o}dicas   $\mu_x$   son las medidas SRB-like. Como $\mu$ satisface la F\'{o}rmula de Pesin para la entrop\'{\i}a, entonces toda componente erg\'{o}dica $\mu_x$ de $\mu$ tambi\'{e}n satisface tal f\'{o}rmula (cf. \cite{Keller}). Luego, del resultado de Sun y Tian concluimos que toda medida SRB-like para $f$, satisface la F\'{o}rmula de Pesin para la Entrop\'{\i}a, si $f$ es un difeomorfismo de Anosov de clase $C^1$ que preserve una medida equivalente a la medida de Lebesgue.

\vspace{.3cm}

  M\'{a}s en general, sin necesidad de asumir hip\'{o}tesis de $C^1$-genericidad ni de existencia de medida invariante absolutamente continua respecto a la medida de Lebesgue, en \cite{CatCerEnrSubmitted} se prueba que para todo $f$ de Anosov   de clase $C^1$, toda medida SRB-like  $\mu$ satisface la F\'{o}rmula de Pesin para la Entrop\'{\i}a. Pero no sabemos si es necesario que adem\'{a}s $\mu$ sea medida de Gibbs, ni   que sea erg\'{o}dica.
  \end{example}
\vspace{.3cm}





\printindex

\end{document}